\documentclass[12pt,a4paper]{report}
\usepackage[latin1]{inputenc}
\usepackage{aeguill}
\usepackage{amsfonts}
\usepackage{amssymb}
\usepackage{xspace}
\usepackage{amsmath}
\usepackage{amsthm}
\usepackage{dsfont}
\usepackage{mathabx}
\usepackage[all,cmtip]{xy}
\usepackage{enumitem}
\usepackage[mathlines]{lineno}
\usepackage{amsmath}
\usepackage{xcolor}
\usepackage[colorlinks=true,linkcolor=blue,pagebackref,breaklinks]{hyperref}
\usepackage[hmargin={1.25in,1.25in},vmargin={1.5in,1.5in}]{geometry}
\usepackage{fancyhdr}

\newtheorem{thm}{Theorem}[section]
\newtheorem{prop}{Proposition}[section]
\newtheorem{df}{Definition}[section]
\newtheorem{lem}{Lemma}[section]
\newtheorem{cor}{Corollary}[section]
\newtheorem{dflem}{Definition-Lemma}[section]
\newtheorem{hyp}{Hypothesis}[section]
\newtheorem{ex}{Example}[section]
\newtheorem{conj}{Conjecture}[section]
\newtheorem{rem}{Remark}[section]

\newenvironment{dem}{\paragraph{Proof}}
{\begin{flushright}$\Box$\end{flushright}}

\newcommand{\N}{\mathbb{N}}
\newcommand{\Z}{\mathbb{Z}}
\newcommand{\Q}{\mathbb{Q}}
\newcommand{\C}{\mathbb{C}}
\newcommand{\Ss}{\mathbb{S}}
\newcommand{\R}{\mathbb{R}}
\newcommand{\LL}{\mathcal{L}}
\newcommand{\F}{\mathbb{F}}

\newcommand{\AK}{\mathbb{A}_{K}}
\newcommand{\AF}{\mathbb{A}_{F}}
\newcommand{\ALL}{\mathbb{A}_{\mathcal{L}}}

\newcommand{\AL}{\mathbb{A}_{L}}

\newcommand{\AFFF}{\mathbb{A}_{\mathcal{F}}}

\newcommand{\AFf}{\mathbb{A}_{F,f}}
\newcommand{\ALf}{\mathbb{A}_{L,f}}

\newcommand{\AFp}{\mathbb{A}_{F^{+}}}
\newcommand{\AFsub}{\mathbb{A}_{F_{0}}}
\newcommand{\AQ}{\mathbb{A}_{\mathbb{Q}}}
\newcommand{\AQf}{\mathbb{A}_{\mathbb{Q},f}}
\newcommand{\lieG}{\mathfrak{g}}
\newcommand{\lieK}{\mathfrak{k}}
\newcommand{\lieP}{\mathfrak{p}}
\newcommand{\lieH}{\mathfrak{h}}

\makeatletter 
\@addtoreset{equation}{section}
\makeatother  

\title{Special values of automorphic $L$-functions for $GL_{n}\times GL_{n'}$ over CM fields, factorization and functoriality of arithmetic automorphic periods}
\author{Jie LIN}
\date{\today}

\begin{document}
\maketitle

\section*{R\'esum\'e}
Michael HARRIS a défini les périodes arithmétiques automorphes pour certaines représentations cuspidales de $GL_{n}$ sur corps quadratiques imaginaires en 1997. Il a aussi montré que les valeurs critiques de fonctions L automorphes pour $GL_{n}\times GL_{1}$ peuvent être interprétées en termes de ces périodes. Dans la thèse, ses travaux sont généralisés sous deux aspects. D'abord, les périodes arithmétiques automorphes ont été définies pour tous corps CM. On montre aussi que ces périodes factorisent comme produits des périodes locales sur les places infinies. De plus, on montre que les valeurs critiques de fonctions L automorphes pour $GL_{n}\times GL_{n'}$ peuvent être interprétées en termes de ces périodes dans beaucoup de cas. Par conséquent on montre que les périodes sont fonctorielles pour l'induction automorphe et changement de base cyclique.\\
On aussi définit des périodes motiviques si le motif est restreint d'un corps CM au corps des nombres rationnels. On peut calculer la période de Deligne pour le produit tensoriel de deux tels motifs. On voit directement que nos résultats automorphes sont compatibles avec la conjecture de Deligne pour les motifs.

\section*{Abstract}
Michael HARRIS defined the arithmetic automorphic periods for certain cuspidal representations of $GL_{n}$ over quadratic imaginary fields in his Crelle paper 1997. He also showed that critical values of automorphic L-functions for $GL_{n}\times GL_{1}$ can be interpreted in terms of these arithmetic automorphic periods. In the thesis, we generalize his results in two ways. Firstly, the arithmetic automorphic periods have been defined over general CM fields. We also prove that these periods factorize as products of local periods over infinity places. Secondly, we show that critical values of automorphic $L$ functions for $GL_{n}\times GL_{n'}$ can be interpreted in terms of these automorphic periods in many situations. Consequently we show that the automorphic periods are functorial for automorphic induction and cyclic base change. \\
We also define certain motivic periods if the motive is restricted from a CM field to the field of rational numbers. We can calculate Deligne's period for tensor product of two such motives. We see directly that our automorphic results are compatible with Deligne's conjecture for motives.

\tableofcontents

\chapter*{Introduction}
\addcontentsline{toc}{chapter}{Introduction}{}
\text{}

Special values of $L$-functions play an important role in the Langlands program. Numerous conjectures predict that special values of $L$-functions reflect arithmetic properties of geometric objects. Most of these conjectures are still open and difficult to attack. \\

At the same time, concrete results on the special values of $L$-functions appear more and more in automorphic settings. For example, in \cite{harris97}, M. Harris constructed complex invariants called arithmetic automorphic periods and showed that the special values of automorphic $L$-function for $GL_{n}*GL_{1}$ could be interpreted in terms of these invariants. \\

We generalize his results in two ways. Firstly, the arithmetic automorphic periods have been defined over general CM fields. Secondly, we show that special values of arithmetic automorphic periods for $GL_{n}*GL_{n'}$ can be interpreted in terms of these arithmetic automorphic periods in many situations. In fact, we have found a concise formula for such critical values. This is our first main automorphic result. One possible application is to construct $p$-adic $L$-functions.\\

We remark that we have not finished the proof for $GL_{n}*GL_{1}$ over general CM fields in the current article. We shall do it later. We have assumed Conjecture \ref{special value for similitude unitary group} throughout the text. This is one important ingredient for automorphic results over general CM fields. \\

The results over quadratic imaginary field follow from the ideas in \cite{harrismotivic} and some technical calculation. Over general CM fields, one can still follow such arguments and get formulas for critical values in terms of arithmetic automorphic periods. But these formulas are ugly and complicated. In fact, we don't know how to write down a formula adapted to most cases. However, if one can show that the arithmetic automorphic periods can be factorized as products of local periods over infinite places, then the generalization to CM fields is straight forward.  \\

The factorization of arithmetic automorphic periods was actually a conjecture of Shimura (c.f. \cite{shimura83}, \cite{shimura88}). One possible way to show this is to define local periods geometrically and prove that special values of $L$-functions can be interpreted in terms of local periods. This was done by M. Harris for Hilbert modular forms in \cite{harrisCMperiod}. But it is extremely difficult to generalize his arguments to $GL_{n}$. Instead, we show that there are relations between arithmetic automorphic periods. These relations lead to a factorization which is our second main automorphic result. \\

We remark that the factorization is not unique. We show that there is a natural way to factorize such that the local periods are functorial for automorphic induction and base change. This is our third main automorphic result. We believe that local periods are also functorial for endoscopic transfer. We will try to prove this in the near future.\\

Although our local periods are not defined geometrically,  they must have geometric meanings. This may be done by defining certain geometric invariants and show that they are related to our local periods with the help of special values of $L$-functions. It is likely to show that our local periods are equal to the geometric invariants defined in \cite{harrisCMperiod} for Hilbert modular forms in this way.\\

\bigskip

On the other hand, Deligne's conjecture related critical values for motives over $\Q$ and Deligne's period (c.f. \cite{deligne79}). When the motive is the restriction to $\Q$ of the tensor product of two motives over a CM field, we may calculate Deligne's period in terms of motivic periods defined in \cite{harrisadjoint}. The formula was first given in \cite{harrisadjoint} when the motives are self-dual. We have dropped the self-dual condition here.\\

If the two motives are associated to automorphic representations of $GL_{n}$ and $GL_{n'}$ respectively, we may define motivic periods which are analogues of the arithmetic automorphic periods. We get a formula of Deligne's period in terms of these motivic periods. Our main motivic result says that our formula for automorphic $L$-functions are at least formally compatible with Deligne's conjecture.\\

\bigskip
\paragraph{Theorems:}

\text{}\\

Let $K$ be a quadratic imaginary field and $F\supset K$ be a CM field of degree $d$ over $K$. We fix an embedding $K\hookrightarrow \C$. Let $\Sigma_{F;K}$ be the set of embeddings $\sigma:F\hookrightarrow \C$ such that $\sigma\mid_{K}$ is the fixed embedding.\\

Let $E$ be a number field. Let $\{ a(\sigma) \} _{\sigma\in Aut(\C/K)}$,  $\{ b(\sigma) \} _{\sigma\in Aut(\C/K)}$ be two families of complex numbers. 
Roughly speaking, we say $a\sim_{E;K} b$ if $a=b$ up to multiplication by elements in $E^{\times}$ and equivariant under $G_{K}$-action.\\

Let $\Pi$ be a cuspidal cohomological representation of $GL_{n}(\AF)$ which has definable arithmetic automorphic periods (c.f. Definition \ref{definable periods}). In particular, we know that $\Pi_{f}$ is defined over a number field $E(\Pi)$. For any $I:\Sigma_{F;K} \rightarrow \{0,1,\cdots,n\}$, we may define the arithmetic automorphic periods $P^{(I)}(\Pi)$ as the Petersson inner product of a rational vector in a certain cohomology space associated to a unitary group of infinity sign $I$. It is a non zero complex number well defined up to multiplication by elements in $E(\Pi)^{\times}$. \\

We assume that Conjecture \ref{special value for similitude unitary group} is true. Our second main automorphic result mentioned above is as follows (c.f. Theorem \ref{special factorization theorem}):
\begin{thm}\label{second automorphic}
If conditions in Theorem \ref{complete theorem} are satisfied, in particular, if $\Pi$ is regular enough, then there exists some complex numbers $P^{(s)}(\Pi,\sigma)$ unique up to multiplication by elements in $(E(\Pi))^{\times}$ such that the following two conditions are satisfied:
\begin{enumerate}
\item $P^{(I)}(\Pi) \sim_{E(\Pi);K} \prod\limits_{\sigma\in\Sigma_{F;K}}P^{(I(\sigma))}(\Pi,\sigma)$ for all $I=(I(\sigma))_{\sigma\in\Sigma_{F;K}}\in \{0,1,\cdots,n\}^{\Sigma_{F;K}}$
\item  and $P^{(0)}(\Pi,\sigma)\sim_{E(\Pi);K} p(\widecheck{\xi_{\Pi}},\overline{\sigma})$
\end{enumerate}
where $\xi_{\Pi}$ is the central character of $\Pi$, $\widecheck{\xi_{\Pi}}:=\xi_{\Pi}^{-1,c}$ and $p(\widecheck{\xi_{\Pi}},\overline{\sigma})$ is the CM period (c.f. Section \ref{section CM periods}).
\end{thm}

We now introduce our first main automorphic result. Let $\Pi'$ be a cuspidal cohomological representation of $GL_{n'}(\AF)$ which has definable arithmetic automorphic periods. For any $\sigma\in\Sigma_{F;K}$, we may define the split indices $sp(j,\Pi;\Pi',\sigma)$ and $sp(k,\Pi';\Pi,\sigma)$ for $0\leq j\leq n$ and $0\leq k\leq n'$ (c.f. Definition \ref{splitindex}). Roughly speaking, we have:

\begin{thm}\label{first automorphic}
If $m\in \Z+\frac{n+n'}{2}$ is critical for $\Pi\times \Pi'$ then 
\begin{equation}L(m,\Pi\times \Pi')\sim_{E(\Pi)E(\Pi');K}(2\pi i)^{nn'md}\prod\limits_{\sigma\in\Sigma_{F;K}}(\prod\limits_{j=0}^{n}P^{(j)}(\Pi,\sigma)^{sp(j,\Pi;\Pi',\sigma)}
\prod\limits_{k=0}^{n'}P^{(k)}(\Pi',\sigma)^{sp(k,\Pi';\Pi,\sigma)})\nonumber
\end{equation} in the following cases:
\begin{enumerate}
\item $n'=1$ and $m$ is bigger than the central value.
\item $n>n'$ and $m\geq 1/2$, both $\Pi$ and $\Pi'$ are conjugate self-dual and the pair $(\Pi,\Pi')$ is in good position (c.f. Definition \ref{definition good position}).
\item $m=1$, both $\Pi$ and $\Pi'$ are conjugate self-dual and the pair $(\Pi,\Pi')$ is regular enough.
\end{enumerate}
\end{thm}

Our third main automorphic result says that the periods are functorial for automorphic induction and base change. Roughly speaking, we have:

\begin{thm}\label{third automorphic}
 (a) Let $\mathcal{F}/F$ be a cyclic extension of CM fields of degree $l$ and $\Pi_{\mathcal{F}}$ be a cuspidal representation of $GL_{n}(\AFFF)$. We write $AI(\Pi_{\mathcal{F}})$ for the automorphic induction of $\Pi_{\mathcal{F}}$. We assume both $AI(\Pi_{\mathcal{F}})$ and $\Pi_{\mathcal{F}}$ have definable arithmetic automorphic periods.

Let $I_{F}\in \{0,1,\cdots,nl\}^{\Sigma_{F;K}}$. We may define $I_{\mathcal{F}}\in\{0,1,\cdots, n\}^{\Sigma_{\mathcal{F};K}}$ as in Lemma \ref{AI infinity sign}. Or locally let $0\leq s\leq nl$ be an integer and $s(\cdot)$ be as in Definition \ref{AI infinity sign local}. We have:
\begin{eqnarray}
P^{(I_{F})}(AI(\Pi_{\mathcal{F}}))\sim_{E(\Pi_{\mathcal{F}});K} P^{(I_{\mathcal{F}})}(\Pi_{\mathcal{F}}) \nonumber \\
\text{or locally } P^{(s)}(AI(\Pi_{\mathcal{F}},\tau) \sim_{E(\Pi_{\mathcal{F}});K} \prod\limits_{\sigma \mid \tau} P^{(s(\sigma))}(\Pi_{\mathcal{F}},\sigma).\nonumber
\end{eqnarray}
(b) Let $\pi_{F}$ be a cuspidal representation of $GL_{n}(\AF)$. We write $BC(\pi_{F})$ for its strong base change to $\mathcal{F}$. We assume that both $\pi_{F}$ and $BC(\pi_{F})$ have definable arithmetic automorphic periods.

Let $I_{F}\in \{0,1,\cdots,n\}^{\Sigma_{F;K}}$. We write $I_{\mathcal{F}}$ the composition of $I_{F}$ and the restriction of complex embeddings of $\mathcal{F}$ to $F$.

We then have:
\begin{eqnarray}
&&P^{(I_{\mathcal{F}})}(BC(\pi_{F})) \sim_{E(\pi_{F});K} p^{I_{F}}(\pi_{F})^{l} \nonumber\\
&\text{or locally}& P^{(s)}(BC(\pi_{F}),\sigma)^{l} \sim_{E(\pi_{F});K}P^{(s)}(\pi_{F},\sigma\mid_{F})^{l}.\nonumber
\end{eqnarray}
Consequently, we know 
\begin{equation}
P^{(s)}(BC(\pi_{F}),\sigma) \sim_{E(\pi_{F})} \lambda^{(s)}(\pi_{F},\sigma)P^{(s)}(\pi_{F},\sigma\mid_{F}).\nonumber
\end{equation}
where $\lambda^{(s)}(\pi_{F},\sigma)$ is an algebraic number whose $l$-th power is in $E(\pi_{F})^{\times}$.

\end{thm}

We now introduce the motivic results. Let $M$, $M'$ be motives over $F$ with coefficients in $E$ and $E'$ of rank $n$ and $n'$ respectively. We assume that $M\otimes M'$ has no $(\omega/2,\omega/2)$-class. We may define motivic periods $Q^{(t)}(M,\sigma)$ for $0\leq t\leq n$ and $\sigma\in\Sigma_{F;K}$. We can calculate Deligne's period of $Res_{F/\Q}(M\otimes M')$ in terms of these periods. If $M$ and $M'$ are motives associated to $\Pi$ and $\Pi'$, Deligne's conjecture is equivalent to the following conjecture:
\begin{conj}
If $m\in \Z+\frac{n+n'}{2}$ is critical for $\Pi\times \Pi'$ then 
\begin{eqnarray}
&L(m,\Pi\times \Pi')=L(m+\frac{n+n'-2}{2},M\otimes M')&\nonumber\\
&\sim_{E(\Pi)E(\Pi');K}(2\pi i)^{mnn'd}\prod\limits_{\sigma\in\Sigma_{F;K}} (\prod\limits_{j=0}^{n}Q^{(j)}(M,\sigma)^{sp(j,\Pi;\Pi',\sigma)}
\prod\limits_{k=0}^{n'}Q^{(k)}(M',\sigma)^{sp(k,\Pi';\Pi,\sigma)})\nonumber&
\end{eqnarray}

\end{conj}

We see that it is compatible with Theorem \ref{first automorphic}. The main point of the proof is to fix proper basis. Deligne's period is defined by rational basis. The basis that we have fixed are not rational. But they are rational up to unipotent transformation matrices. We can still use such basis to calculate determinant. 

\bigskip

\paragraph{Idea of the proof for automorphic results:}
\text{}

Blasius has shown that special values of $L$-functions for Hecke character are related to CM periods. The proof of our automorphic results involve this fact  and the following three main ingredients:

\text{}\\
\textbf{Ingredient A} is Theorem \ref{n*1}. If follows from Conjecture \ref{special value for similitude unitary group}. It says that if $\chi$ is a Hecke characters then critical values $L(m,\Pi\otimes \chi)$ can be written in terms of the arithmetic automorphic periods of $\Pi$ and CM periods of $\chi$.

\text{}\\
\textbf{Ingredient B} is Theorem $3.9$ of \cite{harrismotivic}. It says that if $\Pi^{\#}$ is a certain automorphic representation of $GL_{n-1}(\AF)$ such that $(\Pi,\Pi^{\#})$ is in good position then critical values $L(m,\Pi\otimes \Pi')$ are products of the Whittaker period $p(\Pi)$, $p(\Pi^{\#})$ and an archimedean factor. The advantage of the results in \cite{harrismotivic} is that we don't need $\Pi^{\#}$ to be cuspidal. This gives us large freedom to choose $\Pi^{\#}$.

\text{}\\
\textbf{Ingredient C} is a calculation of Whittaker period $p(\Pi^{\#})$ when $\Pi^{\#}$ is the Langlands sum of cuspidal representations $\Pi_{1}, \cdots, \Pi_{l}$. Following the idea in \cite{mahnkopf05} and \cite{harrismotivic}, we know $p(\Pi^{\#})$ equals to product of $p(\Pi_{i})$ and the value at identity of a certain Whittaker function. Shahidi's calculation in \cite{shahidi10} shows that the latter is related to $\prod\limits_{1\leq i<j\leq l}L(1,\Pi_{i}\times \Pi_{j}^{c})$.\\

The proof of the case (a) in Theorem \ref{third automorphic} is relatively simple. It is enough to take suitable algebraic Hecke character $\eta$ of $F$ and calculate $L(m,AI(\Pi_{\mathcal{F}})\otimes \eta) = L(m,\Pi_{\mathcal{F}}\otimes \eta\circ N_{\AFFF^{\times}/\AF^{\times}})$  by ingredient A. 

The idea for the case (b) is similar. But we have to show that the arithmetic automorphic periods of $BC(\pi_{F})$ are $Gal_{\mathcal{F}/F}$-invariant. This is due to the fact that $BC(\pi_{F})$ itself is $Gal_{\mathcal{F}/F}$-invariant.\\

 We now explain the proof for Theorem \ref{second automorphic} and Theorem \ref{first automorphic}.\\
 
\textbf{Step 0:} determine when a function can factorize through each factor. \\

For example, let $X$, $Y$ be two sets and $f$ be a map from $X\times Y$ to $\C^{\times}$. Then there exists functions $g:X\rightarrow \C^{\times}$ and $h:Y \rightarrow \C^{\times}$ such that $f(x,y)=g(x)h(y)$ for any $x\in X$ and $y\in Y$ if and only if $f(x,y)f(x',y')=f(x,y')f(x',y)$ for any $x,x'\in X$ and $y,y'\in Y$. Therefore, to show that the arithmetic automorphic periods factorize is equivalent to show that there are certain relations between these periods.
 
 \bigskip
 
\text{}\\
\textbf{Step 1:} interpret $p(\Pi)$ in terms of arithmetic automorphic periods. \\

The idea is the same as in \cite{harrismotivic}. We take $\Pi^{\#}$ to be the Langlands sum of Hecke characters $\chi_{1},\cdots,\chi_{n-1}$. We have $L(m,\Pi\times \Pi^{\#})=\prod\limits_{1\leq i\leq n-1}L(m,\Pi\otimes \chi_{i})$.

Ingredient B says that the left hand side equals to the product of $p(\Pi)$, $p(\Pi^{\#})$ and an archimedean factor. Ingredient C tells us that $p(\Pi^{\#})$ is almost $\prod\limits_{1\leq i<j\leq l}L(1,\chi_{i}\times \chi_{j}^{c})$ which equals to product of CM periods by Blasius's result. Therefore, the left hand side equals to product of $p(\Pi)$, the CM periods of $\chi_{i}$ and an archimedean factor.

We may calculate the right hand side by Ingredient $A$. We get that the right hand side equals to product of the arithmetic automorphic periods of $\Pi$, the CM periods of $\chi_{i}$ and a power of $2\pi i$.

Comparing both sides, we will see unsurprisingly that the CM periods of $\chi_{i}$ in two sides coincide. We will get a formula for $p(\Pi)$ in terms of arithmetic automorphic periods. Varying the Hecke characters $\chi_{i}$, we get different formulas for $p(\Pi)$ in terms of arithmetic automorphic periods. We then deduce relations between arithmetic automorphic periods. The factorization property then follows.

We remark that the above procedure can only treat the case when $I(\sigma)\neq 0$ or $n$ for all $\sigma$. The proof for general case is more tricky (see section \ref{complete case}).

 \bigskip
 
\text{}\\
\textbf{Step 2:} repeat step 1 with suitable $\Pi^{\#}$. 

For example, if $n>n'$ and the pair $(\Pi,\Pi')$ is in good position, we may take $\Pi^{\#}$ to be the Langlands sum of $\Pi'$ and some Hecke characters $\chi_{1},\chi_{2},\cdots,\chi_{l}$ where $l=n-n'-1$ such that $(\Pi,\Pi^{\#})$ is in good position. We have $L(m,\Pi\times \Pi^{\#})=L(m,\Pi\times \Pi')\prod\limits_{1\leq i\leq l}L(\Pi\otimes\chi_{i})$.

Again, we calculate the left hand side by ingredient B and ingredient C. We apply step $1$ to $p(\Pi)$ and $p(\Pi')$ and we will get that the left hand side equals to product of arithmetic automorphic periods for $\Pi$ and $\Pi'$ and CM periods of $\chi_{i}$.

We then apply ingredient A to $L(\Pi\otimes\chi_{i})$ and compare both sides. We will get a formula for $L(m,\Pi\times \Pi')$.

For the case where $m=1$, we may take $\Pi^{\#}$ to be the Langlands sum of $\Pi$ and $\Pi'^{c}$. We know that $L(1,\Pi\times \Pi')$ then appears in the calculation of $p(\Pi^{\#})$ by ingredient $C$.

 \bigskip
 
\text{}\\
\textbf{Step 3:} Simplify the archimedean factors.

Once we get a formula of $L(m,\Pi\times \Pi')$ in terms of arithmetic automorphic periods, we may replace $\Pi$ and $\Pi'$ by representations which are automorphic inductions of Hecke characters. Blasius's result says that $L(m,\Pi\times \Pi')$ is equivalent to the product of a power of $2\pi i$ and some CM periods. On the other hand, the arithmetic automorphic periods of $\Pi$ are related to CM periods by Theorem \ref{third automorphic}. We shall deduce that the archimedean factor is equivalent to a power of $2\pi i$ if $\Pi$ and $\Pi'$ are induced from Hecke characters. We can finish the proof by noticing that such representations can have any infinity type.

\bigskip

\paragraph{Plan for the text:}
\text{}\\

In chapter $1$ we introduce our basic notation, in particular, the split index.\\

In chapter $2$ we introduce the base change theory for similitude unitary groups which will help us understanding the descending condition in the definition of arithmetic automorphic periods.\\

We summarize some results on rational structures in Chapter $3$. They play an important role in the proof. In particular, the ingredients $B$ and $C$ are introduced in the second half of this chapter.\\

In chapter $4$ we construct the arithmetic automorphic periods. We generalize the construction of \cite{harris97} to general CM fields. \\

Chapter $5$ contains the details for ingredient $A$. We remark that we have made a hypothesis here (c.f. Conjecture \ref{special value for similitude unitary group}). We will prove it in a forthcoming paper.\\

The motivic results are contained in Chapter $6$. This chapter is independent of others. We show that our main automorphic results are compatible with Deligne's conjecture for motives.\\

We prove the factorization of arithmetic automorphic periods in Chapter $7$ (c.f. Theorem \ref{second automorphic}). This result itself is very important. It is also the crucial step to generalize our results to CM fields. \\

In chapter $8$ we prove that the global and local arithmetic periods are functorial for automorphic induction and base change (c.f. Theorem \ref{third automorphic}). This is a direct corollary of the ingredient $A$ in chapter $5$ and the factorization property in chapter $7$.\\

In chapter $9$ we claim our main conjecture which is an automorphic analogue of Deligne's conjecture. We also claim our main theorem there, namely, Theorem \ref{first automorphic}. Moreover, in the last section of this chapter, we explain why the generalization from quadratic imaginary fields to CM fields is direct by the factorization property. \\

The last two chapters contain the details of the proof for Theorem \ref{first automorphic}. The calculation is not difficult but technique.

\chapter{Notation}
\section{Basic notation}
\text{}

We fix an algebraic closure $\overline{\Q}\hookrightarrow \C$ of $\Q$ and $K\hookrightarrow \overline{\Q}$ a quadratic imaginary field. We denote by $\iota$ the complex conjugation of the fixed embedding $K\hookrightarrow \overline{\Q}$.

We denote by $c$ the complex conjugation on $\C$. Via the fixed embedding $\overline{\Q}\hookrightarrow \C$, it can be considered as an element in $Gal(\overline{\Q}/\Q)$.  

For any number field $L$, let $\AL$ be the adele ring of $L$ and $\ALf$ be the finite part of $\AL$. We denote by $\Sigma_{L}$ the set of embeddings from $L$ to $\overline{\Q}$. If $L$ contains $K$, we write $\Sigma_{L;K}$ for the subset of $\Sigma_{L}$ consisting of elements which is the fixed embedding $K\hookrightarrow \overline{\Q}$ when restricted to $K$.\\

Throughout the text, we fix $\psi$ an algebraic Hecke character of $K$ with infinity type $z^{1}\overline{z}^{0}$ such that $\psi\psi^{c}=||\cdot||_{\AK}$ (see Lemma $4.1.4$ of \cite{CHT} for its existence). It is easy to see that the restriction of $||\cdot||_{\AK}^{\frac{1}{2}}\psi$ to $\AQ^{\times}$ is the quadratic character associated to the extension $K/\Q$ by the class field theory. Consequently our construction is compatible with that in \cite{harrismotivic}.\\

Let $F^{+}$ be totally real field of degree $d$ over $\Q$. We define $F:=F^{+}K$ a CM field. We take $\psi_{F}$ an algebraic Hecke character of $F$ with infinity type $z^{1}$ at each $\sigma\in \Sigma$ such that $\psi_{F}\psi_{F}^{c}=||\cdot||_{\AF}$.\\

For $z\in \C$, we write $\bar{z}$ for its complex conjugation. For $\sigma\in \Sigma_{F}$, we define $\bar{\sigma}:= \sigma^{c}$ the complex conjugation of $\sigma$.

Let $\eta$ be a Hecke character of $F$. We define $\widecheck\eta=\eta^{-1,c}$ and $\widetilde{\eta}=\eta/\eta^{c}$ two Hecke characters of $F$.

Let $n$ be an integer greater or equal to $2$.

\begin{df} Let $N$ be an integer and $\Pi$ be an automorphic representation of $GL_{n}(\AF)$. Let $\sigma$ be an element in $\Sigma_{F;K}$. We denote the infinity type of $\Pi$ at $\sigma$ by $(z^{a_{i}(\sigma)}\overline{z}^{a_{i}(\sigma)'})_{1\leq i\leq n}$. We may assume that $a_{1}(\sigma)\geq a_{2}(\sigma)\geq \cdots\geq a_{n}(\sigma)$ for all $\sigma\in\Sigma_{F;K}$. The representation  $\Pi$ will be called:
\begin{enumerate}
\item \textbf{pure of weight $\omega(\Pi)$} if $a_{i}(\sigma)+a_{i}(\sigma)'=-\omega(\Pi)$ for all $1\leq i\leq n$ and all $\sigma$;
\item \textbf{algebraic} if $a_{i}(\sigma),a_{i}(\sigma)'\in \Z+\frac{n-1}{2}$ for all $1\leq i\leq n$ and all $\sigma$;
\item \textbf{cohomological} if there exists $W$ an irreducible algebraic finite dimensional representation of $GL_{n}(F\otimes_{\Q}\R)$ such that $H^{*}(\lieG_{\infty},F_{\infty};\Pi\otimes W)\neq 0$ (see section \ref{cohomology} for more details);
\item \textbf{regular} if it is pure and $a_{i}(\sigma)-a_{i+1}(\sigma)\geq 1$ for all $1\leq i\leq n-1$ and all $\sigma$.
\item \textbf{$N$-regular} if it is pure and $a_{i}(\sigma)-a_{i+1}(\sigma)\geq N$ for all $1\leq i\leq n-1$ and all $\sigma$.
\end{enumerate}
\end{df}

Finally, let $E\supset K$ be a number field. We now define the relation $\sim_{E;K}$.

Let $\{ a(\sigma) \} _{\sigma\in Aut(\C/K)}$,  $\{ b(\sigma) \} _{\sigma\in Aut(\C/K)}$ be two families of complex numbers.

\begin{df}
We say $a\sim_{E;K} b$ if one of the following conditions is verified:
\begin{enumerate}[label=(\roman*)]
\item $a(\sigma)=0$ for all $\sigma\in Aut(\C/K)$, 
\item $b(\sigma)=0$ for all $\sigma\in Aut(\C/K)$, or
\item $a(\sigma)\neq 0$, $b(\sigma)\neq 0$ for all $\sigma$ and there exists $t\in E^{\times}$ such that $a(\sigma)=\sigma(t)b(\sigma)$ for all $\sigma\in Aut(\C/K)$.

\end{enumerate}

\end{df}

\begin{rem}\label{transitive}
\begin{enumerate}
\item Note that this relation is symmetric but not transitive. More precisely, if $a\sim_{E;K}b$ and $a\sim_{E;K}c$, we do not know whether $b\sim_{E;K} c$ in general unless the condition $a\neq 0$ is provided.
\item If $a\sim_{E;K} b$ with $b(\sigma)\neq 0$, we see that $\frac{a(\sigma)}{b(\sigma)}$ is contained in the Galois closure of $E$. In particular, it is an algebraic number.

\item If moreover, $a(\sigma)=a(\sigma')$ for $\sigma,\sigma'\in Aut(\C)$ such that $\sigma|_{E}=\sigma'|_{E}$, we can then define $a(\sigma)$ for $\sigma\in \Sigma_{\mathcal{F};K}$ by taking any $\widetilde{\sigma}\in Aut(\C)$, a lifting of $\sigma$, and define $a(\sigma):=a(\widetilde{\sigma})$. 
We identify $\C^{\Sigma_{\mathcal{F};K}}$ with $E\otimes_{K}\C$. We consider $A:=(a(\sigma))_{\sigma\in \Sigma_{\mathcal{F};K}}$ as an element in $E\otimes_{K}\C$. 

We assume the same condition for $b$ and define $B$ for $b$ similarly. It is easy to verify that $a\sim_{E;K}b$ if and only if one of the three conditions is verified: $A=0$, $B=0$, or $B\in (E\otimes_{K}\C)^{\times}$ and $AB^{-1}\in E^{\times}\subset (E\otimes_{K}\C)^{\times}$.

We remark that our results will be in this case.

\end{enumerate}

\end{rem}

\begin{lem}\label{anotherdef}
We assume $b(\sigma)\neq 0$ for all $\sigma$. We also assume that $a(\sigma)=a(\sigma')$ and $b(\sigma)=b(\sigma')$ if $\sigma|_{E}=\sigma'|_{E}\text{ for any }\sigma,\sigma'\in Aut(\C)$.
We have $a\sim_{E;K} b$ if and only if
\begin{equation}\nonumber
\tau\left(\cfrac{a(\sigma)}{b(\sigma)}\right)=\cfrac{a(\tau\sigma)}{b(\tau\sigma)}
\end{equation}
 for all $\tau\in Aut(\C/K)$ and $\sigma\in \Sigma_{\mathcal{F};K}$.

\end{lem}

\section{Split index and good position for automorphic pairs}
\begin{df}(\textbf{Split Index})\label{splitindex}

Let $n$ and $n'$ be two positives integers.

Let $\Pi$ and $\Pi'$ be two regular pure representations of $GL_{n}(\AF)$ and $GL_{n'}(\AF)$ respectively. Let $\sigma$ be an element of $\Sigma_{F;K}$. We denote  the infinity type of $\Pi$ and $\Pi'$ at $\sigma$ by $(z^{a_{i}(\sigma)}\overline{z}^{-\omega(\Pi)-a_{i}(\sigma)})_{1\leq i\leq n}$, $a_{1}(\sigma)>a_{2}(\sigma)>\cdots>a_{n}(\sigma)$ and $(z^{b_{j}(\sigma)}\overline{z}^{-\omega(\Pi')-b_{j}(\sigma)})_{1\leq j\leq n'}$, $b_{1}(\sigma)>b_{2}(\sigma)>\cdots>b_{n'}(\sigma)$ respectively. We assume that $a_{i}(\sigma)+b_{j}(\sigma)\neq -\frac{\omega(\Pi)+\omega(\Pi')}{2}$ for all $1\leq i\leq n$ all $1\leq j\leq n'$ and all $\sigma$.

We split the sequence $(a_{1}(\sigma)>a_{2}(\sigma)>\cdots>a_{n}(\sigma))$ with the numbers

 \begin{equation}\nonumber
 -\frac{\omega(\Pi)+\omega(\Pi')}{2}-b_{n'}(\sigma)>-\frac{\omega(\Pi)+\omega(\Pi')}{2}-b_{n'-1}(\sigma)>\cdots>-\frac{\omega(\Pi)+\omega(\Pi')}{2}-b_{1}(\sigma).\end{equation}
 
 This sequence is split into $n'+1$ parts. We denote the length of each part by $sp(0,\Pi':\Pi ,\sigma),sp(1,\Pi';\Pi,\sigma),\cdots,sp(n',\Pi';\Pi,\sigma)$, and call them the \textbf{split indices}.
 
\bigskip

\end{df}
\begin{lem}\label{split lemma}
Let $n$, $n'$, $\Pi$ and $\Pi'$ be as in the above definition. Let $\sigma$ be an element in $\Sigma_{F;K}$. Let $\eta$ be an algebraic Hecke character of $\AF$. Let $0\leq j\leq n'$ be an integer. We have the following formulas:

\begin{enumerate}
\item $\sum\limits_{i=0}^{n'}sp(i,\Pi';\Pi,\sigma)=n=dim(\Pi)$.
\item $sp(j,\Pi';\Pi,\sigma)=sp(n'-j,\Pi'^{c};\Pi^{c},\sigma)=sp(n'-j,\Pi'^{\vee};\Pi^{\vee},\sigma)$.
\item For any $t,s\in \R$, $sp(j,\Pi'\otimes ||\cdot||_{\AK}^{t};\Pi,\sigma)=sp(j,\Pi';\Pi\otimes ||\cdot||_{\AK}^{s},\sigma)=sp(j,\Pi';\Pi,\sigma)$ .
\item $sp(j,\Pi'\otimes\eta; \Pi,\sigma)=sp(j,\Pi';\Pi\otimes \eta,\sigma)$ and $sp(j,\Pi'\otimes\eta^{c};\Pi,\sigma)=sp(j,\Pi'\otimes\eta^{-1};\Pi,\sigma)$. Similarly, $sp(j,\Pi';\Pi\otimes \eta^{c},\sigma)=sp(j,\Pi';\Pi\otimes \eta^{-1},\sigma)$.
\end{enumerate}
\end{lem}

The first two points of the above lemma are direct. For the remaining, we only need to notice that calculating the split index is nothing but comparing $a_{i}(\sigma)+b_{j}(\sigma)$ with $-\frac{\omega(\Pi)+\omega(\Pi')}{2}$.

\bigskip

\begin{ex}
\begin{enumerate}
\item If $F^{+}=\Q$, $n=5$, $n'=4$, $\omega(\Pi)=\omega(\Pi')=0$ and 

\begin{equation}\nonumber
-b_{4}>\textbf{a}_{1}>\textbf{a}_{2}>-b_{3}>-b_{2}>\textbf{a}_{3}>\textbf{a}_{4}>-b_{1}>\textbf{a}_{5},\end{equation}

we have $sp(0,\Pi';\Pi)=0$, $sp(1,\Pi';\Pi)=2$, $sp(2,\Pi';\Pi)=0$, $sp(3,\Pi';\Pi)=2$ and $sp(4,\Pi';\Pi)=1$. We verify that $sp(0,\Pi';\Pi)+sp(1,\Pi';\Pi)+sp(2,\Pi';\Pi)+sp(3,\Pi';\Pi)+sp(4,\Pi';\Pi)=5=dim(\Pi)$ as expected by the previous lemma.

\item If $F^{+}=\Q$, $n'=n-1$, $\omega(\Pi)=\omega(\Pi')=0$ and $\textbf{a}_{1}>-b_{n-1}>\textbf{a}_{2}>-b_{n-2}>\cdots>\textbf{a}_{n-1}>-b_{1}>\textbf{a}_{n}$, we have $sp(j,\Pi';\Pi)=1$ for all $0\leq j\leq n-1$.

Moreover, $sp(k,\Pi;\Pi')=1$ for all $1\leq k\leq n-1$, $sp(0,\Pi;\Pi')=0$ and $sp(n,\Pi;\Pi')=0$.

We verify that $\sum\limits_{j=0}^{n'}sp(j,\Pi';\Pi)=n=dim(\Pi)$ and $\sum\limits_{j=0}^{n}sp(j,\Pi;\Pi')=n-1=dim(\Pi')$ as expected.

\end{enumerate}
\end{ex}

\bigskip

\begin{df}\label{definition good position}
We assume that $n>n'$. Let $\Pi$ and $\Pi'$ be as before. We say the pair$(\Pi,\Pi')$ is \textbf{in good position} if for any $\sigma\in \Sigma_{F;K}$, the $n'$ numbers  \begin{equation}\nonumber
 -\frac{\omega(\Pi)+\omega(\Pi')}{2}-b_{n'}(\sigma)>-\frac{\omega(\Pi)+\omega(\Pi')}{2}-b_{n'-1}(\sigma)>\cdots>-\frac{\omega(\Pi)+\omega(\Pi')}{2}-b_{1}(\sigma).\end{equation}
 lie in different gaps between $(a_{1}(\sigma)>a_{2}(\sigma)>\cdots>a_{n}(\sigma))$.
 
 It is equivalent to saying that $sp(i,\Pi';\Pi,\sigma)\neq 0$ for all $0\leq i\leq n'$ and $\sigma\in\Sigma_{F;K}$. In particular, if $n'=n-1$, we know $(\Pi,\Pi')$ is in good position if and only if $sp(i,\Pi';\Pi,\sigma)=1$ for all $i$ and $\sigma$.
\end{df}

\chapter{Unitary groups and base change}
\section{Unitary groups}
\text{}

In this section, let $\mathcal{L}$ be an arbitrary number field and $\mathcal{F}/\mathcal{L}$ be a quadratic extension of number fields. 

Let $U_{0}$ be the quasi-split unitary group over $\mathcal{L}$ of dimension $n$ with respect to the extension $\mathcal{F}/\mathcal{L}$. We want to know the local behavior of inner forms of $U_{0}$. More generally, we will answer the following question:

Let $G_{0}$ be a connected reductive group over $\mathcal{L}$. If we are given $G_{(v)}$, an inner form of $G_{0,v}$ over $\mathcal{L}_{v}$ for each place $v$ of $\mathcal{L}$, when does $G$, an inner form of $G_{0}$ over $\mathcal{L}$ such that $G_{v}\cong G_{(v)}$ for all $v$, exist?

The answer is given in section $2$ of \cite{clozelIHES}. We recall some results there. We also refer to section $1.2$ \cite{harrislabesse} for further details in the unitary group case. 

\bigskip

The isomorphism classes of inner forms are classified by Galois cohomology. Let $v$ be a place of $\mathcal{L}$. Let $L=\mathcal{L}$ or $\mathcal{L}_{v}$. There exists a bijection between the set of isomorphism classes of inner forms of $G_{0,L}$ and $H^{1}(L,G_{0}^{ad})$. Therefore, the global inner form exists if and only if the element in $\bigoplus\limits_{v}H^{1}(\mathcal{L}_{v},G_{0}^{ad})$ corresponding to the local datum is in the image of $H^{1}(\mathcal{L},G_{0}^{ad}) \rightarrow \bigoplus\limits_{v}H^{1}(\mathcal{L}_{v},G_{0}^{ad})$. 

We remark that if $L$ is local then the quasi-split class corresponds to the trivial element of $H^{1}(L,G_{0}^{ad})$.

\bigskip

If we can calculate this Galois cohomology, then everything is done. Otherwise Kottwitz has given an alternate choice as follows.

For $H$ a connected reductive group over $L$, we define $A(H)=A(H/L):=$ the dual of $\pi_{0}(Z(\hat{H})^{G_{L}}$ where $\hat{H}$ is the neutral component of the dual group of $H$.

Let $A=A(G_{0}^{ad}/\mathcal{L})$ and $A_{v}=A(G_{0}^{ad}/\mathcal{L}_{v})$. 

\begin{prop}
 There exists a natural map $H^{1}(\mathcal{L}_{v},G_{0}^{ad})\rightarrow A_{v}$. Moreover, it is an isomorphism when $v$ is finite.
\end{prop}

The above proposition gives a morphism $\bigoplus\limits_{v}H^{1}(\mathcal{L}_{v},G_{0}^{ad})\rightarrow \bigoplus\limits_{v}A_{v}\rightarrow A$ where the latter is given by restriction. 

\begin{thm}
 The following sequence is exact: 
 \begin{equation}\nonumber
 H^{1}(\mathcal{L},G_{0}^{ad})\rightarrow \bigoplus\limits_{v}H^{1}(\mathcal{L}_{v},G_{0}^{ad})\rightarrow A.\end{equation}  
 In other words, the image of the first map equals the kernel of the second map.
\end{thm}
By this theorem, our question turns to determine the kernel of the second map.
\bigskip

Let us now focus on unitary groups, namely when $G_{0}=U_{0}$. Clozel has calculated $A_{v}$ in the case when $\mathcal{L}$ is totally real and $\mathcal{F}$ is a quadratic imaginary extension over $\mathcal{L}$. We call it the \textbf{CM case}. This is enough for our purpose. Let us list some facts from \cite{clozelIHES}:
\begin{itemize}
 \item If $n$ is odd, then $A=0$. In other words, any local datum $(U_{(v)})_{v}$ which is quasi-split at almost all places come from a global unitary group.
 
 \item If $n$ is even, then
 \begin{enumerate}
  \item $A\cong \Z/2$.
  \item $A_{v}\cong \Z/2$ if $v$ is finite and inert. The non quasi-split unitary group corresponds to the non trivial element of $\Z/2$. The map $A_{v}\rightarrow A$ is identity if we identify both groups with $\Z/2$.   
  \item $A_{v}\cong \Z/n$ if $v$ is finite and split. The element corresponding to the unitary group of a division algebra generates $A_{v}$. The map $A_{v}\rightarrow A$ is the mod $2$ map from $\Z/n$ to $\Z/2$. 
  \item The real unitary group $U(p,q)$ has image $(p-q)/2 \mod 2$ in $A$.
  \end{enumerate}
 \end{itemize}

\begin{rem}
 \begin{enumerate}
  \item The idea of the proof for the last point is to consider the surjective map $H^{1}(\R,T)\rightarrow H^{1}(\R,G_{0})$ where $T\subset G_{0}$ is the maximal elliptic torus over $\R$.
 
 \end{enumerate}

\end{rem}

 The above calculation leads to the following theorem:
 
 \begin{thm}
 Let $F=F^{+}K$
 Let $I$ be as before.
 Let $q$ be a finite place of $\Q$ inert in $F^{+}$ and split in $F$. There exists a Hermitian space $V_{I}$ of dimension $n$ over $F$ with respect to $F/F^{+}$ such that the unitary group $U=U(V_{I})$ over $F^{+}$ associated to $V$ satisfies:
 \begin{itemize}
  \item At each $\sigma\in \Sigma$, $U$ is of sign $(n-I(\sigma),I(\sigma))$.
  \item For $v\neq q$, a finite place of $F^{+}$, $U_{v}$ is unramified;
  \item If $n$ is even and $\sum\limits_{\sigma\in \Sigma}\cfrac{n-2I(\sigma)}{2}\nequiv 0 \mod 2$, then $U_{q}$ is a division algebra. Otherwise $U_{q}$ is also unramified.

 \end{itemize} 
 
\end{thm}
 
 We denote by $U_{I}$ the restriction of $U$ from $F^{+}$ to $\Q$ and $GU_{I}$ the rational similitude group associated to $V_{I}$, namely, for any $\Q$-algebra $R$, 
  \begin{equation}\label{definition of nu(g)}  
  GU_{I}(R)=\{g\in GL(V_{I}\otimes_{\Q}R)|(gv,gw)=\nu(g)(v,w), 
  \nu(g)\in R^{*}\}.\end{equation}

\section{General base change}
\text{}

Let $G$ and $G'$ be two connected quasi-split reductive algebraic groups over $\Q$. Let $\widehat{G}$ be the complex dual group of $G$. The Galois group $G_{\Q}:=Gal(\overline{\Q}/\Q)$ acts on $\widehat{G}$. We define the \textbf{$L$-group} of $G$ by $\text{}^{L}G:=\widehat{G} \rtimes G_{\Q}$ and we define $\text{}^{L}G'$ similarly. A group homomorphism between two $L$-groups $\text{}^{L}G\rightarrow \text{}^{L}G'$ is called an \textbf{$L$-morphism} if it is continuous, its restriction to $\widehat{G}$ is analytic and it is compatible with the projections of $\text{}^{L}G$ and $\text{}^{L}G'$ to $G_{\Q}$. If there exists an $L$-morphism $\text{}^{L}G\rightarrow \text{}^{L}G'$, the \textbf{Langlands' principal of functoriality} predicts a correspondence from automorphic representations of $G(\AQ)$ to automorphic representations of $G'(\AQ)$ (c.f. section $26$ of \cite{arthurtraceformula}). More precisely, we wish to associate an $L$-packet of automorphic representations of  $G(\AQ)$ to that of $G'(\AQ)$.

\bigskip

Locally, we can specify this correspondence for unramified representations. Let $v$ be a finite place of $\Q$  such that $G$ is unramified at $v$. We fix $\Gamma_{v}$ a maximal compact hyperspecial subgroup of $G_{v}:=G(\Q_{v})$.  By definition, for $\pi_{v}$ an admissible representation of $G_{v}$, we say $\pi_{v}$ is \textbf{unramified} (with respect to $\Gamma_{v}$) if it is irreducible and $dim \pi_{v}^{\Gamma_{v}} > 0$. One can show that $\pi_{v}^{\Gamma_{v}}$ is actually one dimensional since $\pi_{v}$ is irreducible.\\

Denote $H_{v}:=\mathcal{H}(G_{v},\Gamma_{v})$ the Hecke algebra consisting of compactly supported continuous functions from $G_{v}$ to $\C$ which are $\Gamma_{v}$ invariants on both sides. We know $H_{v}$ acts on $\pi_{v}$ and preserves $\pi_{v}^{\Gamma_{v}}$ (c.f. \cite{minguez}). Since $\pi_{v}^{\Gamma_{v}}$ is one-dimensional, every element in $H_{v}$ acts as a multiplication by a scalar on it. In other words, $\pi_{v}$ thus determines a character of $H_{v}$. This gives a map from the set of unramified representations of $G_{v}$ to the set of characters of $H_{v}$ which is in fact a bijection (c.f. section $2.6$ of \cite{minguez}).

\bigskip

We can moreover describe the structure of $H_{v}$ in a simpler way. Let $T_{v}$ be a maximal torus of $G_{v}$ contained in a Borel subgroup of $G_{v}$. We denote by $X_{*}(T_{v})$ the set of cocharacters of $T_{v}$ defined over $\Q_{v}$. The Satake transform identifies the Hecke algebra $H_{v}$ with the polynomial ring $\C[X_{*}(T_{v})]^{W_{v}}$ where $W_{v}$ is the Weyl group of $G_{v}$ (c.f. section $1.2.4$ of \cite{harristakagi}).\\

Let $G'$ be a connected quasi-split reductive group which is also unramified at $v$. We can define $\Gamma'_{v}$, $H'_{v}:=\mathcal{H}(G'_{v},\Gamma'_{v})$ and $T'_{v}$ similarly. An $L$-morphism $\text{}^{L}G\rightarrow \text{}^{L}G'$ induces a morphism $\widehat{T_{v}}\rightarrow \widehat{T'_{v}}$ and hence a map $T'_{v}\rightarrow T_{v}$. The latter gives a morphism from $\C[X_{*}(T'_{v})]^{W'_{v}}$ to $\C[X_{*}(T_{v})]^{W_{v}}$ and thus a morphism from $H'_{v}$ to $H_{v}$. Its dual hence gives a map from the set of unramified representations of $G_{v}$ to that of $G'_{v}$. This is the \textbf{local Langlands's principal of functoriality} for unramified representations.

\bigskip

In this article, we are interested in a particular case of the Langlands' functoriality: the cyclic base change. Let $K/\Q$ be a cyclic extension, for example $K$ is a quadratic imaginary field. Let $G$ be a connected quasi-split reductive group over $\Q$. Let $G'=Res_{K/\Q}G_{K}$. We know $\widehat{G'}$ equals to $\widehat{G}^{[K:\Q]}$. The diagonal embedding is then a natural $L$-morphism $\text{}^{L}G \rightarrow \text{}^{L}G'$. The corresponding functoriality is called the base change.\\

More precisely, let $v$ be a place of $\Q$ and $w$ a place of $K$ over $v$. The local Langlands's principal of functoriality gives a map from the set of unramified representations of $G(\Q_{v})$ to that of $G(K_{w})$. We call this map the base change with respect to $K_{w}/\Q_{v}$.\\

Let $\pi$ be an admissible irreducible representation of $G(\AQ)$. We say $\Pi$, a representation of $G(\AK)$, is a \textbf{(weak) base change} of $\pi$ if for almost all $v$, a finite place of $\Q$, such that $\pi$ is unramified at $v$ and all $w$, a place of $K$ over $v$, $\Pi_{w}$ is the base change of $\pi_{v}$. In this case, we say $\Pi$ \textbf{descends to $\pi$} by base change. 

\bigskip

\begin{rem}
The group $Aut(K)$ acts on $G(\AK)$. This induces an action of $Aut(K)$ on automorphic representations of $G(\AK)$. For $\sigma\in Aut(K)$ and $\Pi$ an automorphic representation of $G(\AK)$, we write $\Pi^{\sigma}$ to be the image of $\Pi$ under the action of $\sigma$. It is easy to see that if $\Pi$ is a base change of $\pi$, then $\Pi^{\sigma}$ is one as well. In particular, we have $\Pi^{\sigma}_{w}\cong \Pi_{w}$ for almost every finite place $w$ of $K$. So if we have the strong multiplicity one theorem for $G(\AK)$, we can conclude that every representation in the image of base change is $Aut(K)$-stable.

\end{rem}

\bigskip

\section{Base change for unitary groups and similitude unitary groups}

Recall that $U_{I}(\AK)\cong GL_{n}(\AF)$. The following result on base change comes from Theorem $1.7$ of \cite{endoscopicfour}. We also refer to Corollary $2.5.9$ of \cite{Mokendoscopic} for the quasi-split case.

\bigskip

\begin{prop}\label{base change for unitary group}
\textbf{Base change for unitary group}

Let $\Pi$ be a cuspidal conjugate self-dual and cohomological representation of $GL_{n}(\AF)$. If $n$ is odd then $\Pi$ descends to a cohomological representation of $U_{I}(\AQ)$ unconditionally. If $n$ is even then it descends if $\Pi_{q}$ descends locally.
\end{prop}

We have an exact sequence $1\rightarrow U_{I}\rightarrow GU_{I}\rightarrow \mathbb{G}_{m}\rightarrow 1$ which is split over $K$. Indeed, by Galois descent, it is enough to define $\theta_{I}$, a Galois automorphism on $U_{I,F}\times \mathbb{G}_{m,K}$ such that the subgroup of $U_{I,K}\times \mathbb{G}_{m,K}$ fixed by $\theta_{I}$ is isomorphic to $GU_{I}$. We now define $\theta_{I}$ as follows:\\

For $R$ a $\Q$-algebra, note that $(U_{I,K}\times \mathbb{G}_{m,K})(R)\cong GL(V_{I}\otimes_{\Q}R)\times (K\otimes_{\Q}R)$. We define
\begin{equation}\nonumber
\theta_{I}: GL(V_{I}\otimes_{\Q}R)\times (K\otimes_{\Q}R)\rightarrow GL(V_{I}\otimes_{\Q}R)\times (K\otimes_{\Q}R)
\end{equation}
by sending $(g,z)$ to $((g^{*})^{-1}\bar{z},\bar{z})$ where $g^{*}$ is the adjoint of $g$ with respect to the Hermitian form. It is easy to verify that $\theta_{I}$ satisfies the condition mentioned above.\\

We then have that $GU_{I,K}\cong U_{I,K}\times \mathbb{G}_{m,K}$. In particular, $GU(\AK)\cong GL_{n}(\AF)\times \AK^{\times}$. For $\Pi$ a cuspidal representation of $GL_{n}(\AF)$ and $\xi$ a Hecke character of $K$, $\Pi\otimes \xi$ defines a cuspidal representation of $GU(\AK)$. Conversely, by the tensor product theorem, every irreducible admissible automorphic representation of $GU(\AK)$ is of the form $\Pi\otimes \xi$. The following Lemma is shown in \cite{harristaylor} $VI.2.10$ and \cite{clozelramanujan} Lemma $2.2$.

\bigskip

\begin{lem}\label{existence of xi}
If $\Pi$ is algebraic and conjugate self-dual, then there exists $\xi$, an algebraic Hecke character of $\AK$ such that $\Pi\times\xi$ is $\theta_{I}$-stable.
\end{lem}

\begin{dem}
It is easy to verify that $\Pi\times\xi$ is $\theta_{I}$-stable if and only if $\cfrac{\xi(\overline{z})}{\xi(z)}=\xi_{\Pi}(z)$ for any $z\in \AK^{\times}$ where $\xi_{\Pi}$ is the central character of $\Pi$.\\

We define $U$ the torus over $\Q$ such that $U(\Q)=ker\{Norm: K^{\times}\rightarrow \Q^{\times}\}$. We have $U(\AQ)=\{\cfrac{\overline{z}}{z}\mid z\in\AK^{\times}\}$ by Hilbert 90.\\

Again by Hilbert 90, we have an exact sequence:
\begin{equation}\nonumber
1\rightarrow \Q^{\times}\backslash \AQ^{\times} \rightarrow K^{\times}\backslash \AK^{\times} \rightarrow U(\Q)\backslash U(\AQ) \rightarrow 1
\end{equation}
where the last map sends $z$ to $\cfrac{\overline{z}}{z}$.\\

Therefore, such $\xi$ exists if and only if $\xi_{\Pi}$ is trivial on $\AQ^{\times}$.\\

Since $\Pi$ is conjugate self-dual, we know $\xi_{\Pi}$ is trivial on $Norm_{\AK/\AQ}(\AK^{\times})$. By class field theory, $\Q^{\times}N_{\AK/\AQ}(\AK^{\times})$ has index $2$ in $\AQ^{\times}$. It remains to show that $\xi_{\Pi}$ is trivial at any element $t$ in $\AQ^{\times}-\Q^{\times}N_{\AK/\AQ}(\AK^{\times})$.\\

We consider $t\in \AQ^{\times}$ such that $t=1$ at all finite places and $t=-1$ at the infinity place. It is not in $\Q^{\times}N_{\AK/\AQ}(\AK^{\times})$.\\

Since $\Pi$ is algebraic, we know $\xi_{\Pi}$ has infinity type $z^{a}\overline{z}^{-a}$ with $a\in\Z$. In particular, we have $\xi_{\Pi}(t)=1$ as expected.
\end{dem}

The following proposition follows from Proposition \ref{base change for unitary group}. The idea is the same with Theorem $VI.2.9$ in \cite{harristaylor}.

\begin{prop}\label{base change for similitude unitary group}
Let $\Pi$ be an algebraic automorphic representation of $GL_{n}(\AF)\cong U_{I}(\AK)$ which descends to $U_{I}(\AF)$. 
If $\xi$ is an algebraic Hecke character of $K$ as in Lemma \ref{existence of xi} then $\Pi\otimes \xi$ descends to an automorphic representation of $GU_{I}(\AF)$.
Moreover, if $\Pi$ is cohomological then its descending is also cohomological.\\

In particular, let $\Pi$ be a cuspidal conjugate self-dual and cohomological representation of $GL_{n}(\AF)$. We assume moreover that $\Pi_{q}$ descends locally if $n$ is even. Then there always exists $\xi$ such that $\Pi\otimes \xi$ descends to an irreducible cohomological automorphic representation of $GU_{I}(\AF)$.

\end{prop}

\chapter{Rational structures and Whittaker periods}

\text{}

In this chapter, we will recall some results on the rationality of certain algebraic automorphic representations and also the rationality of the associated Whittaker models and cohomology spaces.\\

In particular, we will define the Whittaker period and present a way to calculate the Whittaker period in certain cases.

\section{Rational structures on certain automorphic representations} 
\text{}

Let $F$ is an arbitrary number field and $n$ be a positive integer.

Let $\Pi$ be an automorphic representation of $GL_{n}(\AF)$. \\

We denote by $V$ the representation space for $\Pi_{f}$. For $\sigma\in Aut(\C)$, we define another $GL_{n}(\AFf)$-representation $\Pi_{f}^{\sigma}$ to be $V\otimes_{\C,\sigma}\C$. Let $\Q(\Pi)$ be the subfield of $\C$ fixed by $\{\sigma\in Aut(\C) \mid  \Pi_{f}^{\sigma} \cong \Pi_{f}\}$. We call it the \textbf{rationality field} of $\Pi$.

\bigskip

For $E$ a number field, $G$ a group and $V$ a $G$-representation over $\C$, we say $V$ has an \textbf{$E$-rational structure} if there exists an $E$-vector space $V_{E}$ endowed with an action of $G$ such that $V=V_{E}\otimes_{E}\C$ as representation of $G$. We call $V_{E}$ an $E$-rational structure of $V$. \\

We denote by $\mathcal{A}lg(n)$ the set of algebraic automorphic representations of $GL_{n}(\AF)$ which are isobaric sums of cuspidal representations as in section $1$ of \cite{clozelaa}.

\bigskip

\begin{thm}{(Th\'eor\`eme $3.13$ in \cite{clozelaa})}

Let $\Pi$ be a regular representation in $\mathcal{A}lg(n)$. We have that:
\begin{enumerate} 
\item $\Q(\Pi)$ is a number field. 
\item $\Pi_{f}$ has a $\Q(\Pi)$-rational structure unique up to homotheties. 
\item For all $\sigma\in Aut(\C)$, $\Pi_{f}^{\sigma}$ is the finite part of a regular representation in $\mathcal{A}lg(n)$. It is unique up to isomorphism by the strong multiplicity one theorem. We denote it by $\Pi^{\sigma}$.
\end{enumerate}

\end{thm}

\bigskip

\begin{rem}
Let $n=n_{1}+n_{2}+\cdots+n_{k}$ be a partitian of positive integers and $\Pi_{i}$ be regular representations in $\mathcal{A}lg(n_{i})$ for $1\leq i\leq k$ respectively.\\

The above theorem implies that, for all $1\leq i\leq k$, the rational field $\Q(\Pi_{i})$ is a number field.\\

Let $\Pi=(\Pi_{1}||\cdot||_{\AK}^{\frac{1-n_{1}}{2}} \boxplus \Pi_{2}||\cdot||_{\AK}^{\frac{1-n_{2}}{2}}\boxplus\cdots \boxplus \Pi_{k}||\cdot||_{\AK}^{\frac{1-n_{k}}{2}})||\cdot||_{\AK}^{\frac{n-1}{2}}$ be the normalized isobaric sum of $\Pi_{i}$. It is still algebraic.\\

We can see from definition that $\Q(\Pi)$ is the compositum of $\Q(\Pi_{i})$ with $1\leq i\leq k$. Moreover, if $\Pi$ is regular, we know from the above theorem that $\Pi$ has a $\Q(\Pi)$-rational structure.
\end{rem}

\section{Rational structures on the Whittaker model}
\text{}
Let $\Pi$ be a regular representation in $\mathcal{A}lg(n)$ and then its rationality field $\Q(\Pi)$ is a number field.\\

We fix a nontrivial additive character $\phi$ of $\AF$. Since $\Pi$ is an isobaric sum of cuspidal representations, it is generic. Let $W(\Pi_{f})$ be the Whittaker model associated to $\Pi_{f}$ (with respect to $\phi_{f}$). It consists of certain functions on $GL_{n}(\AF\text{}_{,f})$ and is isomorphic to $\Pi_{f}$ as $GL_{n}(\AF\text{}_{,f})$-modules.\\

Similarly, we denote the Whittaker model of $\Pi$ (with respect to) $\phi$ by $W(\Pi)$.\\

\begin{df}{\textbf{Cyclotomic character}}

There exists a unique homomorphism $\xi: Aut(\C) \rightarrow \widehat{\Z}^{\times}$ such that for any $\sigma\in Aut(\C)$ and any root of unity $\zeta$, $\sigma(\zeta)=\zeta^{\xi(\sigma)}$, called the cyclotomic character.
\end{df}

For $\sigma\in Aut(\C)$, we define $t_{\sigma}\in (\widehat{\Z} \otimes_{\Z} \mathcal{O}_{F})^{\times}  =\widehat{\mathcal{O}_{F}}^{\times}$ to be the image of $\xi(\sigma)$ by the embedding $(\widehat{\Z})^{\times} \hookrightarrow (\widehat{\Z} \otimes_{\Z} \mathcal{O}_{F})^{\times} $. We define $t_{\sigma,n}$ to be the diagonal matrix $diag(t_{\sigma}^{-n+1},t_{\sigma}^{-n+2},\cdots,t_{\sigma}^{-1},1)\in GL_{n}(\AF\text{}_{,f})$ as in section $3.2$ of \cite{raghuramshahidi}. \\

For $w\in W(\Pi_{f})$, we define a function $w^{\sigma}$ on $GL_{n}(\AF\text{}_{,f})$ by sending $g\in GL_{n}(\AF\text{}_{,f})$ to $\sigma(w(t_{\sigma,n}g))$. For classical cusp forms, this action is just the $Aut(\C)$-action on Fourier coefficients.

\begin{prop}(Lemma $3.2$ of \cite{raghuramshahidi} or Proposition $2.7$ of \cite{harrismotivic})

The map $w\mapsto w^{\sigma}$ gives a $\sigma$-linear $GL_{n}(\AF\text{}_{,f})$-equivariant isomorphism from $W(\Pi_{f})$ to $W(\Pi^{\sigma}_{f})$.

For any extension $E$ of $\Q(\Pi_{f})$, we can define an $E$-rational structure on $W(\Pi_{f})$ by taking the $Aut(\C/E)$-invariants.

Moreover, the $E$-rational structure is unique up to homotheties.

\end{prop}

\begin{dem}
The first part is well-known (see the references in \cite{raghuramshahidi}). \\

For the second part, the original proof in \cite{raghuramshahidi} works for cuspidal representations. The key point is to find a nonzero global invariant vector. It is equivalent to finding a nonzero local invariant vector for every finite place. Then Theorem $5.1(ii)$ of \cite{JPSS} is involved as in \cite{harrismotivic}.\\

The last part follows from the one-dimensional property of the invariant vector which is the second part of Theorem $5.1(ii)$ of \cite{JPSS}.
\end{dem}

\section{Rational structures on cohomology spaces and comparison of rational structures}\label{cohomology}
\text{}
Let $\Pi$ be a regular representation in $\mathcal{A}lg(n)$. The Lie algebra cohomology of $\Pi$ has a rational structure. It is described in section $3.3$ of \cite{raghuramshahidi}. We give a brief summary here.\\

Let $Z$ be the center of $GL_{n}$. Let $\lieG_{\infty}$ be the Lie algebra of $GL_{n}(\R\otimes_{\Q}F)$. Let $S_{real}$ be the set of real places of $F$, $S_{complex}$ be the set of complex places of $F$ and $S_{\infty}=S_{real}\cup S_{complex}$ be the set of infinite places of $F$.\\

For $v\in S_{real}$, we define $K_{v}:=Z(\R)O_{n}(\R)\subset GL_{n}(F_{v})$. For $v\in S_{complex}$, we define $K_{v}:=Z(\C)U_{n}(\C)\subset GL_{n}(F_{v})$. We denote by $K_{\infty}$ the product of $K_{v}$ with $v\in S_{\infty}$, and by $K_{\infty}^{0}$ the topological connected component of $K_{\infty}$.

\bigskip

We fix $T$ the maximal torus of $GL_{n}$ consisting of diagonal matrices and $B$ the Borel subgroup of $G$ consisting of upper triangular matrices. For $\mu$ a dominant weight of $T(\R\otimes_{\Q}F)$ with respect to $B(\R\otimes_{\Q}F)$, we can define $W_{\mu}$ an irreducible representation of $GL_{n}(\R\otimes_{\Q}F)$ with highest weight $\mu$.\\

From the proof of Th\'eor\`eme $3.13$ \cite{clozelaa}, we know that there exists a dominant algebraic weight $\mu$, such that $H^{*}(\lieG_{\infty},K_{\infty}^{0};\Pi_{\infty}\otimes W_{\mu})\neq 0$. \\

Let $b$ be the smallest degree such that $H^{b}(\lieG_{\infty},K_{\infty}^{0};\Pi_{\infty}\otimes W_{\mu})\neq 0$. We have an explicit formula for $b$ in \cite{raghuramshahidi}. More precisely, we set $r_{1}$ and $r_{2}$ the numbers of real and complex embeddings of $F$ respectively. We have $b=r_{1}\left[\frac{n^{2}}{4}\right]+r_{2}\frac{n(n-1)}{2}$.

\bigskip

We can decompose this cohomology group via the action of $K_{\infty}/K_{\infty}^{0}$. There exists a character $\epsilon$ of $K_{\infty}/K_{\infty}^{0}$ described explicitly in \cite{raghuramshahidi} such that:
\begin{enumerate}
\item The isotypic component $H^{b}(\lieG_{\infty},K_{\infty}^{0};\Pi_{\infty}\otimes W_{\mu})(\epsilon)$ is one dimensional.

\item For fixed $w_{\infty}$, a generator of $H^{b}(\lieG_{\infty},K_{\infty}^{0};\Pi_{\infty}\otimes W_{\mu})(\epsilon)$, we have a $GL_{n}(\AF\text{}_{,f})$-equivariant isomorphisms:
\begin{eqnarray}\label{rational isomorphism}
W(\Pi_{f}) &\xrightarrow{\sim}& W(\Pi_{f})\otimes H^{b}(\lieG_{\infty},K_{\infty}^{0};\Pi_{\infty}\otimes W_{\mu})(\epsilon)\nonumber\\
&\xrightarrow{\sim}& H^{b}(\lieG_{\infty},K_{\infty}^{0};W(\Pi)\otimes W_{\mu})(\epsilon)\nonumber\\
&\xrightarrow{\sim}&  H^{b}(\lieG_{\infty},K_{\infty}^{0};\Pi \otimes W_{\mu})(\epsilon)
\end{eqnarray}
where the first map sends $w_{f}$ to $w_{f}\otimes w_{\infty}$ and the last map is given by the isomorphism $W(\Pi)\xrightarrow{\sim} \Pi$.  

\item The cohomology space $H^{b}(\lieG_{\infty},K_{\infty}^{0};\Pi \otimes W_{\mu})(\epsilon)$ is related to the cuspidal cohomology if $\Pi$ is cuspidal and to the Eisenstein cohomology if $\Pi$ is not cuspidal. In both cases, it is endowed with a $\Q(\Pi)$-rational structure (see \cite{raghuramshahidi} for cuspidal case and \cite{harrismotivic} for non cuspidal case).

\end{enumerate}

We denote by $\Theta_{\Pi_{f},\epsilon,w_{\infty}}$ the $GL_{n}(\AFf)$-isomorphism given in (\ref{rational isomorphism}) 
\begin{equation}W(\Pi_{f})\xrightarrow{\sim}  H^{b}(\lieG_{\infty},K_{\infty}^{0};\Pi \otimes W_{\mu})(\epsilon)\nonumber.
\end{equation}
 Both sides have a $\Q(\Pi)$-rational structure. In particular, the preimage of the rational structure on the right hand side gives a rational structure on $W(\Pi_{f})$. But the rational structure on $W(\Pi_{f})$ is unique up to homotheties. Therefore, there exists a complex number $p(\Pi_{f},\epsilon,w_{\infty})$ such that the new map $\Theta^{0}_{\Pi_{f},\epsilon,w_{\infty}}=p(\Pi_{f},\epsilon,w_{\infty})^{-1}\Theta_{\Pi_{f},\epsilon,w_{\infty}}$ preserves the rational structure on both sides. It is easy to see that this number $p(\Pi_{f},\epsilon,w_{\infty})$ is unique up to multiplication by elements in $\Q(\Pi)^{\times}$.

\bigskip

Finally, we observe that the $Aut(\C)$-action preserves rational structures on both the Whittaker models and cohomology spaces. We can adjust the numbers $p(\Pi_{f}^{\sigma},\epsilon^{\sigma},w_{\infty}^{\sigma})$ for all $\sigma\in Aut(\C)$ by elements in $\Q(\Pi)^{\times}$ such that the following diagram commutes:

\[\xymatrixcolsep{12pc}
\xymatrix{
W(\Pi_{f}) \ar[d]^{\sigma} \ar[r]^{p(\Pi_{f},\epsilon,w_{\infty})^{-1}\Theta_{\Pi_{f},\epsilon,w_{\infty}}} &H^{b}(\lieG_{\infty},K_{\infty}^{0};\Pi \otimes W_{\mu})(\epsilon)\ar[d]^{\sigma}\\
W(\Pi_{f}^{\sigma}) \ar[r]^{p(\Pi_{f}^{\sigma},\epsilon^{\sigma},w_{\infty}^{\sigma})^{-1}\Theta_{\Pi_{f}^{\sigma},\epsilon^{\sigma},w_{\infty}^{\sigma}}}          &H^{b}(\lieG_{\infty},K_{\infty}^{0};\Pi^{\sigma} \otimes W_{\mu}^{\sigma})(\epsilon^{\sigma})}
\]
The proof is the same as the cuspidal case in \cite{raghuramshahidi}.\\

In the following, we fix $\epsilon$, $w_{\infty}$ and we define the \textbf{Whittaker period} $p(\Pi):=p(\Pi_{f},\epsilon,w_{\infty})$. For any $\sigma\in Aut(\C)$, we define $p(\Pi^{\sigma}):=p(\Pi_{f}^{\sigma},\epsilon^{\sigma},w_{\infty}^{\sigma})$. It is easy to see that $p(\Pi^{\sigma})=p(\Pi)$ for $\sigma\in Aut(\C/\Q(\Pi))$. 

Moreover, the elements $(p(\Pi^{\sigma}))_{\sigma\in Aut(\C)}$ are well defined up to $\Q(\Pi)^{\times}$ in the following sense: if $(p'(\Pi^{\sigma}))_{\sigma\in Aut(\C)}$ is another family of complex numbers such that $p'(\Pi^{\sigma})^{-1}\Theta_{\Pi_{f}^{\sigma},\epsilon^{\sigma},w_{\infty}^{\sigma}}$ preserves the rational structure and the above diagram commutes, then there exists $t\in \Q(\Pi)^{\times}$ such that $p'(\Pi^{\sigma})=\sigma(t) p(\Pi^{\sigma})$ for any $\sigma\in Aut(\C)$. This also follows from the one dimensional property of the invariant vector. The argument is the same as the last part of the proof of Definition/Proposition $3.3$ in \cite{raghuramshahidi}.

\section{Shahidi's calculation on Whittaker periods}

Let us assume $F$ is a CM field in the following sections of this chapter. In this case, $K_{\infty}$ itself is connected and hence we may omit the index $\epsilon$.\\

let $l$ be a positive integer. Let $n_{1}$, $n_{2}$, $\cdots$ ,$n_{l}$ be positive integers such that $n=n_{1}+n_{2}+\cdots+n_{l}$. \\

Let $\Pi_{1}$, $\Pi_{2}$, $\cdots$, $\Pi_{l}$ be regular cohomological conjugate self-dual automorphic representations of $GL_{n_{1}}$, $GL_{n_{2}}$, $\cdots$, $GL_{n_{l}}$ respectively. We assume that they are Langlands sum of cuspidal representations.\\

We write $P\leqslant GL_{n}$ for the maximal parabolic group of type $(n_{1},n_{2},\cdots,n_{l})$ and $B\leqslant P$ be the corresponding Borel subgroup.\\

Let $\Pi$ be the Langlands sum of $\Pi_{1},\Pi_{2},\cdots, \Pi_{l}$. 
\begin{prop}\label{Shahidi Whittaker pre}
There exists a non zero complex number $\Omega_{(n_{1},n_{2},\cdots,n_{l})}(\Pi_{\infty})$ depending on $\Pi_{\infty}$ and the parabolic type of $P$ which is unique up to elements of $E(\Pi)^{\times}$ such that:
\begin{equation}\label{Whittaker decomposition pre}
p(\Pi)\sim_{E(\Pi_{1})\cdots E(\Pi_{l});K} \Omega_{(n_{1},n_{2},\cdots,n_{l})}(\Pi_{\infty})\prod\limits_{1\leq i\leq l}p(\Pi_{i})\prod\limits_{1\leq i<j\leq l}L(1,\Pi_{i}\times \Pi_{j}^{\vee})
\end{equation}
\end{prop} 

The constructions and ideas come from \cite{mahnkopf05} and \cite{harrismotivic}. We give a sketch of the proof here and will include the details in a forthcoming paper.

\paragraph{Sketch of the proof:}
For simplicity, we may assume that $l=2$. 

For $i=1$ or $2$, we denote by $V_{i}$ the representation space for $\Pi_{i}$ consisting of cusp forms.

We denote by $W_{i,f}$ the Whittaker model for $\Pi_{i,f}$. We write $\mu_{i}$ for the cohomological type of $\Pi_{i}$. We fix $\omega_{i,\infty}$, a generator of $H^{b_{n_{i}}}(\lieG_{i,\infty},K_{i,\infty};\Pi_{i,\infty}\otimes W_{\mu_{i}})$ and we write $\Theta_{i}$ for $\Theta_{\Pi_{i,f},\omega_{i}}$.

We write $H^{b_{i}}(\Pi_{i}\otimes W_{\mu_{i}}):=H^{b_{i}}(\lieG_{i,\infty},K_{i,\infty};\Pi_{i} \otimes W_{\mu_{i}})$ for simplicity.

We use similar notation for $\Pi=\Pi_{1}\boxplus \Pi_{2}$.

We claim that there exists a commutative diagram as follows (c.f. (1.3) of \cite{mahnkopf05}):
\[\xymatrixcolsep{10pc}
\xymatrix{
W_{1,f}\otimes W_{2,f} \ar[d]^{\mathcal{F}^{loc}} \ar[r]^{\Theta_{1}\otimes \Theta_{2}} & H^{b_{1}}( \Pi_{1}\otimes W_{\mu_{1}})\otimes H^{b_{2}}(\Pi_{2} \otimes W_{\mu_{2}}) \cong V_{1,f}\otimes V_{2,f} \ar[d]^{Eis}\\
W(\Pi_{f}) \ar[r]^{\Theta}          & H^{b}(\Pi\otimes W_{\mu})\supset V_{f}
}
\]
We now introduce the maps which appear in the diagram:
\begin{itemize}
\item The map $\mathcal{F}^{loc}$ is an explicit map defined locally in \cite{mahnkopf05}.

\item The map $Eis$ is rational and defined by the theory of Eisenstein series which sends $V_{1,f}\otimes V_{2,f}$ to $V_{f}$. (c.f. Section $1.1$ in \cite{mahnkopf05}).

\item The isomorphism $V_{i,f}\cong H^{b_{i}}( \Pi_{i}\otimes W_{\mu_{i}})$ is rational. The composition of this isomorphism and $\Theta_{i}$ is just the isomorphism between the cuspidal forms and the corresponding Whittaker functions.

\item The theory of Eisenstein cohomology (c.f. \cite{harder90}) gives a rational embedding of $V_{f}$ in $H^{b}(\Pi\otimes W_{\mu})$.

More precisely, we write $S_{n}$ for $GL_{n}(F)\backslash GL_{n}(\AF)/ K_{\infty}$. We denote by $\overline{S}_{n}$ the Borel-Serre compactification of $S_{n}$. We write $\partial \overline{S}_{n}$ for its boundary and $\partial_{B}\overline{S}_{n}$ for the face corresponding to the Borel subgroup $B$.

We know $H^{b}(\Pi\otimes W_{\mu})$ embeds rationally in $H^{b}(\overline{S}_{n},\mathcal{E}_{\mu})$ (c.f. Section $3$ of \cite{harrismotivic} and Section $1.2$ of \cite{harder90}) where $\mathcal{E}_{\mu}$ is a sheaf on $S_{n}$ defined by $\mu$. We restrict the latter to the face $\partial_{B}\overline{S}_{n}$ and get a rational map $H^{b}(\Pi\otimes W_{\mu})\rightarrow H^{b}(\partial_{B}\overline{S}_{n},\mathcal{E}_{\mu})$ which admits a rational section (c.f. Proposition $5.2$ of \cite{harrismotivic}).

We may decompose $H^{b}(\partial_{B}\overline{S}_{n},\mathcal{E}_{\mu})$ as in Theorem $4.2$ of \cite{mahnkopf05} and we see that $V_{f}=Ind_{P}(V_{1,f}\otimes V_{2,f})$ is a rational direct summand of $H^{b}(\Pi\otimes W_{\mu})\rightarrow H^{b}(\partial_{B}\overline{S}_{n},\mathcal{E}_{\mu})$ as $GL_{n}(\AFf)$-module. Here we should take $w$ to be the longest element in the Weyl group and $s$ to be trivial in the \textit{loc.cit}.
\end{itemize}

We take $g_{1}\otimes g_{2}\in W_{1,f}\otimes W_{2,f}$ to be a rational element. We write $f_{1}\otimes f_{2}$ for the image of $g_{1}\otimes g_{2}$ under the map $\Theta_{1}\otimes \Theta_{2}$. We denote by $F:=Eis(f_{1}\otimes f_{2})$. We write $W=W(F)$ to be the corresponding Whittaker function (with respect to a fixed additive Hecke character).

From the diagram it is easy to see that $p_{1}^{-1}p_{2}^{-1}pW$ is a rational element and therefore:
\begin{equation}
p(\Pi)\sim_{E(\Pi_{1})E(\Pi_{2});K} p(\Pi_{1})p(\Pi_{2})W(Id)^{-1}.
\end{equation}

Shahidi's calculation (c.f. Theorem $7.1.2$ of \cite{shahidi10} and Corollary $5.7$ of \cite{harrismotivic}) implies that:
\begin{equation}\nonumber
W(Id)^{-1} \sim_{E(\Pi);K} W_{\infty}(Id_{\infty})^{-1}\prod\limits_{w \text{ ramified places}} W_{w}(id_{w})^{-1}
\prod\limits_{1\leq i<j\leq l}L(1,\Pi_{i}\times \Pi_{j}^{\vee}).
\end{equation}
By the arguments in Corollary $5.7$ of \cite{harrismotivic}) we may choose $g_{1}$, $g_{2}$ such that $W_{w}(id_{w})$ is rational for each ramified place. At last, we conclude the proof by setting $\Omega_{(n_{1},n_{2})}(\Pi_{\infty}):=W_{\infty}(Id_{\infty})^{-1}$. We can read from the construction that it depend onlys on $\Pi_{\infty}$ and the parabolic data.
\begin{flushright}$\Box$\end{flushright}

\bigskip

For any partition $n_{1}+\cdots+n_{l}=n$, we may take $\Pi_{1},\cdots,\Pi_{l}$ as above such that $\Pi_{\infty}=(\Pi_{1}\boxplus \cdots\boxplus \Pi_{l})_{\infty}$. Hence $\Omega_{(n_{1},\cdots,n_{l})}(\Pi_{\infty})$ is well defined. In particular, $\Omega_{(1,1,\cdots,1)}(\Pi_{\infty})$ is well-defined. We denote it by $\Omega(\Pi_{\infty})$. We remark that this is the same archimedean factor appeared in Corollary $5.7$ of \cite{harrismotivic}.

\section{First discussions on archimedean factors}
We will discuss the archimedean factors $\Omega_{(n_{1},n_{2},\cdots,n_{l})}(\Pi_{\infty})$ defined in the last section.\\

One first observation is that 
\begin{equation}\label{first observation}
\Omega_{(n)}(\Pi_{\infty})\sim_{E(\Pi);K} 1.
\end{equation} 
This can be read directly from Equation (\ref{Whittaker decomposition pre}).\\

We observe that if $\Pi_{1}$ is also a Langlands sum of automorphic representations, we can furthermore decompose $p(\Pi_{1})$. We will get relations between the archimedean factors. In fact, we have:

\begin{lem}
If $\Pi$ is Langlands sum of $\Pi_{1}$, $\Pi_{2}$, $\cdots$, $\Pi_{l}$ then we have:
\begin{equation}
\Omega_{(n_{1},n_{2},\cdots,n_{l})}(\Pi_{\infty})\sim_{E(\Pi);K} \cfrac{\Omega(\Pi_{\infty})}{\prod\limits_{1\leq i\leq n}\Omega(\Pi_{i,\infty})}
\end{equation}
\end{lem}

\begin{dem}
We endow $P(n)$, the set of partitions of $n$, with the dictionary order. More precisely, if $(n_{1},\cdots,n_{l})$ and $(n'_{1},\cdots,n'_{l'})$ are two partitions of $n$, we say $(n_{1},\cdots,n_{l})< n'_{1},\cdots,n'_{l'})$ if there exists an integer $s\leq min\{l,l'\}$ such that $n_{i}=n'_{i}$ for $i<s$ and $n_{s}<n'_{s}$. The set $P(n)$ then becomes a totally ordered set. \\

We shall prove the lemma by induction on $n$. For each level $n$, we shall prove by induction on $P(n)$.\\

\textbf{(1) Basis:} When $n=1$, we know both sides are equivalent to $1$ by equation (\ref{first observation}).\\

\textbf{(2) Inductive step:} We assume that the lemma is true for $n_{1}+\cdots+n_{l}=n-1$ with $n\geq 2$. We shall prove it for $n_{1}+\cdots+n_{l}=n$ by induction on $P(n)$.\\

\textbf{(2.1) Basis:} The smallest element in $P(n)$ is $(1,1,\cdots,1)$. In this case, we have $\Omega_{(1,1,\cdots,1)}(\Pi_{\infty})\sim_{E(\Pi);K}\Omega(\Pi_{\infty})$ by definition. Moreover, $\Omega(\Pi_{i,\infty})\sim_{E(\Pi);K} \Omega_{(1)}(\Pi_{i,\infty})\sim_{E(\Pi);K} 1$ by equation (\ref{first observation}) for all $i$. The lemma then follows.\\

\textbf{(2.2) Inductive step:} Let $(n_{1},\cdots,n_{l})\neq (1,1,\cdots,1)\in P(n)$. We assume that the lemma holds for all elements in $P(n)$ smaller than $(n_{1},\cdots,n_{l})$. We now prove the lemma for $(n_{1},\cdots,n_{l})$.

Since $(n_{1},\cdots,n_{l})\neq (1,1,\cdots,1)$, there exists an integer $i$ such that $n_{i}\geq 2$. We take the smallest $i$ with this property and denote it by $t$.

We take positive integers $n^{'}_{t},n^{*}_{t}$ such that $n'_{t}+n^{*}_{t}=n_{t}$. For example, we may take $n'_{t}=n_{t}-1$ and $n^{*}_{t}=1$. We take $\Pi'_{t}$ and $\Pi^{*}_{t}$ to be cohomological conjugate self-dual regular representations of $GL_{n'_{t}}(\AF)$ and $GL_{n^{\*}_{t}}(\AF)$ respectively such that $\Pi_{t,\infty}$ is the same with the infinity type of the Langlands sum of $\Pi'_{t}$ and $\Pi^{*}_{t}$.

Let $\Pi^{\#}$ be the Langlands sum of $\Pi_{1},\cdots,\Pi_{t-1},\Pi'_{t},\Pi^{*}_{t},\Pi_{t+1},\cdots,\Pi_{l}$. We apply Proposition \ref{Shahidi Whittaker pre} to $(\Pi_{1},\cdots,\Pi_{t-1},\Pi'_{t},\Pi^{*}_{t},\Pi_{t+1},\cdots,\Pi_{l})$ and get:
\begin{eqnarray}\nonumber
p(\Pi^{\#})&\sim_{E(\Pi);K} &\Omega_{(n_{1},\cdots,n_{t-1},n'_{t},n^{*}_{t},n_{t+1}\cdots,n_{l})}(\Pi_{\infty})[\prod\limits_{i\neq t}p(\Pi_{i})]p(\Pi'_{t})p(\Pi^{*}_{t})\prod\limits_{ i<j, i\neq t,j\neq t}L(1,\Pi_{i}\times \Pi_{j}^{\vee})
\\\nonumber&&
\prod\limits_{i< t}(L(1,\Pi_{i}\times \Pi'_{t}\text{}^{\vee})L(1,\Pi_{i}\times \Pi_{t}^{*,\vee}))\prod\limits_{j>t}(L(1,\Pi'_{t}\times \Pi_{j}^{\vee})L(1,\Pi^{*}_{t}\times \Pi_{j}^{\vee}))L(\Pi'_{t}\times \Pi_{t}^{*,\vee})
\end{eqnarray}

Similarly, we apply Proposition \ref{Shahidi Whittaker pre} to $(\Pi_{1},\Pi_{t-1},\Pi'_{t}\boxplus \Pi^{*}_{t}, \Pi_{t+1},\cdots,\Pi_{l})$ and then to $(\Pi'_{t},\Pi^{*}_{t})$. We will get:
\begin{eqnarray}\nonumber
p(\Pi^{\#})&\sim_{E(\Pi);K}& \Omega_{(n_{1},\cdots,n_{t-1},n'_{t}+n^{*}_{t},n_{t+1)}\cdots,n_{l})}(\Pi_{\infty})[\prod\limits_{i\neq t}p(\Pi_{i})]p(\Pi'_{t}\boxplus \Pi^{*}_{t})\times \\\nonumber&&
\prod\limits_{ i<j, i\neq t,j\neq t}L(1,\Pi_{i}\times \Pi_{j}^{\vee})\prod\limits_{i<t}L(1,\Pi_{i}\times (\Pi'_{t}\boxplus \Pi^{*}_{t} )^{\vee})\prod\limits_{j> t}L(1,(\Pi'_{t}\boxplus \Pi^{*}_{t} )\times \Pi_{j}^{\vee})\\\nonumber
&\sim_{E(\Pi);K} &\Omega_{(n_{1},\cdots,n_{t-1},n'_{t}+n^{*}_{t},n_{t+1},\cdots,n_{l})}(\Pi_{\infty})\Omega_{(n'_{t},n^{*}_{t})}((\Pi'_{t}\boxplus \Pi^{*}_{t})_{\infty})[\prod\limits_{i\neq t}p(\Pi_{i})]p(\Pi'_{t})p(\Pi^{*}_{t})\times \\\nonumber&&
\prod\limits_{ i<j, i\neq t,j\neq t}L(1,\Pi_{i}\times \Pi_{j}^{c})\prod\limits_{i< t}(L(1,\Pi_{i}\times \Pi'_{t}\text{}^{\vee})L(1,\Pi_{i}\times \Pi_{t}^{*,\vee}))\times \\\nonumber
&&\prod\limits_{j>t}(L(1,\Pi'_{t}\times \Pi_{j}^{\vee})L(1,\Pi^{*}_{t}\times \Pi_{j}^{\vee}))L(\Pi'_{t}\times \Pi_{t}^{*,\vee}).
\end{eqnarray}

Recall that $n'_{t}+n^{*}_{t}=n_{t}$, we obtain that:
\begin{equation}
\Omega_{(n_{1},\cdots,n_{t-1},n'_{t},n^{*}_{t},n_{t+1}\cdots,n_{l})} (\Pi_{\infty})\sim_{E(\Pi;K)}\Omega_{(n_{1},\cdots,n_{t-1},n_{t},n_{t+1},\cdots,n_{l})}(\Pi_{\infty})\Omega_{(n'_{t},n^{*}_{t})}((\Pi'_{t}\boxplus \Pi^{*}_{t})_{\infty}) 
\end{equation}

Since $(n_{1},\cdots,n_{t-1},n'_{t},n^{*}_{t},n_{t+1}\cdots,n_{l})<(n_{1},n_{2},\cdots,n_{l})$, we may apply the hypothesis of the induction step $(2.2)$ and get:
\begin{equation}
\Omega_{(n_{1},\cdots,n_{t-1},n'_{t},n^{*}_{t},n_{t+1}\cdots,n_{l})} (\Pi_{\infty})\sim_{E(\Pi;K)}\cfrac{\Omega(\Pi_{\infty}) }{\prod\limits_{i\neq t}\Omega(\Pi_{i,\infty})\Omega(\Pi'_{t,\infty})\Omega(\Pi^{*}_{t,\infty})}.
\end{equation}

If $n_{t}=n$ then $l=1$ and both sides of the equation of the lemma are equivalent to $1$ by equation (\ref{first observation}). Hence we may assume that $n_{t}<n$. Therefore, the hypothesis of the induction step (2) is satisfied by $(\Pi'_{t},\Pi^{*}_{t})$. We get:
\begin{equation}
\Omega_{(n'_{t},n^{*}_{t})}((\Pi'_{t}\boxplus \Pi^{*}_{t})_{\infty}) 
 \sim_{E(\Pi;K)}\cfrac{\Omega(\Pi_{t,\infty}) }{\Omega(\Pi'_{t,\infty})\Omega(\Pi^{*}_{t,\infty})}.
\end{equation}

Comparing the above three equations, we finally deduce that the lemma is true for $(\Pi_{1},\cdots,\Pi_{n})$ and complete the proof.
\end{dem}

\begin{cor}\label{Shahidi Whittaker}
If $\Pi$ is Langlands sum of $\Pi_{1}$, $\Pi_{2}$, $\cdots$, $\Pi_{l}$ then we have:
\begin{equation}\label{Whittaker decomposition}
p(\Pi)\sim_{E(\Pi);K} \cfrac{\Omega(\Pi_{\infty})}{\prod\limits_{1\leq i\leq n}\Omega(\Pi_{i,\infty})}
\prod\limits_{1\leq i\leq l}p(\Pi_{i})\prod\limits_{1\leq i<j\leq l}L(1,\Pi_{i}\times \Pi_{j}^{\vee})
\end{equation}
\end{cor} 

\section{Special values of tensor products in terms of Whittaker periods, after Grobner-Harris}
 Let $\Pi$ be a regular cuspidal cohomological representation of $GL_{n}(\AF)$. Let $\Pi^{\#}$ be a regular automorphic cohomological representation of $GL_{n-1}(\AF)$ which is the Langlands sum of cuspidal representations. Equivalently, it is a regular element in $\mathcal{A}lg(n-1)$.
 
The arguments in section of \cite{harrismotivic} go over word for word and give the following result:
  
 \begin{prop}\label{Whittaker period theorem CM}
 We assume that $(\Pi,\Pi^{\#})$ is in good position.
 
 There exists a complex number $p(m,\Pi_{\infty},\Pi^{\#}_{\infty})$ which depends on $m,\Pi_{\infty}$ and $\Pi^{\#}_{\infty}$ well defined up to $(E(\Pi)E(\Pi^{\#}))^{\times}$ such that for $m\in \N$ with $m+\cfrac{1}{2}$ critical for $\Pi\times \Pi^{\#}$, we have
\begin{equation}\label{CM Whittaker period}
L(\cfrac{1}{2}+m,\Pi\times \Pi^{\#})\sim_{E(\Pi)E(\Pi^{\#});K} p(m,\Pi_{\infty},\Pi^{\#}_{\infty})p(\Pi)p(\Pi^{\#})
\end{equation}
where $p(\Pi)$ and $p(\Pi^{\#})$ are the Whittaker periods of $\Pi$ and $\Pi^{\#}$ respectively.  \end{prop}
 
We remark that we don't need $\Pi^{\#}$ to be cuspidal here. The above conditions are sufficient to guarantee that a certain Eisenstein series is holomorphic. \\

Moreover, we remark that the good position condition is necessary so that a certain intertwining operator exists.\\

We shall give a proof of this proposition in a separate article.

\chapter{CM periods and arithmetic automorphic periods}
\section{CM periods}\label{section CM periods}
Let $(T,h)$ be a Shimura datum where $T$ is a torus defined over $\Q$ and $h:Res_{\C/\R}\mathbb{G}_{m,\C}\rightarrow G_{\R}$ a homomorphism satisfying the axioms defining a Shimura variety. Such pair is called a \textbf{special} Shimura datum. Let $Sh(T,h)$ be the associated Shimura variety and $E(T,h)$ be its reflex field.\\

Let $(\gamma,V_{\gamma})$ be a one-dimensional algebraic representation of $T$ (the representation $\gamma$ is denoted by $\chi$ in \cite{harrisappendix}). We denote by $E(\gamma)$ a definition field for $\gamma$. We may assume that $E(\gamma)$ contains $E(T,h)$. Suppose that $\gamma$ is motivic (see \textit{loc.cit} for the notion). We know that $\gamma$ gives an automorphic line bundle $[V_{\gamma}]$ over $Sh(T,h)$ defined over $E(\gamma)$. Therefore, the complex vector space $H^{0}(Sh(T,h),[V_{\gamma}])$ has an $E(\gamma)$-rational structure, denoted by $M_{DR}(\gamma)$ and called the De Rham rational structure.

\bigskip

On the other hand, the canonical local system $V_{\gamma}^{\triangledown}\subset [V_{\gamma}]$ gives another $E(\gamma)$-rational structure $M_{B}(\gamma)$ on $H^{0}(Sh(T,h),[V_{\gamma}])$, called the Betti rational structure.\\

We now consider $\chi$ an algebraic Hecke character of $T(\AQ)$ with infinity type $\gamma^{-1}$ (our character $\chi$ corresponds to the character $\omega^{-1}$ in \textit{loc.cit}). Let $E(\chi)$ be the number field generated by the values of $\chi$ on $T(\AQ \text{}_{,f})$ over $E(\gamma)$. We know $\chi$ generates a one-dimensional complex subspace of $H^{0}(Sh(T,h),[V_{\gamma}])$ which inherits two $E(\chi)$-rational structures, one from $M_{DR}(\gamma)$, the other from $M_{B}(\gamma)$. Put $p(\chi,(T,h))$ the ratio of these two rational structures which is well defined modulo $E(\chi)^{\times}$.

\bigskip

\begin{rem}\label{CM complex conjugation}
If we identify $H^{0}(Sh(T,h),[V_{\gamma}])$ with the set $\{f\in \mathbb{C}^{\infty}(T(\Q)\backslash T(\AQ),\C\mid f(tt_{\infty}))=\gamma^{-1}(t_{\infty})f(t), t_{\infty}\in T(\R), t\in T(\AQ)\}$, then $\chi$ itself is in the rational structure inherits from $M_{B}(\gamma)$. See discussion from $A.4$ to $A.5$ in \cite{harrisappendix}.
\end{rem}

\bigskip

Suppose that we have two tori $T$ and $T'$ both endowed with a Shimura datum $(T,h)$ and $(T',h')$. Let $u:(T',h')\rightarrow (T,h)$ be a map between the Shimura data. Let $\chi$ be an algebraic Hecke character of $T(\AQ)$. We put $\chi':=\chi\circ u$ an algebraic Hecke character of $T'(\AQ)$. Since both the Betti structure and the De Rham structure commute with the pullback map on cohomology, we have the following proposition:

\begin{prop}\label{propgeneral}
Let $\chi$, $(T,h)$ and $\chi'$, $(T',h')$ be as above. We have:
\begin{equation}\nonumber
p(\chi,(T,h)) \sim_{E(\chi)} p(\chi',(T',h'))
\end{equation}
\end{prop}

\begin{rem}
In Proposition $1.4$ of \cite{harrisCMperiod}, the relation is up to $E(\chi);E(T,h)$ where $E(T,h)$ is a number field associated to $(T,h)$. Here we consider the action of $G_{\Q}$ and can thus obtain a relation up to $E(\chi)$ (see the paragraph after Proposition $1.8.1$ of \textit{loc.cit}).
\end{rem}

\bigskip

For $F$ a CM field and $\Psi$ a subset of $\Sigma_{F}$ such that $\Psi\cap \iota\Psi=\varnothing$, we can define a Shimura datum $(T_{F},h_{\Psi})$ where $T_{F}:=Res_{F/\Q}\mathbb{G}_{m,F}$ is a torus and $h_{\Psi}:Res_{\C/\R}\mathbb{G}_{m,\C} \rightarrow T_{F,\R}$ is a homomorphism such that over $\sigma\in \Sigma_{F}$, the Hodge structure induced on $F$ by $h_{\Psi}$ is of type $(-1,0)$ if $\sigma\in \Psi$, of type $(0,-1)$ if $\sigma\in \iota\Psi$, and of type $(0,0)$ otherwise. \\

Let $\chi$ be a motivic critical character of a CM field $F$. By definition, $p_{F}(\chi,\Psi)=p(\chi,(T_{F},h_{\Psi}))$ and we call it a \textbf{CM period}. Sometimes we write $p(\chi,\Psi)$ instead of $p_{F}(\chi,\Psi)$ if there is no ambiguity concerning the base field $F$.

\begin{ex}
We have $p(||\cdot||_{\AK},1)\sim_{\Q} (2\pi i)^{-1}.$ See  $(1.10.9)$ on page $100$ of \cite{harris97}. 
\end{ex}

\bigskip

Let $\theta \in Gal(F/\Q)$. We know $\theta$ induces an action on $\Sigma_{F}$ by composition with $\theta$. Moreover, $\theta$ acts on $\AF^{\times}$ and hence acts on the set of Hecke characters of $F$.\\

The CM periods have many good properties. We list below some of them which will be useful in the future.

\begin{prop}\label{propCM}

Let $F$ be a CM field. Let $F_{0}\subset F$ be a sub CM field.

Let $\eta$ be a motivic critical Hecke character of $F_{0}$, $\chi$, $\chi_{1}$, $\chi_{2}$ be motivic critical Hecke characters of $F$.

Let $\tau\in \Sigma_{F}$ be an embedding of $F$ into $\bar{\Q}$ and $\Psi$ be a subset of $\Sigma_{F}$ such that $\Psi\cap \Psi^{c}=\varnothing$. We take $\Psi=\Psi_{1}\sqcup\Psi_{2}$ a partition of $\Psi$. 

Let $\theta$ be an element in $Gal(F/\Q)$. We then have:
\begin{eqnarray}\label{charmulti}
&p(\chi_{1}\chi_{2}), \Psi) \sim _{E(\chi_{1})E(\chi_{2})} p(\chi_{1}, \Psi)p(\chi_{2}^{\sigma}, \Psi).&
\\ \label{separateCMtype}
&p(\chi, \Psi_{1}\sqcup\Psi_{2})\sim _{E(\chi) } p(\chi, \Psi_{1})p(\chi,\Psi_{2}).&
\\
&p(\chi^{\theta}, \Psi^{\theta}) \sim _{E(\chi)} p(\chi, \Psi).&\\
&p_{F}(\eta\circ N_{\AF/\AFsub},\tau ) \sim _{E(\eta)}  p_{F_{0}}(\eta, \tau|_{F_{0}}) .&
\end{eqnarray}

In particular, if we take $\theta=c$ the complex conjugation, we have:
\begin{equation}
p(\chi, \Psi^{c}) \sim _{E(\chi)} p(\chi^{c}, \Psi).\end{equation}

\end{prop}

\begin{rem}
The first three formulas come from Proposition $1.4$, Corollary $1.5$ and Lemma $1.6$ in \cite{harrisCMperiod}. The last formula is a variation of the Lemma $1.8.3$ in \textit{loc.cit}. The idea was explained in the proof of Proposition $1.4$ in \textit{loc.cit}. We sketch the proof here.
\end{rem}

\begin{proof}
All the equations in Proposition \ref{propCM} come from Proposition \ref{propgeneral} by certain maps between Shimura data as follows:
\begin{enumerate}
\item The diagonal map $ (T_{F},h_{\Psi})\rightarrow (T_{F}\times T_{F},h_{\Psi}\times h_{\Psi})$ pulls $(\chi_{1},\chi_{2})$ back to $\chi_{1}\chi_{2}$.
\item The multiplication map $T_{F} \times T_{F} \rightarrow T_{F}$ sends $h_{\Psi_{1}}$, $h_{\Psi_{2}}$ to $h_{\Psi_{1}\sqcup \Psi_{2}}$.
\item The Galois action $\theta: H_{F}\rightarrow H_{F}$ sends $h_{\Psi}$ to $h_{\Psi^{\theta}}$.
\item The norm map $(T_{F},h_{\{\tau\}})\rightarrow (T_{F_{0}},h_{\{\tau|_{F_{0}}\}})$ pulls $\eta$ back to $\eta\circ N_{\AF/\mathbb{A}_{F,0}}$.
\end{enumerate}
\end{proof}

\bigskip
The special values of an $L$-function for a Hecke character over a CM field can be interpreted in terms of CM periods. The following theorem is proved by Blasius. We state it as in Proposition $1.8.1$ in \cite{harrisCMperiod} where $\omega$ should be replaced by $\check{\omega}:=\omega^{-1,c}$ (for this erratum, see the notation and conventions part on page $82$ in the introduction of \cite{harris97}),

\begin{thm}\label{blasius}
Let $F$ be a CM field and $F^{+}$ be its maximal totally real subfield. Put $n$ the degree of $F^{+}$ over $\Q$.

Let $\chi$ be a motivic critical algebraic Hecke character of $F$ and $\Phi_{\chi}$ be the unique CM type of $F$ which is compatible with $\chi$.

 Let $D_{F^{+}}$ be the absolute discriminant of $F^{+}$.
We assume that $D_{F^{+}}^{1/2}\in E(\chi)$ for simplicity.

For $m$ a critical value of $\chi$ in the sense of Deligne (c.f. Lemma \ref{critical}), we have
\begin{equation}\nonumber
L(\chi,m) \sim_{E(\chi)} (2\pi i)^{mn}p(\check{\chi},\Phi_{\chi})
\end{equation}
equivariant under action of $G_{\Q}$\end{thm}

\begin{rem}\label{criticalforcharacter}
\text{}
\begin{enumerate}
\item Let $\{\sigma_{1},\sigma_{2},\cdots,\sigma_{n}\}$ be any CM type of $F$. Let $(\sigma_{i}^{a_{i}}\overline{\sigma}_{i}^{-w-a_{i}})_{1\leq i\leq n}$ denote the infinity type of $\chi$ with $w=w(\chi)$. We may assume $a_{1}\geq a_{2}\geq \cdots \geq a_{n}$. We define $a_{0}:=+\infty$ and $a_{n+1}:=-\infty$ and define $k:=max\{0\leq i\leq n\mid a_{i}>-\cfrac{w}{2}\}$. An integer $m$ is critical for $\chi$ if and only if

\begin{equation}\label{criticalvalue1}
       max(-a_{k}+1,w+1+a_{k+1})\leq m \leq min(w+a_{k},-a_{k+1})
\end{equation}
     (c.f. Lemma \ref{critical}).
\item $D_{F^{+}}^{1/2}$ is well defined up to multiplication by $\pm 1$. More generally, if $\{z_{1},z_{2},\cdots,z_{n}\}$ is any $\Q$-base of $L$, then $det(\sigma_{i}(z_{j}))_{1\leq i,j\leq n}\sim_{\Q} D_{F^{+}}^{1/2}$.

 \end{enumerate}
\end{rem}

\section{Construction of cohomology spaces}
\text{}
Let $\Sigma=\Sigma_{F;K}$ in the current and the following chapters.
Fix an index $I$ as before. Write $s_{\sigma}:=I(\sigma)$ and $r_{\sigma}:=n-I(\sigma)$ for all $\sigma\in\Sigma$.\\

Denote $\Ss:=Res_{\C/\R}\mathbb{G}_{m}$. Recall that $GU_{I}(\R)$ is isomorphic to a subgroup of $\prod\limits_{\sigma\in \Sigma}GU(r_{\sigma},s_{\sigma})$ defined by the same \textit{similitude}. We can define a homomorphism $h_{I}:\Ss(\R) \rightarrow GU_{I}(\R)$ by sending $z\in \C$ to $\left(\begin{pmatrix}
zI_{r_{\sigma}} & 0\\
0 & \bar{z}I_{s_{\sigma}}
\end{pmatrix}\right)_{\sigma\in \Sigma}$. \\

Let $X_{I}$ be the $GU_{I}(\R)$-conjugation class of $h_{I}$. We know $(GU_{I},X_{I})$ is a Shimura datum with reflex field $E_{I}$ and dimension $2\sum\limits_{\sigma\in \Sigma} r_{\sigma}s_{\sigma}$. The Shimura variety associated to $(GU_{I},X_{I})$ is denoted by $Sh_{I}$.

\bigskip

Let $K_{I,\infty}$ be the centralizer of $h_{I}$ in $GU_{I}(\R)$. Via the inclusion $GU_{I}(\R)\hookrightarrow \prod\limits_{\sigma\in \Sigma}GU(r_{\sigma},s_{\sigma})\subset \R^{+,\times}\prod\limits_{\sigma\in\Sigma} U(n,\C)$, we may identify $K_{I,\infty}$ with \begin{equation}\nonumber   \{  (\mu,\begin{pmatrix}
 u_{r_{\sigma}} & 0\\
0 &  v_{s_{\sigma}}
\end{pmatrix}_{\sigma\in \Sigma}) \mid u_{r_{\sigma}}\in U(r_{\sigma},\C), v_{s_{\sigma}}\in U(s_{\sigma},\C), \mu\in \R^{+,\times}\}\end{equation}  
where $U(r,\C)$ is the standard unitary group of degree $r$ over $\C$. Let $H_{I}$ be the subgroup of $K_{I,\infty}$ consisting of the diagonal matrices in $K_{I,\infty}$. Then it is a maximal torus of $GU_{I}(\R)$. Denote its Lie algebra by $\lieH_{I}$.

\bigskip

We observe that $H_{I}(\R) \cong \R^{+,\times}\times \prod\limits_{\sigma\in\Sigma} U(1,\C)^{n}$. Its algebraic characters are of the form 
\begin{equation}
\nonumber   (w,(z_{i}(\sigma))_{\sigma\in\Sigma,1\leq i\leq n})\mapsto w^{\lambda_{0}} \prod\limits_{\sigma\in\Sigma} \prod\limits_{i=1}^{n}z_{i}(\sigma)^{\lambda_{i}(\sigma)}
\end{equation}  
 where $(\lambda_{0},(\lambda_{i}(\sigma))_{\sigma\in\Sigma,1\leq i\leq n})$ is a $(nd+1)$-tuple of integers with $\lambda_{0} \equiv \sum\limits_{\sigma\in\Sigma}\sum\limits_{i=1}^{n}\lambda_{i}(\sigma)$ (mod $2$).

 \bigskip
 
 Recall that $GU_{I}(\C)\cong \C^{\times}\prod_{\sigma\in\Sigma}GL_{n}(\C)$. We fix $B_{I}$ the Borel subgroup of $GU_{I,\C}$ consisting of upper triangular matrices. The highest weights of finite-dimensional irreducible representations of $K_{I,\infty}$ are tuples $\Lambda=(\Lambda_{0},(\Lambda_{i}(\sigma))_{\sigma\in\Sigma,1\leq i\leq n})$ such that $\Lambda_{1}(\sigma)\geq \Lambda_{2}(\sigma)\geq \cdots\geq \Lambda_{r_{\sigma}}(\sigma)$, $ \Lambda_{r_{\sigma}+1}(\sigma) \geq \cdots \geq  \Lambda_{n}(\sigma)$ for all $\sigma$ and $\Lambda_{0} \equiv \sum\limits_{\sigma\in\Sigma}\sum\limits_{i=1}^{n}\Lambda_{i}(\sigma)$ (mod $2$). 
 
\bigskip

We denote the set of such tuples by $\Lambda(K_{I,\infty})$. Similarly, we write $\Lambda(GU_{I})$ for the set of the highest weights of finite-dimensional irreducible representations of $GU_{I}$. It consists of tuples $\lambda=(\lambda_{0},(\lambda_{i}(\sigma))_{\sigma\in\Sigma,1\leq i\leq n})$ such that $\lambda_{1}(\sigma)\geq \lambda_{2}(\sigma) \geq  \cdots \lambda_{n}(\sigma)$ for all $\sigma$ and $\lambda_{0} \equiv \sum\limits_{\sigma\in\Sigma}\sum\limits_{i=1}^{n}\lambda_{i}(\sigma)$ (mod $2$). 
 
 \bigskip
 
We take $\lambda\in\Lambda(GU_{I})$ and $\Lambda\in \Lambda(K_{I,\infty})$.

Let $V_{\lambda}$  and $V_{\Lambda}$ be the corresponding representations.  We define a local system over $Sh_{I}$:
\begin{equation}\nonumber   W_{\lambda}^{\triangledown}:=\lim_{\overleftarrow{K}} GU_{I}(\Q)\backslash V_{\lambda}\times X\times GU_{I}(\AQf)/K\end{equation}  
and an automorphic vector bundle over $Sh_{I}$
\begin{equation}\nonumber   E_{\Lambda}:=\lim_{\overleftarrow{K}} GU_{I}(\Q)\backslash V_{\Lambda}\times GU_{I}(\R) \times GU_{I}(\AQf)/KK_{I,\infty}\end{equation}  
where $K$ runs over open compact subgroup of $GU_{I}(\AQf)$.

\bigskip

The automorphic vector bundles $E_{\Lambda}$ are defined over the reflex field $E$. \\

The local systems $W_{\lambda}^{\triangledown}$ are defined over $K$. The Hodge structure of the cohomology space $H^{q}(Sh_{I},W_{\lambda}^{\triangledown})$ is not pure in general. But the image of $H_{c}^{q}(Sh_{I},W_{\lambda}^{\triangledown})$ in $H^{q}(Sh_{I},W_{\lambda}^{\triangledown})$ is pure of weight $q-c$. We denote this image by $\bar{H}^{q}(Sh_{I},W_{\lambda}^{\triangledown})$. \\

Note that all cohomology spaces have coefficients in $\C$ unless we specify its rational structure over a number field.

\bigskip

\section{The Hodge structures}
\text{}
The results in section $2.2$ of \cite{harrismotives} give a description of the Hodge components of $\bar{H}^{q}(Sh_{I},W_{\lambda}^{\triangledown})$.\\

Denote by $R^{+}$ the set of positive roots of $H_{I,\C}$ in $GU_{I}(\C)$ and by $R_{c}^{+}$ the set of positive compact roots. Define $\alpha_{j,k}=(0,\cdots,0,1,0,\cdots,0,-1,0,\cdots,0)$ for any $1\leq j<k\leq n$. We know $R^{+}=\{(\alpha_{j_{\sigma},k_{\sigma}})_{\sigma\in \Sigma}\mid 1\leq j_{\sigma}<k_{\sigma}\leq n\}$ and 
$R_{c}^{+}=\{(\alpha_{j_{\sigma},k_{\sigma}})_{\sigma\in \Sigma}\mid j_{\sigma}<k_{\sigma}\leq r_{\sigma} \text{ or } r_{\sigma}+1\leq j_{\sigma}<k_{\sigma}\}$.\\

Let $\rho=\cfrac{1}{2}\sum\limits_{\alpha\in R^{+}}\alpha=\left(\left(\cfrac{n-1}{2},\cfrac{n-3}{2},\cdots,-\cfrac{n-1}{2}\right)\right)_{\sigma}$. \\

Let $\lieG$, $\lieK$ and $\lieH$ be Lie algebras of $GU_{I}(\R)$, $K_{I,\infty}$ and $H(\R)$. Write $W$ for the Weyl group $W(\lieG_{\C},\lieH_{\C})$ and $W_{c}$ for the Weyl group $W(\lieK_{\C},\lieH_{\C})$. We can identify $W$ with $\prod\limits_{\sigma\in\Sigma}\mathfrak{S}_{n}$ and $W_{c}$ with $\prod\limits_{\sigma\in\Sigma}\mathfrak{S}_{r_{\sigma}}\times \mathfrak{S}_{s_{\sigma}}$ where $\mathfrak{S}$ refers to the standard permutation group. For $w\in W$, we write the length of $w$ by $l(w)$.

\bigskip

Let $W^{1}:=\{w\in W| w(R^{+})\supset R_{c}^{+}\}$ be a subset of $W$. By the above identification, $(w_{\sigma})_{\sigma}\in W^{1}$ if and only if  $w_{\sigma}(1)<w_{\sigma}(2)<\cdots<w_{\sigma}(r_{\sigma})$ and $w_{\sigma}(r_{\sigma}+1)<\cdots<w_{\sigma}(n)$  One can show that $W^{1}$ is a set of coset representatives of shortest length for $W_{c}\backslash W$. 

\bigskip

Moreover, for $\lambda$ a highest weight of a representation of $GU_{I}$, one can show easily that $w*\lambda:=w(\lambda+\rho)-\rho$ is the highest weight of a representation of $K_{I,\infty}$. More precisely, if $\lambda=(\lambda_{0},(\lambda_{i}(\sigma))_{\sigma\in\Sigma,1\leq i\leq n})$, then $w*\lambda=(\lambda_{0},((w*\lambda)_{i}(\sigma))_{\sigma\in\Sigma,1\leq i\leq n})$ with $(w*\lambda)_{i}(\sigma)=\lambda_{w_{\sigma}(i)}(\sigma)+\cfrac{n+1}{2}-w_{\sigma}(i)-(\cfrac{n+1}{2}-i)=\lambda_{w_{\sigma}(i)}(\sigma)-w_{\sigma}(i)+i$.

\bigskip

\begin{rem}\label{hodgetype}
The results of \cite{harrismotives} tell us that there exists
\begin{equation}
\bar{H}^{q}(Sh_{I},W_{\lambda}^{\triangledown})\cong \bigoplus\limits_{w\in W^{1}}\bar{H}^{q;w}(Sh_{I},W_{\lambda}^{\triangledown})
\end{equation}
a decomposition as subspaces of pure Hodge type $(p(w,\lambda),q-c-p(w,\lambda))$. We now determine the Hodge number $p(w,\lambda)$.

\bigskip

We know that $w*\lambda$ is the highest weight of a representation of $K_{I,\infty}$. We denote this representation by $(\rho_{w*\lambda},W_{w*\lambda})$. We know that $\rho_{w*\lambda} \circ h_{I}|_{\Ss(\R)}: \Ss(\R)\rightarrow K_{I,\infty} \rightarrow GL(W_{w*\lambda})$ is of the form $z\mapsto z^{-p(w,\lambda)}\bar{z}^{-r(w,\lambda)}I_{W_{w*\lambda}}$ with $p(w,\lambda)$, $r(w,\lambda)\in \Z$. The first index $p(w,\lambda)$ is the Hodge type mentioned above. 

Recall that the map 
\begin{eqnarray}
h_{I}|_{\Ss(\R)}: \Ss(\R)&\rightarrow& K_{I,\infty}\subset \R^{+,\times}\times U(n,\C)^{\Sigma}\\ z&\mapsto& \left(|z|,
\begin{pmatrix}
\cfrac{z}{|z|}I_{r_{\sigma}} & 0\\\nonumber
0 & \cfrac{\bar{z}}{|z|}I_{s_{\sigma}}
\end{pmatrix}_{\sigma\in\Sigma}\right)
\end{eqnarray}
and the map 
\begin{eqnarray}
\nonumber \rho_{w*\lambda}: K_{I,\infty} &\rightarrow& GL(W_{w*\lambda}) \\\nonumber (w,diag(z_{i}(\sigma))_{\sigma\in\Sigma,1\leq i\leq n})&\mapsto & w^{\lambda_{0}}\prod\limits_{\sigma\in\Sigma} \prod\limits_{i=1}^{n}z_{i}(\sigma)^{(w*\lambda)_{i}(\sigma)}
\end{eqnarray}
 where $diag(z_{1},z_{2},\cdots,z_{n})$ means the diagonal matrix of coefficients $z_{1},z_{2},\cdots,z_{n}$. 
 
\bigskip

Therefore we have:

\begin{eqnarray}
z^{-p(w,\lambda)}\bar{z}^{-r(w,\lambda)} &=&|z|^{\lambda_{0}}\prod\limits_{\sigma\in\Sigma}(\prod\limits_{1\leq i\leq r_{\sigma}}(\cfrac{z}{|z|})^{(w*\lambda)_{i}(\sigma)}\prod\limits_{r_{\sigma}+1\leq i\leq n}(\cfrac{\overline{z}}{|z|})^{(w*\lambda)_{i}(\sigma)}\nonumber\\
&=& (z^{\frac{1}{2}}\overline{z}^{\frac{1}{2}})^{\lambda_{0}-\sum\limits_{\sigma\in\Sigma}\sum\limits_{1\leq i\leq n}(w*\lambda)_{i}(\sigma)}z^{\sum\limits_{\sigma\in\Sigma}\sum\limits_{1\leq i\leq r_{\sigma}}(w*\lambda)_{i}(\sigma)}
\overline{z}^{\sum\limits_{\sigma\in\Sigma}\sum\limits_{r_{\sigma}+1\leq i\leq n}(w*\lambda)_{i}(\sigma)} \nonumber
\end{eqnarray}

Since $(w*\lambda)_{i}(\sigma)=\lambda_{w_{\sigma}(i)}(\sigma)-w_{\sigma}(i)+i$ and then $\sum\limits_{\sigma\in\Sigma}\sum\limits_{1\leq i\leq n}(w*\lambda)_{i}(\sigma)=\sum\limits_{\sigma\in\Sigma}\sum\limits_{1\leq i\leq n}\lambda_{i}(\sigma)$, we obtain that:
\begin{eqnarray}\label{smallest}
p(w,\lambda) &=& \cfrac{\sum\limits_{\sigma\in\Sigma}\sum\limits_{1\leq i\leq n}\lambda_{i}(\sigma)-\lambda_{0}}{2}-\sum\limits_{\sigma\in\Sigma}\sum\limits_{1\leq i\leq r_{\sigma}}(w*\lambda)_{i}(\sigma)\nonumber\\
&=& \cfrac{\sum\limits_{\sigma\in\Sigma}\sum\limits_{1\leq i\leq n}\lambda_{i}(\sigma)-\lambda_{0}}{2}-\sum\limits_{\sigma\in\Sigma}\sum\limits_{1\leq i\leq r_{\sigma}}(\lambda_{w_{\sigma}(i)}(\sigma)-w_{\sigma}(i)+i)
\end{eqnarray}
\end{rem}

\bigskip

The method of toroidal compactification gives us more information on $\bar{H}^{q;w}(Sh_{I},W_{\lambda}^{\triangledown})$. We take $j:Sh_{I} \hookrightarrow \widetilde{Sh}_{I}$ to be a smooth toroidal compactification. Proposition $2.2.2$ of \cite{harrismotives} tells us that the following results do not depend on the choice of the toroidal compactification. 

\bigskip

The automorphic vector bundle $E_{\Lambda}$ can be extended to $\widetilde{Sh}_{I}$ in two ways: the canonical extension $E_{\Lambda}^{can}$ and the sub canonical extension $E_{\Lambda}^{sub}$ as explained in \cite{harrismotives}. Define:
\begin{equation}\nonumber   \bar{H}^{q}(Sh_{I},E_{\Lambda})=Im(H^{q}(\widetilde{Sh}_{I},E_{\Lambda}^{sub})\rightarrow H^{q}(\widetilde{Sh}_{I},E_{\Lambda}^{can})).\end{equation}

\bigskip

\begin{prop}
There is a canonical isomorphism
\begin{equation}\nonumber   \bar{H}^{q;w}(Sh_{I},W_{\lambda}^{\triangledown})\cong \bar{H}^{q-l(w)}(Sh_{I},E_{w*\lambda})\end{equation}  
\end{prop}

\bigskip

Let $D=2\sum\limits_{\sigma\in\Sigma}r_{\sigma}s_{\sigma}$ be the dimension of the Shimura variety. We are interested in the cohomology space of degree $D/2$. Proposition $2.2.7$ of \cite{harris97} also works here:

\bigskip

\begin{prop}
The space $\bar{H}^{D/2}(Sh_{I},W_{\lambda}^{\triangledown})$ is naturally endowed with a $K$-rational structure, called the de Rham rational structure and noted by $\bar{H}_{DR}^{D/2}(Sh_{I},W_{\lambda}^{\triangledown})$. This rational structure is endowed with a $K$-Hodge filtration $F^{\cdot}\bar{H}_{DR}^{D/2}(Sh_{I},W_{\lambda}^{\triangledown})$ pure of weight $D/2-c$ such that
\begin{equation}\nonumber   F^{p}\bar{H}_{DR}^{D/2}(Sh_{I},W_{\lambda}^{\triangledown})/F^{p+1}\bar{H}_{DR}^{D/2}(Sh_{I},W_{\lambda}^{\triangledown})\otimes_{K} \C \cong \bigoplus\limits_{w\in W^{1},p(w,\lambda)=p} \bar{H}^{D/2;w}(Sh_{I},W_{\lambda}^{\triangledown}).\end{equation}  
Moreover, the composition of the above isomorphism and the canonical isomorphism \begin{equation}\nonumber
\bar{H}^{D/2;w}(Sh_{I},W_{\lambda}^{\triangledown})\cong \bar{H}^{D/2-l(w)}(Sh_{I},E_{w*\lambda})
\end{equation} is rational over $K$.
\end{prop}

\bigskip

\paragraph{Holomorphic part:}

Let $w_{0}\in W^{1}$ defined by $w_{0}(\sigma)(1,2,\cdots,r_{\sigma};r_{\sigma+1},\cdots,n)_{\sigma\in\Sigma}=(s_{\sigma+1},\cdots,n;1,2,\cdots,s_{\sigma})$ for all $\sigma\in\Sigma$. It is the only longest element in $W^{1}$. Its length is $D/2$. 

We have a $K$-rational isomorphism \begin{equation}\bar{H}^{D/2;w_{0}}(Sh_{I},W_{\lambda}^{\triangledown})\cong \bar{H}^{0}(Sh_{I},E_{w_{0}*\lambda}).
\end{equation}

We can calculate the Hodge type of $\bar{H}^{D/2;w_{0}}(Sh_{I},W_{\lambda}^{\triangledown})$ as in Remark \ref{hodgetype}.

By definition we have \begin{equation}\label{twist by w0}
w_{0}*\lambda=(\lambda_{0},(\lambda_{s_{\sigma}+1}(\sigma)-s_{\sigma},\cdots,\lambda_{n}(\sigma)-s_{\sigma};\lambda_{1}(\sigma)+r_{\sigma},\cdots,\lambda_{s_{\sigma}}(\sigma)+r_{\sigma})_{\sigma\in\Sigma}).
\end{equation} By the discussion in Remark \ref{hodgetype}, the Hodge number 

\begin{equation}\nonumber   
p(w_{0},\lambda)=\cfrac{\sum\limits_{\sigma\in\Sigma}\sum\limits_{1\leq i\leq n}\lambda_{i}(\sigma)-\lambda_{0}+D}{2}-\sum\limits_{\sigma\in\Sigma}(\lambda_{s_{\sigma}+1}(\sigma)+\cdots+\lambda_{n}(\sigma)).\end{equation}  

From equation (\ref{smallest}), it is easy to deduce that $p(w_{0},\lambda)$ is the only largest number among $\{p(w,\lambda)\mid w\in W^{1}\}$. Therefore \begin{equation}F^{p(w_{0},\lambda)} (Sh_{I},W_{\lambda}^{\triangledown})\otimes_{K} \C \cong \bar{H}^{0}(Sh_{I},E_{w_{0}*\lambda}).
\end{equation} Moreover, as mentioned in the above proposition, we know that the above isomorphism is $K$-rational.\\

We call $\bar{H}^{D/2;w_{0}}(Sh_{I},W_{\lambda}^{\triangledown})\cong\bar{H}^{0}(Sh_{I},E_{w_{0}*\lambda})$ the \textbf{holomorphic part} of the Hodge decomposition of $\bar{H}^{D/2}(Sh_{I},W_{\lambda}^{\triangledown})$. It is isomorphic to the space of holomorphic cusp forms of type $(w_{0}*\lambda)^{\vee}$.\\

\paragraph{Anti-holomorphic part:}\text{}
The only shortest element in $W^{1}$ is the identity with the smallest Hodge number 
\begin{equation}\nonumber   
p(id,\lambda)=\cfrac{\sum\limits_{\sigma\in\Sigma}\sum\limits_{1\leq i\leq n}\lambda_{i}(\sigma)-\lambda_{0}}{2}-\sum\limits_{\sigma\in\Sigma}(\lambda_{1}(\sigma)+\cdots+\lambda_{r_{\sigma}}(\sigma)).
\end{equation} 
We call $\nonumber   \bar{H}^{D/2;id}(Sh_{I},W_{\lambda}^{\triangledown})\cong \bar{H}^{D/2}(Sh_{I},E_{\lambda})$
the \textbf{anti-holomorphic part} of the Hodge decomposition of $\bar{H}^{D/2}(Sh_{I},W_{\lambda}^{\triangledown})$.

\section{Complex conjugation}\label{Complex conjugation}

We specify some notation first. \\

Let $\lambda=(\lambda_{0},(\lambda_{1}(\sigma)\geq \lambda_{2}(\sigma)\geq \cdots \geq \lambda_{n}(\sigma))_{\sigma\in\Sigma}) \in \Lambda(GU_{I})$ as before. We define $\lambda^{c}:=(\lambda_{0},(-\lambda_{n}(\sigma)\geq -\lambda_{n-1}(\sigma)\geq \cdots \geq -\lambda_{1}(\sigma))_{\sigma\in\Sigma})$ and $\lambda^{\vee}:=(-\lambda_{0},(-\lambda_{n}(\sigma)\geq -\lambda_{n-1}(\sigma)\geq \cdots \geq -\lambda_{1}(\sigma))_{\sigma\in\Sigma})$. They are elements in $\Lambda(GU_{I})$. Moreover, the representation $V_{\lambda^{c}}$ is the complex conjugation of $V_{\lambda}$ and the representation $V_{\lambda^{\vee}}$ is the dual of $V_{\lambda}$ as $GU_{I}$-representation.\\

Similarly, for $\Lambda=(\Lambda_{0},(\Lambda_{1}(\sigma)\geq \cdots \geq \Lambda_{r_{\sigma}}(\sigma),\Lambda_{r_{\sigma}+1}(\sigma)\geq \cdots \geq \Lambda_{n}(\sigma))_{\sigma\in\Sigma})\in \Lambda(K_{I,\infty})$, we define $\Lambda^{*}:=(-\Lambda_{0},(-\Lambda_{r_{\sigma}}(\sigma)\geq \cdots \geq -\Lambda_{1}(\sigma),-\Lambda_{n}\geq\cdots\geq -\Lambda_{r_{\sigma}+1})_{\sigma\in\Sigma})$. 

We know $V_{\Lambda^{*}}$ is the dual of $V_{\Lambda}$ as $K_{I}$-representation. We sometimes write the latter as $\widecheck{V_{\Lambda}}$.

\bigskip

We define $I^{c}$ by $I^{c}(\sigma)=n-I(\sigma)$ for all $\sigma\in\Sigma$. We know $V_{I^{c}}=-V_{I}$ and $GU_{I^{c}}\cong GU_{I}$. The complex conjugation gives an anti-holomorphic isomorphism $X_{I}\xrightarrow{\sim} X_{I^{c}}$. This induces a $K$-antilinear isomorphism \begin{equation}
\bar{H}^{D/2}(Sh_{I},W_{\lambda}^{\triangledown})\xrightarrow{\sim}  \bar{H}^{D/2}(Sh_{I^{c}},W_{\lambda^{c}}^{\triangledown}).
\end{equation}
In particular, it sends holomorphic (resp. anti-holomorphic) elements with respect to $(I,\lambda)$ to those respect to $(I^{c},\lambda^{c})$. If we we denote by $w_{0}^{c}$ the longest element related to $I^{c}$ then we have $K$-antilinear rational isomorphisms
\begin{eqnarray}
c_{DR}:  &\bar{H}^{0}(Sh_{I},E_{w_{0}*\lambda})\xrightarrow{\sim}  \bar{H}^{0}(Sh_{I^{c}},E_{w_{0}^{c}*\lambda^{c}})&\\
& \bar{H}^{D/2}(Sh_{I},E_{\lambda})\xrightarrow{\sim}  \bar{H}^{D/2}(Sh_{I^{c}},E_{\lambda^{c}})&.
\end{eqnarray} 

\bigskip

 The Shimura datum $(GU_{I},h)$ induces a Hodge structure of wights concentrated in $\{(-1,1),(0,0),(1,-1)\}$ which corresponds to the Harish-Chandra decomposition induced by $h$ on the Lie algebra:
$\lieG=\lieK_{\C}\oplus \lieP^{+}\oplus \lieP^{-}$.\\

Let $\mathfrak{P}=\lieK_{\C}\oplus \lieP^{-}$. Let $\mathcal{A}$ (resp. $\mathcal{A}_{0}$, $\mathcal{A}_{(2)}$) be the space of automorphic forms (resp. cusp forms, square-integrable forms) on $GU_{I}(\Q)\backslash GU_{I}(\AQ)$.  \\

We have inclusions for all $q$:
\begin{eqnarray}\nonumber   H^{q}(\lieG,K_{I,\infty};\mathcal{A}_{0}\otimes V_{\lambda})\subset \bar{H}^{q}(Sh_{I},V_{\lambda}^{\triangledown})\subset H^{q}(\lieG,K_{I,\infty};\mathcal{A}_{(2)}\otimes V_{\lambda})\\
\nonumber   H^{q}(\mathfrak{P},K_{I,\infty};\mathcal{A}_{0}\otimes V_{\Lambda})\subset \bar{H}^{q}(Sh_{I},E_{\Lambda})\subset H^{q}(\mathfrak{P},K_{I,\infty};\mathcal{A}_{(2)}\otimes V_{\Lambda}).
\end{eqnarray}

The complex conjugation on the automorphic forms induces a $K$-antilinear isomorphism:
\begin{equation} \label{definition cB}
c_{B}: \bar{H}^{0}(Sh_{I},E_{w_{0}*\lambda}) \xrightarrow{\sim}   \bar{H}^{D/2}(Sh_{I},E_{\lambda^{\vee}})
\end{equation}

More precisely, we summarize the construction in \cite{harris97} as follows.

\paragraph{Automorphic vector bundles:}
\text{}

We recall some facts on automorphic vector bundles first. We refer to page $113$ of \cite{harris97} and \cite{harrisarithmeticvectorbundle1} for notation and further details.\\

Let $(G,X)$ be a Shimura datum such that its special points are all CM points. Let $\widecheck{X}$ be the compact dual symmetric space of $X$. There is a surjective functor from the category of $G$-homogeneous vector bundles on $\widecheck{X}$ to the category of automorphic vector bundles on $Sh(G,X)$. This functor is compatible with inclusions of Shimura data as explained in the second part of Theorem $4.8$ of \cite{harrisarithmeticvectorbundle1}. It is also rational over the reflex field $E(G,X)$. \\

Let $\mathcal{E}$ be an automorphic vector bundle on $Sh(G,X)$ corresponding to $\mathcal{E}_{0}$, a $G$-homogeneous vector bundle on $\widecheck{X}$. Let $(T,x)$ be a special pair of $(G,X)$, i.e. $(T,x)$ is a sub-Shimura datum of $(G,X)$ with $T$ a maximal torus defined over $\Q$. Since the functor mentioned above is compatible with inclusions of Shimura datum, we know that the restriction of $\mathcal{E}$ to $Sh(T,x)$ corresponds to the restriction of $\mathcal{E}_{0}$ to $\check{x}\in \widecheck{X}$ by the previous functor. Moreover, by the construction, the fiber of $\mathcal{E}\mid_{Sh(T,x)}$ at any point of $Sh(T,x)$ is identified with the fiber of $\mathcal{E}_{0}$ at $\check{x}$. The $E(\mathcal{E})\cdot E(T,x)$-rational structure on the fiber of $\mathcal{E}_{0}$ at $\check{x}$ then defines a rational structure of $\mathcal{E}\mid_{Sh(T,x)}$ and called the \textbf{canonical trivialization} of $\mathcal{E}$ associated to $(T,x)$.

\bigskip

\paragraph{Complex conjugation on automorphic vector bundles:}
\text{}

Let $\mathcal{E}$ be as in page $112$ of \cite{harris97} and $\overline{\mathcal{E}}$ be its complex conjugation. The key step of the construction is to identify $\overline{\mathcal{E}}$ with the dual of $\mathcal{E}$ in a rational way. \\

More precisely, we recall Proposition $2.5.8$ of the \textit{loc.cit} that there exists a non-degenerate $G(\AQf)$-equivariant paring of real-analytic vector bundle $\mathcal{E}\otimes \overline{\mathcal{E}}\rightarrow \mathcal{E}_{\nu}$ such that its pullback to any CM point is rational with respect to the canonical trivializations. \\

We now explain the notion $\mathcal{E}_{\nu}$. Let $h\in X$ and $K_{h}$ be the stabilizer of $h$ in $G(\R)$. We know $\mathcal{E}$ is associated to an irreducible complex representation of $K_{h}$, denoted by $\tau$ in the \textit{loc.cit}. The complex conjugation of $\tau$ can be extended as an algrebraic representation of $K_{h}$, denoted by $\tau'$. We know $\tau'$ is isomorphic to the dual of $\tau$ and then there exists $\nu$, a one-dimensional representation $K_{h}$, such that a $K_{h}$-equivariant rational paring $V_{\tau}\otimes V_{\tau'}\rightarrow V_{\nu}$ exists. We denote by $\mathcal{E}_{\nu}$
the automorphic vector bundle associated to $V_{\nu}$.

\bigskip

In our case, we have $(G,X)=(GU_{I},X_{I})$, $h=h_{I}$ and $K_{h}=K_{I,\infty}$. Let $\tau=\Lambda=w_{0}*\lambda$ and $\mathcal{E}=E_{\Lambda}$. As explained in the last second paragraph before Corollary $2.5.9$ in the \textit{loc.cit}, we may identify the holomorphic sections of $V_{\Lambda}$ with holomorphic sections of the dual of $\overline{V_{\Lambda}}$. The complex conjugation then sends the latter to the anti-holomorphic sections of $\widecheck{V_{\Lambda}}=V_{\Lambda^{*}}$. The latter can be identified with harmonic (0,d)-forms  with values in $\mathbb{K}\otimes E_{\Lambda^{\*}}$ where $\mathbb{K}=\Omega^{D/2}_{Sh_{I}}$ is the canonical line bundle of $Sh_{I}$. \\

By $2.2.9$ of \cite{harris97} we have $\mathbb{K}=E_{(0,(-s_{\sigma},\cdots,-s_{\sigma},r_{\sigma},\cdots,r_{\sigma})_{\sigma\in\Sigma})}$ where the number of $-s_{\sigma}$ in the last term is $r_{\sigma}$. Therefore, complex conjugation gives an isomorphism:
\begin{equation}
c_{B}: \bar{H}^{0}(Sh_{I},E_{\Lambda}) \xrightarrow{\sim}   \bar{H}^{D/2}(Sh_{I},E_{\Lambda^{*}+0,(-s_{\sigma},\cdots,-s_{\sigma},r_{\sigma},\cdots,r_{\sigma})_{\sigma\in\Sigma})}).
\end{equation}

Recall equation (\ref{twist by w0}) that 
\begin{equation}\nonumber\Lambda=w_{0}*\lambda=(\lambda_{0},(\lambda_{s_{\sigma}+1}(\sigma)-s_{\sigma},\cdots,\lambda_{n}(\sigma)-s_{\sigma};\lambda_{1}(\sigma)+r_{\sigma},\cdots,\lambda_{s_{\sigma}}(\sigma)+r_{\sigma})_{\sigma\in\Sigma}).
\end{equation} 

We have \begin{equation}
\Lambda^{*}=(-\lambda_{0},(-\lambda_{n}(\sigma)+s_{\sigma},\cdots,-\lambda_{s_{\sigma}+1}(\sigma)+s_{\sigma};-\lambda_{s_{\sigma}}(\sigma)-r_{\sigma},\cdots,-\lambda_{1}(\sigma)+r_{\sigma})_{\sigma\in\Sigma}).
\end{equation}

Therefore, $\Lambda^{*}+(0,(-s_{\sigma},\cdots,-s_{\sigma},r_{\sigma},\cdots,r_{\sigma})_{\sigma\in\Sigma})=\lambda^{\vee}$. We finally get equation (\ref{definition cB}).

\bigskip

Similarly, if we start from the anti-holomorphic part, we will get a $K$-antilinear isomorphism which is still denoted by $c_{B}$:
\begin{equation} 
c_{B}: \bar{H}^{D/2}(Sh_{I},E_{\lambda}) \xrightarrow{\sim}   \bar{H}^{0}(Sh_{I},E_{w_{0}*\lambda^{\vee}})
\end{equation} which sends anti-holomorphic elements with respect to $\lambda$ to holomorphic elements for $\lambda^{\vee}$.

\section{The rational paring}
Let $\Lambda\in \Lambda(K_{I,\infty})$. We write $V=V_{\Lambda}$ in this section for simplicity.  As in section $2.6.11$ of \cite{harris97}, we denote by $\C_{\Lambda}$ the corresponding highest weight space. We know $\Lambda^{*}:=\Lambda^{\#}-(2\Lambda_{0},(0))$ is the tuple associated to $\widecheck{V}$, the dual of this $K_{I}$ representation. We denote by $\C_{-\Lambda}$ the lowest weight of $\widecheck{V}$.\\

The restriction from $V$ to $\C_{\Lambda}$ gives an isomorphism
\begin{equation}
Hom_{K_{I,\infty}}(V,\mathcal{C}^{\infty}(GU_{I}(\F)\backslash GU_{I}(\AF))) \xrightarrow{\sim} Hom_{H}(\C_{\Lambda},\mathcal{C}^{\infty}(GU_{I}(\F)\backslash GU_{I}(\AF))_{V}) 
\end{equation}
where $\mathcal{C}^{\infty}(GU_{I}(\F)\backslash GU_{I}(\AF))_{V}$ is the $V$-isotypic subspace of $\mathcal{C}^{\infty}(GU_{I}(\F)\backslash GU_{I}(\AF))$.

Similarly, we have \begin{equation}\label{trivialization}
Hom_{K_{I,\infty}}(\widecheck{V},\mathcal{C}^{\infty}(GU_{I}(\F)\backslash GU_{I}(\AF))) \xrightarrow{\sim} Hom_{H}(\C_{-\Lambda},\mathcal{C}^{\infty}(GU_{I}(\F)\backslash GU_{I}(\AF))_{\widecheck{V}}) 
\end{equation}

Proposition $2.6.12$ of \cite{harris97} says that up to a rational factor the perfect paring
\begin{equation}
Hom_{H}(\C_{\Lambda},\mathcal{C}^{\infty}(GU_{I}(\F)\backslash GU_{I}(\AF))_{V}) \times Hom_{H}(\C_{-\Lambda},\mathcal{C}^{\infty}(GU_{I}(\F)\backslash GU_{I}(\AF))_{\widecheck{V}}) 
\end{equation}
given by integration over the diagonal equals to restriction of the canonical paring (c.f. $(2.6.11.4)$ of \cite{harris97})
 \begin{eqnarray}
 &&Hom_{K_{I,\infty}}(V,\mathcal{C}^{\infty}(GU_{I}(\F)\backslash GU_{I}(\AF))) \times Hom_{K_{I,\infty}}(\widecheck{V},\mathcal{C}^{\infty}(GU_{I}(\F)\backslash GU_{I}(\AF))) 
\nonumber\\
&\rightarrow& Hom_{K_{I,\infty}}(V\otimes \widecheck{V},\mathcal{C}^{\infty}(GU_{I}(\F)\backslash GU_{I}(\AF)))\nonumber\\
&\rightarrow& Hom_{K_{I,\infty}}(\C,\mathcal{C}^{\infty}(GU_{I}(\F)\backslash GU_{I}(\AF)))\nonumber\\
&\rightarrow& \C.
 \end{eqnarray}

 \bigskip

 We identify $\Gamma^{\infty}(Sh_{I},E_{\Lambda})$ with $Hom_{GU_{I}K_{I,\infty}}(\widecheck{V},\mathcal{C}^{\infty}(GU_{I}(\F)\backslash GU_{I}(\AF)))$ and  regard the latter as subspace of $Hom_{K_{I,\infty}}(\widecheck{V},\mathcal{C}^{\infty}(GU_{I}(\F)\backslash GU_{I}(\AF)))$.\\
 
 The above construction gives a $K$-rational perfect paring between holomorphic sections of $E_{\Lambda}$ and anti-holomorphic sections of $E_{\Lambda^{*}}$.\\
  
 If $\Lambda=w_{0}*\lambda$, as we have seen in Section \ref{Complex conjugation} that the anti-holomorphic sections of $E_{\Lambda^{*}}$ can be identified with harmonic $(0,d)$-forms with values in $E_{\lambda^{\vee}}$.\\
 
 We therefore obtain a $K$-rational perfect paring 
 \begin{equation}
 \Phi=\Phi^{I,\lambda}: \bar{H}^{0}(Sh_{I},E_{w_{0}*\lambda}) \times \bar{H}^{D/2}(Sh_{I},E_{\lambda^{\vee}})\rightarrow \C. 
 \end{equation}
 In other words, there is a rational paring between the holomorphic elements for $(I,\lambda)$ and anti-holomorphic elements for $(I,\lambda^{\vee})$.
  
  \bigskip
 
It is easy to see that the isomorphism $Sh_{I}\xrightarrow{\sim} Sh_{I^{c}}$ commutes with the above paring and hence:

 \begin{lem}\label{period stable under complex conjugation}
For any $f\in  \bar{H}^{0}(Sh_{I},E_{w_{0}*\lambda}) $ and $g\in \bar{H}^{D/2}(Sh_{I},E_{\lambda^{\vee}})$, we have $\Phi^{I,\lambda}(f,g)=\Phi^{I^{c},\lambda^{c}}(c_{DR}f,c_{DR}g)$.
 \end{lem}
  
  The next lemma follows from Corollary $2.5.9$ and Lemma $2.8.8$ of \cite{harris97}.
  
 \begin{lem}
 Let $0\neq f\in \bar{H}^{0}(Sh_{I},E_{w_{0}*\lambda})$. We have
 $\Phi(f,c_{B}f)\neq 0$.
 
 More precisely, if we consider $f$ as an element in $Hom_{K_{I,\infty}}(\widecheck{V},\mathcal{C}^{\infty}(GU_{I}(\F)\backslash GU_{I}(\AF)))$ then by (\ref{trivialization}) and the fixed trivialization of $\C_{-w_{0}*\lambda}$, we may consider $f$ as an element in $\mathcal{C}^{\infty}(GU_{I}(\F)\backslash GU_{I}(\AF)))$. We have:
 \begin{equation}
 \Phi(f,c_{B}f) =\pm i^{\lambda_{0}}\int\limits_{GU_{I}(\Q)Z_{GU_{I}}(\AQ)\backslash GU_{I}(\AQ)} f(g)\overline{f}(g)||\nu(g)||^{c}dg.
 \end{equation}
 Recall that $\nu(\cdot)$ is the similitude defined in (\ref{definition of nu(g)}).
 \end{lem}

\bigskip
Similarly, if we start from anti-holomorphic elements, we get a paring:
\begin{equation}\Phi^{-}=\Phi^{I,\lambda,-}: \bar{H}^{D/2}(Sh_{I},E_{\lambda}) \times \bar{H}^{0}(Sh_{I},E_{w_{0}*\lambda^{\vee}})\rightarrow \C. \end{equation}
We use the script $\text{}^{-}$ to indicate that is anti-holomorphic. It is still $c_{DR}$ stable. For $0\neq f^{-} \in \bar{H}^{D/2}(Sh_{I},E_{\lambda}) $, we also know that $\Phi^{-}(f^{-},c_{B}f^{-})\neq 0$.

\section{Arithmetic automorphic periods}

Let $\pi$ be an irreducible cuspidal representation of $GU_{I}(\AQ)$ defined over a number field $E(\pi)$. We may assume that $E(\pi)$ contains the quadratic imaginary field $K$.\\

We assume that $\pi$ is cohomological with type $\lambda$, i.e. $H^{*}(\lieG,K_{I,\infty};\pi\otimes W_{\lambda})\neq 0$.\\

For $M$ a $GU_{I}(\AQf)$-module, define the $K$-rational $\pi_{f}$-isotypic components of $M$ by
\begin{equation}\nonumber
M^{\pi}:=Hom_{GU_{I}(\AFf)}(Res_{E(\pi)/K}(\pi_{f}),M)=\bigoplus\limits_{\tau \in \Sigma_{E(\pi);K}} Hom(\pi_{f}^{\tau},M).
\end{equation} 

Therefore, if $M$ has a $K$-rational structure then $M^{\pi}$ also has a $K$-rational structure.\\

As in section \ref{Complex conjugation}, we have inclusions:
\begin{equation}\nonumber    H^{q}(\mathfrak{P},K_{I,\infty};\mathcal{A}_{0}^{\pi}\otimes V_{\Lambda})\subset \bar{H}^{q}(Sh_{I},E_{\Lambda})^{\pi}\subset H^{q}(\mathfrak{P},K_{I,\infty};\mathcal{A}_{(2)}^{\pi}\otimes V_{\Lambda}).\end{equation}

Under these inclusions, $c_{B}$ sends $\bar{H}^{0}(Sh_{I},E_{w_{0}*\lambda})^{\pi}$ to $\bar{H}^{D/2}(Sh_{I},E_{\lambda^{\vee}})^{\pi^{\vee}}$.

These inclusions are compatible with those $K$-rational structures and then induce $K$-rational parings
\begin{eqnarray}\Phi^{\pi}: \bar{H}^{0}(Sh_{I},E_{w_{0}*\lambda})^{\pi} \times \bar{H}^{D/2}(Sh_{I},E_{\lambda^{\vee}})^{\pi^{\vee}}\rightarrow \C\\
 \text{and }
\Phi^{-,\pi}: \bar{H}^{D/2}(Sh_{I},E_{\lambda})^{\pi} \times \bar{H}^{0}(Sh_{I},E_{w_{0}*\lambda^{\vee}})^{\pi^{\vee}}\rightarrow \C. \end{eqnarray}

\bigskip

\begin{df}
Let $\beta$ be a non zero $K$-rational element of $\bar{H}^{0}(Sh_{I},E_{w_{0}*\lambda})^{\pi}$. We define the \textbf{holomorphic arithmetic automorphic period associated to } $\beta$ by $P^{(I)}(\beta,\pi):=(\Phi^{\pi(}\beta^{\tau},c_{B}\beta^{\tau}))_{\tau\in\Sigma_{E(\pi);K}}$. It is an element in $(E(\pi)\otimes_{K}\C)^{\times}$.\\

Let $\gamma$ be a non zero $K$-rational element of $\bar{H}^{D/2}(Sh_{I},E_{\lambda})^{\pi}$. We define the \textbf{anti-holomorphic arithmetic automorphic period associated to} $\gamma$ by $P^{(I),-}(\gamma,\pi):=(\Phi^{-,\pi}(\gamma^{\tau},c_{B}\gamma^{\tau}))_{\tau\in\Sigma_{E(\pi);K}}$. It is an element in $(E(\pi)\otimes_{K}\C)^{\times}$.

\end{df}

\bigskip

\begin{dflem}
Let us assume now $\pi$ is tempered and $\pi_{\infty}$ is discrete series representation. In this case, $\bar{H}^{0}(Sh_{I},E_{w_{0}*\lambda})^{\pi}$ is a rank one $E(\pi)\otimes_{K} \C$-module (c.f. \cite{endoscopicfour}). \\

We define the \textbf{holomorphic arithmetic automorphic period} of $\pi$ by $P^{(I)}(\pi):= P^{(I)}(\beta,\pi)$ by taking $\beta$ any non zero rational element in $\bar{H}^{0}(Sh_{I},E_{w_{0}*\lambda})^{\pi}$. It is an element in $(E(\pi)\otimes_{K}\C)^{\times}$ well defined up to $E(\pi)^{\times}$-multiplication.\\

We define $P^{(I),-}(\pi)$ the \textbf{anti-holomorphic arithmetic automorphic period} of $\pi$ similarly.
\end{dflem} 

\bigskip

\begin{lem}\label{pair to ratio}
We assume that $\pi$ is tempered and $\pi_{\infty}$ is discrete series representation. Let $\beta$ be a non zero rational element in $\bar{H}^{0}(Sh_{I},E_{w_{0}*\lambda})^{\pi}$ and $\beta^{\vee}$ be a non zero rational element in $\bar{H}^{0}(Sh_{I},E_{\lambda}^{\vee})^{\pi^{\vee}}$. 

We have $c_{B}(\beta)\sim_{E(\pi)} P^{(I)}(\pi) \beta^{\vee}$.
\end{lem}
\begin{dem}
It is enough to notice that $\Phi^{\pi}(\beta,\beta^{\vee})\in E(\pi)^{\times}$.
\end{dem}

\begin{lem}\label{lemma inverse}
If $\pi$ is tempered and $\pi_{\infty}$ is discrete series representation then we have:
\begin{enumerate}
\item $P^{(I^{c})}(\pi^{c})\sim_{E(\pi);K} P^{(I)}(\pi)$.
\item $P^{(I)}(\pi^{\vee})*P^{(I),-}(\pi)\sim_{E(\pi);K}1$.
\end{enumerate}
\end{lem}
\begin{dem}
The first part comes from Lemma \ref{period stable under complex conjugation} and the fact that $c_{DR}$ preserves rational structures. 

For the second part, recall that the following two parings are actually the same:
\begin{eqnarray}\Phi^{\pi^{\vee}}: \bar{H}^{0}(Sh_{I},E_{w_{0}*\lambda^{\vee}})^{\pi^{\vee}} \times \bar{H}^{D/2}(Sh_{I},E_{\lambda})^{\pi}\rightarrow \C\\
 \text{and }
\Phi^{-,\pi}: \bar{H}^{D/2}(Sh_{I},E_{\lambda})^{\pi} \times \bar{H}^{0}(Sh_{I},E_{w_{0}*\lambda^{\vee}})^{\pi^{\vee}}\rightarrow \C. \end{eqnarray}

We take $\beta$ a rational element in $\bar{H}^{0}(Sh_{I},E_{w_{0}*\lambda^{\vee}})^{\pi^{\vee}}$ and $\gamma$ a rational element in $\bar{H}^{D/2}(Sh_{I},E_{\lambda})^{\pi}$. We may assume that $\Phi^{\pi^{\vee}}(\beta^{\tau},\gamma^{\tau})=\Phi^{-,\pi}(\gamma^{\tau},\beta^{\tau})=1$ for all $\tau\in \Sigma_{E(\pi);K}$.

By definition $p^{(I)}(\pi^{\vee})=(\Phi^{\pi^{\vee}}(\beta^{\tau},c_{B}\beta^{\tau}))_{\tau\in\Sigma_{E(\pi);K}}$. Since $\bar{H}^{D/2}(Sh_{I},E_{\lambda})^{\pi}$ is a rank one $E(\pi)\otimes \C$-module, there exists $C \in (E(\pi)\otimes_{K} \C)^{\times}$ such that $(c_{B}\beta^{\tau})_{\tau\in\Sigma_{E(\pi);K}}=C (\gamma^{\tau})_{\tau\in\Sigma_{E(\pi);K}}$. Therefore $p^{I}(\pi^{\vee})=C (\Phi^{\pi^{\vee}}(\beta^{\tau},\gamma^{\tau}))_{\tau\in\Sigma_{E(\pi);K}}=C$.

On the other hand, since $c_{B}^{2}=Id$, we have $(c_{B}\gamma^{\tau})_{\tau\in\Sigma_{E(\pi);K}}=C^{-1}(\beta^{\tau})_{\tau\in\Sigma_{E(\pi);K}}$. We can deduce that $p^{(I),-}(\pi)=C^{-1}$ as expected.
\end{dem}
\begin{df}
We say $I$ is \textbf{compact} if $U_{I}(\C)$ is. In other words, $I$ is compact if and only if $I(\sigma)=0$ or $n$ for all $\sigma\in\Sigma$.
\end{df}

\begin{cor}\label{unitary similitude product 1}
If $I$ is compact then $P^{(I)}(\pi) \sim_{E(\pi);K}P^{(I),-}(\pi)$.
We have $P^{(I)}(\pi^{\vee})*P^{(I)}(\pi)\sim_{E(\pi);K}1$.
\end{cor}
\begin{dem}
If $I$ is compact, then $w_{0}=Id$. The anti-holomorphic part and holomorphic part are the same. We then have $P^{(I)}(\pi) \sim_{E(\pi);K}P^{(I),-}(\pi)$. The last assertion comes from Lemma \ref{lemma inverse}.
\end{dem}

\chapter{Special values of automorphic $L$-functions: the start point}
\section{Special values of automorphic $L$-functions for similitude unitary group}
The method in \cite{harris97} should work for a general CM field. We state the predicted formula in this section.

Let $\pi$ be a  tempered representation of $GU_{I}(\AQ)$ such that $\pi_{\infty}$ is discrete series representation. In particular, the holomorphic arithmetic automorphic periods $P^{(I)}(\pi)$ is well defined.

We assume that $\pi$ is cohomological with type $\lambda=(\lambda_{0},(\lambda_{1}(\sigma)\geq \lambda_{2}(\sigma)\cdots\geq \lambda_{n}(\sigma))_{\sigma\in\Sigma})$. 

We say $\lambda$ or $\pi$ is $2$-regular if $\lambda_{i}(\sigma)-\lambda_{i+1}(\sigma)\geq 1$ for all $1\leq i\leq n-1$ and all $\sigma\in\Sigma$.

Let $\chi$ be an algebraic Hecke character of $\AF^{\times}$ with infinity type $(z^{-k(\sigma)})_{\sigma\in\Sigma}$. Let $\alpha$ be an algebraic Hecke character of $\AF^{\times}$ with infinity type $(z^{\kappa})_{\sigma\in\Sigma}$.

\begin{df}\label{motivic critical triple}
We set $\lambda_{0}(\sigma)=+\infty$ and $\lambda_{n+1}(\sigma)=-\infty$. We say $m\in \Z$ is critical for $M(\pi,\chi,\alpha)$ (c.f. \cite{harris97}) if  for all $\sigma\in\Sigma$,
\begin{eqnarray}
\nonumber \lambda_{r_{\sigma}+1}(\sigma)+k(\sigma)+s_{\sigma}-\kappa\leq m\leq -\lambda_{r_{\sigma}+1}(\sigma)-k(\sigma)+r_{\sigma}\\
\text{ and } -\lambda_{r_{\sigma}}(\sigma)-k(\sigma)+r_{\sigma}\leq m\leq \lambda_{r_{\sigma}}(\sigma)+k(\sigma)+s_{\sigma}-\kappa.\nonumber
\end{eqnarray}
\end{df}

\bigskip

This definition generalize the condition in Lemma $3.3.7$ of the \textit{loc.cit}. In the \textit{loc.cit}, it is assumed that $\mu$ is self-dual. In general, the index $\Lambda$ in that Lemma should be $\Lambda(\mu^{c};r,s)$.\\

We assume the following conjecture throughout the text.
\begin{conj}\label{special value for similitude unitary group}
Let $\pi$ be a  tempered representation of $GU_{I}(\AQ)$ such that $\pi_{\infty}$ is discrete series representation and cohomological with type $\lambda$. In particular, the holomorphic arithmetic automorphic periods $P^{(I)}(\pi)$ is well defined. If an integer $m\geq \cfrac{n-\kappa}{2}$ is critical for $M(\pi,\chi,\alpha)$ then a certain rational differential operator exists and we have
\begin{eqnarray}\nonumber
&L^{mot,S}(m,\pi\otimes\chi,St,\alpha)\sim_{E(\pi)E(\chi)E(\alpha);K} &\\&
(2\pi i)^{(m-\frac{n-1}{2})nd}(2\pi)^{-\lambda_{0}} P^{(I)}(\pi)\prod\limits_{\sigma\in\Sigma} (p(\widetilde{\chi}\alpha,\sigma)^{-I(\sigma)}p(\widetilde{\chi}\alpha,\overline{\sigma})^{-n+I(\sigma)})&.
\end{eqnarray}
Recall that $\widetilde{\chi}=\cfrac{\chi}{\chi^{c}}$ is a Hecke character of $\AF$.
\end{conj}

\begin{rem}
\begin{enumerate}
\item If $m\geq \cfrac{n-\kappa}{2}$ is critical for $M(\pi,\chi,\alpha)$, we have $2\lambda_{r_{\sigma}+1}-2(r_{\sigma}+1)\leq \kappa-2k(\sigma)-n-2$ and $2\lambda_{r_{\sigma}}-2r_{\sigma}\geq \kappa-2k(\sigma)-n$. We see that $r_{\sigma}=\max\{r\mid 2\lambda_{n'}-2r\geq \kappa-2k(\sigma)-n\}$. 

\item We didn't state the CM periods in the above conjecture as in Theorem $3.5.13$ of \cite{harris97}. Instead, the current form appears in middle steps of the proof for Theorem $3.5.13$. We refer to equation $(2.9.12)$ or the third line in page $138$ of the \textit{loc.cit}.
\end{enumerate}
\end{rem}

\bigskip

Let us examine the condition in the above conjecture. After simple calculation, we see that such $m$ always exists. Moreover, if $\lambda_{r_{\sigma}}>\lambda_{r_{\sigma}+1}$, we may have $m\geq \cfrac{n-\kappa+1}{2}$. In this case, we know $L^{mot,S}(m,\pi\otimes\chi,St,\alpha)$ does not vanish.\\

Let $GU$ and $GU'$ be two rational similitude group associated to two unitary groups over $F$ with respect to $F/F^{+}$ of dimension $n$. We know $GU'$ is an inner form of $GU$ and thus they are isomorphic to each other at almost all primes. \\

Let $\pi$ and $\pi'$ be automorphic representations of $GU(\AQ)$ and $GU'(\AQ)$ respectively. We say $\pi$ is \textbf{nearly equivalent} to $\pi'$ if they are isomorphic to each other at almost all primes. In particular, they have the same value of partial $L$-functions. \\

We then deduce that:
 
\begin{cor}
The arithmetic automorphic periods $P^{(I)}(\pi)$ and $P^{(I),-}(\pi)$ depends only on the nearly equivalence class of $\pi$ if $\pi$ is $2$-regular.
\end{cor}

 \section{Special values of automorphic $L$-functions for $GL_{n}\times GL_{1}$}
 
 Let $\Pi$ be a cuspidal automorphic representation of $GL_{n}(\AF)$ which is regular conjugate self-dual, cohomological. We assume moreover that $\Pi_{q}$ descends locally if $n$ is even. By Lemma \ref{existence of xi}, there exists an Hecke character of $K$, denoted by $\xi$, such that $\Pi^{\vee} \otimes\xi$ is $\theta_{I}$-stable. By Proposition \ref{base change for similitude unitary group}, we know $\Pi^{\vee}\otimes \xi$ descends to $\pi$, an automorphic cohomological representation of $GU_{I}(\AQ)$ with is tempered and discrete series at the infinity place. In particular, the arithmetic automorphic period of $\pi$ can be defined.
 
\bigskip

 \begin{df}\label{arithmetic period conjugate self-dual case}
 We fix one Hecke character $\xi$ as above. We denote its infinity type by $z^{u}\overline{z}^{v}$. We define the \textbf{arithmetic automorphic period} for $\Pi$ by $P^{(I)}(\Pi,\xi):=(2\pi)^{u+v}P^{(I)}(\pi)$. 
 
If $\Pi$ is $2$-regular then $P^{(I)}(\Pi,\xi)$ does not depend on the choice of $\xi$ up to elements in $E(\Pi)$. This is a corollary of Theorem \ref{main theorem CM}. Therefore we may define $P^{(I)}(\Pi):=P^{(I)}(\Pi,\xi)$ for any fixed $\xi$ in this case.
  \end{df}
 
\bigskip

\begin{lem}
Let $\Pi$ be as in the above theorem. 
\begin{enumerate}
\item For all $I$ we have $P^{(I)}(\Pi)\sim_{E(\Pi);K} P^{(I^{c})}(\Pi^{c})$.
\item If $I$ is compact then $P^{(I)}(\Pi^{\vee})*P^{(I)}(\Pi)\sim_{E(\Pi);K} 1$. 
In particular, since $\Pi$ is conjugate self-dual, we have $P^{(I^{c})}(\Pi)*P^{(I)}(\Pi)\sim_{E(\Pi);K} 1$ in this case.
\end{enumerate}
\end{lem}
\begin{dem}
The first part comes from the first part of Lemma \ref{lemma inverse}. It is also a direct corollary from Theorem \ref{main theorem CM} below.

The second one follows from Corollary \ref{unitary similitude product 1}. Let $\xi$ be the auxiliary Hecke character with infinity type $z^{u}\overline{z}^{v}$ as above and $\pi$ be the representation of $GU_{I}(\AQ)$ with base change $\Pi^{\vee}\otimes \xi$. We know $\pi^{\vee}$ is a representation of $GU_{I}(\AQ)$ with base change $\Pi\otimes \xi^{-1}$.

By definition $P^{(I)}(\Pi)\sim_{E(\Pi);K} (2\pi)^{u+v}P^{(I)}(\pi)$ and $P^{(I)}(\Pi^{\vee})\sim_{E(\Pi);K} (2\pi)^{-u-v}P^{(I)}(\pi^{\vee})$. Thus $P^{(I)}(\Pi^{\vee})*P^{(I)}(\Pi)\sim_{E(\Pi);K} P^{(I)}(\pi^{\vee})*P^{(I)}(\pi)\sim_{E(\Pi);K} 1$.

\end{dem}

\paragraph{Critical points:}
 Let $n$, $n'$ be two integers. Let $\Pi$ and $\Pi'$ be algebraic automorphic representations of $GL_{n}(\AF)$ and $GL_{n'}(\AF)$ with pure infinity type $(z^{a_{i}(\sigma)}\overline{z}^{-\omega(\Pi)-a_{i}(\sigma)})_{1\leq i\leq n}$ and $(z^{a'_{j}(\sigma)}\overline{z}^{-\omega(\Pi)-a'_{j}(\sigma)})_{1\leq j\leq n'}$ respectively. \\
 
 We assume the existence of motives $M(\Pi)$ and $M(\Pi')$ associated to $\Pi$ and $\Pi'$. Let $m\in\Z+\cfrac{n+n'}{2}$. We say $m$ is \textbf{critical} for $\Pi\otimes \Pi'$ if $m+\cfrac{n+n'-2}{2}$ is critical for $M(\Pi)\otimes M(\Pi')$ in the sense of Deligne (c.f. \cite{deligne79} or Chapter \ref{chapter motive}). \\
 
 If $a_{i}(\sigma)+a'_{j}(\sigma)\neq \cfrac{-\omega(\Pi)-\omega(\Pi')}{2}$ for all $i$, $j$ and $\sigma$ then critical points always exist. In this case, we have an explicit description for them (c.f. $(1.3.1)$ of \cite{deligne79}). More precisely, $m$ is critical if and only if for all $i$, $j$, $\sigma$, if $-a_{i}(\sigma)-a'_{j}(\sigma)>\cfrac{\omega(\Pi)+\omega(\Pi')}{2}$ then $\omega(\Pi)+\omega(\Pi')+1+a_{i}(\sigma)+a'_{j}(\sigma)\leq m\leq -a_{i}(\sigma)-a'_{j}(\sigma)$, if $-a_{i}(\sigma)-a'_{j}(\sigma)<\cfrac{\omega(\Pi)+\omega(\Pi')}{2}$ then $1-a_{i}(\sigma)-a'_{j}(\sigma)\leq m\leq \omega(\Pi)+\omega(\Pi')+a_{i}(\sigma)+a'_{j}(\sigma)$. Roughly speaking, $m$ should be closer to the central point than any of the $a_{i}(\sigma)+a'_{j}(\sigma)$.\\
 
If $a_{i}(\sigma)+a'_{j}(\sigma)= \cfrac{-\omega(\Pi)-\omega(\Pi')}{2}$ for some $i$, $j$ and $\sigma$ then there is no critical points (c.f. Lemma $1.7.1$ of \cite{harris97}).

\bigskip

The following theorem follows directly from Conjecture \ref{special value for similitude unitary group}.

 \begin{thm}\label{main theorem CM}\label{n*1}
 Let us assume that Conjecture \ref{special value for similitude unitary group} is true. Let $\Pi$ be a regular, conjugate self-dual, cohomological, cuspidal automorphic representation of $GL_{n}(\AF)$ which descends to $U_{I}(\AFp)$ for any $I$ (c.f. Proposition \ref{base change for unitary group}
). We denote the infinity type of $\Pi$ at $\sigma\in \Sigma$ by $(z^{a_{i}(\sigma)}\overline{z}^{-a_{i}(\sigma)})_{1\leq i\leq n}$. 

 Let $\eta$ be an algebraic Hecke character of $F$ with infinity type $z^{a(\sigma)}\overline{z}^{b(\sigma)}$ at $\sigma\in \Sigma$. We know that $a(\sigma)+b(\sigma)$ is a constant independent of $\sigma$, denoted by $-\omega(\eta)$.
 
 We suppose that $a(\sigma)-b(\sigma)+2a_{i}(\sigma)\neq 0$ for all $1\leq i\leq n$ and $\sigma\in \Sigma$. We define $I:=I(\Pi,\eta)$ to be the map on $\Sigma$ which sends $\sigma\in\Sigma$ to $I(\sigma):=\#\{i:a(\sigma)-b(\sigma)+2a_{i}(\sigma)<0\}$. 
 
Let $m\in \Z+\cfrac{n-1}{2}$. If $m\geq \cfrac{1+\omega(\eta)}{2}$ is critical for $\Pi\otimes \eta$, we have:

\begin{equation}\label{main result CM}
L(m,\Pi\otimes \eta) \sim_{E(\Pi)E(\eta);K} (2 \pi i)^{mnd} P^{(I(\Pi,\eta))}(\Pi) \prod\limits_{\sigma\in\Sigma}p(\widecheck{\eta},\sigma)^{I(\sigma)}p(\widecheck{\eta},\overline{\sigma})^{n-I(\sigma)}.
\end{equation}
 \end{thm}

\bigskip

 \begin{lem}
Let $\chi$, $\alpha$ be as in Conjecture \ref{special value for similitude unitary group}. Let $\pi$ be a representation of $GU_{I}$ with base change $\Pi^{c}\times\xi$ for certain auxiliary $\xi$. We set $\eta^{c}=\widetilde{\chi}\alpha$. Let $m\in\Z+\cfrac{n+1}{2}$. Then $m$ is critical for $\Pi\otimes\eta$ if and only if $m+\cfrac{n-1}{2}$ is critical for $M(\pi,\chi,\alpha)$.
\end{lem}

\begin{dem}

Since $\eta^{c}=\widetilde{\chi}\alpha$, we have $b(\sigma)=-k(\sigma)+\kappa$ and $a(\sigma)=k(\sigma)$. Note that $-\omega(\eta)=a(\sigma)+b(\sigma)=\kappa$.

We write the cohomology type of $\pi$ by $(\lambda_{0},(\lambda_{1}(\sigma)\geq \cdots\lambda_{n}(\sigma)))$. We order $a_{i}(\sigma)$ in decreasing order. The cohomology type of $\pi$ is then $(-u-v; (a_{1}(\sigma)-\cfrac{n-1}{2}\geq a_{2}(\sigma)-\cfrac{n-3}{2}\geq \cdots\geq a_{n}(\sigma)+\cfrac{n-1}{2})_{\sigma\in\Sigma})$. This gives $\lambda_{0}=-u-v$ and $\lambda_{i}(\sigma)=a_{i}(\sigma)+i-\cfrac{n+1}{2}$.

Let $m\in \Z+\cfrac{n-1}{2}$. By the above discussion, $m$ is critical for $\Pi\otimes\eta$ if and only if $a_{i}(\sigma)-b(\sigma)+1\leq m\leq -a_{i}(\sigma)-a(\sigma)$ if $a(\sigma)-b(\sigma)+2a_{i}(\sigma)<0$, $1-a_{i}(\sigma)-a(\sigma)\leq m\leq a_{i}(\sigma)-b(\sigma)$ if $a(\sigma)-b(\sigma)+2a_{i}(\sigma)>0$. Since  $r_{\sigma}:=\max\{i:a(\sigma)-b(\sigma)+2a_{i}(\sigma)>0\}$, we deduce that $m$ is critical for $\Pi\otimes \eta$ if and only if $a_{r_{\sigma}+1}(\sigma)-b(\sigma)+1\leq m\leq -a_{r_{\sigma}+1}(\sigma)-a(\sigma)$ and $1-a_{r_{\sigma}}(\sigma)-a(\sigma)\leq m\leq a_{r_{\sigma}}(\sigma)-b(\sigma)$. It is easy to see that these two equations are the same with those in Definition \ref{motivic critical triple}.

\end{dem}

  \paragraph{Proof of Theorem \ref{main theorem CM}:}

We can always choose $\chi$ and $\alpha$ as in Conjecture \ref{special value for similitude unitary group} such that $\eta^{c}=\widetilde{\chi}\alpha$. 

As $m\geq \cfrac{1+\omega(\eta)}{2}$ we have $m+\cfrac{n-1}{2}\geq \cfrac{n+\omega(\eta)}{2}$.
Moreover, the above lemma implies that $m+\cfrac{n-1}{2}$ is critical for $M(\pi,\chi,\alpha)$ and then Conjecture \ref{special value for similitude unitary group} applies, namely:
 \begin{eqnarray}\label{CM proof step 1}
&L(m+\cfrac{n-1}{2},\Pi^{\vee}\otimes \widetilde{\chi}\alpha)\sim_{E(\pi)E(\chi)E(\alpha);K} &\\&\nonumber
(2\pi i)^{mnd}(2\pi)^{-\lambda_{0}} P^{(I(\Pi,\eta))}(\pi)\prod\limits_{\sigma\in\Sigma} (p(\widetilde{\chi}\alpha,\sigma)^{-s_{\sigma}}p(\widetilde{\chi}\alpha,\overline{\sigma})^{-n+s_{\sigma}})&.
\end{eqnarray}

Recall by definition that $P^{(I(\Pi,\eta))}(\Pi,\xi)=(2\pi)^{-\lambda_{0}} P^{(I(\Pi,\eta))}(\pi)$. Moreover, $\eta^{c}=\widetilde{\chi}\alpha$ and then $p(\widetilde{\chi}\alpha,\sigma)\sim p(\eta^{c},\sigma)\sim p(\widecheck{\eta},\sigma)^{-1}$ and $p(\widetilde{\chi}\alpha,\overline{\sigma})\sim p(\widecheck{\eta},\overline{\sigma})^{-1}$. We deduce that the right hand side of (\ref{CM proof step 1}) is equivalent to $(2\pi i)^{mnd}(P^{(I(\Pi,\eta))}(\Pi,\xi)\prod\limits_{\sigma\in\Sigma} (p(\widecheck{\eta},\sigma)^{s_{\sigma}}p(\widecheck{\eta},\overline{\sigma})^{n-s_{\sigma}})$.

We end our proof by the fact that  $L(m,\Pi\otimes\eta)=L(m,\Pi^{c}\otimes \eta^{c}) =L(m,\Pi^{\vee}\otimes \widetilde{\chi}\alpha) =L^{mot}(m+\cfrac{n-1}{2},\pi\otimes\chi,St,\alpha)$.

\begin{flushright}$\Box$\end{flushright}
 
 \section{Arithmetic automorphic periods for  conjugate self-dual representations}

If we consider $\Pi^{*}:=\Pi\otimes \eta$, it is not conjugate self-dual but $\Pi^{*}\otimes \eta^{-1}$ is. We call such representations \textbf{ conjugate self-dual}. We want to generalize the definition for arithmetic automorphic period to such representations.

We firstly generalize the definition for $I(\Pi,\eta)$ in Theorem \ref{n*1}.

\begin{dflem}\label{index I }
Let $\Pi^{*}$ be an algebraic regular representation of $GL_{n}(F)$ with infinity type $(z^{a_{i}(\sigma)}\overline{z}^{b_{i}(\sigma)})_{1\leq i\leq n}$ at $\sigma\in \Sigma$. Let $\eta$ be am algebraic Hecke character of $F$ with infinity type $z^{a(\sigma)}\overline{z}^{b(\sigma)}$ at $\sigma\in \Sigma$.
We assume that $a(\sigma)-b(\sigma)+a_{i}(\sigma)-b_{i}(\sigma)\neq 0$ for all $\sigma$ and $i$.

We define $I(\Pi^{*},\eta)$ to be the map sending $\sigma\in\Sigma$ to $\#\{i:a(\sigma)-b(\sigma)+a_{i}(\sigma)-b_{i}(\sigma)<0\}$.

It is easy to see that $I(\Pi^{*},\eta_{1}\eta_{2})=I(\Pi^{*}\otimes\eta_{1},\eta_{2})$ for any $\eta_{1}$, $\eta_{2}$.
\end{dflem}

\bigskip

We can now define the arithmetic automorphic periods.
\begin{dflem}\label{definable periods}\label{general period definition}
\text{}

We say a $3$-regular cohomological cuspidal automorphic representation $\Pi^{*}$ of $GL_{n}(\AF)$ has \textbf{definable arithmetic automorphic periods} if there exists an algebraic Hecke character $\eta$ of $F$ such that $\Pi^{*}\otimes \eta^{-1}$ descends to unitary groups of any sign. In particular, $\Pi^{*}\otimes \eta^{-1}$ is conjugate self-dual.

In this case, for any sign $I$, i.e. a map from $\Sigma$ to $\{0,1,\cdots,n\}$, we define the \textbf{arithmetic automorphic period} for $\Pi^{*}$ by $P^{(I)}(\Pi^{*}):=P^{(I)}(\Pi^{*}\otimes \eta^{-1})\prod\limits_{\sigma\in\Sigma}p(\widecheck{\eta},\sigma)^{I(\sigma)}p(\widecheck{\eta},\overline{\sigma})^{n-I(\sigma)}.$ 

This definition does not depend on the choice of $\eta$ and hence is compatible with Definition \ref{arithmetic period conjugate self-dual case} if $\Pi^{*}$ itself is conjugate self-dual. \\

\end{dflem}

\begin{dem}
The last part comes from Theorem \ref{n*1}. In fact, for any $I$, let $\chi$ be an algebraic Hecke character such that $I(\Pi^{*},\chi)=I$. Since $\Pi^{*}$ is $3$-regular, we may choose $\chi$ such that there exists $m\geq 1+\cfrac{\omega(\eta)+\omega(\chi)}{2}$ critical for $\Pi^{*}\otimes \chi$.\\

Let $\eta$ be an algebraic Hecke character such that $\Pi:=\Pi^{*}\otimes \eta^{-1}$ is conjugate self-dual, we have $\Pi^{*}\otimes \chi= \Pi \otimes (\eta\chi)$.\\

Since $I(\Pi,\eta\chi)=I(\Pi^{*},\chi)=I$, Theorem \ref{n*1} gives that:
\begin{eqnarray}\nonumber
&&L(m,\Pi^{*}\otimes \chi)=L(m,\Pi\otimes (\eta\chi)) \\\nonumber&\sim_{E(\Pi)E(\eta)E(\chi);K}& (2 \pi i)^{mnd} P^{(I)}(\Pi) \prod\limits_{\sigma\in\Sigma}p(\widecheck{\eta\chi},\sigma)^{I(\sigma)}p(\widecheck{\eta\chi},\overline{\sigma})^{n-I(\sigma)} \\\nonumber
&\sim_{E(\Pi)E(\eta)E(\chi);K}& (2 \pi i)^{mnd} P^{(I)}(\Pi) \prod\limits_{\sigma\in\Sigma}p(\widecheck{\eta},\sigma)^{I(\sigma)}p(\widecheck{\eta},\overline{\sigma})^{n-I(\sigma)}\prod\limits_{\sigma\in\Sigma}p(\widecheck{\chi},\sigma)^{I(\sigma)}p(\widecheck{\chi},\overline{\sigma})^{n-I(\sigma)}
\end{eqnarray}
with both sides non zero.\\

In particular, $P^{(I)}(\Pi^{*}\otimes \eta^{-1}) \prod\limits_{\sigma\in\Sigma}p(\widecheck{\eta},\sigma)^{I(\sigma)}p(\widecheck{\eta},\overline{\sigma})^{n-I(\sigma)}$\begin{equation}\nonumber
\sim_{E(\Pi^{*});K} (2\pi i)^{-mnd} (\prod\limits_{\sigma\in\Sigma}p(\widecheck{\chi},\sigma)^{I(\sigma)}p(\widecheck{\chi},\overline{\sigma})^{n-I(\sigma)})^{-1}L(m,\Pi^{*}\otimes \chi)
\end{equation} which does not depend on the choice of $\eta$.

\end{dem}

\begin{rem}
Let $\Pi$ be a $3$-regular cohomological cuspidal automorphic representation of $GL_{n}(\AF)$ which is conjugate self-dual after tensoring an algebraic Hecke character $\eta$. Let $q$ be a prime number inert in $F^{+}$ and split in $F$. If $(\Pi\otimes \eta)_{q}$ descends locally then $\Pi$ has definable arithmetic automorphic periods. In particular, this holds true if $\Pi_{q}$ is in discrete series.
\end{rem}
\bigskip

We read from the above proof that Theorem \ref{n*1} can be rewritten as follows:

\begin{thm}\label{n*1 essential}
Let $\Pi$ be an algebraic automorphic representation of $GL_{n}(\AF)$ which has definable arithmetic automorphic periods. Let $\eta$ be an algebraic Hecke character as in Definition \ref{index I }. We write $I:=I(\Pi,\eta)$.

Let $m\in \Z+\cfrac{n-1}{2}$ be critical for $\Pi\otimes \eta$. If $m\geq \cfrac{1+\omega(\Pi)+\omega(\eta)}{2}$ then
\begin{equation}
L(m,\Pi\otimes \eta) \sim_{E(\Pi)E(\eta);K} (2 \pi i)^{mnd} P^{(I(\Pi,\eta))}(\Pi) \prod\limits_{\sigma\in\Sigma}p(\widecheck{\eta},\sigma)^{I(\sigma)}p(\widecheck{\eta},\overline{\sigma})^{n-I(\sigma)}.
\end{equation}
Moreover, there always exists such $m$ with both sides non zero. 
\end{thm}

\bigskip

\begin{rem}
The last part comes from the fact that $\Pi$ is $3$-regular. The $3$-regular condition is not needed to define the arithmetic automorphic periods in general. We assume it here to guarantee that Definition \ref{definable periods} does not depend on the choice of $\eta$. One can replace this condition by a weaker one on the non vanishing property for certain $L$-functions.
\end{rem}

\chapter{Motives and Deligne's conjecture}\label{chapter motive}
\section{Motives over $\Q$}
In this article, a \textbf{motive} simply means a pure motive for absolute Hodge cycles in the sense of Deligne \cite{deligne79}. \\

More precisely, a motive over $\Q$ with coefficients in a number field $E$ is given by its Betti realization $M_{B}$, its de Rham realization $M_{DR}$ and its $l$-adic realization $M_{l}$ for all prime numbers $l$ where $M_{B}$ and $M_{DR}$ are finite dimensional vector space over $E$, $M_{l}$ is a finite dimensional vector space over $E_{l}:=E\otimes_{\Q} \Q_{l}$ endowed with:
\begin{itemize}
\item $I_{\infty}: M_{B}\otimes \C \xrightarrow{\sim}  M_{DR}\otimes \C $ as $E\otimes_{\Q}\C$-module;
\item $I_{l}: M_{B}\otimes \Q_{l}\xrightarrow{\sim} M_{l}$ as $E\otimes _{\Q} \Q_{l}$-module.

\end{itemize}

\bigskip

From the isomorphisms above, we see that $dim_{E}M_{B}=dim_{E}M_{DR}=dim_{E_{l}}M_{l}$ and this is called the \textbf{rank} of $M$. We need moreover:
\begin{enumerate}
\item An $E$-linear involution (infinite Frobenius) $F_{\infty}$ on $M_{B}$ and a Hodge decomposition $M_{B}\otimes \C=\bigoplus\limits_{p,q\in \Z}M^{p,q}$ as $E\otimes \C$-module such that $F_{\infty}$ sends $M^{p,q}$ to $M^{q,p}$.\\

For $w$ an integer, we say $M$ is \textbf{pure of weight} $w$ if $M^{p,q}=0$ for $p+q\neq w$. Throughout this paper, all the motives are assumed to be pure. We assume also that $F_{\infty}$ acts on $M^{p,p}$ as a scalar for all $p\in \Z$.\\

We say $M$ is \textbf{regular} if $dim M^{p,q}\leq 1$ for all $p,q\in\Z$.

\item An $E$-rational Hodge filtration on $M_{DR}$: $\cdots \supset M^{i}\supset M^{i+1}\supset \cdots$ which is compatible with the Hodge structure on $M_{B}$ via $I_{\infty}$, i.e.,
$$I_{\infty}(\bigoplus\limits_{p\geq i}M^{p,q})=M^{i}\otimes \C.$$
\item A Galois action of $G_{\Q}$ on each $M_{l}$ such that $(M_{l})_{l}$ forms a compatible system of $l$-adic representations $\rho_{l}:G_{\Q} \longrightarrow GL(M_{l})$. More precisely, for each prime number $p$, let $I_{p}$ be the inertia subgroup of a decomposition group at $p$ and $F_{p}$ the geometric Frobenius of this decomposition group. We have that for all $l\neq p$, the polynomial  $det(1-F_{p}|M_{l}^{I_{p}})$ has coefficients in $E$ and is independent of the choice of $l$.
We can then define $L_{p}(s,M):=det(1-p^{-s}F_{p}|M_{l}^{I_{p}})^{-1}\in E(p^{-s})$ for whatever $l\neq p$.
\end{enumerate}

\bigskip

For any fixed embedding $\sigma: E\hookrightarrow \C$, we may consider $L_{p}(s,M,\sigma)$ as a complex valued function. We define $L(s,M,\sigma)=\prod\limits_{p}L_{p}(s,M,\sigma)$. It converges for $Re(s)$ sufficiently large. It is conjectured that the $L$-function has analytic continuation and functional equation on the whole complex plane.\\

We can also define $L_{\infty}(s,M)$, the infinite part of the $L$-function, as in chapter $5$ of \cite{deligne79}.\\

Deligne has defined the critical values for $M$ as follows:

\begin{df}
We say an integer $m$ is \textbf{critical} for $M$ if neither $L_{\infty}(M,s)$ nor $L_{\infty}(\check{M},1-s)$ has a pole at $s=m$ where $\check{M}$ is the dual of $M$. We call $m$ a \textbf{critical value} of $M$.
\end{df}

\bigskip

\begin{rem}\label{infinitytype}
The notion $L_{\infty}(s,M)$ implicitly indicates that the infinity type of the $L$-function does not depend on the choice of $\sigma:E\hookrightarrow \C$. More precisely, for every $\sigma:E\hookrightarrow \C$, put $M_{B,\sigma}:=M_{B}\otimes_{E,\sigma}\C$. \\

We then have $M_{B}\otimes \C=\bigoplus\limits_{\sigma:E\hookrightarrow \C}M_{B,\sigma}$. Since $M^{p,q}$ is stable by $E$, each $M_{B,\sigma}$ inherits a Hodge decomposition $M_{B,\sigma}=\bigoplus M_{B,\sigma}^{p,q}$. We may define $L_{\infty}(s,M,\sigma)$ with help of the Hodge decomposition of $M_{B}\otimes_{E,\sigma}\C$. It is a product of $\Gamma$ factors which depend only on $\dim M_{B,\sigma}^{p,q}$ and the action of $F_{\infty}$ on $M_{B,\sigma}^{p,p}$. The latter is independent of $\sigma$ since we have assumed that $F_{\infty}$ acts on $M^{p,p}$ by a scalar.\\

It remains to show that $\dim M_{B,\sigma}^{p,q}$ is also independent of $\sigma$. In fact, since $M$ is pure, $M^{p,q}$ can be reconstructed from the Hodge filtration $M^{i}$. Hence $M^{p,q}=\bigoplus M_{B,\sigma}^{p,q}$ is a free $E\otimes\C$-module. One can show $M_{B,\sigma}^{p,q}=M^{p,q}\otimes_{E,\sigma}\C$ and hence $dim M_{B,\sigma}^{p,q}$ is independent of $\sigma$.

\end{rem}

\bigskip

If $F_{\infty}$ acts as a scalar at $M^{p,p}$ for every $p$ then Deligne's period can be defined. We will only treat the case when $M$ has no $(p,p)$ class. Therefore, Deligne's period can always be defined. \\

\begin{df}
Let $M$ be a motive over $\Q$ of weight $\omega$ which has no $(\omega/2,\omega/2)$ class. We denote by $M_{B}^{+}$ the subspace of $M_{B}$ fixed by $F_{\infty}$. We denote by $F^{+}(M):=F^{\omega/2}(M)$ a subspace of $M_{DR}$. It is easy to see that $I_{\infty}^{-1}(F^{+}(M)\otimes \C)$ equals to $\bigoplus\limits_{p>q}M^{p,q}$.\\

The comparison isomorphism then induces an isomorphism:
\begin{equation}
M_{B}^{+}\otimes \C \hookrightarrow M_{B}\otimes\C \xrightarrow{\sim} M_{DR}\otimes \C\rightarrow (M_{DR}/F^{+}(M))\otimes \C. \end{equation}
Deligne's period $c^{+}(M)$ is defined to be the determinant of the above isomorphism with respect to fixed $E(M)$-bases of $M_{B}^{+}$ and $M_{DR}/F^{+}(M)$. It is well defined up to $E(M)^{\times}$,
\end{df}

\bigskip

Deligne has conjectured in \cite{deligne79} that:
\begin{conj}\label{deligne conjecture}
If $0$ is critical for $M$ (see the \textit{loc.cit} for the definition of critical), then $L(0,M,\sigma)\sim_{E(M)} c^{+}(M)$.
\end{conj}

\bigskip

More generally, tensoring $M$ by the Tate motive $\Q(m)$ (c.f. \cite{deligne79} chapter $1$), we obtained a new motive $M(m)$. We remark that $L(s,M(m),\sigma)=L(s+m,M,\sigma)$. The following conjecture is a corollary of the previous conjecture:

\begin{conj}
If $m$ is critical for $M$, then
\begin{equation}
L(m,M,\sigma)\sim_{E(M)} (2\pi i)^{d^{+}m}c^{+}(M) 
\end{equation}
where $d^{+}=dim_{E(M)}(M_{B}^{+})$.
\end{conj}

\bigskip

Deligne has given a criteria to determine whether $0$ is critical for $M$ (see $(1.3.1)$ of \cite{deligne79}). We observe that $n$ is critical for $M$ if and only if $0$ is critical for $M(n)$. Thus we can rewrite the criteria of Deligne for arbitrary $n$. In the case where $M^{p,p}=0$ for all $p$, this criteria becomes rather simple.\\

We first define the \textbf{Hodge type} of $M$ by the set $T=T(M)$ consisting of pairs $(p,q)$ such that $M^{p,q}\neq 0$. Since $M$ is pure, there exists an integer $w$ such that $p+q=w$ for all $(p,q)\in T(M)$. We remark that if $(p,q)$ is an element of $T(M)$, then $(q,p)$ is also contained $T(M)$.\\

 \begin{lem}\label{critical}
Let $M$ be a pure motive of weight $w$. We assume that for all $(p,q)\in T(M)$, $p\neq q$ which is equivalent to that $p\neq \cfrac{w}{2}$.\\

Let $p_{1}<p_{2}<\cdots<p_{n}$ be some integers such that $$T(M)=\{(p_{1},q_{1}),(p_{2},q_{2}),\cdots, (p_{n},q_{n})\}\cup \{(q_{1},p_{1}),(q_{2},p_{2}),\cdots, (q_{n},p_{n})\}$$ where $q_{i}=w-p_{i}$ for all $1\leq i\leq n$.\\

We set $p_{0}=-\infty$ and $p_{n+1}=+\infty$. Denote by $k:=max\{0\leq i \leq n\mid  p_{i}< \cfrac{w}{2} \}$. We have that $m$ is critical for $M$ if and only if
$$\max(p_{k}+1,w+1-p_{k+1})\leq m\leq  \min(w-p_{k},p_{k+1}).$$
In particular, critical value always exists in the case where $p_{i}\neq q_{i}$ for all $i$.
 \end{lem}

\begin{dem}

The Hodge type of $M(m)$ is $\{(p_{i}-m,w-p_{i}-m)\mid 1\leq i\leq n\}\cup \{(w-p_{i}-m,p_{i}-m)\mid 1\leq i\leq n\}$. By Deligne's criteria, $0$ is critical for $M$ if and only if for all $i$, either $p_{i}-m\leq -1$ and $w-p_{i}-m\geq 0$, or $p_{i}-m\geq 0$ and $w-p_{i}-m\leq -1$. Hence the set of critical values for $M$ are  $\bigcap\limits_{1\leq i\leq n} ([w+1-p_{i}, p_{i}]\cup [p_{i}+1,w-p_{i}])$.

For $i\leq k$, $p_{i}<\cfrac{w}{2}$ and then $p_{i}<w+1-p_{i}$. Therefore $\bigcap\limits_{1\leq i\leq k} ([w+1-p_{i}, p_{i}]\cup [p_{i}+1,w-p_{i}])=\bigcap\limits_{1\leq i\leq k}  [p_{i}+1,w-p_{i}]=[p_{k}+1,w-p_{k}]$. Similarly we have $\bigcap\limits_{k< i\leq n} ([w+1-p_{i}, p_{i}]\cup [p_{i}+1,w-p_{i}])=\bigcap\limits_{k< i\leq n} [w+1-p_{i}, p_{i}]=[w+1-p_{k+1},p_{k+1}]$.

We deduce, at last, that the set of critical values for $M$ is $[\max(p_{k}+1,w+1-p_{k+1}), \min(w-p_{k},p_{k+1})]$. It is easy to verify that the latter set is non empty.
\end{dem}

\section{Motivic periods over quadratic imaginary fields}\label{definition of motivic period}

Recall that $K$ is a quadratic imaginary field with fixed embedding $K\hookrightarrow \overline{\Q}$. Let $E$ be a number field.\\

Let $M$ be a regular motive over $K$ (with respect to the fixed embedding) with coefficients in $E$ of dimension $n$ pure of weight $\omega(M)$.
\\

Recall that $M_{B}$, the Betti realization of $M$, is a finite dimensional $E$-vector space. The infinite Frobenius gives an E-linear isomorphism $F_{\infty}: M_{B} \xrightarrow{\sim} M_{B}^{c}.$\\

Since $M$ is regular of dimension $n$, we can write its Hodge type by $(p_{i},q_{i}=\omega(M)-p_{i})_{1\leq i\leq n}$ with $p_{1}>p_{2}>\cdots>p_{n}$. The Betti realization $M_{B}$ has a Hodge decomposition $M_{B}\otimes_{\Q}\C=\bigoplus\limits_{i=1}^{n}M^{p_{i},q_{i}}$ as $E\otimes_{\Q}\C$-modules.\\

We write the Hodge type of $M^{c}$ by $(p_{i}^{c},q_{i}^{c}=\omega(M)-p_{i}^{c})_{1\leq i\leq n}$ with $p_{i}^{c}=q_{n+1-i}=\omega(M)-p_{n+1-i}$. Note that the Hodge numbers $p_{i}^{c}$ are still in decreasing order. We know $M_{B}^{c}\otimes_{\Q}\C=\bigoplus\limits_{i=1}^{n}(M^{c})^{p_{i}^{c},q_{i}^{c}}$ and $F_{\infty}$ induces $E$-linear isomorphisms: $M^{p_{i},q_{i}}\xrightarrow{\sim} (M^{c})^{p^{c}_{n+1-i},q^{c}_{n+1-i}}$.\\

The De Rham realization $M_{DR}$ is also a finite dimensional $E$-linear space endowed with a Hodge filtration $M_{DR}=M^{p_{n}}\supset M^{p_{n-1}}\supset\cdots\supset M^{p_{1}}$. The comparison isomorphism:
\begin{equation}
I_{\infty}: M_{B}\otimes_{\Q} \C\xrightarrow{\sim} M_{DR}\otimes_{\Q}\C
\end{equation}
induces compatibility isomorphisms on the Hodge decomposition of $M_{B}$ and the Hodge filtration on $M_{DR}$.\\

 More precisely, for each $1\leq i\leq n$, $I_{\infty}$ induces an isomorphism:
\begin{equation}\label{comparison Hodge De Rham}
I_{\infty}: \bigoplus\limits_{p_{j}\geq p_{i}}M^{p_{j},q_{j}}=\bigoplus\limits_{j\leq i}M^{p_{j},q_{j}}\xrightarrow{\sim} M^{p_{i}}\otimes_{\Q}\C
\end{equation}

\bigskip

\begin{df}
For any fixed $E$-bases of $M_{B}$ and $M_{DR}$, we can extend them to $E\otimes_{\Q}\C$ bases of $M_{B}\otimes \C$ and $M_{DR}\otimes \C$ respectively. We define $\delta^{Del}(M)$ to be the determinant of $I_{\infty}$ with respect to the fixed $E$-rational bases, called the \textbf{determinant period}. It is an element in $(E\otimes\C)^{\times}$well defined up to multiplication by elements in $E^{\times}\subset (E\otimes\C)^{\times}$.
\end{df} This is an analogue of Deligne's period $\delta$ defined in ($1.7.3$) of \cite{deligne79}.

\bigskip

Let us now fix some bases. We take $(e_{i})_{1\leq i\leq n}$ an $E$-base of $M_{B}$. Since $F_{\infty}$ is $E$-linear on $M_{B}$, we know $(e^{c}_{i}:=F_{\infty}e_{i})_{1\leq i\leq n}$ forms an $E$-base of $M_{B}^{c}$.\\

From (\ref{comparison Hodge De Rham}) we see that $I_{\infty}$ induces an isomorphism $M^{p_{i},q_{i}}\xrightarrow{\sim} (M^{p_{i}})\otimes_{\Q}\C/ (M^{p_{i-1}})\otimes_{\Q}\C$ for any $1\leq i\leq n$. Here we set $M^{p_{0}}=\{0\}$. Let $\omega_{i}$ be a non zero element in $M^{p_{i},q_{i}}$ such that the image of $\omega_{i}$ by the above isomorphism is in $M^{p_{i}} (\text{mod }(M^{p_{i-1}})\otimes_{\Q}\C)$. In other words, $I_{\infty}(\omega_{i})$ is equivalent to an element in $M^{p_{i}}$ modulo $(M^{p_{i-1}})\otimes_{\Q}\C$. \\

Since $(\omega_{i})_{1\leq i\leq n}$ forms an $E\otimes\C$-base of $M_{B}\otimes \C$, we know $(I_{\infty}(\omega_{i}))_{1\leq i\leq n}$ forms an $E\otimes\C$-base of $M_{DR}\otimes \C$. This base is not rational, i.e. is not contained in $M^{DR}$. But by the above construction, it can pass to a rational base of $M_{DR}\otimes_{\Q}\C$ with a unipotent matrix for change of basis. Since the determinant of a unipotent matrix is always one, we can use this base for calculating $\delta^{Del}(M)$. \\

We define $\omega^{c}_{i}\in (M^{c})^{p_{i}^{c},q_{i}^{c}}$ similarly. We will use $(I_{\infty}(\omega^{c}_{i}))_{1\leq i\leq n}$ as an $E\otimes_{Q} \C$ base of $M_{DR}\otimes \C$ to calculate motivic period henceforth.\\

\begin{df}
Since $M^{p_{i}^{c},q_{i}^{c}}$ is a rank one free $E\otimes\C$-module, we know there exists numbers $Q_{i}(M)\in (E\otimes\C)^{\times}$, $1\leq i\leq n$ such that
$F_{\infty}\omega_{i}=Q_{i}(M)\omega^{c}_{n+1-i}$ for all $i$. These numbers in $(E\otimes \C)^{\times}$ are called \textbf{motivic period} and well defined up to $E^{\times}$.
\end{df}

\bigskip

Since $F_{\infty}^{2}=Id$, we have $F_{\infty}\omega^{c}_{n+1-i}=Q_{i}(M)^{-1}\omega_{i}$. We deduce that:

\begin{lem}\label{period for Mc}
For all $1\leq i\leq n$, $Q_{i}(M^{c})\sim_{E(M);K} Q_{n+1-i}(M)^{-1}$.
\end{lem}

\bigskip

We write $\omega_{a}=\sum\limits_{i=1}^{n}A_{ia}e_{i}$, $\omega_{t}^{c}=\sum\limits_{i=1}^{n}A_{it}^{c}e^{c}_{i}$ for all $1\leq a,t\leq n$.\\

We know $\delta^{Del}(M)=det(A_{ia})_{1\leq i,a \leq n}$ and  $\delta^{Del}(M^{c})=det(A^{c}_{it})_{1\leq i,t \leq n}$. This implies that $\bigwedge_{i=1}^{n}\omega_{i}=\delta^{Del}(M)\bigwedge_{i=1}^{n}e_{i}$.\\

We denote by $det(M)$ the determinant motive of $M$ as in section $1.2$ of \cite{harrisadjoint}. We know $I_{\infty}(\bigwedge_{i=1}^{n}\omega_{i})$ is an $E$-base of $det(M)_{DR}$ and $\bigwedge_{i=1}^{n}e_{i}$ is an $E$-base of $det(M)_{B}$. Moreover, $F_{\infty}(\bigwedge_{i=1}^{n}\omega_{i})=\prod\limits_{1\leq i\leq n}Q_{i}(M)\bigwedge_{i=1}^{n}\omega_{i}^{c}$. \\

We deduce that:

\begin{lem}\label{determinant motive lemma}
\begin{eqnarray}
\delta^{Del}(M) \sim_{E(M);K} \delta^{Del}(det(M))\\
Q_{1}(det(M))\sim_{E(M);K} \prod\limits_{i=1}^{n}Q_{i}(M)
\end{eqnarray}
\end{lem}

\bigskip

\begin{rem}
The determinant period $\delta^{Del}(M)$ is inverse of the period $\delta$ defined in \cite{harrisadjoint}. In fact, the period $\delta(M)$ is defined by equation $(1.2.4)$ of \cite{harrisadjoint}, namely, $\bigwedge_{i=1}^{n}e_{i}=\delta(M)\bigwedge_{i=1}^{n}\omega_{i}$. Therefore $\delta(M)\sim_{E(M);K}\delta^{Del}(det(M))^{-1}\sim_{E(M);K}\delta^{Del}(M)^{-1}$.

\end{rem}

\bigskip

\begin{lem}\label{delta c}
For all motive $M$ as above, we have:
\begin{equation}\nonumber
\delta^{Del}(M^{c})\sim_{E(M);K}(\prod\limits_{1\leq i\leq n}Q_{i}^{-1})\delta^{Del}(M)
\end{equation}
\end{lem}
\begin{dem}
This follows directly from equation (\ref{coefficient relation}). \\

One can also prove this with help of Lemma \ref{determinant motive lemma}. In fact, by Lemma \ref{determinant motive lemma}, we may assume that $n=1$. We take $\omega\in M_{DR}$, $\omega^{c} \in M^{c}_{DR}$ and $e\in M_{B}$. Then $\omega=\delta^{Del}(M)e$ and $\omega^{c}=\delta^{Del}(M^{c})e^{c}$ where $e^{c}=F_{\infty}e$.\\

By definition of  motivic period, we have $F_{\infty}\omega=Q_{1}(M)\omega^{c}$ and then $\omega^{c}=Q_{1}(M)^{-1}F_{\infty}\omega=Q_{1}(M)^{-1}F_{\infty}(\delta^{Del}(M)e)=Q_{1}(M)^{-1}\delta^{Del}(M)e^{c}$. It follows that $\delta^{Del}(M^{c})\sim_{E(M);K}Q_{1}(M)^{-1}\delta^{Del}(M)$ as expected.
\end{dem}

\bigskip

\begin{ex}\textbf{Tate motive}
\text{}

Let $\Z(1)_{K}$ be the extension of $\Z(1)$ from $\Q$ to $K$. It is a motive with coefficients in $K$. As in section $3.1$ of \cite{deligne79}, $\Z(1)_{K,B}=H_{1}(\mathbb{G}_{m,K})\cong K$ and $\Z(1)_{K,DR}$ is the dual of $H^{1}_{DR}(\mathbb{G}_{m,K})$ with generator $\cfrac{dz}{z}$. Therefore the comparison isomorphism $\Z(1)_{K,B}\otimes \C\cong K\otimes \C\rightarrow \Z(1)_{K,DR}\otimes \C\cong K\otimes \C$ sends $K$ to $\oint\cfrac{dz}{z}K=(2\pi i )K$. We have $\delta^{Del}(\Z(1)_{K})\sim_{K;K} 2\pi i$. \\

In general, let $M$ be a motive over $K$ with coefficients in $E$. We have 
\begin{equation}\label{Tate twist}
\delta^{Del}(M(n))\sim_{E(M);K} (2\pi i)^{n}\delta^{Del}(M).
\end{equation}
\end{ex}

\bigskip

\begin{rem}
All the determinants and the coefficients we consider here are elements in $(E\otimes_{\Q}\C)^{\times}$.
\end{rem}

\section{Deligne's conjecture for tensor product of motives}
Let $E$ and $E'$ be two number fields.\\

Let $M$ be a regular motive over $K$ (with respect to the fixed embedding) with coefficients in $E$ of dimension $n$ pure of weight $\omega(M)$.
Let $M'$ be a regular motive over $K$ with coefficients in $E'$ of dimension $n'$  pure of weight $\omega(M')$.\\

We denote by $R(M\otimes M')$ the restriction from $K$ to $\Q$ of the motive $M\otimes M'$. It is a motive of weight $\omega:=\omega(M)+\omega(M')$ with Betti realization $M_{B}\otimes M'_{B}\oplus M^{c}_{B}\otimes M'^{c}_{B}$ and De Rham realization $M_{DR}\otimes M'_{DR}\oplus M^{c}_{DR}\otimes M'^{c}_{DR}$.\\

We denote the Hodge type of $M$ by $(p_{i},\omega(M)-p_{i})_{1\leq i\leq n}$ with $p_{1}>\cdots>p_{n}$ and the Hodge type of $M'$ by $(r_{j},\omega(M')-r_{j})_{1\leq j\leq n}$ with $r_{1}>r_{2}>\cdots>r_{n'}$. As before, we define $p_{i}^{c}=\omega(M)-p_{n+1-i}$ and $r_{j}^{c}=\omega(M')-r_{n'+1-j}$. There are indices for Hodge type of $M^{c}$ and $M'^{c}$ respectively.\\

We assume that $R(M\otimes M')$ has no $(w/2,w/2)$ class. In other words, $p_{a}+r_{b}\neq \frac{\omega}{2}$ and then $p_{t}^{c}+r_{u}^{c}\neq \frac{\omega}{2}$ for all $1\leq a,t\leq n$, $1\leq b,u\leq n'$.\\

As in the above section, we take $(e_{i})_{1\leq j\leq n}$ an $E$-base of $M_{B}$ and define $(e^{c}_{i}:=F_{\infty}e_{i})_{1\leq i\leq n}$ which is an $E$-base of $M^{c}_{B}$. Similarly, we take $(f_{j})_{1\leq j\leq n'}$ an $E'$-base of $M'_{B}$ and define $f^{c}_{j}:=F_{\infty}f_{j}$ for $1\leq j\leq n'$.\\

We also take $\omega_{i}\in M^{p_{i},\omega(M)-p_{i}}$, $(\omega^{c}_{i})\in (M^{c})^{p^{c}_{i},\omega(M)-p^{c}_{i}}$ for $1\leq i\leq n$ as in previous section and $\eta_{j}\in M^{r_{j},\omega(M')-r_{j}}$, $\eta^{c}_{j}\in (M'^{c})^{r^{c}_{j},\omega(M')-r^{c}_{j}}$ for $1\leq j\leq n'$ similarly.\\

\bigskip

Recall the motive periods are complex numbers $Q_{i}$, $1\leq i\leq n$ and $Q'_{j}$, $1\leq j\leq n'$ such that
\begin{equation}\label{Qi}\tag{P}
F_{\infty}\omega_{i}=Q_{i}\omega^{c}_{n+1-i}, F_{\infty}\mu_{j}=Q'_{j}\mu^{c}_{n'+1-j}.
\end{equation}

The aim of this section is to calculate the Deligne's period for $R(M\otimes M')$ in terms of motivic periods.

\begin{rem}
If we define a paring $(M_{B}\otimes\C)\otimes (M_{B}\otimes\C)\rightarrow \C$ such that $<\omega_{i},\omega_{n+1-i}^{c}>=1$ and $<\omega_{i},\omega_{n+1-j}^{c}>=0$ for $j\neq i$ then $Q_{i}=<\omega_{i},F_{\infty}\omega_{i}>$. 
\end{rem}

\bigskip

Let $M^{\#}=R(M\otimes M')$. It is a motive over $\Q$. We are going to calculate $c^{+}(M^{\#})$.

We define $A=\{(a,b) \mid p_{a}+r_{b}> \frac{\omega}{2}\}$ and $T=\{(t,u) \mid p^{c}_{t}+r^{c}_{u}> \frac{\omega}{2}\}=\{(t,u) \mid p_{n+1-t}+r_{n'+1-u}< \frac{\omega}{2}\}$. 

\bigskip

\begin{rem} 
Keeping in mind that \begin{equation}\label{AT}
(t,u)\in T\text{ if and only if }(n+1-t,n'+1-u)\notin A.
\end{equation}
\end{rem}

\bigskip

\begin{prop}\label{deligne period tensor product}
Let $M$, $M'$ be motive over $K$ with coefficients in $E$ and $E'$ respectively. We assume that $M\otimes M'$ has no $(\omega/2,\omega/2)$-class. 
We then have
\begin{eqnarray}c^{+}(R(M\otimes M')) &\sim_{E(M)E(M');K} &\left(\prod_{(t,u)\notin T(M,M')}Q_{n+1-t}(M)^{-1}Q_{n'+1-u}(M')^{-1}\right)\delta^{Del}(M\otimes M')\nonumber\\&\sim_{E(M)E(M');K} &\left(\prod_{(t,u)\in A(M,M')}Q_{t}(M)^{-1}Q_{u}(M')^{-1}\right)\delta^{Del}(M\otimes M')
\end{eqnarray} 

\end{prop}

\begin{dem}

For simplification of notation, we identify $\omega_{i}\in M_{B}\otimes \C$ and $I_{\infty}(\omega_{i})\in M_{DR}\otimes \C$ and similarly, we identify $\omega^{c}$, $\mu_{j}$, $\mu_{j}^{c}$ with their image under $I_{\infty}$ in the following.\\

We fixe bases for $M^{+}_{B}$ and $M^{\#}_{DR}/F^{+}(M^{\#})$ now. For $M^{+}_{B}$, we know $(e_{i}\otimes f_{j}+e^{c}_{i}\otimes f^{c}_{j})_{1\leq i\leq n, 1\leq j\leq n'}$ forms an $EE'$-base. For $M^{\#}_{DR}/F^{+}(M^{\#})$, as in the above section, we consider $\mathcal{B}:= (\omega_{a}\otimes \mu_{b}, \omega^{c}_{t}\otimes \mu^{c}_{u}(\text{mod } F^{+}(M^{\#})) \mid (a,b)\notin A, (t,u)\notin T)$ as an $E\otimes\C$ base of $(M^{\#}_{DR}/F^{+}(M^{\#}))\otimes\C$ which is not rational but can change to a rational base with a unipotent matrix for change of basis. Therefore we can use this base to calculate Deligne's period.\\

If $(a,b)\notin A$ then $(n+1-a,n'+1-b)\in T$ by (\ref{AT}). Along with (\ref{Qi}), we know that $ F_{\infty}(\omega_{a}\otimes \mu_{b})=Q_{a}Q'_{b}\omega_{n+1-a}^{c}\mu_{n'+1-b}^{c} \in F^{+}(M^{\#})\otimes \C$.  Similarly, $F_{\infty}( \omega^{c}_{t}\otimes \mu^{c}_{u}) \in F^{+}(M^{\#})\otimes \C$ for all $(t,u)\notin T$.\\

Note that $F_{\infty}$ is an endomorphism on $M^{\#}(B)\otimes \C$ and $M^{\#}_{DR}\otimes \C$. For any $\phi\in M^{\#}(B)\otimes \C$ or $M^{\#}(DR)\otimes \C$, we write $(1+F_{\infty})\phi:=\phi+F_{\infty}(\phi)$.\\

Recall that $(M^{\#}_{DR}/F^{+}(M^{\#}))\otimes \C \cong (M^{\#}_{DR}\otimes \C)/(F^{+}(M^{\#})\otimes \C)$. Thus $ \mathcal{B}=((1+F_{\infty})\omega_{a}\otimes \mu_{b},(1+F_{\infty}) \omega^{c}_{t}\otimes \mu^{c}_{u}(\text{mod } F^{+}(M^{\#})\otimes\C) \mid (a,b)\notin A, (t,u)\notin T)$.\\

We write $\omega_{a}=\sum\limits_{i=1}^{n}A_{ia}e_{i}$, $\omega_{t}^{c}=\sum\limits_{i=1}^{n}A_{it}^{c}e^{c}_{i}$, $\mu_{b}=\sum\limits_{j=1}^{n'}B_{jb}f_{j}$, $\mu_{u}=\sum\limits_{j=1}^{n'}B^{c}_{ju}f^{c}_{j}$ for all $1\leq a,t\leq n$ and $1\leq b,u\leq n'$. \\

We then have 
\begin{eqnarray}\nonumber &&(1+F_{\infty})\omega_{a}\mu_{b}=(1+F_{\infty})\sum\limits_{i,j}A_{ia}B_{jb}e_{i}\otimes f_{j}=\sum\limits_{i,j}A_{ia}B_{jb}(e_{i}\otimes f_{j}+e^{c}_{i}\otimes f^{c}_{j})\\
&\text{and }& (1+F_{\infty})\omega^{c}_{t}\omega^{c}_{u}=(1+F_{\infty})\sum\limits_{i,j}A^{c}_{it}B^{c}_{ju}e_{i}^{c}\otimes f_{j}^{c}=\sum\limits_{i,j}A^{c}_{it}B^{c}_{ju}(e_{i}\otimes f_{j} +e_{i}^{c}\otimes f_{j}^{c})\nonumber.
\end{eqnarray}

Up to multiplication by elements in $(EE')^{\times}$, the Deligne's period then equals the determinant of the matrix 
\begin{equation}\nonumber Mat_{1}:=\left(A_{ia}B_{jb}, A^{c}_{it}B^{c}_{ju}\right)\end{equation} with $1\leq i\leq n, 1\leq j\leq n'$, $(a,b)\notin A$, $(t,u)\notin T$.\\

\bigskip

By the relation \ref{Qi}, we have $F_{\infty}\omega_{n+1-t}=Q_{n+1-t}\omega^{c}_{t}$. We get 
\begin{equation}
\sum\limits_{i=1}^{n}A_{i,n+1-t}e^{c}_{i}=Q_{n+1-t}\omega^{c}_{t}=Q_{n+1-t}\sum\limits_{i=1}^{n}A_{it}^{c}e^{c}_{i} \nonumber\end{equation}
Therefore, for all $i,j$, we obtain,
\begin{equation}\label{coefficient relation}
A_{it}^{c}=(Q_{n+1-t})^{-1}A_{i,n+1-t}, B_{ju}^{c}=(Q'_{n'+1-u})^{-1}B_{j,n'+1-u} .
\end{equation}

We then deduce that $A^{c}_{it}B^{c}_{ju}=(Q_{n+1-t})^{-1}(Q'_{n'+1-u})^{-1}A_{i,n+1-t}B_{j,n'+1-u}$. \\

Thus the Deligne's period:
\begin{equation}\nonumber
c^{+}(R(M\otimes M')) \sim_{E(M^{\#});K} det(Mat_{1})=\prod\limits_{(t,u)\notin T}((Q_{n+1-t})^{-1}(Q'_{n'+1-u})^{-1})\times det(Mat_{2})
\end{equation} where $Mat_{2}=\left(A_{ia}B_{jb},A_{i,n+1-t,j,n'+1-u}\right)$ with $1\leq i \leq n, 1\leq j\leq n'$, $(a,b)\notin A$ and $(t,u)\notin T$.

\bigskip

Recall that $(t,u)\notin T$ if and only if $(n+1-t,n'+1-u)\in A$. Therefore the index $(n+1-t,n'+1-u)$ above runs over the pairs in $A$. We see that $Mat_{2}=(A_{ia}B_{jb})$ with both $(i,j)$ and $(a,b)$ runs over all the pair in $\{1,2,\cdots,n\}\times \{1,2,\cdots,n'\}$. It is noting but $(A_{ia})\otimes (B_{jb})$.\\

Let us back to the definition of $A_{ia}$. It is the coefficients with respect to the chosen rational bases of the map $M_{B}\otimes \C\rightarrow M_{DR}\otimes \C$. Therefore $(A_{ia})\otimes (B_{jb})$ is the coefficient matrix of the comparison isomorphism $(M\otimes M')_{B}\otimes \C\rightarrow (M\otimes M')_{DR}\otimes \C$. We then get $det((A_{ia})\otimes (B_{jb}))=\delta^{Del}(M\otimes M')$ which terminates the proof.

\end{dem}

\section{Motivic periods for automorphic representations over quadratic imaginary fields}

\paragraph{Hecke character case:}
\text{}
Let $\eta$ be an arbitrary algebraic Hecke character of $K$ with infinity type $z^{a}\overline{z}^{b}$. We assume that $a\neq b$.\\

Let $M(\eta)$ be the motive associated to $\eta$ (c.f. \cite{deligne79} section $8$.) It is of Hodge type $(-a,-b)$. \\

For the motivic period for $M(\eta)$, we use $\eta$ to indicate $M(\eta)$ for simplification. For example. $\delta^{Del}(\eta):=\delta^{Del}(M(\eta))$.\\

On one hand, by Blasius's result, $c^{+}(R(M(\eta)))\sim_{E(\eta);K}p(\widecheck{\eta},1)$ if $a<b$; $c^{+}(R(M(\eta)))\sim_{E(\eta);K}p(\widecheck{\eta},\iota)$ if $a>b$.\\

On the other hand, by Proposition \ref{deligne period tensor product}, we have 
\begin{equation}
c^{+}(R(M(\eta)))\sim_{E(\eta);K} \prod\limits_{t\in A}Q_{t}(\eta)^{-1}\delta^{Del}(\eta)
\end{equation}
where $A=\{1\}$ if $-a>-b$ and $A=\emptyset$ if $-a<-b$.\\

\bigskip

Let us assume $a<b$ first. We have $Q_{1}(\eta)^{-1}\delta^{Del}(\eta)\sim_{E(\eta):K} p(\widecheck{\eta},1)$. We apply the above to $\eta^{c}$ and get $\delta^{Del}(M(\eta^{c}))\sim_{E(\eta):K} p(\widecheck{\eta^{c}},\iota)\sim_{E(\eta):K} p(\widecheck{\eta},1)$. \\

Notice that there is there is a rational paring: $M(\eta)\times M(\eta^{c})\rightarrow M(\eta_{0})(a+b)$ where $\eta_{0}$ is a Dirichlet character over $\AQ$ such that $\eta\eta^{c}=(\eta_{0}\circ N_{\AK/AQ}) ||\cdot||_{\AK}^{a+b}$. We obtain that \begin{equation}
\delta^{Del}(\eta)\times \delta^{Del}(\eta^{c})\sim_{E(\eta);K}\delta^{Del}(\eta_{0})(2\pi i)^{a+b}\end{equation}
by equation (\ref{Tate twist}).\\

We deduce by Lemma \ref{propCM} that 
\begin{eqnarray}
&& \delta^{Del}(\eta)\times \delta^{Del}(\eta^{c})\sim_{E(\eta);K}\mathcal{G}(\eta_{0})(2\pi i)^{a+b}\nonumber\\
&\sim_{E(\eta);K}&p(\eta_{0}\circ N_{\AK/\AQ},1)^{-1}p(||\cdot||_{\AK}^{a+b},1)^{-1}\sim_{E(\eta);K}p((\eta_{0}\circ N_{\AK/AQ}) ||\cdot||_{\AK}^{a+b},1)^{-1}\nonumber\\
&\sim_{E(\eta);K}&p(\eta\eta^{c},1)^{-1}\sim_{E(\eta);K}p(\widecheck{\eta},1)p(\widecheck{\eta^{c}},1)\nonumber
\end{eqnarray}

Therefore,  \begin{equation}\label{Q1 Hecke}
\delta^{Del}(\eta)\sim_{E(\eta);K} p(\widecheck{\eta^{c}},1)\text{ and then } Q_{1}(\eta)\sim_{E(\eta);K} \cfrac{p(\widecheck{\eta^{c}},1)}{p(\widecheck{\eta},1)} \sim_{E(\eta);K} p(\cfrac{\eta^{c}}{\eta},1).
\end{equation} 
If $a>b$, we follow the above procedure and can see easily the last two formulas are still true.\\

\paragraph{Conjugate self-dual case:}
\text{}
Let $\Pi$ be a regular cuspidal cohomological conjugate self-dual representation of $GL_{n}(\AK)$. We denote the infinity type of $\Pi$ by $(z^{a_{i}}\overline{z}^{-a_{i}})_{1\leq i\leq n}$ with $a_{1}>a_{2}>\cdots>a_{n}$. \\

We assume that there exists a motive $M$ over $K$ associated to $\Pi$ with coefficients in $E(M)=E(\Pi)$.\\

We know $M$ is a reguler motive of Hodge type $(-a_{i}+\frac{n-1}{2},a_{i}+\frac{n-1}{2})_{1\leq i\leq n}$. \\

We define $p_{i}:=-a_{n+1-i}+\frac{n-1}{2}$, $q_{i}:=n-1-p_{i}=a_{n+1-i}++\frac{n-1}{2}$, $p_{i}^{c}=a_{i}+\frac{n-1}{2}=q_{n+1-i}$, $q_{i}^{c}:=n-1-p_{i}^{c}$ for all $1\leq i\leq n$. \\

We define $\omega_{i}\in M^{p_{i},q_{i}}$, $\omega_{i}^{c}\in (M^{c})^{p_{i}^{c},q_{i}^{c}}$ as in the previous sections.\\

If $\Pi$ is conjugate self-dual, then $M(\Pi)$ is polarized. The polarization on the De Rham realization induces an $E(M)$-rational perfect parings $<,>$:
\begin{equation}
M(\Pi)^{p_{i},q_{i}}\otimes (M(\Pi)^{c})^{p^{c}_{n+1-i},q^{c}_{n+1-i}} \rightarrow E(M)(1-n)\cong E(M) \nonumber .\\
\end{equation}
We may assume that $<\omega_{i},\omega^{c}_{n+1-i}>=1$ by adjusting $\omega^{c}_{n+1-i}$ with multiplication by elements in $E(M)^{\times}$.\\

Let $1\leq i\leq n$. We write $Q_{i}(\Pi):=Q_{i}(M(\Pi))$ as we did for $\eta$. The motivic period $Q_{i}(\Pi)$ then equals $<w_{i},F_{\infty}w_{n+1-i}>$ up to multiplication by elements in $E(\Pi)^{\times}$. \\

\paragraph{ conjugate self-dual case:}
\text{}
In the general case, we write $\Pi=\Pi'\otimes \eta$ with $\Pi'$ conjugate self-dual and $\eta$ be an algebraic Hecke character of $\AK$. \\

We take $\omega_{i}^{0}\in M(\Pi')^{p_{i}(\Pi'),q_{i}(\Pi')}$ as before. Let $\omega$ be a base of $M(\eta)_{DR}$ and $\omega^{c}$ be a base of $M(\eta^{c})_{DR}$. We know $F_{\infty}(\omega)=Q_{1}(\eta)\omega^{c}$ up to multiplication by elements in $E(\eta)^{\times}$.\\

Then $(\omega_{i}:=\omega_{i}^{0}\otimes \omega)_{1\leq i\leq n}\in M^{p_{i},q_{i}}$ which is equivalent to a rational element in $F^{p_{i}(M)}$ modulo $F^{p_{i-1}}(M)\otimes \C$. We have similar properties for $(\omega_{i}^{c}:=\omega_{i}^{0,c}\otimes \omega^{c})_{1\leq i\leq n}$. \\

Moreover,
\begin{equation}\label{motive main relation}
F_{\infty}(\omega_{i}\otimes \omega)=Q_{i}(\Pi')Q_{1}(\eta)(\omega'_{n+1-i}\otimes \omega').
\end{equation}

The motivic period for $\Pi$ then equals $Q_{i}(\Pi):=Q_{i}(\Pi')Q_{1}(\eta)$ for $1\leq i\leq n$ up to multiplication by elements in $E(\Pi)^{\times}$.

\section{Deligne's conjecture for automorphic pairs over quadratic imaginary fields}

Let $\Pi$ (resp. $\Pi'$) be a regular cuspidal cohomological  conjugate self-dual representation of $GL_{n}(\AK)$ (resp. $GL_{n'}(\AK)$). We denote the infinity type of $\Pi$ (resp. $\Pi'$) by $(z^{a_{i}}\overline{z}^{-\omega(\Pi)-a_{i}})_{1\leq i\leq n}$ (resp. $(z^{b_{j}}\overline{z}^{-\omega(\Pi)-b_{j}})_{1\leq j\leq n'}$) with $a_{1}>a_{2}>\cdots>a_{n}$ (resp. $b_{1}>b_{2}>\cdots>b_{n'}$). We suppose $\Pi\times \Pi'$ is regular, i.e. $a_{i} +b_{j}\neq -\frac{\omega(\Pi)+\omega(\Pi')}{2}$ for all $i,j$.\\

We assume that there exists a motive $M$ (resp. $M'$) over $K$ associated to $\Pi$ (resp. $\Pi'$) with coefficients in $E(M)$ (resp. $E(M')$). \\

We know $M$ (resp. $M'$) is reguler motive of Hodge type $(-a_{i}+\frac{n-1}{2},a_{i}+\omega(\Pi)+\frac{n-1}{2})_{1\leq i\leq n}$ (resp. $(-b_{j}+\frac{n'-1}{2},b_{j}+\omega(\Pi')+\frac{n-1}{2})_{1\leq j\leq n'}$ ). \\

We define $p_{i}:=-a_{n+1-i}+\frac{n-1}{2}$ and $p_{i}^{c}=\omega(\Pi)+n-1-p_{n+1-i}=a_{i}+\omega(\Pi)+\frac{n-1}{2}$ for all $1\leq i\leq n$. We define like this to guarantee that $p_{1}>p_{2}>\cdots>p_{n}$ and $p_{1}^{c}>p_{2}^{c}>\cdots>p_{n}^{c}$ as in the previous sections. Similarly, we define $r_{j}:=-b_{n'+1-j}+\frac{n'-1}{2}$, $r_{j}^{c}:=b_{j}+\omega(\Pi')+\frac{n'-1}{2}$ for all $1\leq j\leq n'$.\\

Proposition \ref{deligne period tensor product} implies that:
\begin{prop} The Deligne's period of $(M\otimes M')$ satisfies:
\begin{equation}\label{deligne for pair rough}
c^{+}(R(M(\Pi)\otimes M(\Pi')))\sim_{E(M)E(M');K}\left(\prod_{(t,u)\in A(M,M')}Q_{t}(\Pi)^{-1}Q_{u}(\Pi')^{-1}\right)\delta^{Del}(M\otimes M').
\end{equation} 
where the set $A(M,M')=\{(t,u)\mid p_{t}+r_{u}> \frac{\omega(\Pi)+\omega(\Pi')+(n-1)+(n'-1)}{2}\}$.
\end{prop}

\bigskip 

Recall $p_{t}=-a_{n+1-t}+\frac{n-1}{2}$ and $r_{u}:=-b_{n'+1-u}+\frac{n'-1}{2}$. We obtain that $A(M,M')=\{(t,u)\mid p_{t}+r_{u}> \frac{\omega(\Pi)+\omega(\Pi')+(n-1)+(n'-1)}{2}\}
=\{(t,u)\mid a_{n+1-t}+b_{n'+1-u}< -\frac{\omega(\Pi)+\omega(\Pi')}{2}\}$.\\

Therefore, 
\begin{equation}
\prod_{(t,u)\in A(M,M')}Q_{t}(\Pi)^{-1}=\prod\limits_{t=1}^{n}Q_{t}(\Pi)^{-\#\{u\mid b_{n+1-u}<-a_{n+1-t}-\frac{\omega(\Pi)+\omega(\Pi')}{2}\}}
\end{equation}
In this section, we define $sp(j):=sp(j,\Pi;\Pi')$ for $0\leq j\leq n$ and $sp'(k):=sp(k,\Pi';\Pi)$ for $0\leq k\leq n'$. Recall $sp(j)$ are the lengths of different parts of $b_{1}>b_{2}>\cdots>b_{n'}$ separated by $-a_{n}-\frac{\omega(\Pi)+\omega(\Pi')}{2}>-a_{n-1}-\frac{\omega(\Pi)+\omega(\Pi')}{2}>\cdots>-a_{1}-\frac{\omega(\Pi)+\omega(\Pi')}{2}$.

Therefore $\#\{u\mid b_{n+1-u}<-a_{n+1-t}-\frac{\omega(\Pi)+\omega(\Pi')}{2}\}=\#\{u\mid b_{u}<-a_{n+1-t}-\frac{\omega(\Pi)+\omega(\Pi')}{2}\}=sp(t)+sp(t+1)+\cdots+sp(n)$, we have $\prod_{(t,u)\in A(M,M')}Q_{t}(\Pi)^{-1}$
\begin{eqnarray}
&=&Q_{1}(\Pi)^{-sp(1)-sp(2)-\cdots-sp(n)}Q_{2}(\Pi)^{-sp(2)-sp(3)-\cdots-sp(n)}\cdots Q_{n}(\Pi)^{-sp(n)}\nonumber\\
&=& [Q_{1}(\Pi)^{-1}]^{sp(1)}[Q_{1}(\Pi)^{-1}Q_{2}(\Pi)^{-1}]^{sp(2)}\cdots[Q_{1}(\Pi)^{-1}Q_{2}(\Pi)^{-1}\cdots Q_{n}(\Pi)^{-1}]^{sp(n)}\nonumber
\end{eqnarray}

We define $Q_{\leq j}(\Pi)=Q_{1}(\Pi)^{-1}Q_{2}(\Pi)^{-1}\cdots Q_{j}(\Pi)^{-1}$ for $1\leq j\leq n$ and $Q_{\leq 0}(\Pi)=1$. We define $Q_{\leq k}(\Pi')$ similarly for $0\leq k\leq n'$.\\

We have obtained that 
\begin{equation}\label{Q_leq}
\prod_{(t,u)\in A(M,M')}Q_{t}(\Pi)^{-1}Q_{u}(M')^{-1}=\prod\limits_{j=0}^{n}Q_{\leq j}(\Pi)^{sp(j)}\prod\limits_{k=0}^{n'}Q_{\leq k}(\Pi')^{sp'(k)}.\end{equation}

\bigskip

We define $\Delta(M(\Pi))=\Delta(\Pi):=(2\pi i)^{\frac{n(n-1)}{2}}\delta^{Del}(\Pi)$. In fact, let $\xi_{\Pi}$ be the central character of $\Pi$. Since $\Lambda^{n}M(\Pi)\equiv M(\xi_{\Pi})(-\frac{n(n-1)}{2})$, we have $\delta^{Del}(\Pi)=\delta^{Del}(\Lambda^{n}M(\Pi))=\delta^{Del}(\xi_{\Pi})(2\pi i)^{-\frac{n(n-1)}{2}}$. Therefore, 
\begin{equation}\label{Delta}
\Delta(\Pi)\sim_{E(\Pi):K}\delta^{Del}(\xi_{\Pi}).
\end{equation} 

We define $\Delta(\Pi')$ similarly. We have $\delta^{Del}(M\otimes M')=\delta^{Del}(\Pi)^{n'}\delta^{Del}(\Pi')^{n}=(2\pi i)^{-nn'(n+n'-2)}\Delta(\Pi)^{n'}\Delta(\Pi')^{n}$. Since $\sum\limits_{i=0}^{n}sp(i)=n'$ and $\sum\limits_{i=0}^{n}sp(i)=n'$, we know:
\begin{equation}\label{Delta}
\delta^{Del}(M\otimes M')=(2\pi i)^{-\frac{nn'(n+n'-2)}{2}}\prod\limits_{j=0}^{n}\Delta(\Pi)^{sp(j)}\prod\limits_{k=0}^{n'}\Delta(\Pi')^{sp'(k)}.
\end{equation}

At last, we define for all $0\leq j\leq n$ that\begin{equation}\label{Q^j}
Q^{(j)}(\Pi)=Q_{\leq j}(\Pi)\times \Delta(\Pi)\sim_{E(\Pi);K} Q_{1}(\Pi)^{-1}\cdots Q_{j}(\Pi)^{-1}\delta^{Del}(\xi_{\Pi})
\end{equation}

We define $Q^{(k)}(\Pi')$ for $1\leq k\leq n'$ similarly. Comparing (\ref{deligne for pair rough}) with (\ref{Q_leq}) and (\ref{Delta}), we have
\begin{equation}
c^{+}(R(M\otimes M'))\sim_{E(\Pi)E(\Pi');K}(2\pi i)^{\frac{-nn'(n+n'-2)}{2}}\prod\limits_{j=0}^{n}Q^{(j)}(\Pi)^{sp(j,\Pi;\Pi')}
\prod\limits_{k=0}^{n'}Q^{(k)}(\Pi')^{sp(k,\Pi';\Pi)}\nonumber
\end{equation}

\bigskip

We can now state Deligne's conjecture for automorphic pairs:

\begin{conj}\label{Deligne automorphic}
Let $n$ and $n'$ be two positive integers. Let $\Pi$ and $\Pi'$ be regular cohomological cuspidal representation of $GL_{n}(\AK)$ and $GL_{n'}(\AK)$ respectively which are  conjugate self-dual. We suppose that $\Pi\otimes\Pi'$ is regular.

We assume that there exists motives $M$ and $M'$ over $K$ associated to $\Pi$ and $\Pi'$ respectively.

Let $m\in \Z+\frac{n+n'}{2}$ be critical for $\Pi\otimes \Pi'$. It is equivalent to saying that $m+\frac{n+n'-2}{2}$ is critical for $M\otimes M'$. Deligne's conjecture predicts that:
\begin{eqnarray}
&L(m,\Pi\times \Pi')=L(m+\frac{n+n'-2}{2},M\otimes M')&\nonumber\\
&\sim_{E(\Pi)E(\Pi');K}(2\pi i)^{nn'm}\prod\limits_{j=0}^{n}Q^{(j)}(\Pi)^{sp(j,\Pi;\Pi')}
\prod\limits_{k=0}^{n'}Q^{(k)}(\Pi')^{sp(k,\Pi';\Pi)}\nonumber&
\end{eqnarray}

\end{conj}

\section{The picture for general CM fields}

Let $F$ be a CM field containing $K$ and $F^{+}$ be the maximal totally real subfield of $F$. 

Let $M$ be a motive over $F$ with coefficients in $E(M)$ of dimension $n$ and pure of weight $\omega(M)$.\\

For each $\sigma\in \Sigma_{F}$, we may define the motivic period $\delta^{Del}(M,\sigma)$ and $Q_{i}(M,\sigma)$ as in Section \ref{definition of motivic period}. We write the Hodge type of $M$ at $\sigma$ by $(p_{i}(\sigma),q_{i}(\sigma))_{1\leq i\leq n}$.\\

We take $M'$ another motive over $F$ with coefficients in $E(M')$ of dimension $n'$ and pure of weight $\omega(M')$. Similarly, we write the Hodge type of $M'$ at $\sigma$ by $(r_{j}(\sigma),s_{j}(\sigma))_{1\leq i\leq n'}$. \\

We assume that $p_{i}(\sigma)+r_{j}(\sigma)\neq \cfrac{\omega(M)+\omega(M')}{2}$ for all $\sigma$, $i$, $j$. \\

Define $A(M,M')(\sigma)=\{(a,b) \mid p_{a}(\sigma)+r_{b}(\sigma)> \frac{\omega}{2}\}$. \\

For any CM type  $\Phi$ of $F$, we have $Res_{F/\Q}(M\otimes M')=(\bigoplus\limits_{\sigma\in\Phi}M^{\sigma}\otimes M'^{\sigma})\oplus (\bigoplus\limits_{\sigma\in\Phi}M^{\sigma^{c}}\otimes M'^{\sigma^{c}})$.\\

The proof of Proposition \ref{deligne period tensor product} can be easily generalized to the CM field case and we get:

\begin{prop}\label{tensor product motive general}
Let $M$, $M'$ be motive over $F$ with coefficients in $E$ and $E'$ respectively. We assume that $M\otimes M'$ has no $(\omega/2,\omega/2)$-class. We have:
\begin{eqnarray}
&c^{+}(Res_{F/\Q}(M\otimes M')) \sim_{E(M)E(M');K}&\\
& \prod\limits_{\sigma\in\Psi} \left(\prod\limits_{(t,u)\in A(M,M')(\sigma)}Q_{t}(M,\sigma)^{-1}Q_{u}(M',\sigma)^{-1}\right)\prod\limits_{\sigma\in\Psi}\delta^{Del}(M\otimes M',\sigma)\nonumber&.
\end{eqnarray} 
\end{prop}

\bigskip

Let us assume that $M$ and $M'$ are motives associated to certain representations $\Pi$ and $\Pi'$ respectively. We still write the motivic period $Q_{i}(M(\Pi))$ as $Q_{i}(\Pi)$ for simplicity.\\

We define $Q_{\leq j}(\Pi,\sigma):=\prod\limits_{i=1}^{j}Q_{i}(\Pi,\sigma)$. We have (c.f. equation (\ref{Q_leq}))
\begin{equation}
\prod_{(t,u)\in A(M,M')(\sigma)}Q_{t}(\Pi)^{-1}Q_{u}(M')^{-1}=\prod\limits_{j=0}^{n}Q_{\leq j}(\Pi,\sigma)^{sp(j,\Pi;\Pi', \sigma)}\prod\limits_{k=0}^{n'}Q_{\leq k}(\Pi',\sigma)^{sp'(k,\Pi';\Pi,\sigma)}.\end{equation}

Recall that $\Lambda^{n}M(\Pi)\cong M(\xi_{\Pi})(-\frac{n(n-1)}{2})$. We have \begin{equation}
\delta^{Del}(\Pi,\sigma) \sim_{E(\Pi);K} (2\pi i)^{-\frac{n(n-1)}{2}}\delta^{Del}(\xi_{\Pi},\sigma).
\end{equation}

We define $\Delta(\Pi,\sigma):=(2\pi i)^{\frac{n(n-1)}{2}}\delta^{Del}(\Pi,\sigma)=\delta^{Del}(\xi_{\Pi},\sigma)$ and $Q^{(j)}(\Pi,\sigma):=Q_{\leq j}(\Pi,\sigma)\Delta(\Pi,\sigma)$ as before. We have:
\begin{eqnarray}
c^{+}(Res_{F/\Q}(M\otimes M'))\sim_{E(\Pi)E(\Pi');K}(2\pi i)^{\frac{-nn'd(n+n'-2)}{2}}\times \nonumber\\
\prod\limits_{\sigma\in\Psi} \prod\limits_{j=0}^{n}Q^{(j)}(\Pi,\sigma)^{sp(j,\Pi;\Pi',\sigma)}
\prod\limits_{k=0}^{n'}Q^{(k)}(\Pi',\sigma)^{sp(k,\Pi';\Pi,\sigma)}\nonumber
\end{eqnarray}

In particular, if we take $\Psi=\Sigma_{F;K}$, Deligne's conjecture predicts that:
\begin{conj}\label{Deligne automorphic CM}
Let $n$ and $n'$ be two positive integers.
Let $\Pi$ and $\Pi'$ be regular cohomological cuspidal representation of $GL_{n}(\AF)$ and $GL_{n'}(\AF)$ respectively which are  conjugate self-dual. We suppose that $\Pi\otimes\Pi'$ is regular.

We assume that there exists motives $M$ and $M'$ over $F$ associated to $\Pi$ and $\Pi'$ respectively.

If $m\in \Z+\frac{n+n'}{2}$ is critical for $\Pi\times \Pi'$ then 
\begin{eqnarray}
&L(m,\Pi\times \Pi')=L(m+\frac{n+n'-2}{2},M\otimes M')&\nonumber\\
&\sim_{E(\Pi)E(\Pi');K}(2\pi i)^{mnn'd}\prod\limits_{\sigma\in\Sigma_{F;K}} (\prod\limits_{j=0}^{n}Q^{(j)}(\Pi,\sigma)^{sp(j,\Pi;\Pi',\sigma)}
\prod\limits_{k=0}^{n'}Q^{(k)}(\Pi',\sigma)^{sp(k,\Pi';\Pi,\sigma)})\nonumber&
\end{eqnarray}
\end{conj}

\section{Motivic periods for Hecke characters over CM fields}

Let $\eta,\eta'$ be two algebraic Hecke character of $F$. We assume that $\eta\eta'$ is critical and then is compatible with a CM type $\Psi(\eta\eta')$. Proposition \ref{tensor product motive general} can be rewritten as 
\begin{eqnarray}\label{Hecke character pair CM}
&c^{+}(Res_{F/\Q}(M(\eta)\otimes M(\eta'))) \sim_{E(\eta)E(\eta');K}&\\
&\prod\limits_{\sigma\in\Psi} \left(\prod\limits_{\sigma\in \Psi\cap \Psi(\eta\eta')}Q_{1}(\eta,\sigma)^{-1}Q_{1}(\eta',\sigma)^{-1}\right)\prod\limits_{\sigma\in\Psi}\delta^{Del}(M\otimes M',\sigma)\nonumber&.
\end{eqnarray} 

Fix any Hecke character $\eta$. We may take $\eta'$ such that $\Psi(\eta\eta')=\Psi^{c}$.

Equation (\ref{Hecke character pair CM}) implies that
\begin{equation}
c^{+}(Res_{F/\Q})(M(\eta)\otimes M(\eta'))\sim_{E(\eta)E(\eta');K}\prod\limits_{\sigma\in\Psi}\delta^{Del}(\eta,\sigma) \prod\limits_{\sigma\in\Psi}\delta^{Del}(\eta',\sigma).
\end{equation}

On the other hand, by Blasius's result, we have:
\begin{eqnarray}
c^{+}(Res_{F/\Q}(M(\eta)\otimes M(\eta')))\sim_{E(\eta)E(\eta');K} p(\widecheck{\eta\eta'},\Psi^{c})\nonumber\\
\sim_{E(\eta)E(\eta');K}\prod\limits_{\sigma\in\Psi}p(\widecheck{\eta},\sigma^{c}) \prod\limits_{\sigma\in\Psi}p(\widecheck{\eta'},\sigma^{c}).\nonumber
\end{eqnarray}

Let $\eta'$ vary. We get for any $CM$ type $\Psi$ that:
\begin{equation}\label{delta CM}
\prod\limits_{\sigma\in\Psi}\delta^{Del}(\eta',\sigma)\sim_{E(\eta);K} \prod\limits_{\sigma\in\Psi}p(\widecheck{\eta},\sigma^{c}).
\end{equation}

Let $\Psi$ vary now. It is easy to deduce that there exists $\zeta_{d}$, an $d$-th root of unity, such that $\delta^{Del}(\eta',\sigma)\sim_{E(\eta);K} \zeta_{d}p(\widecheck{\eta},\sigma^{c})$ for all $\sigma\in\Sigma$.\\

For simplicity, we assume that $E(\Pi)$ contains all $d$-th roots of unity then we get $\delta^{Del}(\eta',\sigma)\sim_{E(\eta);K} p(\widecheck{\eta},\sigma^{c})$ for all $\sigma\in\Sigma$.\\

We can now calculate $Q_{1}(\eta,\sigma)$. Let $\sigma_{0}$ be in $\Psi$. We take $\eta'$ such that $\Psi_{\eta\eta'}=\{\sigma\}\cup (\Psi-\{\sigma_{0}\})^{c}$. 
Equation (\ref{Hecke character pair CM}) implies that
\begin{eqnarray}&\nonumber
c^{+}(Res_{F/\Q})(M(\eta)\otimes M(\eta'))\sim_{E(\eta)E(\eta');K} Q_{1}(\eta,\sigma)^{-1}Q_{1}(\eta',\sigma)^{-1}\times& \\\nonumber&\prod\limits_{\sigma\in(\Psi-\{\sigma_{0}\})}\delta^{Del}(\eta,\sigma) \prod\limits_{\sigma\in(\Psi-\{\sigma_{0}\})}\delta^{Del}(\eta',\sigma)\times \delta^{Delta}(\eta,\sigma_{0})\delta^{Delta}(\eta',\sigma_{0}).&
\end{eqnarray}

Blasius's result implies that 
\begin{equation}
c^{+}(Res_{F/\Q})(M(\eta)\otimes M(\eta'))\sim_{E(\eta)E(\eta');K}(-1)^{\epsilon(\Psi)} \prod\limits_{\sigma\in(\Psi-\{\sigma_{0}\})}p(\widecheck{\eta\eta'},\sigma^{c})\times p(\widecheck{\eta\eta'},\sigma_{0}).
\end{equation}

Along with equation (\ref{delta CM}), we obtain that \begin{equation}
Q_{1}(\eta,\sigma_{0})\sim_{E(\eta);K} p(\widecheck{\eta},\sigma_{0}^{c})p(\widecheck{\eta},\sigma_{0})^{-1}\sim_{E(\eta);K} p(\cfrac{\eta^{c}}{\eta},\sigma_{0}).
\end{equation}

\chapter{Factorization of arithmetic automorphic periods and a conjecture}
We want to show that the arithmetic automorphic periods can be factorized as products of local periods over infinite places. We may assume that $\Pi$ is conjugate self-dual in this Chapter. The essential conjugate self-dual case then follows by Definition \ref{definable periods} and the fact that the CM periods is factorable.

\section{Basic lemmas}
 Let $X$, $Y$ be two sets and $Z$ be a multiplicative abelian group. We will apply the result of this section to $Z=\C^{\times}/E^{\times}$ where $E$ is a proper number field.
 
\begin{lem}
 Let $f$ be a map from $X\times Y$ to $Z$. The following two statements are equivalent:
 \begin{enumerate}
 \item There exists two maps $g:X\rightarrow Z$ and $h:Y\rightarrow Z$ such that $f(x,y)=g(x)h(y)$ for all $(x,y)\in X\times Y$.
 \item For all $x,x'\in X$ and $y,y'\in Y$, we have $f(x,y)f(x',y')=f(x,y')f(x',y)$.
 \end{enumerate} 
 Moreover, if the above equivalent statements are satisfied, the maps $g$ and $h$ are unique up to scalars.
 
 \end{lem}
 
 \begin{dem}
 The direction that $1$ implies $2$ is trivial. Let us prove the inverse. 
 We fix any $y_{0}\in Y$ and define $g(x):=f(x,y_{0})$ for all $x\in X$. We then fix any $x_{0}\in X$ and define $h(y):=\cfrac{f(x_{0},y)}{g(x_{0})}=\cfrac{f(x_{0},y)}{f(x_{0},y_{0})}$. 
 
For any $x\in X$ and $y\in Y$, Statement $2$ tells us that $f(x,y)f(x_{0},y_{0})=f(x,y_{0})f(x_{0},y)$. Therefore $f(x,y)=f(x,y_{0})\times\cfrac{f(x_{0},y)}{f(x_{0},y_{0})}=g(x)h(y)$ as expected.

  \end{dem}

 Let $n$ be a positive integer and $X_{1},\cdots,X_{n}$ be some sets. Let $f$ be a map from $X_{1}\times X_{2}\times \cdots\times X_{n}$ to $Z$.\\
 
 The following corollary can be deduced from the above Lemma by induction on $n$.
 
  \begin{cor}\label{factorization lemma}
 
The following two statements are equivalent:
 \begin{enumerate}
 \item There exists some maps $f_{k}:X_{k}\rightarrow Z$ for $1\leq k\leq n$ such that $f(x_{1},x_{2},\cdots,x_{n})=\prod\limits_{1\leq k\leq n} f_{k}(x_{k})$ for all $x_{k}\in X_{k}$, $1\leq k\leq n$.
 \item Given any $x_{j},x'_{j}\in X_{j}$ for each $1\leq j\leq n$, we have 
$
f(x_{1},x_{2},\cdots,x_{n})\times f(x'_{1},x'_{2},\cdots,x'_{n})=f(x_{1},\cdots,x_{k-1},x'_{k},x_{k+1},x_{n})\times f(x'_{1},\cdots x'_{k-1},x_{k},x'_{k+1},\cdots,x'_{n})
$ for any $1\leq k\leq n$.
 \end{enumerate} 
 Moreover, if the above equivalent statements are satisfied then for any $\lambda_{1},\cdots,\lambda_{n}\in Z$ such that $\lambda_{1}\cdots\lambda_{n}=1$, we have another factorization $f(x_{1},\cdots,x_{n})=\prod\limits_{1\leq i\leq n}(\lambda_{i}f_{i})(x_{i})$. Each factorization of $f$ is of the above form.\\

We fix $a_{i}\in X_{i}$ for each i and $c_{1},\cdots,c_{n}\in Z$ such that $f(a_{1},\cdots,a_{n})=c_{1}\cdots c_{n}$. If the above equivalent statements are satisfied then there exists a unique factorization such that $f_{i}(a_{i})=c_{i}$.

  \end{cor}
  \bigskip

\begin{rem}  \label{remark factorization}
If $\#X_{k}\geq 3$ for all $k$, it is enough to verify the condition in statement $2$ of the above corollary in the case $x_{j}\neq x'_{j}$ for all $1\leq j\leq n$.

In fact, when $\#X_{k}\geq 3$ for all $k$, for any $1\leq j\leq n$ and any $y_{j},y'_{j}\in X_{j}$, we may take $x_{j} \in X_{j}$ such that $x_{j}\neq y_{j}$, $x_{j}\neq y_{j}'$.

We fix any $1\leq k\leq n$. If statement $2$ is verified when $x_{j}\neq x'_{j}$ for all $j$ then for any $y_{k}\neq y'_{k}$, we have
\begin{eqnarray}
&&f(y_{1},y_{2},\cdots,y_{n})f(y'_{1},y'_{2},\cdots,y'_{n}) f(x_{1},x_{2},\cdots,x_{n}) \nonumber\\
&=&f(y_{1},y_{2},\cdots,y_{n})f(y'_{1},\cdots y'_{k-1},x_{k},y'_{k+1},\cdots,y'_{n})
 f(x_{1},\cdots,x_{k-1},y'_{k},x_{k+1},\cdots x_{n}) \nonumber\\
 &=&f(y_{1},y_{2},\cdots,y_{n})f(x_{1},\cdots,x_{k-1},y'_{k},x_{k+1},\cdots x_{n}) f(y'_{1},\cdots y'_{k-1},x_{k},y'_{k+1},\cdots,y'_{n})
 \nonumber\\
 &=& f(y_{1},\cdots,y_{k-1},y'_{k},y_{k+1},\cdots,y_{n})f(x_{1},\cdots,x_{k-1},y_{k},x_{k+1},\cdots,x_{n})\times \nonumber\\
 &&f(y'_{1},\cdots y'_{k-1},x_{k},y'_{k+1},\cdots,y'_{n})\nonumber\\
 &=& f(y_{1},\cdots,y_{k-1},y'_{k},y_{k+1},\cdots,y_{n})f(y'_{1},\cdots y'_{k-1},y_{k},y'_{k+1},\cdots,y'_{n})f(x_{1},x_{2},\cdots,x_{n}).\nonumber
\end{eqnarray}
We have assumed $y_{k}\neq y'_{k}$ to guarantee that each time we apply the formula in Statement $2$, the coefficients satisfy $x_{j}\neq x'_{j}$ for all $1\leq j\leq n$.\\

Therefore $f(y_{1},y_{2},\cdots,y_{n})\times f(y'_{1},y'_{2},\cdots,y'_{n})=f(y_{1},\cdots,y_{k-1},y'_{k},y_{k+1},\cdots,y_{n})\times f(y'_{1},\cdots y'_{k-1},y_{k},y'_{k+1},\cdots,y'_{n})$ if $y_{k}\neq y'_{k}$. If $y_{k}=y'_{k}$, this formula is trivially true. \\

We conclude that we can weaken the condition in Statement $2$ of the above Corollary to $x_{j}\neq x'_{j}$ for all $1\leq j\leq n$ when $\#X_{k}\geq 3$ for all $k$. We will verify this weaker condition in the application to the factorization of arithmetic automorphic periods.
\end{rem}
  
\section{Formula for the Whittaker period: even dimensional}
 Let $\Pi$ be a regular cuspidal representation of $GL_{n}(\AF)$ as in Theorem \ref{main theorem CM} with infinity type $(z^{a_{i}(\sigma)}\overline{z}^{-a_{i}(\sigma)})_{1\leq i\leq n}$ at $\sigma\in \Sigma$. We may assume that $a_{1}(\sigma)>a_{2}(\sigma)>\cdots>a_{n}(\sigma)$ for all $\sigma\in\Sigma$.
 
Recall that we say $\Pi$ is \textbf{$N$-regular} if $a_{i}(\sigma)-a_{i+1}(\sigma)\geq N$ for all $1\leq i\leq n-1$ and $\sigma\in\Sigma$.\\

 For $1\leq u \leq n-1$, let $\chi_{u}$ be an algebraic conjugate self-dual Hecke character of $F$ with infinity type $z^{k_{u}(\sigma)}\overline{z}^{-k_{u}(\sigma)}$ at $\sigma\in \Sigma$.\\

Let us first consider the case $n$ even. In this case, $a_{i}(\sigma)\in \Z+\cfrac{1}{2}$ for all $1\leq i\leq n$ and all $\sigma\in \Sigma$. We assume the following hypothesis:

\begin{hyp}\label{hyp infinity type}\textbf{Even dimensional}

For all $\sigma\in\Sigma$, the numbers $\{k_{u}(\sigma)\mid 1\leq u\leq n-1\}$ lie in the $n-1$ gaps between $-a_{n}(\sigma)>-a_{n-1}(\sigma)>\cdots>-a_{1}(\sigma)$.
\end{hyp}

  We define $\Pi^{\#}$ to be the Langlands sum of $\chi_{u}$, $1\leq u\leq n-1$. It is an algebraic regular automorphic representation of $GL_{n-1}(\AF)$. It follows by the above hypothesis that $(\Pi,\Pi^{\#})$ is in good position. By Proposition \ref{Whittaker period theorem CM} we have
  
\begin{equation}\label{CM Whittaker period}
L(\cfrac{1}{2}+m,\Pi\times \Pi^{\#})\sim_{E(\Pi)E(\Pi^{\#});K} p(\Pi)p(\Pi^{\#})p(m,\Pi_{\infty},\Pi^{\#}_{\infty})
\end{equation}
where $p(m,\Pi_{\infty},\Pi^{\#}_{\infty})$ is a complex number which depends on $m,\Pi_{\infty}$ and $\Pi^{\#}_{\infty}$. 
 
 \paragraph{Simplification of $L(\cfrac{1}{2}+m,\Pi\times \Pi^{\#})$:}
 
 \text{}\\
 
 Since $\Pi^{\#}$ is the Langlands sum of $\chi_{u}$, $1\leq u\leq n-1$, we have \begin{equation}\nonumber L(\cfrac{1}{2}+m,\Pi\times \Pi^{\#}) =\prod\limits_{1\leq u\leq n-1}L(\cfrac{1}{2}+m,\Pi\times \chi_{u}).
 \end{equation}
  We then apply Theorem \ref{main result CM} to the right hand side and get:
 \begin{eqnarray}
 &&L(\cfrac{1}{2}+m,\Pi\times \Pi^{\#}) =\prod\limits_{1\leq u\leq n-1}L(\cfrac{1}{2}+m,\Pi\times \chi_{u})\nonumber\\
 &\sim_{E(\Pi)E(\Pi^{\#});K}&\prod\limits_{1\leq u\leq n-1}[(2\pi i)^{d(m+\frac{1}{2})n}P^{(I(\Pi,\chi_{u}))}(\Pi)\prod\limits_{\sigma\in\Sigma}p(\widecheck{\chi_{u}},\sigma)^{I_{u}(\sigma)} p(\widecheck{\chi_{u}},\overline{\sigma})^{n-I_{u}(\sigma)} ]\nonumber
 \end{eqnarray}
Here we write $I_{u}$ for $I(\Pi,\chi_{u})$. In particular, $I_{u}(\sigma)=\#\{i\mid-a_{i}(\sigma)> k_{u}(\sigma)\}$ for $\sigma\in\Sigma$.\\
 
 Note that $\chi_{u}$ is conjugate self-dual, we have $p(\widecheck{\chi_{u}},\overline{\sigma})\sim_{E(\Pi^{\#});K}p(\widecheck{\chi_{u}^{c}},\sigma)\sim_{E(\Pi^{\#});K}p(\widecheck{\chi_{u}^{-1}},\sigma)\sim_{E(\Pi^{\#});K}p(\widecheck{\chi_{u}},\sigma)^{-1}$. We deduce that:
 
 \begin{equation}\label{the left hand side step 1}
L(\cfrac{1}{2}+m,\Pi\times \Pi^{\#}) \sim_{E(\Pi)E(\Pi^{\#});K} (2\pi i)^{d(m+\frac{1}{2})n(n-1)}\prod\limits_{1\leq u\leq n-1} [P^{(I_{u})}(\Pi)\prod\limits_{\sigma\in\Sigma}p(\widecheck{\chi_{u}},\sigma)^{2I_{u}(\sigma)-n}  ]\end{equation}
 
 \bigskip

  \paragraph{Calculate $p(\Pi^{\#})$:}
 By Proposition \ref{Shahidi Whittaker}, there exists a constant $\Omega(\Pi^{\#}_{\infty})\in \C^{\times}$ well defined up to $E(\Pi^{\#})^{\times}$ such that 
 \begin{equation}
 p(\Pi^{\#})\sim_{E(\Pi^{\#});K} \Omega(\Pi^{\#}_{\infty})\prod\limits_{1\leq u<v\leq n-1}L(1,\chi_{u}\chi_{v}^{-1}).
 \end{equation}
 
 By Blasius's result, we have:
 \begin{equation}\nonumber
 L(1,\chi_{u}\chi_{v}^{-1}) \sim_{E(\Pi^{\#});K} (2\pi i)^{d}\prod\limits_{\sigma\in\Sigma} p(\widecheck{\chi_{u}\chi_{v}^{-1}},\sigma')
 \end{equation}
 
 If $k_{u}(\sigma)<k_{v}(\sigma)$ we have $\sigma'=\sigma$ and $p(\widecheck{\chi_{u}\chi_{v}^{-1}},\sigma') \sim_{E(\chi_{u});K} p(\widecheck{\chi_{u}},\sigma) p(\widecheck{\chi_{v}},\sigma)^{-1}$.\\

Otherwise we have $\sigma'=\overline{\sigma}$ and
$
p(\widecheck{\chi_{u}\chi_{v}^{-1}},\sigma')\sim_{E(\chi_{u});K}  p(\widecheck{\chi_{u}},\sigma)^{-1}p(\widecheck{\chi_{v}},\sigma)
$.\\

 Therefore, we have the Whittaker period $p(\Pi^{\#})$
 \begin{equation}\label{the right hand side step 1}
  \sim_{E(\Pi^{\#});K}(2\pi i)^{\frac{d(n-1)(n-2)}{2}}\Omega(\Pi^{\#}_{\infty}) \prod\limits_{1\leq u\leq n-1}\prod\limits_{\sigma\in \Sigma} p(\widecheck{\chi_{u}},\sigma)^{\#\{v\mid k_{v}(\sigma)>k_{u}(\sigma)\}-\#\{v\mid k_{v}(\sigma)<k_{u}(\sigma)\}}
 \end{equation}
 
 We know $\#\{v\mid k_{v}(\sigma)<k_{u}(\sigma)\}=n-2-\#\{v\mid k_{v}(\sigma)>k_{u}(\sigma)\}$.
 
 Moreover, by Hypothesis \ref{hyp infinity type}, we have $\#\{v\mid k_{v}(\sigma)>k_{u}(\sigma)\}=\#\{i\mid -a_{i}(\sigma)>k_{u}(\sigma)\}-1=I_{u}(\sigma)-1$. 
 Therefore, 
 \begin{equation}\label{the right hand side step 2}
 \#\{v\mid k_{v}(\sigma)>k_{u}(\sigma)\}-\#\{v\mid k_{v}(\sigma)<k_{u}(\sigma)\}=2I_{u}(\sigma)-n
 \end{equation}
 
 \bigskip

 We compare equations (\ref{CM Whittaker period}), (\ref{the left hand side step 1}), (\ref{the right hand side step 1}) and (\ref{the right hand side step 2}). If $L(\cfrac{1}{2}+m,\Pi\times \Pi^{\#})\neq 0$, we obtain that:
 
 \begin{equation}
 (2\pi i)^{d(m+\frac{1}{2})n(n-1)}\prod\limits_{1\leq u\leq n-1} P^{(I_{u})}(\Pi)\sim_{E(\Pi)E(\Pi^{\#});K}  (2\pi i)^{\frac{d(n-1)(n-2)}{2}}p(\Pi)\Omega(\Pi^{\#}_{\infty}) p(m,\Pi_{\infty},\Pi^{\#}_{\infty}).
 \end{equation}
 
 Hence we have
 \begin{equation}\nonumber
p(\Pi)\sim_{E(\Pi)E(\Pi^{\#};K)} (2\pi i)^{d(m+\frac{1}{2})n(n-1)-\frac{d(n-1)(n-2)}{2}}\Omega(\Pi^{\#}_{\infty})^{-1} p(m,\Pi_{\infty},\Pi^{\#}_{\infty})^{-1}\prod\limits_{1\leq u\leq n-1} P^{(I_{u})}(\Pi).
 \end{equation}
 
 If we take $Z(m,\Pi_{\infty},\Pi'_{\infty})=(2\pi i)^{d(m+\frac{1}{2})n(n-1)-\frac{d(n-1)(n-2)}{2}}\Omega(\Pi^{\#}_{\infty})^{-1} p(m,\Pi_{\infty},\Pi^{\#}_{\infty})^{-1}$ then $p(\Pi)\sim_{E(\Pi)E(\Pi^{\#};K)} Z(m,\Pi_{\infty},\Pi'_{\infty})\prod\limits_{1\leq u\leq n-1} P^{(I_{u})}(\Pi).$ We see that $Z(m,\Pi_{\infty},\Pi'_{\infty})$ depends only on $\Pi_{\infty}$. \\
 
 We may define:
 \begin{equation}\label{definition of Z}
 Z(\Pi_{\infty}):=Z(m,\Pi_{\infty},\Pi'_{\infty})=(2\pi i)^{d(m+\frac{1}{2})n(n-1)-\frac{d(n-1)(n-2)}{2}}\Omega(\Pi^{\#}_{\infty})^{-1} p(m,\Pi_{\infty},\Pi^{\#}_{\infty})^{-1}.
 \end{equation}
  It is well defined up to elements in $E(\Pi)^{\times}$.\\

 We deduce that:
 \begin{equation}\label{Whittaker period whole formula}
 p(\Pi)\sim_{E(\Pi)E(\Pi^{\#};K)} Z(\Pi_{\infty})\prod\limits_{1\leq u\leq n-1} P^{(I_{u})}(\Pi).
  \end{equation}
 
 \section{Formula for the Whittaker period: odd dimensional}\label{Whittaker period odd dimensional}

Let $n$ be an odd positive integer. We keep the notation in the above section We have $a_{i}(\sigma)\in \Z$ for all $1\leq i\leq n$ and all $\sigma\in \Sigma$. We assume that:

\begin{hyp}\label{hyp infinity type odd}\textbf{Odd dimensional}
For all $\sigma\in\Sigma$, the numbers $\{k_{u}(\sigma)+\cfrac{1}{2}
\mid 1\leq u\leq n-1\}$ lies in the $n-1$ gaps between $-a_{n}(\sigma)>-a_{n-1}(\sigma)>\cdots>-a_{1}(\sigma)$.
\end{hyp}
 
Recall that $\psi_{F}$ is an algebraic Hecke character of $F$ with infinity type $z^{1}$ at each $\sigma\in \Sigma$ such that $\psi_{F}\psi_{F}^{c}=||\cdot||_{\AF}$. We take $\Pi^{\#}$ to be the Langlands sum of $\chi_{u}\psi_{F}||\cdot||_{\AF}^{-\frac{1}{2}}$, $1\leq u\leq n-1$. It is an algebraic regular automorphic representation of $GL_{n-1}(\AF)$. The conditions of \ref{Whittaker period theorem CM} hold. \\
 
 We repeat the above process for $\Pi$ and $\Pi^{\#}$ and get
  \begin{eqnarray}\nonumber
&&L(\cfrac{1}{2}+m,\Pi\times \Pi^{\#}) \sim_{E(\Pi)E(\Pi^{\#});K} (2\pi i)^{dmn(n-1)}\prod\limits_{1\leq u\leq n-1} [P^{(I(\Pi,\chi_{u}\psi_{F}))}\prod\limits_{\sigma\in\Sigma}p(\widecheck{\chi_{u}},\sigma)^{2I_{u}(\sigma)-n}  ]\times\\
&&\prod\limits_{\sigma\in\Sigma}(p(\widecheck{\psi_{F}},\sigma)^{\sum_{1\leq u\leq n-1}I_{u}(\sigma)}p(\widecheck{\psi_{F}}^{c},\sigma)^{\sum_{1\leq u\leq n-1}(n-I_{u}(\sigma))})
\end{eqnarray}
 where $I_{u}:=I(\Pi,\chi_{u}\psi_{F})$ with $I_{u}(\sigma)=\#\{i\mid-a_{i}(\sigma)> k_{u}(\sigma)+\cfrac{1}{2}\}$.
 
 \bigskip
 
It is easy to verify that Hypothesis \ref{hyp infinity type} or \ref{hyp infinity type odd} is equivalent to the following hypothesis:
\begin{hyp}\label{hyp infinity type 2}
For all $\sigma\in\Sigma$, the $(n-1)$ numbers $I_{u}(\sigma)$, $1\leq u\leq n-1$, run over the numbers $1,2,\cdots,n-1$.
\end{hyp}
 
 We see $\prod\limits_{1\leq u\leq n-1}I_{u}(\sigma)=\cfrac{n(n-1)}{2}$ and $\sum_{1\leq u\leq n-1}(n-I_{u}(\sigma))=\cfrac{n(n-1)}{2}$. \\
 
 We then have 
 \begin{eqnarray}
&& \prod\limits_{\sigma\in\Sigma}(p(\widecheck{\psi_{F}},\sigma)^{\sum_{1\leq u\leq n-1}I_{u}(\sigma)}p(\widecheck{\psi_{F}}^{c},\sigma)^{\sum_{1\leq u\leq n-1}(n-I_{u}(\sigma))}
)\nonumber\\
&\sim_{E(\psi_{F});K}& \prod\limits_{\sigma\in\Sigma}p(\widecheck{\psi_{F}\psi_{F}^{c}},\sigma)^{\frac{n(n-1)}{2}} \sim_{E(\psi_{F});K}  \prod\limits_{\sigma\in\Sigma}p(||\cdot||_{\AF}^{-1},\sigma)^{\frac{n(n-1)}{2}} \sim_{E(\psi_{F});K} (2\pi i)^{\frac{dn(n-1)}{2}}\nonumber.
 \end{eqnarray}
 
 We verify that the equation (\ref{the right hand side step 1}) and (\ref{the right hand side step 2}) remain unchanged. We can see that equation (\ref{Whittaker period whole formula}) still holds here.
 
 \bigskip
 
 Let us start from the numbers $I_{u}(\sigma)$. If we are given some numbers $I_{u}(\sigma)$, $\sigma\in\Sigma$, $1\leq u\leq n-1$, such that Hypothesis \ref{hyp infinity type 2} is satisfied, we can always choose $k_{u}(\sigma)\in \Z$ such that $I_{u}(\sigma)=\#\{i\mid-a_{i}(\sigma)> k_{u}(\sigma)\}$ if $n$ is even, $I_{u}(\sigma)=\#\{i\mid-a_{i}(\sigma)> k_{u}(\sigma)+\cfrac{1}{2}\}$ if $n$ is odd.\\
 
 We may then take $\chi_{u}$, $1\leq u\leq n-1$ with infinity type $z^{k_{u}(\sigma)}\overline{z}^{-k_{u}(\sigma)}$ at $\sigma\in\Sigma$. Equation  (\ref{Whittaker period whole formula}) tells us that 
 \begin{equation}
 p(\Pi)
\sim_{E(\Pi);K} Z(\Pi_{\infty})  \prod\limits_{1\leq u\leq n-1} P^{(I_{u})}(\Pi)
 \end{equation}
 provided a non vanishing condition of the $L$-function, for example, if $\Pi$ is 3-regular.\\
  
 \begin{thm}\label{main step for factorization}
 Let $I_{u}(\sigma)$, $1\leq u\leq n-1$, $\sigma\in\Sigma$ be some integers such that Hypothesis \ref{hyp infinity type 2} is verified. There exists a complex number $Z(\Pi_{\infty})$ such that if a non vanishing condition of a global $L$-function is verified, in particular, if $\Pi$ is 3-regular, then:
 \begin{equation}
 p(\Pi)\sim_{E(\Pi);K} Z(\Pi_{\infty})  \prod\limits_{1\leq u\leq n-1} P^{(I_{u})}(\Pi).
 \end{equation}
 
 \end{thm}

\section{Factorization of arithmetic automorphic periods: restricted case}

  We consider the function $\prod\limits_{\sigma\in\Sigma} \{0,1,\cdots,n\}\rightarrow \C^{\times}/E(\Pi)^{\times}$ which sends  $(I(\sigma))_{\sigma\in\Sigma}$ to $P^{(I)}(\Pi)$.
  
 The motivic calculation predicts that:
 \begin{conj}\label{factorization conjecture}
 There exists some non zero complex numbers $P^{(s)}(\Pi,\sigma)$ for all $0\leq s\leq n$ and $\sigma\in\Sigma$ such that $P^{(I)}(\Pi) \sim_{E(\Pi);K} \prod\limits_{\sigma\in\Sigma}P^{(I(\sigma))}(\Pi,\sigma)$ for all $I=(I(\sigma))_{\sigma\in\Sigma}\in \{0,1,\cdots,n\}^{\Sigma}$.
  \end{conj}  
  
  \bigskip
  
  In this section, we will prove the above conjecture restricted to $\{1,2,\cdots,{n-1}\}^{\Sigma}$. More precisely, we will prove that
  
 \begin{thm}\label{restricted theorem}
If $n\geq 4$ and $\Pi$ satisfies a global non vanishing condition, in particular, if $\Pi$ is $3$-regular, then there exists some non zero complex numbers $P^{(s)}(\Pi,\sigma)$ for all $1\leq s\leq n-1$, $\sigma\in\Sigma$ such that $P^{(I)}(\Pi) \sim_{E(\Pi);K} \prod\limits_{\sigma\in\Sigma}P^{(I(\sigma))}(\Pi,\sigma)$ for all $I=(I(\sigma))_{\sigma\in\Sigma}\in \{1,2,\cdots,n-1\}^{\Sigma}$.

 \end{thm}
\begin{dem}
 For all $\sigma\in\Sigma$, let $I_{1}(\sigma)\neq I_{2}(\sigma)$ be two numbers in $\{1,2,\cdots,n-1\}$. We consider $I_{1},I_{2}$ as two elements in $\{1,2,\cdots,n-1\}^{\Sigma}$.
 
 Let $\sigma_{0}$ be any element in $\Sigma$. We define $I'_{1}, I'_{2}\in \{1,2,\cdots,n-1\}^{\Sigma}$ by $I_{1}'(\sigma):=I_{1}(\sigma)$, $I_{2}'(\sigma):=I_{2}(\sigma)$ if $\sigma\neq \sigma_{0}$ and $I'_{1}(\sigma_{0}):=I_{2}(\sigma_{0})$, $I'_{2}(\sigma_{0}):=I_{1}(\sigma_{0})$.
 
  By Remark \ref{remark factorization}, it is enough to prove that 
  \begin{equation}\nonumber
  P^{(I_{1})}(\Pi)P^{(I_{2})}(\Pi)\sim_{E(\Pi);K}P^{(I'_{1})}(\Pi)P^{(I'_{2})}(\Pi).
  \end{equation}

Since $I_{1}(\sigma)\neq I_{2}(\sigma)$  for all $\sigma\in\Sigma$, we can always find $I_{3},\cdots,I_{n-1}\in \{1,2,\cdots,n-1\}^{\Sigma}$ such that for all $\sigma\in\Sigma$, the $(n-1)$ numbers $I_{u}(\sigma)$, $1\leq u\leq n-1$ run over $1,2,\cdots,n-1$. In other words, Hypothesis \ref{hyp infinity type 2} is verified.\\
   
 By Theorem \ref{main step for factorization}, we have 
 \begin{equation}\nonumber p(\Pi)\sim_{E(\Pi);K} Z(\Pi_{\infty})  P^{(I_{1})}(\Pi)P^{(I_{2})}(\Pi)\prod\limits_{3\leq u\leq n-1} P^{(I_{u})}(\Pi).
 \end{equation}
 
 On the other hand, it is easy to see that $I'_{1}$, $I'_{2}$, $I_{3},\cdots,I_{n-1}$ also satisfy Hypothesis \ref{hyp infinity type 2}. Therefore  \begin{equation}\nonumber
 p(\Pi)\sim_{E(\Pi);K} Z(\Pi_{\infty})  P^{(I'_{1})}(\Pi)P^{(I'_{2})}(\Pi)\prod\limits_{3\leq u\leq n-1} P^{(I_{u})}(\Pi).
\end{equation}

We conclude at last $P^{(I_{1})}(\Pi)P^{(I_{2})}(\Pi)\sim_{E(\Pi);K}P^{(I'_{1})}(\Pi)P^{(I'_{2})}(\Pi)$ and then the above theorem follows.
 \end{dem}

  \section{Factorization of arithmetic automorphic periods: complete case}\label{complete case}
In this section, we will prove Conjecture \ref{factorization conjecture} when $\Pi$ is regular enough. More precisely, we have

\begin{thm}\label{complete theorem}
Conjecture \ref{factorization conjecture} is true if $\Pi$ is $2$-regular and satisfies a global non vanishing condition, in particular, if $\Pi$ is $6$-regular.
\end{thm}

\begin{cor}\label{Whittaker period local formula}
If $\Pi$ satisfied the conditions in the above theorem then we have:
\begin{equation}
p(\Pi)
\sim_{E(\Pi);K} Z(\Pi_{\infty}) \prod\limits_{\sigma\in\Sigma} \prod\limits_{1\leq i\leq n-1} P^{(i)}(\Pi,\sigma)
\end{equation}
\end{cor}

\bigskip

 If $n=1$, Conjecture \ref{factorization conjecture} is known as multiplicity of CM periods. We may assume that $n\geq 2$. The set $\{0,1,\cdots,n\}$ has at least $3$ elements and then Remark \ref{remark factorization} can apply.\\
 
  For all $\sigma\in\Sigma$, let $I_{1}(\sigma)\neq I_{2}(\sigma)$ be two numbers in $\{0,1,\cdots,n\}$. We have $I_{1},I_{2}\in \{0,1,2,\cdots,n\}^{\Sigma}$.\\
 
 Let $\sigma_{0}$ be any element in $\Sigma$. We define $I'_{1}, I'_{2}\in \{0,1,2,\cdots,n\}^{\Sigma}$ as in the proof of Theorem \ref{restricted theorem}. \\
 
 It remains to show that \begin{equation}\label{factorization equation}
 P^{(I_{1})}(\Pi)P^{(I_{2})}(\Pi)\sim_{E(\Pi);K}P^{(I'_{1})}(\Pi)P^{(I'_{2})}(\Pi).
 \end{equation}\\
 
 Let us assume that $n$ is odd at first.
 Since $\Pi$ is $2$-regular, we can find $\chi_{u}$ a conjugate self-dual algebraic Hecke character of $F$ such that $I(\Pi,\chi_{u})=I_{u}$ for $u=1,2$. We denote the infinity type of $\chi_{u}$ at $\sigma\in\Sigma$ by $z^{k_{u}(\sigma)}\overline{z}^{-k_{u}(\sigma)}$, $u=1,2$. We remark that $k_{1}(\sigma)\neq k_{2}(\sigma)$ for all $\sigma$ since $I_{1}(\sigma)\neq I_{2}(\sigma)$.\\
 
 Let $\Pi^{\#}$ be the Langlands sum of $\Pi$, $\chi_{1}^{c}$ and $\chi_{2}^{c}$. We write the infinity type of $\Pi^{\#}$ at $\sigma\in\Sigma$ by $(z^{b_{i}(\sigma)}\overline{z}^{-b_{i}(\sigma)})_{1\leq i\leq n+2}$ with $b_{1}(\sigma)>b_{2}(\sigma)>\cdots >b_{n+2}(\sigma)$. The set $\{b_{i}(\sigma),1\leq i\leq n+2\}=\{a_{i}(\sigma),1\leq i\leq n\}\cup \{-k_{1}(\sigma),-k_{2}(\sigma)\}$.\\
 
 Let $\Pi^{\diamondsuit}$ be a cuspidalconjugate self-dual cohomological representation of $GL_{n+3}(\AF)$ with infinity type $(z^{c_{i}(\sigma)}\overline{z}^{-c_{i}(\sigma)})_{1\leq i\leq n+3}$ such that $-c_{n+3}(\sigma)>b_{1}(\sigma)>-c_{n+2}(\sigma)>b_{2}(\sigma)>\cdots >-c_{2}(\sigma)>b_{n+2}(\sigma)>-c_{1}(\sigma)$ for all $\sigma\in\Sigma$. We may assume that $\Pi^{\diamondsuit}$ has definable arithmetic automorphic periods.\\
 
Proposition \ref{Whittaker period theorem CM} is true for $(\Pi^{\diamondsuit},\Pi^{\#})$. Namely,
\begin{equation}
L(\cfrac{1}{2}+m,\Pi^{\diamondsuit}\times \Pi^{\#})\sim_{E(\Pi^{\diamondsuit})E(\Pi^{\#});K} p(\Pi^{\diamondsuit})p(\Pi^{\#})p(m,\Pi^{\diamondsuit}_{\infty},\Pi^{\#}_{\infty}).
\end{equation}

We know \begin{equation}
L(\cfrac{1}{2}+m,\Pi^{\diamondsuit}\times \Pi^{\#})=L(\cfrac{1}{2}+m,\Pi^{\diamondsuit}\times \Pi)L(\cfrac{1}{2}+m,\Pi^{\diamondsuit}\times \chi_{1}^{c})L(\cfrac{1}{2}+m,\Pi^{\diamondsuit}\times \chi_{2}^{c})
\end{equation}

For $u=1$ or $2$, by Theorem \ref{main result CM} and the fact that $\chi_{u}$ is conjugate self-dual, we have
\begin{equation}
L(\cfrac{1}{2}+m,\Pi^{\diamondsuit}\times \chi_{u}) \sim_{E(\Pi^{\diamondsuit})E(\Pi^{\#});K} (2\pi i)^{(\frac{1}{2}+m)d(n+3)} P^{I(\Pi^{\diamondsuit},\chi_{u}^{c})}(\Pi)\prod\limits_{\sigma\in\Sigma} p(\widecheck{\chi_{u}},\sigma)^{-2I(\Pi^{\diamondsuit},\chi_{u}^{c})(\sigma)+(n+3)}.
\end{equation}

Proposition \ref{Shahidi Whittaker} implies that
\begin{equation}
p(\Pi^{\#})\sim_{E(\Pi^{\#});K} \Omega(\Pi^{\#}_{\infty})p(\Pi)L(1,\Pi\otimes \chi_{1})L(1,\Pi\otimes \chi_{2})L(1,\chi_{1}\chi_{2}^{c})
\end{equation} where $\Omega(\Pi^{\#}_{\infty})$ is a non zero complex numbers depend on $\Pi^{\#}_{\infty}$.\\
 
 By Theorem \ref{main result CM} again, for $u=1,2$, we have \begin{equation}
L(1,\Pi\times \chi_{u}) \sim_{E(\Pi^{\#});K} (2\pi i)^{dn} P^{I(\Pi,\chi_{u})}\prod\limits_{\sigma\in\Sigma} p(\widecheck{\chi_{u}},\sigma)^{2I(\Pi,\chi_{u})(\sigma)-n}.
\end{equation}
 
 Moreover, $L(1,\chi_{1}\chi_{2}^{c})\sim_{E(\Pi^{\#});K} (2\pi i)^{d}\prod\limits_{\sigma\in\Sigma}p(\widecheck{\chi_{1}},\sigma)^{t(\sigma)}p(\widecheck{\chi_{2}},\sigma)^{-t(\sigma)}$ where $t(\sigma)=1$ if $k_{1}(\sigma)<k_{2}(\sigma)$, $t(\sigma)=-1$ if $k_{1}(\sigma)>k_{2}(\sigma)$.\\

\begin{lem}
For all $\sigma\in\Sigma$, \begin{eqnarray}\nonumber
-2I(\Pi^{\diamondsuit},\chi_{1}^{c})(\sigma)+(n+3)=2I(\Pi,\chi_{1})(\sigma)-n+t(\sigma),\\\nonumber
-2I(\Pi^{\diamondsuit},\chi_{2}^{c})(\sigma)+(n+3)=2I(\Pi,\chi_{1})(\sigma)-n-t(\sigma).
\end{eqnarray}
\end{lem}
 \begin{dem}
 
 By definition we have 
 \begin{equation}\nonumber I(\Pi^{\diamondsuit},\chi_{1}^{c})(\sigma)=\#\{1\leq i\leq n+3\mid -c_{i}(\sigma)>-k_{1}(\sigma)\}.
 \end{equation}

Recall that $-c_{n+3}(\sigma)>b_{1}(\sigma)>-c_{n+2}(\sigma)>b_{2}(\sigma)>\cdots >-c_{2}(\sigma)>b_{n+2}(\sigma)>-c_{1}(\sigma)$ and $\{b_{i}(\sigma),1\leq i\leq n+2\}=\{a_{i}(\sigma),1\leq i\leq n\}\cup \{-k_{1}(\sigma),-k_{2}(\sigma)\}$.\\

Therefore \begin{eqnarray}\nonumber
I(\Pi^{\diamondsuit},\chi_{1}^{c})(\sigma)&=&\#\{1\leq i\leq n+2\mid b_{i}(\sigma)>-k_{1}(\sigma)\}+1\\ \nonumber &=&
\#\{1\leq i\leq n\mid a_{i}(\sigma)>-k_{1}(\sigma)\}+\mathds{1}_{-k_{2}(\sigma)>-k_{1}(\sigma)}+1.
\end{eqnarray}

By definition we have
\begin{equation}\nonumber
I(\Pi,\chi_{1})(\sigma)=\#\{1\leq i\leq n\mid -a_{i}(\sigma)>k_{1}(\sigma)\}=n-\#\{1\leq i\leq n\mid a_{i}(\sigma)>-k_{1}(\sigma)\}.
\end{equation}

Therefore, $I(\Pi^{\diamondsuit},\chi_{1}^{c})(\sigma)=n-I(\Pi,\chi_{1})(\sigma)+\mathds{1}_{-k_{2}(\sigma)>-k_{1}(\sigma)}+1$. Hence we have
$-2I(\Pi^{\diamondsuit},\chi_{1}^{c})(\sigma))+(n+3)=2I(\Pi,\chi_{1})(\sigma)-n+1-2\mathds{1}_{-k_{2}(\sigma)>-k_{1}(\sigma)}$.\\

It is easy to verify that $1-2\mathds{1}_{-k_{2}(\sigma)>-k_{1}(\sigma)}=t(\sigma)$. The first statement then follows and the second is similar to the first one.

 \end{dem}

We deduce that if $L(\cfrac{1}{2}+m,\Pi^{\diamondsuit}\times \Pi^{\#})\neq 0$, then \begin{eqnarray}\nonumber\label{factorization complete final step}
&&L(\cfrac{1}{2}+m,\Pi^{\diamondsuit}\times \Pi) (2\pi i)^{(1+2m)d(n+3)}P^{I(\Pi^{\diamondsuit},\chi_{1}^{c})}(\Pi^{\diamondsuit})P^{I(\Pi^{\diamondsuit},\chi_{2}^{c})}(\Pi^{\diamondsuit})\sim_{E(\Pi^{\diamondsuit})E(\Pi^{\#});K}\\
&&(2\pi i)^{d(2n+1)}p(\Pi^{\diamondsuit})\Omega(\Pi^{\#}_{\infty})p(m,\Pi_{\infty}^{\#},\Pi^{\#}_{\infty})P^{I(\Pi,\chi_{1})}(\Pi)P^{I(\Pi,\chi_{2})}(\Pi).
\end{eqnarray}

Now let $\chi_{1}'$, $\chi_{2}'$ be two conjugate self-dual algebraic Hecke characters of $F$ such that $\chi'_{1,\sigma}=\chi_{1,\sigma}$ and $\chi'_{2,\sigma}=\chi_{2,\sigma}$ for $\sigma\neq \sigma_{0}$, $\chi'_{1,\sigma_{0}}=\chi_{2,\sigma_{0}}$ and $\chi'_{2,\sigma_{0}}=\chi_{1,\sigma_{0}}$.\\

We take $\Pi^{\#\#}$ as Langlands sum of $\Pi$, $\chi'_{1}\text{}^{c}$ and $\chi'_{2}\text{}^{c}$. Since the infinity type of $\Pi^{\#\#}$ is the same with $\Pi^{\#}$, we can repeat the above process and we see that equation (\ref{factorization complete final step}) is true for $(\Pi^{\diamondsuit},\Pi^{\#\#})$. Observe that most terms remain unchanged. \\

Comparing equation (\ref{factorization complete final step}) for $(\Pi^{\diamondsuit},\Pi^{\#})$ and that for $(\Pi^{\diamondsuit},\Pi^{\#\#})$, we get

\begin{equation}\label{last compare factorization}
\cfrac{P^{I(\Pi^{\diamondsuit},\chi'_{1}\text{}^{c})}(\Pi^{\diamondsuit})P^{I(\Pi^{\diamondsuit},\chi'_{2}\text{}^{c})}(\Pi^{\diamondsuit})}{P^{I(\Pi^{\diamondsuit},\chi_{1}^{c})}(\Pi^{\diamondsuit})P^{I(\Pi^{\diamondsuit},\chi_{2}^{c})}(\Pi^{\diamondsuit})} \sim_{E(\Pi^{\diamondsuit})E(\Pi);K}
\cfrac{P^{I(\Pi,\chi'_{1})}(\Pi)P^{I(\Pi,\chi'_{2})}(\Pi)}{P^{I(\Pi,\chi_{1})}(\Pi)P^{I(\Pi,\chi_{2})}(\Pi)}.
\end{equation}

By construction, $I(\Pi,\chi_{u})=I_{u}$ and $I(\Pi,\chi'_{u})=I'_{u}$ for $u=1,2$. Hence to prove (\ref{factorization equation}), it is enough to show the left hand side of the above equation is a number in $E(\Pi^{\diamondsuit})^{\times}$.\\

There are at least two ways to see this. We observe that $I(\Pi^{\diamondsuit},\chi'_{1}\text{}^{c})(\sigma)=I(\Pi^{\diamondsuit},\chi_{1}\text{}^{c})(\sigma)$, $I(\Pi^{\diamondsuit},\chi'_{2}\text{}^{c})(\sigma)=I(\Pi^{\diamondsuit},\chi_{2}\text{}^{c})(\sigma)$ for $\sigma\neq \sigma_{0}$ and $I(\Pi^{\diamondsuit},\chi'_{1}\text{}^{c})(\sigma_{0})=I(\Pi^{\diamondsuit},\chi_{2}\text{}^{c})(\sigma_{0})$, $I(\Pi^{\diamondsuit},\chi'_{2}\text{}^{c})(\sigma_{0})=I(\Pi^{\diamondsuit},\chi_{1}\text{}^{c})(\sigma_{0})$. Moreover, these numbers are all in $\{1,2,\cdots,(n+3)-1\}$. Theorem \ref{restricted theorem} gives a factorization of the holomorphic arithmetic automorphic periods through each place. In particular, it implies that the left hand side of  (\ref{last compare factorization}) is in $E(\Pi^{\diamondsuit})^{\times}$ as expected.\\

One can also show this by taking $\Pi^{\diamondsuit}$ an automorphic induction of a Hecke character. We can then calculate $L(\cfrac{1}{2}+m,\Pi^{\diamondsuit}\times\chi_{u}^{c})$ in terms of CM periods. Since the factorization of CM periods is clear, we will also get the expected result.

\bigskip

When $n$ is even, we consider $\Pi^{\#}$ the Langlands sum of $\Pi$, $(\chi_{1}\psi_{F}||\cdot||^{-1/2})^{c}$ and $(\chi_{2}\psi_{F}||\cdot||^{-1/2})^{c}$ where $\chi_{1}$, $\chi_{2}$ are two suitable algebraic Hecke characters of $F$. We follow the above steps and will get the factorization in this case. We leave the details to the reader and just remark that as in section \ref{Whittaker period odd dimensional}, some CM periods of $\psi_{F}$ appear but they will be eliminated at the end.

\section{Specify the factorization}\label{specify the factorization}
\text{}
 Let us assume that Conjecture \ref{factorization conjecture} is true. We want to specify one factorization. \\
 
 We denote by $I_{0}$ the map which sends each $\sigma\in\Sigma$ to $0$. By the last part of Corollary \ref{factorization lemma}, it is enough to choose $c(\Pi,\sigma)\in (\C/E(\Pi))^{\times}$ which is $G_{K}$-equivariant such that $P^{(I_{0})}(\Pi)\sim_{E(\Pi);K} \prod\limits_{\sigma\in\Sigma}c(\Pi,\sigma)$. Then there exists a unique factorization of $P^{(\cdot)}(\Pi)$ such that $P^{(0)}(\Pi,\sigma)= c(\Pi,\sigma)$ . We may then define the \textbf{local arithmetic automorphic periods} $P^{(s)}(\Pi,\sigma)$ as an element in $\C^{\times}/ (E(\pi))^{\times}$.\\
 
In this section, we shall prove $P^{(I_{0})}(\Pi)\sim_{E(\Pi);K} p(\widecheck{\xi_{\Pi}},\overline{\Sigma}) \sim_{E(\Pi);K} \prod\limits_{\sigma\in\Sigma} p(\widecheck{\xi_{\Pi}},\overline{\sigma})$. Therefore, we may take $c(\Pi,\sigma)=p(\widecheck{\xi_{\Pi}},\overline{\sigma})$.\\

More generally, we will see that:
\begin{lem}\label{period in compact case}
If $I$ is compact then $P^{(I)}(\Pi)\sim_{E(\Pi);K}  \prod\limits_{I(\sigma)=0} p(\widecheck{\xi_{\Pi}},\overline{\sigma})\times \prod\limits_{I(\sigma)=n} p(\widecheck{\xi_{\Pi}},\sigma)$. 
\end{lem}

This lemma leads to the following theorem:
\begin{thm}\label{special factorization theorem}
If Conjecture \ref{factorization conjecture} is true, in particular, if conditions in Theorem \ref{complete theorem} are satisfied, then there exists some complex numbers $P^{(s)}(\Pi,\sigma)$ unique up to multiplication by elements in $(E(\Pi))^{\times}$ such that the following two conditions are satisfied:
\begin{enumerate}
\item $P^{(I)}(\Pi) \sim_{E(\Pi);K} \prod\limits_{\sigma\in\Sigma}P^{(I(\sigma))}(\Pi,\sigma)$ for all $I=(I(\sigma))_{\sigma\in\Sigma}\in \{0,1,\cdots,n\}^{\Sigma}$,
\item  and $P^{(0)}(\Pi,\sigma)\sim_{E(\Pi);K} p(\widecheck{\xi_{\Pi}},\overline{\sigma})$
\end{enumerate}
where $\xi_{\Pi}$ is the central character of $\Pi$.

Moreover, we know $P^{(n)}(\Pi,\sigma)\sim_{E(\Pi);K} p(\widecheck{\xi_{\Pi}},\sigma)$ or equivalently $P^{(0)}(\Pi,\sigma)\times P^{(n)}(\Pi,\sigma)\sim_{E(\Pi);K} 1$.
\end{thm}

\paragraph{Proof of Lemma \ref{period in compact case}:}
Recall that $D/2=\sum\limits_{\sigma\in\Sigma}I_{\sigma}(n-I_{\sigma})=0$ since $I$ is compact.\\

Let $T$ be the center of $GU_{I}$. We have $T(\R)\cong \{(z_{\sigma})\in (\C^{\times})^{\Sigma}\mid |z_{\sigma}| \text{ does not depend on }\sigma\}$.
We define a homomorphism $h_{T}:\Ss(\R) \rightarrow T(\R)$ by sending $z\in \C$ to $((z)_{I(\sigma)=0},(\overline{z})_{I(\sigma)=n})$. \\

Since $I$ is compact, we see that $h_{I}$ is the composition of $h_{T}$ and the embedding $T\hookrightarrow GU_{I}$. We get an inclusion of Shimura varieties:
$Sh_{T}:=Sh(T,h_{T})\hookrightarrow Sh_{I}=Sh(GU_{I},h_{I})$.\\

Let $\xi$ be a Hecke character of $K$ such that $\Pi^{\vee}\otimes \xi$ descends to $\pi$, a representation of $GU_{I}(\AQ)$, as before. We write $\lambda\in \Lambda(GU_{I})$ the cohomology type of $\pi$. We define $\lambda^{T}:=(\lambda_{0},(\sum\limits_{1\leq i\leq n}\lambda_{i}(\sigma))_{\sigma\in\Sigma})$. Since $\pi$ is irreducible, it acts as scalars when restrict to $T$. This gives $\pi^{T}$, a one dimensional representation of $T(\AQ)$ which is cohomology of type ${\lambda}^{T}$. We denote by $V_{\lambda^{T}}$ the character of $T(\R)$ with highest weight $\lambda^{T}$.\\ 

The automorphic vector bundle $E_{\lambda}$ pulls back to the automorphic vector bundle $[V_{\lambda^{T}}]$ (see \cite{harrisappendix} for notation) on $Sh_{T}$. \\

Let $\beta$ be an element in $\bar{H}^{0}(Sh_{I},E_{\lambda})^{\pi}$. We fix a non zero $E(\pi)$-rational element in $\pi$ and then we can lift $\beta$ to $\phi$, an automorphic form on $GU_{I}(\AQ)$.\\

There is an isomorphism $H^{0}(Sh_{T},[V_{\lambda^{T}}])\xrightarrow{\sim} \{f\in \mathbb{C}^{\infty}(T(\Q)\backslash T(\AQ),\C\mid f(tt_{\infty}))=\pi^{T}(t_{\infty})f(t), t_{\infty}\in T(\R), t\in T(\AQ)\}$ (c.f. \cite{harrisappendix}). We send $\beta$ to the element in $H^{0}(Sh_{T},[V_{\lambda^{T}}])^{\pi^{T}}$ associated to $\phi|_{T(\AQ)}$.\\

We then obtain rational morphisms
\begin{eqnarray}
&\bar{H}^{0}(Sh_{I},E_{\lambda})^{\pi} \xrightarrow{\sim} H^{0}(Sh_{T},[V_{\lambda^{T}}])^{\pi^{T}}&\\
\text{ and similarly } &\bar{H}^{0}(Sh_{I},E_{\lambda^{\vee}})^{\pi^{\vee}} \xrightarrow{\sim} H^{0}(Sh_{T},[V_{\lambda^{T,\vee}}])^{\pi^{T,\vee}}.&
\end{eqnarray}

\bigskip

These morphisms are moreover isomorphisms. In fact, since both sides are one dimensional, it is enough to show the above morphisms are injective. Indeed, if $\phi$, a lifting of an element in $\bar{H}^{0}(Sh_{I},E_{\lambda})^{\pi}$, vanishes at the center, in particular, it vanishes at the identity. Hence it vanishes at $GU_{I}(\AQf)$ since it is an automorphic form. We observe that $GU_{I}(\AQf)$ is dense in $GU_{I}(\Q)\backslash GU_{I}(\AQ)$. We know $\phi=0$ as expected.\\

We are going to calculate the arithmetic automorphic period. Let $\beta$ be rational. We take a rational element $\beta^{\vee}\in \bar{H}^{0}(Sh_{I},E_{\lambda^{\vee}})^{\pi^{\vee}} $ and lift it to an automorphic form $\phi^{\vee}$. We have $ c_{B}(\phi) \sim_{E(\pi);K}P^{(I)}(\pi) \phi^{\vee} $ by Lemma \ref{pair to ratio}.\\

For the torus, by Remark \ref{CM complex conjugation}, we know \begin{equation}\nonumber
\phi^{\vee}|_{T(\AQ)}\sim_{E(\pi);K} p(Sh(T,h_{T}),\pi^{T})^{-1}(\phi|_{T(\AQ)})^{-1}.
\end{equation}

Recall that $c_{B}(\phi)=\pm i^{\lambda_{0}}\overline{\phi}||\nu(\cdot)||^{\lambda_{0}}$. Therefore $(c_{B}(\phi))|_{T(\AQ)}=\pm i^{\lambda_{0}}(\phi|_{T(\AQ)})^{-1}$. We then get \begin{equation}\label{pass to torus}
 i^{\lambda_{0}}P^{(I)}(\pi)\sim_{E(\pi);K}p(Sh(T,h_{T}),\pi^{T}).
\end{equation}

\bigskip

We now set $T^{\#}:=Res_{K/\Q}T_{K}$. We have $T^{\#} \cong Res_{K/\Q}\mathbb{G}_{m}\times Res_{F/\Q}\mathbb{G}_{m}$. In particular, $T^{\#}(\R)\cong \C^{\times} \times (\R\otimes_{\Q}F)^{\times}\cong \C^{\times}\times (\C^{\times})^{\Sigma}$.\\

We define $h_{T^{\#}}:\Ss(\R) \rightarrow T^{\#}(\R)$ to be the composition of $h_{T}$ and the natural embedding $T(\R)\rightarrow T^{\#}(\R)$. We know $h_{T^{\#}}$ sends $z\in\C^{\times}$ to $(z\overline{z},(z)_{I(\sigma)=0}, (\overline{z})_{r(\sigma)=0})$. The embedding $(T,h_{T})\rightarrow (T^{\#},h_{T^{\#}})$ is a map between Shimura datum.\\

We observe that $\pi^{T,\#}:=||\cdot||^{-\lambda_{0}}\times \xi_{\Pi}^{-1}$ is a Hecke character on $T^{\#}$. Its restriction to $T$ is just $\pi^{T}$. By Proposition \ref{propgeneral}, we have $p(Sh(T,h_{T}),\pi^{T})\sim_{E(\pi);K}p(Sh(T^{\#},h_{T^{\#}}),\pi^{T^{\#}})$. \\

By the definition of CM period and Proposition \ref{propCM}, we have
 \begin{equation}
p(Sh(T^{\#},h_{T^{\#}}),\pi^{T^{\#}}) \sim_{E(\pi);K} (2\pi i)^{\lambda_{0}} \prod\limits_{I(\sigma)=0}p(\xi_{\Pi}^{-1},\sigma) \prod\limits_{I(\sigma)=n}p(\xi_{\Pi}^{-1},\overline{\sigma}).
\end{equation}

\bigskip

Since $\xi_{\Pi}$ is conjugate self-dual, we have $p(\xi_{\Pi}^{-1},\overline{\sigma})\sim_{E(\Pi);K} p(\xi_{\Pi},\sigma)$.

By equation (\ref{pass to torus}), we get:
\begin{equation}
 i^{\lambda_{0}}P^{(I)}(\pi)\sim_{E(\pi);K}(2\pi i)^{\lambda_{0}} \prod\limits_{I(\sigma)=0}p(\xi_{\Pi}^{-1},\sigma) \prod\limits_{I(\sigma)=n}p(\xi_{\Pi},\sigma).
\end{equation}
Recall that by definition $P^{(I)}(\Pi)\sim_{E(\Pi);K} (2\pi)^{-\lambda_{0}} P^{(I)}(\pi)$,  we get finally
\begin{eqnarray}\nonumber
P^{(I)}(\Pi)&\sim_{E(\Pi);K}&  \prod\limits_{I(\sigma)=0} p(\xi_{\Pi}^{-1},\sigma)\times \prod\limits_{I(\sigma)=n} p(\xi_{\Pi},\sigma)\\\nonumber
&\sim_{E(\Pi);K}&  \prod\limits_{I(\sigma)=0} p(\widecheck{\xi_{\Pi}},\overline{\sigma})\times \prod\limits_{I(\sigma)=n} p(\widecheck{\xi_{\Pi}},\sigma). 
\end{eqnarray} 
The last formula comes from the fact that $\xi_{\Pi}$ is conjugate self-dual.

\begin{flushright}$\Box$\end{flushright}

We recall that the arithmetic automorphic periods can be defined for essential conjugate self-dual representations. More precisely, let $\Pi$ be conjugate self-dual as in Theorem \ref{special factorization theorem}, let $\eta$ be an algebraic Hecke character. By Definition \ref{definable periods}, we have defined $P^{(I)}(\Pi\otimes \eta)$ as $P^{(I)}(\Pi)\prod\limits_{\sigma\in\Sigma}p(\widecheck{\eta},\sigma)^{I(\sigma)}p(\widecheck{\eta},\overline{\sigma})^{n-I(\sigma)}$. As we showed above that $P^{(I)}(\Pi) \sim_{E(\Pi);K} \prod\limits_{\sigma\in\Sigma}P^{(I(\sigma))}(\Pi,\sigma)$, it is natural to define:

\begin{df}\label{factorization essential case}
We define \textbf{local arithmetic automorphic periods} for  conjugate self-dual representations by \begin{equation}
P^{(s)}(\Pi\otimes \eta,\sigma)=P^{(s)}(\Pi,\sigma)p(\widecheck{\eta},\sigma)^{I(\sigma)}p(\widecheck{\eta},\overline{\sigma})^{n-I(\sigma)}.
\end{equation} 

\end{df}

\begin{rem}
If $s=0$, we see that 
\begin{eqnarray}
P^{(0)}(\Pi\otimes \eta,\sigma)&=&P^{(0)}(\Pi,\sigma)p(\widecheck{\eta},\overline{\sigma})^{n}\sim_{E(\Pi;K)}p(\widecheck{\xi_{\Pi}},\overline{\sigma})p(\widecheck{\eta},\overline{\sigma})^{n}\\\nonumber
&\sim_{E(\Pi;K)}&p(\widecheck{\xi_{\Pi}\eta^{n}},\overline{\sigma})\sim_{E(\Pi;K)}p(\widecheck{\xi_{\Pi\otimes \eta}},\overline{\sigma})
\end{eqnarray}

Therefore, if $\Pi$ has definable arithmetic automorphic periods and regular enough, we still have \begin{enumerate}
\item $P^{(I)}(\Pi) \sim_{E(\Pi);K} \prod\limits_{\sigma\in\Sigma}P^{(I(\sigma))}(\Pi,\sigma)$ for all $I=(I(\sigma))_{\sigma\in\Sigma}\in \{0,1,\cdots,n\}^{\Sigma}$,
\item  and $P^{(0)}(\Pi,\sigma)\sim_{E(\Pi);K} p(\widecheck{\xi_{\Pi}},\overline{\sigma})$.
\end{enumerate}
Moreover, these two properties determine the local periods.

\end{rem}

\bigskip

\begin{rem}\label{remark n=1}
If $n=1$ and $\Pi=\eta$ is a Hecke character, we obtain that:
$
P^{(0)}(\eta,\sigma)\sim_{E(\eta);K} p(\widecheck{\eta},\overline{\sigma})
$
and similarly
$
P^{(1)}(\eta,\sigma)\sim_{E(\eta);K}p(\widecheck{\eta},\sigma).$
\end{rem}

\chapter{Functoriality of arithmetic automorphic periods}
\section{Period relations for automorphic inductions: settings}\label{AI setting}
\text{}

Let $F$ be a CM field containing $K$ as before.

Let $\mathcal{F}/F$ be a cyclic extension of CM fields of degree $l$. 

Let $\Pi_{\mathcal{F}}$ be a cuspidal representation of $GL_{n}(\AFFF)$.\\

By Theorem $6.2$ of \cite{arthurclozel}, there exists $\Pi_{F}$, an automorphic representation of $GL_{nl}(\AF)$ which lifts $\Pi_{\mathcal{F}}$. We assume moreover that $\Pi_{\mathcal{F}}\ncong \Pi_{\mathcal{F}}^{g}$ for all $g\in Gal(\mathcal{F}/F)$ non trivial. We can read from the proof of Theorem $6.2$ in the \textit{loc.cit} that $\Pi_{F}$ is then cuspidal.\\

 We want to compare the arithmetic automorphic periods of $\Pi_{\mathcal{F}}$ and $\Pi_{F}$ if they are defined. For this purpose, we assume that $\Pi_{\mathcal{F}}$ has definable arithmetic automorphic periods as in Definition \ref{definable periods}. In other words, $\Pi_{\mathcal{F}}$ is $3$-regular, cohomological and descends to unitary groups of any sign after tensoring by an algebraic Hecke character. \\
  
 We write the infinity type of $\Pi_{\mathcal{F}}$ as $(z^{a_{i}(\sigma)}\overline{z}^{b_{i}(\sigma)})_{1\leq i\leq n}$ at $\sigma\in\Sigma_{\mathcal{F};K}$. We remark that $a_{i}(\sigma), b_{i}(\sigma)\in \Z+\cfrac{n-1}{2}$. \\
 
 The restriction of embeddings gives a map: \begin{equation}\nonumber
 \Psi_{\mathcal{F}/F}: \Sigma_{\mathcal{F};K} \rightarrow \Sigma_{F;K}.
 \end{equation}
 For $\tau\in \Sigma_{F;K}$, the infinity type of $\Pi_{F}$ at $\tau$ is $(z^{a_{i}(\sigma)}\overline{z}^{b_{i}(\sigma)})_{1\leq i\leq n,\sigma\in \Psi_{\mathcal{F}/F}^{-1}(\tau)}$. We assume in this chapter that for any $\tau$ the $nl$ numbers $a_{i}(\sigma)$, $1\leq i\leq n$ and $\sigma\in\Psi_{\mathcal{F}/F}^{-1}(\tau)$, are different. We assume moreover their differences are at least $3$. Hence $\Pi_{F}$ is also $3$-regular.\\
 
 If $l$ is odd or $n$ is even, we know $\Pi_{F}$ is algebraic and then cohomological. We write $\Pi'_{F}:=\Pi_{F}$ in this case.\\
   
 If $l$ is even and $n$ is odd, $\Pi_{F}$ is no longer algebraic. We define $\Pi'_{F}:=\Pi_{F}||\cdot||_{\AF}^{-1/2}$.  It is then a cuspidal cohomological representation of $GL_{nl}(\AF)$. It is conjugate self-dual after tensoring by an algebraic Hecke character. To see this, we take $\psi_{F}$ an algebraic Hecke character of $F$ with infinity type $z^{1}\overline{z}^{0}$ at each infinity place such that $\psi_{F}\psi_{F}^{c}=||\cdot||_{\AF}$. We remark that the Hecke character $||\cdot||_{\AF}^{-1/2}\otimes \psi_{F}$ is conjugate self-dual.\\
 
 We also assume that $\Pi'_{F}$ descends to unitary groups of any sign after tensoring an algebraic Hecke character. Therefore, $\Pi'_{F}$ has definable arithmetic automorphic periods.\\
 
Let $I_{F}$ be a map from $\Sigma_{F;K}$ to the set $\{0,1,\cdots,nl\}$. We want to relate $P^{(I_{F})}(\Pi'_{F})$ to arithmetic automorphic periods of $\Pi_{\mathcal{F}}$ in the following sections. \\

We take $\eta$ an algebraic Hecke character of $F$ such that $I(\Pi_{F},\eta)=I_{F}$. We take $m$ as in the last part Theorem \ref{n*1 essential}. We assume that Conjecture \ref{special value for similitude unitary group} is true and we have:
\begin{equation}\label{AI step 1}
L(m,\Pi'_{F}\otimes \eta) \sim_{E(\Pi_{F})E(\eta);K} (2\pi i)^{mnld} P^{(I_{F})}(\Pi'_{F})\prod\limits_{\tau\in\Sigma_{F;K}} p(\widecheck{\eta},\tau)^{I_{F}(\tau)}p(\widecheck{\eta},\overline{\tau})^{nl-I_{F}(\tau)}
\end{equation}
with both sides non zero.

\section{Relations of global periods for automorphic inductions}

\paragraph{The case $l$ is odd or $n$ is even:}
 In this case, $\Pi'_{F}=\Pi_{F}$ is the automorphic induction of $\Pi_{\mathcal{F}}$.\\

 We know $L(m,\Pi_{F}\otimes \eta) = L(m,\Pi_{\mathcal{F}}\otimes \eta\circ N_{\AFFF^{\times}/\AF^{\times}}) $. It is easy to see that $m$ is also critical for $\Pi_{\mathcal{F}}\otimes \eta\circ N_{\AFFF^{\times}/\AF^{\times}}$. We can apply Theorem \ref{n*1 essential} to $(\Pi_{\mathcal{F}}, \eta\circ N_{\AFFF^{\times}/\AF^{\times}})$. \\
 
 We write $I_{\mathcal{F}}:=I(\Pi_{\mathcal{F}}, \eta\circ N_{\AFFF^{\times}/\AF^{\times}})$ and get:
\begin{eqnarray}\nonumber
&L(m,\Pi_{F}\otimes \eta)=L(m,\Pi_{\mathcal{F}}\otimes \eta\circ N_{\AFFF^{\times}/\AF^{\times}})  \sim_{E(\Pi_{\mathcal{F}})E(\eta);K}&\\&\label{AI step 2}
 (2\pi i)^{mndl} P^{(I_{\mathcal{F}})}(\Pi_{\mathcal{F}})
 \prod\limits_{\sigma\in\Sigma_{\mathcal{F};K}} p(\widecheck{\eta\circ N_{\AFFF^{\times}/\AF^{\times}}},\sigma)^{I_{\mathcal{F}}(\sigma)}p(\widecheck{\eta\circ N_{\AFFF^{\times}/\AF^{\times}}},\overline{\sigma})^{n-I_{\mathcal{F}}(\sigma)}.&
\end{eqnarray}

\bigskip

We first calculate $I_{\mathcal{F}}=I(\Pi_{\mathcal{F}}, \eta\circ N_{\AFFF^{\times}/\AF^{\times}})$. We write the infinity type of $\eta$ at $\tau\in \Sigma_{F;K}$ by $z^{a(\tau)}\overline{z}^{b(\tau)}$.  \\

For $\sigma\in \Sigma_{\mathcal{F};K}$, the infinity type of $\eta\circ N_{\AFFF/\AF}$ at $\sigma$ is then $z^{a(\Psi_{\mathcal{F}/F}(\sigma))}\overline{z}^{b(\Psi_{\mathcal{F}/F}(\sigma))}$. \\

We have by definition that
\begin{equation}
I_{\mathcal{F}}(\sigma)=\#\{i \mid 1\leq i\leq n, a(\Psi_{\mathcal{F}/F}(\sigma))-b(\Psi_{\mathcal{F}/F}(\sigma))+a_{i}(\sigma)-b_{i}(\sigma)<0\}
\end{equation}

Recall that the infinity type of $\Pi_{F}$ at $\tau$ is $(z^{a_{i}(\sigma)}\overline{z}^{b_{i}(\sigma)})_{1\leq i\leq n,\sigma\in \Psi_{\mathcal{F}/F}^{-1}(\tau),1\leq i\leq n}$. We have:
\begin{equation}
I_{F}(\tau)=I(\Pi_{F},\eta)(\tau)=\#\{(i,\sigma)\mid 1\leq i\leq n,  \sigma\in \Psi_{\mathcal{F}/F}^{-1}(\tau), a(\tau)-b(\tau)+a_{i}(\sigma)-b_{i}(\sigma)<0\}.
\end{equation}

We observe that $I_{\mathcal{F}}$ is uniquely determined by $I_{F}$. More precisely, it is easy to show that:

\begin{lem} \label{AI infinity sign}
The integer $I_{\mathcal{F}}(\sigma)$ is the number of elements in $\{a_{i}(\sigma)\mid 1\leq i\leq n\}$ which is one of the $I_{F}(\tau)$-th smallest numbers in the set $\{a_{i}(\sigma')\mid 1\leq i\leq n, \sigma'\in \Psi_{\mathcal{F}/F}^{-1}(\tau)\}$ where $\tau=\Psi_{\mathcal{F}/F}(\sigma)$. 
\end{lem}

Moreover, it is easy to see that \begin{equation}\label{AI sum index}
\sum\limits_{\sigma\in \Psi_{\mathcal{F}/F}^{-1}(\tau)}I_{\mathcal{F}}(\sigma)=I_{F}(\tau).
\end{equation} \\

By Proposition \ref{propCM}, we get 
\begin{eqnarray}\nonumber
\prod\limits_{\sigma\in\Sigma_{\mathcal{F};K}}p(\widecheck{\eta\circ N_{\AFFF^{\times}/\AF^{\times}}},\sigma)^{I_{\mathcal{F}}(\sigma)}\sim_{E(\eta);K} \prod\limits_{\sigma\in\Sigma_{\mathcal{F};K}}p(\widecheck{\eta},\Psi_{\mathcal{F}/F}(\sigma))^{I_{\mathcal{F}}(\sigma)}\\
\sim_{E(\eta);K} \prod\limits_{\tau\in\Sigma_{F;K}}p(\widecheck{\eta},\tau)^{\sum\limits_{\sigma\in \Psi_{\mathcal{F}/F}^{-1}(\tau)}I_{\mathcal{F}}(\sigma)}\sim_{E(\eta);K} \prod\limits_{\tau\in\Sigma_{F;K}}p(\widecheck{\eta},\tau)^{I_{F}(\tau)}
\end{eqnarray}

Similarly, we have 
\begin{equation}
\prod\limits_{\sigma\in\Sigma_{\mathcal{F};K}}p(\widecheck{\eta\circ N_{\AFFF^{\times}/\AF^{\times}}},\overline{\sigma})^{n-I_{\mathcal{F}}(\sigma)} \sim_{E(\eta);K} \prod\limits_{\tau\in\Sigma_{F;K}}p(\widecheck{\eta},\tau)^{nl-I_{F}(\tau)}.
\end{equation}

Comparing the above two equations with equations (\ref{AI step 1}) and (\ref{AI step 2}), we deduce that:
\begin{equation}
P^{(I_{F})}(\Pi_{F})\sim_{E(\Pi_{\mathcal{F}});K} P^{(I_{\mathcal{F}})}(\Pi_{\mathcal{F}}).
\end{equation}

\bigskip

\paragraph{The case $l$ is even and $n$ is odd:}

In this case $\Pi_{F}$ is no longer algebraic and we consider $\Pi_{F}'=\Pi_{F}\otimes ||\cdot||^{-1/2}$.

We know $L(m,\Pi'_{F}\otimes \eta)=L(m-\cfrac{1}{2},\Pi_{F}\otimes \eta)=L(m-\cfrac{1}{2},\Pi_{\mathcal{F}}\otimes \eta\circ N_{\AFFF^{\times}/\AF^{\times}})$.\\

As in the previous case, we get:
\begin{eqnarray}\nonumber
&L(m-\cfrac{1}{2},\Pi_{F}\otimes \eta)=L(m-\cfrac{1}{2},\Pi_{\mathcal{F}}\otimes \eta\circ N_{\AFFF^{\times}/\AF^{\times}})  &\\&\nonumber \sim_{E(\Pi_{\mathcal{F}})E(\eta);K}
 (2\pi i)^{(m-\frac{1}{2})ndl} P^{(I_{\mathcal{F}})}(\Pi_{\mathcal{F}})
 \prod\limits_{\sigma\in\Sigma_{\mathcal{F};K}} p(\widecheck{\eta\circ N_{\AFFF^{\times}/\AF^{\times}}},\sigma)^{I_{\mathcal{F}}(\sigma)}p(\widecheck{\eta\circ N_{\AFFF^{\times}/\AF^{\times}}},\overline{\sigma})^{n-I_{\mathcal{F}}(\sigma)}&\\\nonumber
 &\sim_{E(\Pi_{\mathcal{F}})E(\eta);K}(2\pi i)^{(m-\frac{1}{2})ndl} P^{(I_{\mathcal{F}})}(\Pi_{\mathcal{F}})
\prod\limits_{\tau\in\Sigma_{F;K}} p(\widecheck{\eta},\tau)^{I_{F}(\tau)}p(\widecheck{\eta},\overline{\tau})^{nl-I_{F}(\tau)}&
\end{eqnarray}

We conclude that:
\begin{equation}
P^{(I_{F})}(\Pi_{F}\otimes ||\cdot||^{-1/2})\sim_{E(\Pi_{\mathcal{F}});K} (2\pi i)^{-\frac{nld}{2}}P^{(I_{\mathcal{F}})}(\Pi_{\mathcal{F}}).
\end{equation}

\section{Relations of local periods for automorphic inductions}
Recall that the arithmetic automorphic periods admit a factorization (c.f. Theorem \ref{special factorization theorem}) $P^{(I)}(\Pi) \sim_{E(\Pi);K} \prod\limits_{\sigma\in\Sigma}P^{(I(\sigma))}(\Pi,\sigma)$ such that 
\begin{equation}
P^{(0)}(\Pi,\sigma)\sim_{E(\Pi);K} p(\xi_{\Pi}^{-1},\sigma)\sim_{E(\Pi);K} p(\xi_{\Pi},\sigma)^{-1}.
\end{equation} We will discuss the functoriality of local periods in this section.\\

Let $\tau$ be an element of $\Sigma_{F;K}$.

It is easy to see from Lemma \ref{AI infinity sign} or equation (\ref{AI sum index}) that if $I_{F}(\tau)=0$ then $I_{\mathcal{F}}(\sigma)=0$ for all $\sigma\in \Sigma_{\mathcal{F};K}$ over $\tau$.\\

Fix any $\tau_{0}\in \Sigma_{F;K}$ and an integer $0\leq s_{0}\leq n$. We define $I_{F}$ such that $I_{F}(\tau_{0})=s_{0}$ and $I_{F}(\tau)=0$ for $\tau\neq \tau_{0}\in \Sigma_{F;K}$.\\

\paragraph{The case $l$ is odd or $n$ is even:}
 
Recall $P^{(I_{F})}(AI(\Pi_{\mathcal{F}}))\sim_{E(\Pi_{\mathcal{F}});K} P^{(I_{\mathcal{F}})}(\Pi_{\mathcal{F}})$ in this case. 
We get \begin{equation}\label{AI local step 1}
P^{(s_{0})}(\Pi_{F},\tau_{0})\prod\limits_{\tau\neq \tau_{0}\in\Sigma_{F;K}} P^{(0)}(\Pi_{F},\tau) \sim_{E(\Pi_{\mathcal{F}});K} \prod\limits_{\sigma_{0} \mid \tau_{0}} P^{(I_{\mathcal{F}}(\sigma_{0}))}(\Pi_{\mathcal{F}},\sigma_{0})  \prod\limits_{\tau\neq \tau_{0}\in\Sigma_{F;K}}\prod\limits_{\sigma\mid \tau} P^{(0)}(\Pi_{\mathcal{F}},\sigma)
\end{equation}\\

For $\tau\neq \tau_{0}$, we have $P^{(0)}(\Pi_{F},\tau) \sim_{E(\Pi_{\mathcal{F}});K} p(\xi_{\Pi_{F}},\tau)^{-1}$. Similarly, we have $P^{(0)}(\Pi_{\mathcal{F}},\sigma)\sim_{E(\Pi_{\mathcal{F}});K} p(\xi_{\Pi_{\mathcal{F}}},\sigma)^{-1}$.\\

Let $g$ be a generator of $Gal(\mathcal{F}/F)$. From the construction in \cite{clozelaa} we know $\Pi_{F}$ has base change $\Pi_{\mathcal{F}}\times \Pi_{\mathcal{F}}^{g}\times \cdots \times \Pi_{\mathcal{F}}^{g^{l-1}}$. In particular, we know $\xi_{\Pi_{F}}\circ N_{\AFFF^{\times}/\AF^{\times}}=\prod\limits_{0\leq i\leq l-1}\xi_{\Pi_{\mathcal{F}}}^{g^{i}}$.\\

We fix any $\sigma_{1}\in\Sigma_{\mathcal{F};K}$ over $\tau$. We know \begin{eqnarray}
p(\xi_{\Pi_{F}},\tau)&\sim_{E(\Pi_{\mathcal{F}});K}&\nonumber p(\xi_{\Pi_{F}}\circ N_{\AFFF^{\times}/\AF^{\times}},\sigma_{1})\\
&\sim_{E(\Pi_{\mathcal{F}});K}&\nonumber \prod\limits_{0\leq i\leq l-1} p(\xi_{\Pi_{\mathcal{F}}}^{g^{i}},\sigma_{1})\\
&\sim_{E(\Pi_{\mathcal{F}});K}&\nonumber \prod\limits_{0\leq i\leq l-1} p(\xi_{\Pi_{\mathcal{F}}},\sigma_{1}^{g^{-i}})\\
&\sim_{E(\Pi_{\mathcal{F}});K}&\nonumber \prod\limits_{\sigma\mid \tau} p(\xi_{\Pi_{\mathcal{F}}},\sigma)
\end{eqnarray}

Equation (\ref{AI local step 1}) then gives:
\begin{equation}\label{AI local step 2}
P^{(s_{0})}(\Pi_{F},\tau_{0}) \sim_{E(\Pi_{\mathcal{F}});K} \prod\limits_{\sigma_{0} \mid \tau_{0}} P^{(I_{\mathcal{F}}(\sigma_{0}))}(\Pi_{\mathcal{F}},\sigma_{0}) 
\end{equation}

We can read from Lemma \ref{AI infinity sign} that $I_{\mathcal{F}}(\sigma_{0})$ depends only on $I_{F}(\tau_{0})=s_{0}$.\\

 It is natural to define:
\begin{df}\label{AI infinity sign local}
Let $0\leq s\leq nl$ be any integer. Let $\tau\in \Sigma_{F;K}$. For any $\sigma\in\Sigma_{\mathcal{F};K}$ over $\tau$, we define $s(\sigma)$ to be the number of elements in $\{a_{i}(\sigma)\mid 1\leq i\leq n\}$ which is one of the $s$-th smallest numbers in the set $\{a_{i}(\sigma')\mid 1\leq i\leq n, \sigma'\in \Psi_{\mathcal{F}/F}^{-1}(\tau)\}$.
\end{df}

Equation (\ref{AI local step 2}) can be rewritten as
\begin{equation}
P^{(s)}(\Pi_{F},\tau) \sim_{E(\Pi_{\mathcal{F}});K} \prod\limits_{\sigma \mid \tau} P^{(s(\sigma))}(\Pi_{\mathcal{F}},\sigma).
\end{equation}

\bigskip

\paragraph{The case $l$ is even and $n$ is odd:}

In this case, we have:
\begin{equation}
\xi_{\Pi_{F}'} \circ N_{\AFFF^{\times}/\AF^{\times}}= \prod\limits_{0\leq i\leq l-1}\xi_{\Pi_{\mathcal{F}}}^{g^{i}} \times ||\cdot||_{\AF}^{-nl/2}.
\end{equation}

We repeat the above procedure and get:
\begin{equation}
P^{(s)}(AI(\Pi_{\mathcal{F}})\otimes ||\cdot||^{-1/2},\tau)\sim_{E(\Pi_{\mathcal{F}});K} (2\pi i)^{-\frac{nl}{2}} \prod\limits_{\sigma \mid \tau} P^{(s(\sigma))}(\Pi_{\mathcal{F}},\sigma).
\end{equation}

We conclude the functoriality of arithmetic automorphic periods for automorphic induction:
\begin{thm}\label{AI theorem}
Let $F\supset K$ be a CM field of degree $d$ over $K$.

Let $\mathcal{F}/F$ be a cyclic extension of CM fields of degree $l$ and $\Pi_{\mathcal{F}}$ be a cuspidal representation of $GL_{n}(\AFFF)$ which has definable arithmetic automorphic periods. 

We assume that $\Pi_{\mathcal{F}}\ncong \Pi_{\mathcal{F}}^{g}$ for all $g \in Gal(\mathcal{F}/F)$ non trivial. We define $AI(\Pi_{\mathcal{F}})$ to be the automorphic induction of $\Pi_{\mathcal{F}}$. It is a cuspidal representation of $GL_{nl}(\AF)$.

We assume that $AI(\Pi_{\mathcal{F}})$ (resp. $AI(\Pi_{\mathcal{F}})\otimes ||\cdot||^{-1/2}$) also has definable arithmetic automorphic periods if $l$ is odd or $n$ is even (resp. if $l$ is even and $n$ is odd) (c.f. Section \ref{AI setting}). 

Let $I_{F}$ be any map from $\Sigma_{F;K}$ to $\{0,1,\cdots,nl\}$. Let $I_{\mathcal{F}}$ be the map from $\Sigma_{\mathcal{F};K}$ to $\{0,1,\cdots, n\}$ determined by $I_{F}$ and $\Pi_{\mathcal{F}}$ as in Lemma \ref{AI infinity sign}. Or locally let $0\leq s\leq nl$ be an integer and $s(\cdot)$ be as in Definition \ref{AI infinity sign local}.

If $l$ is odd or $n$ is even, we have:
\begin{eqnarray}
P^{(I_{F})}(AI(\Pi_{\mathcal{F}}))\sim_{E(\Pi_{\mathcal{F}});K} P^{(I_{\mathcal{F}})}(\Pi_{\mathcal{F}}) \nonumber \\
\text{or locally } P^{(s)}(AI(\Pi_{\mathcal{F}}),\tau) \sim_{E(\Pi_{\mathcal{F}});K} \prod\limits_{\sigma \mid \tau} P^{(s(\sigma))}(\Pi_{\mathcal{F}},\sigma).\nonumber
\end{eqnarray}
Otherwise we have:
\begin{eqnarray}
&&P^{(I_{F})}(AI(\Pi_{\mathcal{F}})\otimes ||\cdot||^{-1/2})\sim_{E(\Pi_{\mathcal{F}});K} (2\pi i)^{-\frac{nld}{2}}P^{(I_{\mathcal{F}})}(\Pi_{\mathcal{F}})\nonumber\\ &\text{ or locally}  &
P^{(s)}(AI(\Pi_{\mathcal{F}})\otimes ||\cdot||^{-1/2},\tau)
\sim_{E(\Pi_{\mathcal{F}});K} (2\pi i)^{-\frac{nl}{2}} \prod\limits_{\sigma \mid \tau} P^{(s(\sigma))}(\Pi_{\mathcal{F}},\sigma).\nonumber
\end{eqnarray}
\end{thm}

\section{Period relations under Galois action}
We are going to prove period relations for base change. Before that, we first prove that arithmetic periods are equivariant under Galois actions.\\

More precisely, let $F\supset K$ be a CM field and $\Pi$ be a cuspidal representation of $GL_{n}(\AF)$ which has definable arithmetic automorphic periods. \\

We fix any $I:\Sigma_{F;K}\rightarrow \{0,1,\cdots,n\}$ and take $\eta$, an algebraic Hecke character of $F$ such that $I(\Pi,\eta)=I$. Assuming Conjecture \ref{special value for similitude unitary group} and Theorem \ref{n*1} gives:
\begin{equation}\label{period relations Galois 1}
L(m,\Pi\otimes\eta) \sim_{E(\Pi)E(\eta);K} P^{(I)}(\Pi) \prod\limits_{\sigma\in\Sigma_{F;K}}p(\widecheck{\eta},\sigma)^{I(\sigma)}p(\widecheck{\eta},\overline{\sigma})^{n-I(\sigma)}
\end{equation}
for a critical point $m$ with both sides non zero.\\

Let $g\in Gal(F/K)$. We observe that $L(s,\Pi\otimes \eta)=L(s,\Pi^{g}\otimes \eta^{g})$. \\

We then get:
\begin{eqnarray}\label{period relations Galois 2}
&&L(m,\Pi\otimes\eta)=L(s,\Pi^{g}\otimes\eta^{g}) \\
&\sim_{E(\Pi)E(\eta);K} &P^{(I(\Pi^{g}))}(\Pi^{g},\eta^{g}) \prod\limits_{\sigma\in\Sigma_{F;K}}p(\widecheck{\eta^{g}},\sigma)^{I(\Pi^{g},\eta^{g})(\sigma)}p(\widecheck{\eta^{g}},\overline{\sigma})^{n-I(\Pi^{g},\eta^{g})(\sigma)} \nonumber
\end{eqnarray}

It is easy to see that $I(\Pi^{g},\eta^{g})(\sigma)=I(\Pi,\eta)(\sigma^{g^{-1}})$. \\

Moreover, by Proposition \ref{propCM} we have $p(\widecheck{\eta^{g}},\sigma)\sim_{E(\eta);K} p(\widecheck{\eta},\sigma^{g^{-1}})$. \\

We obtain:

\begin{eqnarray}
\prod\limits_{\sigma\in\Sigma_{F;K}}p(\widecheck{\eta^{g}},\sigma)^{I(\Pi^{g},\eta^{g})(\sigma)}&\sim_{E(\eta);K}& \prod\limits_{\sigma\in\Sigma_{F;K}}p(\widecheck{\eta},\sigma^{g^{-1}})^{I(\Pi,\eta)(\sigma^{g^{-1}})}\nonumber\\&\sim_{E(\eta);K}&
\prod\limits_{\sigma\in\Sigma_{F;K}}p(\widecheck{\eta},\sigma)^{I(\Pi,\eta)(\sigma)}.
\end{eqnarray}

Similarly, $\prod\limits_{\sigma\in\Sigma_{F;K}}p(\widecheck{\eta^{g}},\overline{\sigma})^{n-I(\Pi^{g},\eta^{g})(\sigma)}\sim_{E(\eta);K}\prod\limits_{\sigma\in\Sigma_{F;K}}p(\widecheck{\eta},\overline{\sigma})^{n-I(\Pi,\eta)(\sigma)}$.

We write $I^{g}:=I(\Pi^{g},\eta^{g})$. Then $I^{g}(\sigma)=I(\sigma^{g^{-1}})$. Compare with equation (\ref{period relations Galois 1}) and equation (\ref{period relations Galois 2}), we deduce that:
\begin{equation}
P^{(I)}(\Pi)\sim_{E(\Pi);K} P^{(I^{g})}(\Pi^{g}).
\end{equation}

We can moreover get relations of local periods. Let us fix $\sigma_{0}\in \Sigma_{F;K}$ and $0\leq s\leq n$ an integer.

We set $I(\sigma_{0})=s$ and $I(\sigma)=0$ for $\sigma\neq \sigma_{0}$. Then $I^{g}(\sigma_{0}^{g})=s$ and $I^{g}(\sigma)=0$ for $\sigma\neq \sigma_{0}^{g}$.

By the results in Section \ref{specify the factorization}, we have 
\begin{equation}
P^{(I)}(\Pi) \sim_{E(\Pi);K} P^{(s)}(\Pi,\sigma_{0})\prod\limits_{\sigma\neq \sigma_{0}} P^{(0)}(\Pi,\sigma) \sim_{E(\Pi);K} P^{(s)}(\Pi,\sigma_{0})\prod\limits_{\sigma\neq \sigma_{0}} p(\xi_{\Pi},\sigma)^{-1}
\end{equation}
and similarly:
\begin{equation}
P^{(I^{g})}(\Pi^{g})  \sim_{E(\Pi);K} P^{(s)}(\Pi^{g},\sigma_{0}^{g})\prod\limits_{\sigma\neq \sigma_{0}^{g}} p(\xi_{\Pi^{g}},\sigma)^{-1}.
\end{equation}

Again by Proposition \ref{propCM}, we have 
\begin{eqnarray}
\prod\limits_{\sigma\neq \sigma_{0}^{g}} p(\xi_{\Pi^{g}},\sigma)^{-1} &\sim_{E(\Pi);K} &\prod\limits_{\sigma\neq \sigma_{0}^{g}} p(\xi_{\Pi}^{g},\sigma)^{-1} \nonumber\\
&\sim_{E(\Pi);K} &\prod\limits_{\sigma\neq \sigma_{0}^{g}} p(\xi_{\Pi},\sigma^{g^{-1}})^{-1}\nonumber\\
&\sim_{E(\Pi);K} &\prod\limits_{\sigma\neq \sigma_{0}} p(\xi_{\Pi},\sigma)^{-1}.
\end{eqnarray}

We conclude that:
\begin{equation}
P^{(s)}(\Pi,\sigma_{0})\sim_{E(\Pi);K} P^{(s)}(\Pi^{g},\sigma_{0}^{g}).
\end{equation}

\begin{thm}
Let $F\supset K$ be a CM field and $\Pi$ be a cuspidal representation of $GL_{n}(\AF)$ which has definable arithmetic automorphic periods. We assume that Conjecture \ref{special value for similitude unitary group} is true. Let $g\in Gal(F/K)$, $\sigma\in\Sigma_{F;K}$ and $0\leq s\leq n$ be an integer. We have:
\begin{eqnarray}
&&P^{(I)}(\Pi)\sim_{E(\Pi);K} P^{(I^{g})}(\Pi^{g})\\
&\text{or locally }& P^{(s)}(\Pi,\sigma)\sim_{E(\Pi);K} P^{(s)}(\Pi^{g},\sigma^{g}).\label{local period relation Galois action}
\end{eqnarray}

\end{thm}

\section{Relations of global periods for base change}
Let $\mathcal{F}/F$ a cyclic extension of CM fields of degree $l$ as before. Let $\pi_{F}$ be a cuspidal representation of $GL_{n}(\AF)$. The strong base change of $\pi_{F}$ exists. We denote it by $\pi_{\mathcal{F}}$ or $BC(\pi_{F})$.\\

By the class field theory, we have $(F^{\times}N_{\AFFF^{\times}/\AF^{\times}}\AFFF^{\times} )\backslash\AF^{\times}\cong Gal(\mathcal{F}/F)$. Since $Gal(\mathcal{F}/F)$ is cyclic, its dual is also cyclic. We fix any generator of $Hom(Gal(\mathcal{F}/F),\C^{\times})$ which gives $\eta_{\mathcal{F}/F}$ a Hecke character of $F$.\\

We remark that $\eta_{\mathcal{F}/F}$ is an order $l$ Hecke character. In particular it has trivial infinity type. It is also unitary and thus conjugate self-dual.\\

We assume that $\pi_{F}\otimes \eta_{\mathcal{F}/F}^{t} \ncong \pi_{F}$ for all $1\leq t\leq l-1$. Then $\Pi_{\mathcal{F}}$ is cuspidal (Th\'eor\`em $4.2$ of \cite{arthurclozel}). We want to compare the arithmetic automorphic periods of $\Pi_{\mathcal{F}}$ to those of $\pi_{F}$ if they are defined.\\

We assume that $\pi_{F}$ has definable arithmetic automorphic periods. In other words, it is $3$-regular, cohomological and descends to unitary groups of any sign after tensoring an algebraic Hecke character. We know $\Pi_{\mathcal{F}}$ is also $3$-regular and cohomological. We assume that $\Pi_{\mathcal{F}}$ also descends to unitary groups of any sign after tensoring an algebraic Hecke character.\\

Let $I_{F}$ be any map from $\Sigma_{F;K}$ to the set $\{0,1,\cdots,n\}$.\\

We take $\eta$ an algebraic Hecke character of $F$ with $I(\pi_{F},\eta)=I_{F}$ such that $(\pi_{F},\eta)$ satisfies conditions in Theorem \ref{n*1 essential}. Let $\eta':=\eta\circ N_{\AFFF^{\times}/\AF^{\times}}$ be the base change of $\eta$ to $\mathcal{F}$. \\

There is a relation between the $L$-function of $\pi_{\mathcal{F}}$ and that of $\pi_{F}$, namely:
\begin{equation}\label{BC step 0}
L(s,\pi_{\mathcal{F}}\otimes \eta') =\prod\limits_{i=0}^{l-1} L(s,\pi_{F}\otimes \eta\eta_{\mathcal{F}/F}^{i}).
\end{equation}

We write $I_{\mathcal{F}}:=I(\pi_{\mathcal{F}},\eta')$. For any $\tau\in\Sigma_{F;K}$ and any $\sigma\in\Sigma_{\mathcal{F};K}$ with $\sigma\mid \tau$, it is easy to see that $I_{\mathcal{F}}(\sigma)=I_{F}(\tau)$. In other words, $I_{\mathcal{F}}$ is the composition of $I_{F}$ and $\Psi_{\mathcal{F}/F}$.\\

We assume that Conjecture \ref{special value for similitude unitary group} is true. We can interpret both sides in terms of arithmetic automorphic periods and CM periods and then deduce period relations.\\

More precisely, for a certain critical point $m$ we have:
\begin{equation}\label{BC step 1}
L(m,\pi_{\mathcal{F}}\otimes \eta') \sim_{E(\pi_{F})E(\eta);K} (2\pi i)^{mnld} P^{(I_{\mathcal{F}})}(\pi_{\mathcal{F}})\prod\limits_{\tau\in\Sigma_{F;K}} \prod\limits_{\sigma\mid \tau} p(\widecheck{\eta'},\sigma)^{I_{\mathcal{F}}(\sigma)}p(\widecheck{\eta'},\overline{\sigma})^{n-I_{\mathcal{F}}(\sigma)}
\end{equation}
with both sides non zero.\\

Since $\eta'=\eta\circ N_{\AFFF^{\times}/\AF^{\times}}$, we have $ p(\widecheck{\eta'},\sigma)= p(\widecheck{\eta},\sigma\mid_{F})$. Moreover, $I_{\mathcal{F}}(\sigma)=I_{F}(\tau)$ if $\sigma\mid \tau$. Therefore, 
\begin{equation}\label{BC step 2}
L(m,\pi_{\mathcal{F}}\otimes \eta') \sim_{E(\pi_{F})E(\eta);K} (2\pi i)^{mnld} P^{(I_{\mathcal{F}})}(\pi_{\mathcal{F}})\prod\limits_{\tau\in\Sigma_{F;K}} p(\widecheck{\eta},\tau)^{lI_{F}(\tau)}p(\widecheck{\eta},\overline{\tau})^{l(n-I_{F}(\tau))}.
\end{equation}

On the other hand, we apply Theorem \ref{n*1 essential} to $(\pi_{F},\eta\eta_{\mathcal{F}/F}^{i})$ and get:
\begin{eqnarray}\label{BC step 3}
&&L(m,\pi_{F}\otimes \eta\eta_{\mathcal{F}/F}^{i}) \\&\sim_{E(\pi_{F})E(\eta);K}& (2\pi i)^{mnd} P^{(I_{F})}(\pi_{F})\prod\limits_{\tau\in\Sigma_{F;K}} p(\widecheck{\eta\eta_{\mathcal{F}/F}^{i}},\tau)^{I_{F}(\tau)}p(\widecheck{\eta\eta_{\mathcal{F}/F}^{i}},\overline{\tau})^{n-I_{F}(\tau)}\nonumber\\
&\sim_{E(\pi_{F})E(\eta);K}& (2\pi i)^{mnd} P^{(I_{F})}(\pi_{F})\prod\limits_{\tau\in\Sigma_{F;K}} p(\widecheck{\eta},\tau)^{I_{F}(\tau)}p(\widecheck{\eta},\overline{\tau})^{n-I_{F}(\tau)}p(\widecheck{\eta_{\mathcal{F}/F}},\tau)^{iI_{F}(\tau)}p(\widecheck{\eta_{\mathcal{F}/F}^{-1}},\tau)^{i(n-I_{F}(\tau))}\nonumber\\
&\sim_{E(\pi_{F})E(\eta);K}& (2\pi i)^{mnd} P^{(I_{F})}(\pi_{F})\prod\limits_{\tau\in\Sigma_{F;K}} p(\widecheck{\eta},\tau)^{I_{F}(\tau)}p(\widecheck{\eta},\overline{\tau})^{n-I_{F}(\tau)}p(\widecheck{\eta_{\mathcal{F}/F}},\tau)^{2iI_{F}(\tau)-in}.\nonumber
\end{eqnarray}

Comparing equation (\ref{BC step 0}), (\ref{BC step 2}) and (\ref{BC step 3}), we get:
\begin{eqnarray}
P^{(I_{\mathcal{F}})}(\pi_{\mathcal{F}}) \sim_{E(\pi_{F});K} p^{I_{F}}(\pi_{F})
\prod\limits_{\tau\in \Sigma_{F;K}}\prod\limits_{i=0}^{l-1}p(\widecheck{\eta_{\mathcal{F}/F}},\tau)^{2iI_{F}(\tau)-in}\nonumber\\
\sim_{E(\pi_{F});K} p^{I_{F}}(\pi_{F})^{l}\prod\limits_{\tau\in \Sigma_{F;K}}
p(\widecheck{\eta_{\mathcal{F}/F}},\tau)^{l(l-1)I_{F}(\tau)-l(l-1)n/2}.
\end{eqnarray}

Since $p(\widecheck{\eta_{\mathcal{F}/F}},\tau)^{l}\sim_{K;K} p(\widecheck{\eta_{\mathcal{F}/F}^{l}},\tau)\sim_{K;K}1$, we have:
\begin{equation}
P^{(I_{\mathcal{F}})}(\pi_{\mathcal{F}}) \sim_{E(\pi_{F});K} p^{I_{F}}(\pi_{F})^{l}\prod\limits_{\tau\in \Sigma_{F;K}}
p(\widecheck{\eta_{\mathcal{F}/F}},\tau)^{-l(l-1)n/2}.
\end{equation}

If $l$ is odd, we have $p(\widecheck{\eta_{\mathcal{F}/F}},\tau)^{-l(l-1)n/2}\sim_{K;K} 1$. Otherwise we assume that $p(\widecheck{\eta_{\mathcal{F}/F}}^{l/2},\tau)\in E(\pi_{F})^{\times}$ for simplicity.

We conclude that:
\begin{equation}\label{BC global}
P^{(I_{\mathcal{F}})}(\pi_{\mathcal{F}}) \sim_{E(\pi_{F});K} p^{I_{F}}(\pi_{F})^{l}.
\end{equation}

\section{Relations of local periods for base change}
There are relations between local periods of $\pi_{F}$ and those of $\pi_{\mathcal{F}}$. 

Let $0\leq s_{0}\leq n$ be an integer. We fix $\tau_{0}\in \Sigma_{F;K}$. We take $I_{F}$ to be the map which sends $\tau_{0}$ to $s_{0}$ and $\tau\neq \tau_{0}$ to $0$. Equation (\ref{BC global}) then becomes:

\begin{equation}
\prod\limits_{\sigma_{0} \mid \tau_{0}}P^{(s_{0})}(\pi_{\mathcal{F}},\sigma_{0}) \prod\limits_{\tau\neq \tau_{0}}\prod\limits_{\sigma\mid \tau} P^{(0)}(\pi_{\mathcal{F}},\sigma) \sim_{E(\pi_{F});K} p^{(s_{0})}(\pi_{F},\tau_{0})^{l}\prod\limits_{\tau\neq \tau_{0}}p^{(0)}(\pi_{F},\tau)^{l}.
\end{equation}

Recall that for $ \sigma\mid \tau$, we have:
\begin{eqnarray}
P^{(0)}(\pi_{\mathcal{F}},\sigma)&\sim_{E(\pi_{F});K} &p(\xi_{\pi_{\mathcal{F}}},\sigma)^{-1}\nonumber\\
&\sim_{E(\pi_{F});K} &p(\xi_{\pi_{F}}\circ N_{\AFFF^{\times}/\AF^{\times}},\sigma)^{-1}\nonumber\\
&\sim_{E(\pi_{F});K} &p(\xi_{\pi_{F}},\tau)^{-1}\nonumber\\
&\sim_{E(\pi_{F});K} & P^{(0)}(\pi_{F},\tau).
\end{eqnarray}

We deduce that:
\begin{equation}\label{local base change step 1}
\prod\limits_{\sigma \mid \tau_{0}}P^{(s)}(\pi_{\mathcal{F}},\sigma)  \sim_{E(\pi_{F});K} P^{(s)}(\pi_{F},\tau_{0})^{l}.\end{equation}

\bigskip

We observe that $\pi_{\mathcal{F}}$ is $Gal(\mathcal{F}/F)$-invariant. The local periods are then $Gal(\mathcal{F}/F)$-invariant. 

Indeed, for any $g\in Gal(L/K)$, we have $\pi_{\mathcal{F}}^{g}\cong \pi_{\mathcal{F}}$. Theorem \ref{local period relation Galois action} implies that 
\begin{equation}
P^{(s)}(\pi_{\mathcal{F}},\sigma)\sim_{E(\pi_{F});K} P^{(s)}(\pi_{\mathcal{F}}^{g},\sigma^{g})\sim_{E(\pi_{F});K} P^{(s)}(\pi_{\mathcal{F}},\sigma^{g}).
\end{equation}
Recall that $Gal(\mathcal{F}/F)$ acts faithfully and transitively on the set $\{\sigma: \sigma\mid \tau\}$. We fix any $\sigma_{0}\mid \tau$ and then:
\begin{equation}\label{local base change step 2}
\prod\limits_{\sigma \mid \tau}P^{(s)}(\pi_{\mathcal{F}},\sigma) =\prod\limits_{g\in Gal(\mathcal{F}/F)}P^{(s)}(\pi_{\mathcal{F}},\sigma_{0}^{g})  \sim_{E(\pi_{F});K}P^{(s)}(\pi_{\mathcal{F}},\sigma_{0})^{l}.
\end{equation}

Comparing equation (\ref{local base change step 1}) and equation (\ref{local base change step 2}), we conclude that for any $\sigma\in \Sigma_{\mathcal{F};K}$:
\begin{equation}
P^{(s)}(\pi_{\mathcal{F}},\sigma)^{l} \sim_{E(\pi_{F});K}P^{(s)}(\pi_{F},\sigma\mid_{F})^{l}.
\end{equation}
Consequently, there exists an algebraic number $\lambda^{(s)}(\pi_{F}.\sigma)$ with $\lambda(\sigma)^{l}\in E(\pi_{F})^{\times}$ such that $P^{(s)}(\pi_{\mathcal{F}},\sigma) \sim_{E(\pi_{F})} \lambda^{(s)}(\pi_{F},\sigma)P^{(s)}(\pi_{F},\sigma\mid_{F}).
$

It is expected that $\lambda^{(s)}(\pi_{F},\sigma)$ is equivariant under Galois action. But we don't know how to prove it at this moment.

We summarize the above results on period relations for base change as follows.

\begin{thm}
Let $\mathcal{F}/F$ be a cyclic extension of CM fields of degree $l$. Let $\pi_{F}$ be a cuspidal representation of $GL_{n}(\AF)$. We denote by  $BC(\pi_{F})$ its strong base change to $\mathcal{F}$.

We assume that $\pi_{F}\otimes \eta_{\mathcal{F}/F}^{t} \ncong \pi_{F}$ for all $1\leq t\leq l-1$ and then $BC(\Pi_{F})$ is cuspidal (Th\'eor\`em $4.2$ of \cite{arthurclozel}). We assume that both $\pi_{F}$ and $BC(\pi_{F})$ have definable arithmetic automorphic periods.

Let $I_{F}$ be any map from $\Sigma_{F;K}$ to $\{0,1,\cdots,n\}$. We write $I_{\mathcal{F}}$ the composition of $I_{F}$ and $\Psi_{\mathcal{F}/F}$.

We then have:
\begin{eqnarray}
&&P^{(I_{\mathcal{F}})}(BC(\pi_{F})) \sim_{E(\pi_{F});K} p^{I_{F}}(\pi_{F})^{l} \\
&\text{or locally}& P^{(s)}(BC(\pi_{F}),\sigma)^{l} \sim_{E(\pi_{F});K}P^{(s)}(\pi_{F},\sigma\mid_{F})^{l}.
\end{eqnarray}
Consequently, we know 
\begin{equation}
P^{(s)}(BC(\pi_{F}),\sigma) \sim_{E(\pi_{F})} \lambda^{(s)}(\pi_{F},\sigma)P^{(s)}(\pi_{F},\sigma\mid_{F}).
\end{equation}
where $\lambda^{(s)}(\pi_{F},\sigma)$ is an algebraic number whose $l$-th power is in $E(\pi_{F})^{\times}$.
\end{thm}

\chapter{An automorphic version of Deligne's conjecture}

\section{A conjecture}
\begin{conj}\label{main conjecture split}
Let $n$ and $n'$ be two positive integers.
Let $\Pi$ and $\Pi'$ be cuspidal representations of $GL_{n}(\AF)$ and $GL_{n'}(\AF)$ respectively which have definable arithmetic automorphic periods.

 Let $m\in \Z+\cfrac{n+n'}{2}$ be critical for $\Pi\otimes \Pi'$. We predict that:
\begin{equation}
L(m,\Pi\times \Pi')\sim_{E(\Pi)E(\Pi');K}(2\pi i)^{nn'md}\prod\limits_{\sigma\in\Sigma_{F;K}}(\prod\limits_{j=0}^{n}P^{(j)}(\Pi,\sigma)^{sp(j,\Pi;\Pi',\sigma)}
\prod\limits_{k=0}^{n'}P^{(k)}(\Pi',\sigma)^{sp(k,\Pi';\Pi,\sigma)}).\nonumber
\end{equation}

\end{conj}

\bigskip

\begin{ex} (Known cases for the above conjecture:)

Let $F=K$ be the quadratic imaginary field. Then the above conjecture is already known in the following cases:

\begin{enumerate} 
\item  $n'=1$ and $m$ is strictly bigger than the central value. This is the main theorem in \cite{harris97}. We keep the notation as in Therem \ref{n*1 essential}. Let $\Pi'=\eta$. It is easy to verify that $sp(0,\Pi;\Pi')=n-s$, $sp(1,\Pi';\Pi)=s$, $sp(i,\Pi;\Pi')=0$ unless $i=s$ and $sp(s,\Pi;\Pi')=1$.

Recall that $P^{(0)}(\eta)\sim p(\widecheck{\eta},\iota)$ and  $P^{(1)}(\eta)\sim p(\widecheck{\eta},1)$. The formula in the above conjecture is the same with the formula in Theorem \ref{n*1 essential}.

\item $n'=n-1$, $\Pi$, $\Pi'$ conjugate self-dual in good position and $m>\frac{1}{2}$ or $m=\frac{1}{2}$ along with a non vanishing condition.

 In this case, we have $-a_{n}>b_{1}>-a_{n-1}>b_{2}>\cdots>b_{n-1}>-a_{1}$. Equivalently we have $sp(k,\Pi';\Pi)=1$ for all $0\leq k\leq n-1$; $sp(j,\Pi;\Pi')=1$ for $1\leq j\leq n-1$ and $=0$ for $j=0$ or $n$. Recall that $P^{(0)}(\Pi')\sim_{E(\Pi);K} P^{(n-1)}(\Pi')^{-1}$. Above conjecture is equivalent to that:
 \begin{equation}
 L(m,\Pi\times \Pi')\sim_{E(\Pi)E(\Pi');K}(2\pi i)^{n(n-1)m}\prod\limits_{j=1}^{n-1}P^{(j)}(\Pi)
\prod\limits_{k=2}^{n-2}P^{(k)}(\Pi')\nonumber
 \end{equation}
This is Theorem $6.11$ of \cite{harrismotivic} and Theorem $5.1$ of \cite{lincomptesrendus}.
\end{enumerate}
\end{ex}

\bigskip

The shall prove the following special cases of the above conjecture in the next chapters:
\begin{thm}\label{theorem split}
We assume that $\Pi$ and $\Pi'$ are $6$-regular if $F\neq K$ to guarantee the factorization of the arithmetic automorphic periods. We know Conjecture \ref{main conjecture split} is true for the following cases:
\begin{enumerate}
\item $n>n'$, the pair $(\Pi, \Pi')$ is in good position (see Definition \ref{definition good position}), $m$ is strictly bigger than the central value, or $m$ equals to the central value along with a non vanishing condition and moreover
\begin{enumerate}[label=(i)]
\item $\Pi$, $\Pi'$ conjugate self-dual if $n\nequiv n' (\text{mod 2})$, 
\item $\Pi$ conjugate self-dual, $\Pi'\otimes \psi_{F}^{-1}$ conjugate self-dual if $n\equiv n' (\text{mod 2})$.
\end{enumerate}
\item Any $n$, $n'$ and any position for $\Pi$, $\Pi'$ , the pair $(\Pi,\Pi')$ is very regular (\ref{r1*r2 very regular}) and moreover:
\begin{enumerate}[label=(i)]
\item $m=1$, $\Pi_{1}$, $\Pi_{2}$ conjugate self-dual if $n\equiv n' (\text{mod 2})$;
\item $m=\frac{1}{2}$, $\Pi_{1}$, $\Pi_{2}\otimes \psi_{F}^{-1}$ conjugate self-dual if $n\nequiv n' (\text{mod 2})$.
\end{enumerate}
\end{enumerate}
\end{thm}


\section{Compatibility with Deligne's conjecture over quadratic imaginary fields}

One see easily that Conjecture \ref{main conjecture split} is formally compatible with Conjecture \ref{Deligne automorphic}, Deligne's conjecture for automorphic pairs. For this, it is enough to compare the arithmetic automorphic period $P^{(j)}(\Pi)$ with the motivic period $Q^{(j)}(\Pi)$ where $\Pi$ is an  conjugate self-dual representation.\\

When $F$ is not $K$, this is difficult since we don't have geometric meanings for our local periods $P^{(s)}(\Pi,\sigma)$. But for the case when $F=K$, this is already discussed in Section $4$ of \cite{harrismotivic}. We now give a detailed explanation here. \\

First, let $\Pi$ be conjugate self-dual. In the construction of the arithmetic automorphic period, we have chosen $\xi$, an algebraic Hecke character of $\AK$, such that $\Pi^{\vee}\otimes \xi$ descends to a similitude unitary group. It is easy to verify that $\xi_{\Pi}=\cfrac{\xi}{\xi^{c}}$ (c.f. Theorem $VI.2.1$ or $VI.2.9$ of \cite{harristaylor}). The arithmetic automorphic period is defined to be the Petterson inner product of a rational class in the bottom stage of the Hodge filtration of a cohomology space related to $\Lambda^{j}M(\Pi^{c})\otimes M(\xi)$. In other words, $P^{(j)}(\Pi)$ is related to $Q_{n-j+1}(\Pi^{c})Q_{n-j+2}(\Pi^{c})\cdots Q_{n}(\Pi^{c})\times Q_{1}(\xi)$.\\

By Lemma \ref{period for Mc} we have $Q_{n-i+1}(\Pi^{c})\sim_{E(M)}Q_{i}(\Pi)^{-1}$ for all $1\leq i\leq n$. \\

By equation (\ref{Q1 Hecke}), we see $Q_{1}(\xi)\sim_{E(\xi);K} p(\cfrac{\xi^{c}}{\xi},1)\sim_{E(\xi);K} p(\xi_{\Pi}^{c},1)\sim_{E(\xi);K} p(\widecheck{\xi_{\Pi}^{c}},1)\sim_{E(\xi_{\Pi});K} \delta^{Del}(\xi_{\Pi})$. We deduce that:
\begin{equation}
Q_{n-j+1}(\Pi^{c})\cdots Q_{n}(\Pi^{c})\times Q_{1}(\xi^{c})
\sim_{E(\Pi)E(\xi);K} Q_{1}(\Pi)^{-1}Q_{2}(\Pi)^{-1}\cdots Q_{j}(\Pi)^{-1}\delta^{Del}(\xi_{\Pi}).
\end{equation}

\bigskip

Recall equation (\ref{Q^j}), the right hand side of the above formula is just $Q^{(j)}(\Pi)$ as expected.\\

\begin{rem}
We can also deduce the above result without passing to the motivic period of $\Pi^{c}$. In fact, we can also consider $P^{(j)}(\Pi)$ as Petterson inner product of a rational class in the bottom degree of a cohomology space related to $\Lambda^{n-j}M(\Pi)\otimes M(\xi^{c})$. It should be related to $Q_{j+1}(\Pi)Q_{j+2}(\Pi)\cdots Q_{n}(\Pi)Q_{1}(\xi^{c})$.\\

Lemma \ref{delta c} implies that
\begin{equation}\nonumber
\delta^{Del}(\xi_{\Pi}^{c})\sim_{E(M);K}(\prod\limits_{1\leq i\leq n}Q_{i}^{-1})\delta^{Del}(\xi_{\Pi}).
\end{equation}
Therefore, 
\begin{eqnarray}
Q^{(j)}(\Pi)&=&Q_{1}(\Pi)^{-1}Q_{2}(\Pi)^{-1}\cdots Q_{j}(\Pi)^{-1}\delta^{Del}(\xi_{\Pi})\nonumber \\ \nonumber&=&Q_{j+1}(\Pi)Q_{j+2}(\Pi)\cdots Q_{n}(\Pi)\delta^{Del}(\xi_{\Pi}^{c}).
\end{eqnarray}
We can deduce the comparison by the fact that $Q_{1}(\xi^{c})\sim_{E(\xi);K} p(\widecheck{\cfrac{\xi}{\xi^{c}}},1)\sim_{E(\xi);K} p(\widecheck{\xi_{\Pi}},1)\sim_{E(\xi);K}  \delta^{Del}(\xi_{\Pi}^{c})$.
\end{rem}
\bigskip

For the general cases, we write $\Pi=\Pi'\otimes \eta$ with $\Pi'$ conjugate self-dual. For the automorphic part, we see from Definition-Lemma (\ref{general period definition}) that \begin{equation}\nonumber P^{(j)}(\Pi)\sim_{E(\Pi);K}P^{(j)}(\Pi')p(\widecheck{\eta},1)^{j}p(\widecheck{\eta},\iota)^{n-j}.
\end{equation}

For the motivic part, we have $Q_{i}(\Pi)=Q_{i}(\Pi')Q_{1}(\eta)$ and $\Delta(\Pi)=\Delta(\Pi')\delta^{Del}(\eta)^{n}$. Therefore $Q^{(j)}(\Pi)=Q^{(j)}(\Pi')Q_{1}(\eta)^{-j}\delta^{Del}(\eta)^{n}$.\\

By (\ref{Q1 Hecke}) again, we see at first that $Q_{1}(\eta)\sim_{E(\eta);K} \cfrac{p(\widecheck{\eta^{c}},1)}{p(\widecheck{\eta},1)} $ and $\delta^{Del}(\eta)\sim_{E(\eta);K} p(\widecheck{\eta^{c}},1)$. We obtain finally:
\begin{equation}\nonumber Q^{(j)}(\Pi)\sim_{E(\Pi);K}Q^{(j)}(\Pi')p(\widecheck{\eta},1)^{j}p(\widecheck{\eta},\iota)^{n-j}.
\end{equation}

We have already related $P^{(j)}(\Pi')$ to $Q^{(j)}(\Pi')$. The relation for the general cases then comes.

\begin{rem}
We believe that the above comparison also works over general CM fields. However, the local periods $P^{(s)}(\Pi,\sigma)$ are not defined geometrically. It is expected that their geometric meaning can be obtained by comparing special values of $L$-functions. 
\end{rem}

\section{Simplify archimedean factors}

We observe that in Conjecture \ref{main conjecture split} the right hand side only concerns arithmetic automorphic periods and a power of $2\pi i$. 
Sometimes we will get a formula of $L(m,\Pi\otimes \Pi')$ which also involves archimedean factors as in Theorem $6.10$ of \cite{harrismotivic}. We need to show that the contribution of these archimedean factors is equivalent to a power of $2\pi i$:\\

\begin{prop}\label{proposition archimedean factor}
Let $\Pi$ and $\Pi'$ be as in Conjecture \ref{main conjecture split}. 
We assume that either the critical value $m$ is strictly bigger then the central value, either it is equal to the central value along with a nonvanishing condition on a certain $L$-function that we shall see in the proof.

If there exists an archimedean factor $a(m,\Pi_{\infty},\Pi'_{\infty})$ depending only on $m$, $\Pi_{\infty}$ and $\Pi'_{\infty}$ such that 
\begin{eqnarray}\label{simplify a}
&L(m,\Pi\times \Pi')\sim_{E(\Pi)E(\Pi');K}&\\\nonumber
&a(m,\Pi_{\infty},\Pi'_{\infty})\prod\limits_{\sigma\in\Sigma_{F;K}}(\prod\limits_{l=0}^{n}P^{(l)}(\Pi,\sigma)^{sp(l,\Pi;\Pi',\sigma)}
\prod\limits_{k=0}^{n'}P^{(k)}(\Pi',\sigma)^{sp(k,\Pi';\Pi,\sigma)})&,
\end{eqnarray}
then we have $a(m,\Pi_{\infty},\Pi'_{\infty})\sim_{E(\Pi);K}(2\pi i)^{nn'md}$. In particular, Conjecture \ref{main conjecture split} then follows.
\end{prop}

\bigskip

 Sometimes it is possible to calculate the archimedean factors directly. 
 
 A simpler way is to take $\Pi$ and $\Pi'$ as representations induced from Hecke characters. Then we may write the left hand side of equation (\ref{simplify a}) in terms of a power of $2\pi i$ and products of CM periods. For the right hand side, note that we have already related the arithmetic automorphic periods of a representation induced from Hecke characters and the CM periods by Theorem \ref{AI theorem}. 
 
 We shall deduce that the archimedean factor $a(m,\Pi_{\infty},\Pi'_{\infty})$ is equivalent to a power of $2\pi i$ if $\Pi$ and $\Pi'$ are induced from Hecke characters. But such representations can have any infinity type. The only non trivial point is that if $\Pi$ is conjugate self-dual then we may take a conjugate self-dual Hecke character such that its automorphic induction has the same infinity type as $\Pi$. We prove this in the lemma below. Hence the above proposition is true for any $\Pi$ and $\Pi'$. This is the idea of the proof of Theorem $5.1$ in \cite{lincomptesrendus}.\\
 
 \begin{lem}\label{chiandchi'}
Let $L\supset F$ be a cyclic extension of CM fields of degree $n$. We assume that $n$ is odd.
If $\Pi$ is a conjugate self-dual representation of $GL_{n}(\AF)$ then there exists $\chi$ a conjugate self-dual algebraic Hecke character of $L$ such that $\Pi_{\infty}\cong AI(\chi)_{\infty}$.

\end{lem}

\begin{dem}
We denote by $L^{+}$ the maximal totally real subfield of $L$.

We may take an algebraic Hecke character $\chi'$ of $L$ such that $\Pi_{\infty}\cong AI(\chi')_{\infty}$ (c.f. Lemma $4.1.1$ and paragraphs before Lemma $4.1.3$ in \cite{CHT}). 

Since $\Pi$ is conjugate self-dual, we see that $\chi'_{\infty}$ is conjugate self-dual. In particular, $\chi'|_{L^{+}}$ is trivial at infinity places. By Lemma $4.1.4$ of \cite{CHT}, we may find $\phi$ an algebraic Hecke character of $L$ with trivial infinity type such that $\phi\phi^{c}=\chi'\chi'\text{}^{c}$. Put $\chi=\chi'\phi^{-1}$. It is then a conjugate self-dual Hecke character with $\Pi_{\infty}\cong AI(\chi)_{\infty}$.

\end{dem}

 We now give the details of the proof for Proposition \ref{proposition archimedean factor}. 
 \begin{dem}
 For simplicity, we assume that both $n$ and $n'$ are odd. For general case, we have to twist $AI(\chi)$ or $AI(\chi')$ by $||\cdot||^{-1/2}\psi_{F}$ as before. The following proof goes through as well.\\
 
 We take $L\supset F$ (resp. $L'\supset F$) a CM field which is a cyclic extension of $F$ of degree $n$ (resp. $n'$). We assume that $L$ and $L'$ are linearly independent over $F$. Let $\LL:=LL'$. It is then a CM field of degree $nn'$ over $F$.
 
 We may take $\chi$ (resp. $\chi'$) an algebraic Hecke character of $L$ such that $\Pi_{\infty}=AI(\chi)_{\infty}$ (resp. $\Pi'_{\infty}=AI(\chi')_{\infty}$) where $AI(\chi)$ (resp. $AI(\chi')$) is the automorphic induction of $\chi$ (resp. $\chi'$) from $L$ (resp. $L'$) to $F$. Moreover, we may assume that $AI(\chi)$ and $AI(\chi')$ are cuspidal and have definable arithmetic automorphic periods.\\
 
 For $\sigma\in\Sigma_{F;K}$, we write $\sigma_{1},\cdots,\sigma_{n}$ for the elements in $\Sigma_{L;K}$ above $\sigma$ and $\sigma'_{1},\cdots,\sigma'_{n'}$ for the elements in $\Sigma_{L';K}$ above $\sigma$. Let $1\leq i\leq n$ and $1\leq j\leq n'$. We write $\sigma_{i,j}$ for the only element in $\Sigma_{\LL;K}$ such that $\sigma_{i,j}\mid_{L}=\sigma_{i}$ and $\sigma_{i,j}\mid_{L'}=\sigma'_{j}$. 
 
 We write $z^{a_{i}(\sigma)}\overline{z}^{-\omega(\Pi)-a_{i}(\sigma)}$ for the infinity type of $AI(\chi)$ at $\sigma_{i}$ and $z^{b_{j}(\sigma)}\overline{z}^{-\omega(\Pi')-b_{j}(\sigma)}$ for the infinity type of $AI(\chi)$ at $\sigma'_{j}$ .
 
Then $\Pi$ has infinity type $(z^{a_{i}(\sigma)}\overline{z}^{-\omega(\Pi)-a_{i}(\sigma)})_{1\leq i\leq n}$ and $\Pi'$ has infinity type $(z^{b_{j}(\sigma)}\overline{z}^{-\omega(\Pi')-b_{j}(\sigma)})_{1\leq j\leq n'}$ at $\sigma$.\\
 
 By equation (\ref{simplify a}), we have
 \begin{eqnarray}\label{archimedean factor step 1}
&L(m,AI(\chi) \times AI(\chi')) \sim_{E(\chi)E(\chi');K}&\\\nonumber
&a(m,\Pi_{\infty},\Pi'_{\infty})\prod\limits_{\sigma\in\Sigma_{F;K}}(\prod\limits_{l=0}^{n}P^{(l)}(AI(\chi),\sigma)^{sp(l,\Pi;\Pi',\sigma)}
\prod\limits_{k=0}^{n'}P^{(k)}(AI(\chi'),\sigma)^{sp(k,\Pi';\Pi,\sigma)})&
\end{eqnarray}

 On one hand, we have $L(m,AI(\chi) \times AI(\chi'))=L(m,(\chi\circ N_{\ALL/\AL})(\chi' \circ N_{\ALL/\mathcal{A}_{L'}}))$. 
 
 We observe that the infinity type of $(\chi\circ N_{\ALL/\AL})(\chi' \circ N_{\ALL/\mathcal{A}_{L'}})$ at $\sigma_{i,j}\in \Sigma_{\LL;K}$ is $z^{a_{i}(\sigma)+b_{j}(\sigma)}\overline{z}^{-\omega(\Pi)-\omega(\Pi')-a_{i}(\sigma)-b_{j}(\sigma)}$. 
 
 We denote $J_{\sigma}:=\{(i,j)\mid a_{i}(\sigma)+b_{j}(\sigma)<-\cfrac{\omega(\Pi)+\omega(\Pi')}{2}\}$.
 
We write $\chi_{\LL}=(\chi\circ N_{\ALL/\AL})(\chi' \circ N_{\ALL/\mathcal{A}_{L'}})$. By Blasius's result,
 \begin{equation}
 L(m,AI(\chi) \times AI(\chi'))\sim_{E(\chi)E(\chi');K} (2\pi i)^{mnn'd} \prod\limits_{\sigma\in \Sigma_{F;K}}\prod\limits_{(i,j)\in J_{\sigma}}p(\widecheck{\chi_{\LL}},\sigma_{i,j})\prod\limits_{(i,j)\notin J_{\sigma}}p(\widecheck{\chi_{\LL}},\overline{\sigma_{i,j}}).
 \end{equation}
 
 We need to assume that we may choose $\chi$ and $\chi'$ such that $L(m,AI(\chi) \times AI(\chi'))\neq 0$. When $m$ is strictly bigger then the central value, this is always true. When $m$ is equal to the central value, we assume this as a hypothesis.\\
 
 Recall that the CM periods are multiplicative and functorial for base change. Hence 
 \begin{eqnarray}\nonumber
 p(\widecheck{\chi_{\LL}},\sigma_{i,j})&\sim_{E(\chi)E(\chi');K} &p(\widecheck{\chi\circ N_{\ALL/\AL}},\sigma_{i,j})p(\widecheck{\chi\circ N_{\ALL/\mathcal{A}_{\mathcal{L'}}}},\sigma_{i,j})\\
 &\sim_{E(\chi)E(\chi');K} & p(\widecheck{\chi},\sigma_{i})p(\widecheck{\chi'},\sigma_{j}).
 \end{eqnarray}

We have deduced that \begin{eqnarray}\label{archimedean factor step 2}
 &L(m,AI(\chi) \times AI(\chi'))\sim_{E(\chi)E(\chi');K}&\\\nonumber
& (2\pi i)^{mnn'd} \prod\limits_{\sigma\in\Sigma_{F;K}} \prod\limits_{1\leq i\leq n}(p(\widecheck{\chi},\sigma_{i})^{s_{i}(\sigma)}p(\widecheck{\chi},\overline{\sigma_{i}})^{n'-s_{i}(\sigma)}) \prod\limits_{1\leq j\leq n'}(p(\widecheck{\chi'},\sigma_{j})^{t_{j}(\sigma)}p(\widecheck{\chi'},\overline{\sigma_{j}})^{n-t_{j}(\sigma)})&
\end{eqnarray}
where $s_{i}(\sigma)=\#\{1\leq j\leq n'\mid (i,j)\in J_{\sigma}\}$ and $t_{j}(\sigma)=\#\{1\leq i\leq n\mid (i,j)\in J_{\sigma}\}$.\\

\bigskip

On the other hand, for $0\leq l\leq n$, we have
\begin{equation}
P^{(l)}(AI(\chi),\sigma) \sim_{E(\Pi_{\mathcal{F}});K} \prod\limits_{1\leq i\leq n} P^{(u_{i}(l))}(\chi,\sigma_{i})
\end{equation}
where $u_{i}(l)=1$ if $a_{i}(\sigma)$ is in the $l$-th smallest numbers in the set $\{a_{i}(\sigma)\mid 1\leq i\leq n\}$ and $u_{i}(l)=0$ otherwise by Definition \ref{AI infinity sign local}.

We order $a_{i}(\sigma)$ and $b_{j}(\sigma)$ in decreasing order. We have $u_{i}(l)=1$ if and only if $i\geq n-l+1$. We get
\begin{equation}
P^{(l)}(AI(\chi),\sigma) \sim_{E(\Pi_{\mathcal{F}});K} \prod\limits_{1\leq i\leq n-l} P^{(0)}(\chi,\sigma_{i})\prod\limits_{n-l+1\leq i\leq n} P^{(1)}(\chi,\sigma_{i})
\end{equation}

Recall that $P^{(0)}(\chi,\sigma_{i})\sim_{E(\chi);K}p(\widecheck{\chi},\overline{\sigma_{i}})$ and $P^{(1)}(\chi,\sigma_{i})\sim_{E(\chi);K}p(\widecheck{\chi},\sigma_{i})$ by Remark \ref{remark n=1}. We obtain that:
\begin{equation}\label{archimedean factor step 3}
P^{(l)}(AI(\chi),\sigma) \sim_{E(\Pi_{\mathcal{F}});K} \prod\limits_{1\leq i\leq n-l} p(\widecheck{\chi},\overline{\sigma_{i}})\prod\limits_{n-l+1\leq i\leq n} p(\widecheck{\chi},\sigma_{i})
\end{equation}

\bigskip

Comparing equations (\ref{archimedean factor step 1}), (\ref{archimedean factor step 2}) and (\ref{archimedean factor step 3}), we observe that it remains to show:
\begin{equation}\label{archimedean temp}
\sum\limits_{l=n-i+1}^{n}sp(l,\Pi;\Pi',\sigma)=s_{i}(\sigma)
\end{equation}
\begin{equation}
\text{and }\sum\limits_{l=0}^{n-i}sp(l,\Pi;\Pi',\sigma)=n'-s_{i}(\sigma).
\end{equation}

Since $\sum\limits_{l=0}^{n}sp(l,\Pi;\Pi',\sigma)=n'$ by Lemma \ref{split lemma}. We see the above two equations are equivalent. We now prove the first one.

Recall by definition that $sp(l,\Pi;\Pi',\sigma)$ is the length of the $l$-th part of the sequence $b_{1}(\sigma)>b_{2}(\sigma)>\cdots>b_{n'}(\sigma)$ split by the numbers $-\cfrac{\omega(\Pi)+\omega(\Pi')}{2}-a_{n}>-\cfrac{\omega(\Pi)+\omega(\Pi')}{2}-a_{n-1}>\cdots>-\cfrac{\omega(\Pi)+\omega(\Pi')}{2}-a_{1}$.

Therefore, $\sum\limits_{l=n-i+1}^{n}sp(l,\Pi;\Pi',\sigma)=\#\{j \mid b_{j}<-\cfrac{\omega(\Pi)+\omega(\Pi')}{2}-a_{i}\}$. This is exactly $s_{i}(\sigma)$ as expected.

 \end{dem}
 
 \begin{rem}\label{a hope}
 Roughly speaking, the above proposition tells us that if we have a formula like equation (\ref{simplify a}) then the archimedean factor must be equivalent to a power of $2\pi i$. If one can show that the CM periods $p(\chi,\tau)$, $\tau\in \Sigma_{L;K}$ is algebraically independent, we can moreover prove that the power of arithmetic automorphic periods must be the split indices. 
 
 More precisely, the following statement is true:\\
 
If there exists an archimedean factor $a(m,\Pi_{\infty},\Pi'_{\infty})$ depending only on $m$, $\Pi_{\infty}$ and $\Pi'_{\infty}$ and integers $b(j,\Pi_{\infty};\Pi'_{\infty},\sigma)$, $c(k,\Pi'_{\infty};\Pi_{\infty},\sigma)$ for $0\leq j\leq n$, $0\leq k\leq n'$ and $\sigma\in\Sigma_{F;K}$ depending on $\Pi_{\infty}$, $\Pi'_{\infty}$ such that 
\begin{eqnarray}\label{simplify b}
&L(m,\Pi\times \Pi')\sim_{E(\Pi)E(\Pi');K}&\\\nonumber
&a(m,\Pi_{\infty},\Pi'_{\infty})\prod\limits_{\sigma\in\Sigma_{F;K}}(\prod\limits_{j=0}^{n}P^{(j)}(\Pi,\sigma)^{b(j,\Pi_{\infty};\Pi'_{\infty},\sigma)}
\prod\limits_{k=0}^{n'}P^{(k)}(\Pi',\sigma)^{c(k,\Pi'_{\infty};\Pi_{\infty},\sigma))})&,
\end{eqnarray}
then we have $b(j,\Pi_{\infty};\Pi'_{\infty},\sigma)=sp(j,\Pi;\Pi',\sigma)$ and $c(k,\Pi'_{\infty};\Pi_{\infty},\sigma)=sp(k,\Pi';\Pi,\sigma)$ provided that the local CM periods are algebraically independent. In particular, Conjecture \ref{main conjecture split} then follows.\\

The proof of the above statement is the same as proof for Proposition \ref{proposition archimedean factor}. We remark that the indices for the arithmetic automorphic periods are determined by equation (\ref{archimedean temp}).\\

This statement is very powerful. Sometimes it is easy to show that there exists a formula in the form of equation (\ref{simplify b}) but difficult to calculate the exact indices. In fact, one can explain in several lines that there exists such formulas for the cases in Theorem \ref{theorem split}. But we have devoted the next two whole chapters to calculate the precise indices. \\

Unfortunately we don't know how to prove the algebraically independency of the CM periods. So the calculation in the next two chapters are inevitable at the moment.

  \end{rem}

 \section{More discussions on the archimedean factors}
 As discussed in the previous section, one can leave the archimedean factors to the end of the proof and show that they contribute as a power of $2\pi i$. \\
 
 In our situation, we happen to be able to calculate the product of the archimedean factors directly. Let us first recall some archimedean factors. \\
 
 Let $\Pi_{\infty}$ (resp. $\Pi^{\#}_{\infty}$) be an algebraic regular generic representation of $GL_{n}(F\otimes_{\Q}\R$ (resp.$GL_{n-1}(F\otimes_{\Q}\R$)). We have defined:
 
 \begin{enumerate}
 \item $\Omega(\Pi_{\infty})$ which appears in the calculation of Whittaker period (c.f. Section $3.4$).
 
 \item $p(m,\Pi_{\infty},\Pi'_{\infty})$ which appears in the calculation of critical values for automorphic representations of $GL_{n}\times GL_{n-1}$ (c.f. Proposition \ref{Whittaker period theorem CM}).
 
 \item $Z(\Pi_{\infty})$ defined in equation (\ref{definition of Z}) by  \begin{equation}\nonumber
 Z(\Pi_{\infty}):=(2\pi i)^{d(m+\frac{1}{2})n(n-1)-\frac{d(n-1)(n-2)}{2}}\Omega(\Pi^{\#}_{\infty})^{-1} p(m,\Pi_{\infty},\Pi^{\#}_{\infty})^{-1}.
 \end{equation}
 \end{enumerate} 
 \bigskip
 
\begin{lem}\label{threeproduct}
The archimedean factors satisfy:
\begin{equation}\nonumber
Z(\Pi_{\infty})\Omega(\Pi'_{\infty})p(m,\Pi_{\infty},\Pi'_{\infty})\sim_{E(\Pi_{\infty})E(\Pi'_{\infty});K}(2\pi i)^{dn(n-1)(m+\frac{1}{2})-\frac{d(n-1)(n-2)}{2}}\end{equation}

 for all $m\geq 0$.
\end{lem}

We now take $\Pi$ and $\Pi^{\#}$ to be cuspidal conjugate self-dual representations of $GL_{n}(\AF)$ and $GL_{n-1}(\AF)$ respectively such that $(\Pi,\Pi^{\#})$ is in good position. We assume that they all have definable arithmetic automorphic periods. 

By Proposition \ref{Whittaker period theorem CM}, we have
\begin{equation}
L(\cfrac{1}{2}+m,\Pi\times \Pi^{\#})\sim_{E(\Pi)E(\Pi^{\#});K} p(m,\Pi_{\infty},\Pi^{\#}_{\infty})p(\Pi)p(\Pi^{\#})
 \end{equation}
 for some critical $\cfrac{1}{2}+m\geq \cfrac{1}{2}$. 

Recall from equation (\ref{Whittaker period local formula}) that:
 \begin{equation}
 p(\Pi)\sim_{E(\Pi)E(\Pi^{\#};K)} Z(\Pi_{\infty})\prod\limits_{\sigma\in\Sigma_{F;K}}\prod\limits_{1\leq i\leq n-1} P^{(i)}(\Pi).
  \end{equation}

We have a similar formula for $\Pi^{\#}$ since $\Pi^{\#}$ is also cuspidal. We then deduce that:
\begin{eqnarray}
&&L(\cfrac{1}{2}+m,\Pi\times \Pi^{\#})\\\nonumber
&\sim_{E(\Pi)E(\Pi^{\#});K}& p(m,\Pi_{\infty},\Pi^{\#}_{\infty})Z(\Pi_{\infty})Z(\Pi^{\#}_{\infty})\prod\limits_{\sigma\in\Sigma_{F;K}}(\prod\limits_{1\leq i\leq n-1} P^{(i)}(\Pi,\sigma)\prod\limits_{1\leq j\leq n-2} P^{(j)}(\Pi^{\#},\sigma))\\\nonumber
&\sim_{E(\Pi)E(\Pi^{\#});K}& p(m,\Pi_{\infty},\Pi^{\#}_{\infty})Z(\Pi_{\infty})Z(\Pi^{\#}_{\infty})\prod\limits_{\sigma\in\Sigma_{F;K}}(\prod\limits_{1\leq i\leq n-1} P^{(i)}(\Pi,\sigma)\prod\limits_{0\leq j\leq n-1} P^{(j)}(\Pi^{\#},\sigma)).
 \end{eqnarray}
 Here we have used the fact that $P^{(0)}(\Pi^{\#},\sigma)P^{(n-1)}(\Pi^{\#},\sigma)\sim_{E(\Pi^{\#});K} 1$ by Theorem \ref{special factorization theorem}.
\bigskip

Proposition \ref{proposition archimedean factor} then gives the following result on the archimedean factors:
\begin{prop}\label{propositiion archimedean n*n-1}
The archimedean factors satisfy:
\begin{equation}\nonumber
p(m,\Pi_{\infty},\Pi^{\#}_{\infty})Z(\Pi_{\infty})Z(\Pi^{\#}_{\infty})\sim_{E(\Pi_{\infty})E(\Pi'_{\infty});K}
(2\pi i)^{d(m+\frac{1}{2})n(n-1)}\end{equation}

provided $m\geq 1$ or $m=0$ along with a non vanishing condition for the central value of a certain $L$-function.
\end{prop}

This is Theorem $5.1$ of \cite{lincomptesrendus} when $F=K$ is a quadratic imaginary field.

\bigskip

Comparing Lemma \ref{threeproduct} and Proposition \ref{propositiion archimedean n*n-1}, we change the notation $\Pi^{\#}$ to $\Pi$ and deduce that:
\begin{cor}\label{divise}
We write $r=n-1$. For $\Pi_{\infty}$ an algebraic and generic representation of $GL_{n'}(F\otimes_{\Q} \R)$, we have

\begin{equation}\nonumber
Z(\Pi_{\infty})\Omega(\Pi_{\infty})^{-1}\sim _{E(\Pi_{\infty});K} (2\pi i)^{\frac{d(n-1)(n-2)}{2}}=(2\pi i)^{\frac{d(r-1)r}{2}}\end{equation}

provided $m\geq 1$ or $m=0$ along with a non vanishing condition for the central value of a certain $L$-function.
\end{cor}

In the following, we assume that $m\geq 1$, or $m=0$ along with a non vanishing condition for $\Pi$.

 \section{From quadratic imaginary fields to general CM fields}
\text{}

We shall prove Theorem \ref{theorem split} in the following two chapters over quadratic imaginary field. The proof only requires little change for general CM fields. This is because the automorphic arithmetic periods and the CM periods are all factorable. We now explain the details for the first case of Theorem \ref{theorem split} in the current section.\\

Let $\Pi$ and $\Pi'$ be cuspidal conjugate self-dual representations of $GL_{n}(\AF)$ and $GL_{n'}(\AF)$ which has definable arithmetic automorphic periods. We assume that $(\Pi,\Pi')$ is in good position and both $\Pi$ and $\Pi'$ are regular enough. For simplicity, we also assume that $n$ is even and $n'$ is odd.

Let $l=n-n'-1$. We take some conjugate self-dual Hecke characters $\chi_{1},\cdots,\chi_{l}$ such that if we write $\Pi^{\#}$ for the Langlands sum of $\Pi'$ and $\chi_{1},\cdots,\chi_{n}$ then $(\Pi,\Pi^{\#})$ is in good position. By the assumption on the parity of $n$ and $n'$ we know that $\Pi^{\#}$ is algebraic.

We may assume that for each $\sigma\in \Sigma_{F;K}$, the first index of the infinity type of $\chi_{i}$ is in decreasing order. Therefore, $I(\Pi,\chi_{i})(\sigma)$ is determined by the infinity type of $\Pi$ and $\Pi'$ at $\sigma$.

As explained in the introduction, the proof requires three main ingredients.

\text{}\\
\textbf{Ingredient A:} Theorem \ref{n*1 essential} says that for certain Hecke characters $\eta$ and critical points $m$ we have:
\begin{equation}
L(m,\Pi\otimes \eta) \sim_{E(\Pi)E(\eta);K} (2 \pi i)^{mnd} P^{(I(\Pi,\eta))}(\Pi) \prod\limits_{\sigma\in\Sigma}p(\widecheck{\eta},\sigma)^{I(\Pi,\eta)(\sigma)}p(\widecheck{\eta},\overline{\sigma})^{n-I(\Pi,\eta)(\sigma)}.
\end{equation}
where $I(\Pi,\eta)(\sigma)$ depends only the infinity type of $\Pi$ and $\eta$ at $\sigma$.

By Theorem \ref{special factorization theorem}, we may rewrite the above equation as:
\begin{equation}\label{from quadratic 1}
L(m,\Pi\otimes \eta) \sim_{E(\Pi)E(\eta);K} (2 \pi i)^{mnd}  \prod\limits_{\sigma\in\Sigma}[P^{(I(\Pi,\eta)(\sigma))}(\Pi,\sigma)p(\widecheck{\eta},\sigma)^{I(\Pi,\eta)(\sigma)}p(\widecheck{\eta},\overline{\sigma})^{n-I(\Pi,\eta)(\sigma)}].
\end{equation}

\text{}\\
\textbf{Ingredient B:} Proposition \ref{Whittaker period theorem CM} says that if $m\geq 0$ and $m+\cfrac{1}{2}$ is critical for $\Pi\times \Pi^{\#}$ then
\begin{equation}\label{from quadratic 2}
L(\cfrac{1}{2}+m,\Pi\times \Pi^{\#})\sim_{E(\Pi)E(\Pi^{\#});K} p(m,\Pi_{\infty},\Pi^{\#}_{\infty})p(\Pi)p(\Pi^{\#})
\end{equation}
where $p(\Pi)$ and $p(\Pi^{\#})$ are the Whittaker periods.

\text{}\\
\textbf{Ingredient C}: Corollary \ref{Shahidi Whittaker} implies that
\begin{equation}\label{from quadratic 3}
p(\Pi^{\#})\sim_{E(\Pi);K} p(\Pi')\cfrac{\Omega(\Pi_{\infty})}{\Omega(\Pi'_{\infty})}
\prod\limits_{1\leq i\leq l}L(1,\Pi'\otimes \chi_{i}^{\vee})\prod\limits_{1\leq i<j\leq l}L(1,\chi_{i}\times \chi_{j}^{\vee})
\end{equation}

Moreover, by Corollary \ref{Whittaker period local formula}, we have
\begin{eqnarray}\label{from quadratic 4}
p(\Pi)
\sim_{E(\Pi);K} Z(\Pi_{\infty}) \prod\limits_{\sigma\in\Sigma} \prod\limits_{1\leq i\leq n-1} P^{(i)}(\Pi,\sigma)\\\label{from quadratic 5}
p(\Pi')
\sim_{E(\Pi');K} Z(\Pi'_{\infty}) \prod\limits_{\sigma\in\Sigma} \prod\limits_{1\leq j\leq n'-1} P^{(j)}(\Pi',\sigma)
\end{eqnarray}

\bigskip
On one hand, note that $L(\cfrac{1}{2}+m,\Pi\times \Pi^{\#})=L(\cfrac{1}{2}+m,\Pi\times \Pi')\prod\limits_{1\leq i\leq l}L(\cfrac{1}{2}+m,\Pi\times \chi_{i})$. We replace $\eta$ by $\chi_{i}$ in equation (\ref{from quadratic 1}) and will get:
\begin{eqnarray}
&L(\cfrac{1}{2}+m,\Pi\times \Pi^{\#})=L(\cfrac{1}{2}+m,\Pi\times \Pi')\times& \\\nonumber&
\prod\limits_{\sigma\in\Sigma_{F;K}}[\prod\limits_{1\leq i\leq l}
P^{(I(\Pi,\chi_{i})(\sigma))}(\Pi,\sigma)p(\widecheck{\chi_{i}},\sigma)^{I(\Pi,\chi_{i})(\sigma)}p(\widecheck{\chi_{i}},\overline{\sigma})^{n-I(\Pi,\chi_{i})(\sigma)}
]&
\end{eqnarray}

\bigskip

On the other hand, apply equations (\ref{from quadratic 3}), (\ref{from quadratic 4}) and (\ref{from quadratic 5}) to the right hand side of equation (\ref{from quadratic 2}), we get:
\begin{eqnarray}
&L(\cfrac{1}{2}+m,\Pi\times \Pi^{\#}) \sim_{E(\Pi)E(\Pi^{\#});K} a_{0}(m,\Pi_{\infty},\Pi'_{\infty})& \times\\\nonumber
&\prod\limits_{\sigma\in\Sigma} [\prod\limits_{1\leq i\leq n-1} P^{(i)}(\Pi,\sigma)\prod\limits_{1\leq j\leq n'-1} P^{(j)}(\Pi',\sigma)]\prod\limits_{1\leq i\leq l}L(1,\Pi'\otimes \chi_{i}^{\vee})\prod\limits_{1\leq i<j\leq l}L(1,\chi_{i}\times \chi_{j}^{\vee})&
\end{eqnarray}
where $a_{0}(m,\Pi_{\infty},\Pi'_{\infty})$ is a non zero complex number depending only on $m$ and the infinity type.

We apply equation (\ref{from quadratic 1}) to $(\Pi',\chi_{i})$ and Blasius's result to $L(1,\chi_{i}\times \chi_{j}^{\vee})$, we get:
\begin{eqnarray}
&L(\cfrac{1}{2}+m,\Pi\times \Pi^{\#}) \sim_{E(\Pi)E(\Pi^{\#});K} a(m,\Pi_{\infty},\Pi'_{\infty})& \times\\\nonumber
&\prod\limits_{\sigma\in\Sigma} [\prod\limits_{1\leq i\leq n-1} P^{(i)}(\Pi,\sigma)\prod\limits_{1\leq j\leq n'-1} P^{(j)}(\Pi',\sigma)\prod\limits_{1\leq i\leq l}P^{(I(\Pi',\chi_{i})(\sigma))}(\Pi',\sigma)p(\chi_{i},\sigma)^{I_{1}(\Pi,\Pi')(\sigma)}
p(\chi_{i},\overline{\sigma})^{I_{2}(\Pi,\Pi')(\sigma)}
]&
\end{eqnarray}
where $a(m,\Pi_{\infty},\Pi'_{\infty})$ is an archimedean factor as before, $I_{1}(\Pi,\Pi')(\sigma)$ and $I_{2}(\Pi,\Pi')(\sigma)$ are two integers which depend only on the infinity type of $\Pi$ and $\Pi'$ at $\sigma$.

\bigskip

The first thing we need to show is that $I_{1}(\Pi,\chi_{i})(\sigma)=I(\Pi,\chi_{i})(\sigma)$ and $I_{2}(\Pi,\chi_{i})(\sigma)=n-I(\Pi,\chi_{i})(\sigma)$. Since we have ordered the first index of the infinity type of $\chi_{i}$ at $\sigma$ in decreasing order, we know that both sides only concern the infinity type of $\Pi$ and $\Pi'$ at the fixed place $\sigma$. So the proof is the same with the quadratic imaginary case.

\bigskip

We then deduce a formula in the following form:
\begin{eqnarray}
&L(m,\Pi\times \Pi')\sim_{E(\Pi)E(\Pi');K}&\\\nonumber
&a(m,\Pi_{\infty},\Pi'_{\infty})\prod\limits_{\sigma\in\Sigma_{F;K}}(\prod\limits_{j=0}^{n}P^{(j)}(\Pi,\sigma)^{b(j,\Pi_{\infty};\Pi'_{\infty},\sigma)}
\prod\limits_{k=0}^{n'}P^{(k)}(\Pi',\sigma)^{c(k,\Pi'_{\infty};\Pi_{\infty},\sigma))})&
\end{eqnarray}
where $b(j,\Pi_{\infty};\Pi'_{\infty},\sigma)$ and $c(k,\Pi_{\infty};\Pi'_{\infty},\sigma)$ are integers which depend only on $j$, $k$ and the infinity type of $\Pi_{\infty}$ and $\Pi'_{\infty}$ at $\sigma$.

\bigskip

If we know the CM periods are algebraically independent then we can finish the proof by Remark \ref{a hope}. Unfortunately this is hard to prove and hence we need to calculate $b(j,\Pi_{\infty};\Pi'_{\infty},\sigma)$ and $c(k,\Pi_{\infty};\Pi'_{\infty},\sigma)$ explicitly. Again, since they only concern infinity type of $\Pi_{\infty}$ and $\Pi'_{\infty}$ at $\sigma$, we may repeat our calculation for the quadratic imaginary field case for the fixed place $\sigma$. We shall see that the indices are just the split indices.

\bigskip

Finally we may show that the archimedean factor $a(m,\Pi_{\infty},\Pi'_{\infty})\sim_{E(\Pi)E(\Pi');K} (2\pi i)^{mnn'd}$ by Proposition \ref{proposition archimedean factor} and complete the proof.

\chapter{Special values of $L$-functions for automorphic pairs over quadratic imaginary fields}

\section{Settings, the simplest case}
 
In the current and the following chapters, let $\Pi$ and $\Pi'$ be conjugate self-dual cuspidal representations of $GL_{n}(\AK)$ and $GL_{n'}(\AK)$ respectively which have definable arithmetic automorphic periods. We will interpret the critical values for $L(s,\Pi\otimes\Pi')$ in terms of arithmetic automorphic periods when $(\Pi,\Pi')$ is in good position (see Definition \ref{definition good position}).\\

We write $(z^{a_{i}}\overline{z}^{-a_{i}})_{1\leq i\leq n}$ for the infinity type of $\Pi$ and $(z^{b_{j}}\overline{z}^{-b_{j}})_{1\leq j\leq n'}$ for the infinity type of $\Pi$. We may order $a_{i}$ and $b_{j}$ such that $a_{1}>a_{2}>\cdots>a_{n}$ and $b_{1}>b_{2}>\cdots>b_{n'}$.\\

 We assume $n$ is even and $n'$ is odd at first. Then the numbers $a_{i}, 1\leq i\leq n$ are half integers and the numbers $b_{j},1\leq j\leq n'$ are integers. \\

We assume the pair $(\Pi,\Pi')$ is \textbf{in good position}, namely, each $b_{i} $ are included in one of the intervals $]-a_{j+1},-a_{j}[$, $1\leq j\leq n-1$ and each such interval contains at most one $b_{i}$. \\

Let $w(1)>w(2)>\cdots>w(n)$ be the integers such that $a_{n-w(i)}>-b_{n'+1-i}>a_{n+1-w(i)}$ for all $1\leq i \leq n'$. More precisely, we have: 
\begin{eqnarray}\label{relationn*r}
\textbf{a}_{1}>\cdots>\textbf{a}_{n-w(1)}>-b_{n'}&>&\nonumber\\
\textbf{a}_{n+1-w(1)}>\cdots>\textbf{a}_{n-w(2)}>-b_{n'-1}&>&\nonumber\\
\cdots&&\nonumber\\
>\cdots>\textbf{a}_{n-w(n'+1-i)}>-b_{i}&>&\nonumber\\
\cdots\nonumber\\
\cdots>\textbf{a}_{n-w(n')}>-b_{1}&>&\nonumber\\
\textbf{a}_{n+1-w(n')}>\cdots>\textbf{a}_{n}.&&
\end{eqnarray}
It is easy to see:
\begin{eqnarray}\label{split number}
sp(0,\Pi';\Pi)=n-w(1), sp(n',\Pi';\Pi)=w(n')\nonumber\\
 sp(j,\Pi';\Pi)=w(j)-w(j+1) \text{ for all  }1\leq j\leq n'-1.
\end{eqnarray}

Hence we have
\begin{equation}\label{spandw}
w(j)=\sum\limits_{k=j}^{n'}sp(k,\Pi';\Pi)  \text{ for all  }1\leq j\leq n'.\end{equation}

We put $l=n-n'-1$. Let $\chi_{1},\chi_{2},\cdots,\chi_{l}$ be conjugate self-dual algebraic Hecke characters of $\AK$ of infinity type $z^{k_{1}}\overline{z}^{-k_{1}},z^{k_{2}}\overline{z}^{-k_{2}},\cdots z^{k_{l}}\overline{z}^{-k_{l}}$ respectively. We assume that $k_{1}>k_{2}>\cdots>k_{l}$ lie in different intervals $]-a_{j+1},-a_{j}[$ which do not contain any of $b_{i}$. 

More precisely, we have 

\begin{eqnarray}\label{mainrelation}
k_{1} > k_{2} >\cdots >k_{w(n')-1}>&b_{1}&>\nonumber \\
>k_{w(n')} >k_{w(n')+1}>\cdots >k_{w(n'-1)-2} >&b_{2}&>\nonumber \\
\cdots& &\nonumber \\
k_{w(n'+2-i)-i+2}>k_{w(n'+2-i)-i+3}>\cdots>k_{w(n'+1-i)-i}>&b_{i}&>\nonumber \\
\cdots& &\nonumber\\
k_{w(2)-n'+2}>k_{w(2)-n'+3}>\cdots k_{w(1)-n'}>&b_{n'}&>\nonumber\\
k_{w(1)-n'+1}>k_{w(1)-n'+2}>\cdots>k_{l}& &
\end{eqnarray}
and the above $l+n'=n-1$ numbers lie in different gaps between the $n$ numbers $-a_{n}>-a_{n-1}>\cdots>-a_{1}$. Note that in this case, the $n-1$ numbers above are integers and the $(a_{i})_{1\leq i\leq n}$ are half integers.\\

Let $\Pi^{\#}$ be the Langlands sum of $\Pi'$ and $\chi_{1},\chi_{2},\cdots,\chi_{l}$. It is a generic cohomological conjugate self-dual automorphic representation of $GL_{n-1}(\AK)$.

Let $m\geq 0$ be an integer. By Proposition \ref{Whittaker period theorem CM}, we know that if $m+\frac{1}{2}$ is critical for $\Pi\times \Pi^{\#}$, then 
\begin{equation}\label{simplest case starting point}
L(\frac{1}{2}+m,\Pi\times \Pi^{\#})\sim_{E(\Pi)E(\Pi^{\#});K} p(\Pi)p(\Pi^{\#})p(m,\Pi_{\infty},\Pi^{\#}_{\infty})\end{equation}

\bigskip

We shall simplify both sides of the above formula. We first calculate the left hand side.\\

We know $L(\frac{1}{2}+m,\Pi\times \Pi^{\#})=L(\cfrac{1}{2}+m, \Pi\times \Pi')\prod\limits_{j=1}^{l}L(\frac{1}{2}+m,\Pi\otimes \chi_{j})$.\\

For each $j$ with $1\leq j\leq l$, we apply Theorem \ref{main theorem CM} to $\Pi\otimes \chi_{j}$ and get:

\begin{equation}\nonumber
L(\frac{1}{2}+m,\Pi\otimes \chi_{j})\sim_{E(\Pi)E(\chi_{j});K} (2\pi i)^{(m+\frac{1}{2})n}P^{(s_{j})}(\Pi)
p(\widecheck{\chi_{j}},1)^{s_{j}}
p(\widecheck{\chi_{j}},\iota)^{n-s_{j}}\end{equation}

\begin{equation}\nonumber
\sim_{E(\Pi)E(\chi_{j});K} (2\pi i)^{(m+\frac{1}{2})n}P^{(s_{j})}(\Pi)
p(\widecheck{\chi_{j}},1)^{2s_{j}-n}\end{equation}

where \begin{equation}\label{sj}
s_{j}=\#\{1\leq i\leq n \mid k_{j}<-a_{i}\}=j+\#\{1\leq i\leq n'\mid b_{i}>k_{j}\}.
\end{equation}

\bigskip

By equation (\ref{mainrelation}), we see that 
\begin{eqnarray}\nonumber
s_{1}=1,s_{2}=2,\cdots, s_{w(n')-1}=w(n')-1,\\
\nonumber
s_{w(n')}=w(n')+1, s_{w(n')+1}=w(n')+2,\cdots, s_{w(n'-1)-2}=w(n'-1)-1,\\
\nonumber
\cdots\\
\nonumber
s_{w(1)-n'+1}=w(1)+1,s_{w(1)-n'+2}=w(1)+2,\cdots, s_{l}=l+n'=n-1.\end{eqnarray}

Shortly, $s_{1}<s_{2}<\cdots<s_{l}$ are the numbers in $\{1,2,\cdots, n-1\}\backslash\{w(n'),w(n'-1),\cdots,w(1))\}$.

We then deduce that:

\begin{equation}\nonumber
L(\frac{1}{2}+m,\Pi\times \Pi^{\#})\sim_{E(\Pi)E(\Pi')E;K}L(\cfrac{1}{2}+m, \Pi\times \Pi')(2\pi i)^{(m+\frac{1}{2})nl}\times\end{equation}

\begin{equation}\nonumber
\prod\limits_{i=1}^{n-1}P^{(i)}(\Pi)\prod\limits_{k=1}^{n'}P^{(w(k))}(\Pi)^{-1}\prod\limits_{j=1}^{l}p(\widecheck{\chi_{j}},1)^{2s_{j}-n}
\end{equation}

where $E$ is the compositum of $E(\chi_{j})$, $1\leq j\leq l$.

\section{Calculate the Whittaker period, the simplest case}
By Corollary \ref{Shahidi Whittaker}, we know that
\begin{equation}\nonumber
p(\Pi^{\#})
\sim_{E(\Pi);K} \Omega(\Pi^{\#}_{\infty})p(\Pi')\Omega(\Pi'_{\infty})^{-1}\prod\limits_{1\leq j\leq l}L(1,\Pi'\otimes \chi_{j}^{c})\prod\limits_{1\leq i<j\leq l}L(1,\chi_{i}\otimes \chi_{j}^{c}).\end{equation}
Recall that $\chi_{j}^{\vee}=\chi_{j}^{c}$ since $\chi_{j}$ is conjugate self-dual.

\paragraph{\textbf{Calculate }$\prod\limits_{1\leq j\leq l}L(1,\Pi'\otimes \chi_{j}^{c})$:} 
\text{}

For $1\leq j\leq l$, applying Theorem \ref{main theorem CM} to $\Pi'\times \chi_{j}^{c}$, we get

\begin{equation}\nonumber
L(1,\Pi'\otimes \chi_{j}^{c})\sim_{E(\Pi')E(\chi_{j});K} (2\pi i)^{n'}P^{(t_{j})}(\Pi')p(\widecheck{\chi_{j}^{c}},1)^{t_{j}}p(\widecheck{\chi_{j}^{c}},\iota)^{n'-t_{j}}\end{equation}

\begin{equation}\nonumber
\sim_{E(\Pi')E(\chi_{j});K} (2\pi i)^{n'}P^{(t_{j})}(\Pi')p(\widecheck{\chi_{j}},1)^{n'-2t_{j}}\end{equation}

where $t_{j}=\#\{1\leq i\leq n'\mid b_{i}-k_{j}<0\}$. 
The last step is due to the fact that $\chi_{j}^{c}=\chi_{j}^{-1}$.\\

It is easy to verify that $1$ is critical for $\Pi'\times \chi_{j}^{c}$ by considering the Hodge type and the original definition by Deligne. Recall that $\Pi'$ is of infinity type $(z^{b_{i}}\overline{z}^{-b_{i}})_{1\leq i\leq n'}$ and $\chi_{j}^{c}$ is of infinity type $z^{-k_{j}}\overline{z}^{k_{j}}$.

Compare with (\ref{sj}), we see that $t_{j}=n'-\#\{1\leq i\leq n'\mid b_{i}>k_{j}\}=n'+j-s_{j}$. Then $n'-2t_{j}=2s_{j}-n'-2j.$

Therefore, we have deduced that:
\begin{equation}\nonumber
\prod\limits_{1\leq j\leq l}L(1,\Pi'\otimes \chi_{j}^{c})\sim_{E(\Pi')E;K} (2\pi i)^{rl}\prod\limits_{j=1}^{l}P^{(t_{j})}(\Pi')\prod\limits_{j=1}^{l}p(\widecheck{\chi_{j}},1)^{2s_{j}-n'-2j}\end{equation}

\bigskip

\paragraph{\textbf{Calculate }$\prod\limits_{1\leq i<j\leq l}L(1, \chi_{i}\otimes\chi_{j}^{c})$:} 
\text{}

For $1\leq i<j\leq l$, since $k_{i}>k_{j}$, we have
\begin{equation}\nonumber
L(1,\chi_{i}\otimes \chi_{j}^{c})\sim_{E(\chi_{j});K} (2\pi i)p(\widecheck{\chi_{i}\chi_{j}^{c}},\iota)\sim_{E(\chi_{j});K} (2\pi i)p(\widecheck{\chi_{i}},1)^{-1}p(\widecheck{\chi_{j}},1)\end{equation}

by Blasius's result.\\

Therefore, we know that \begin{equation}\label{chitimeschi}
\prod\limits_{1\leq i<j\leq l}L(1, \chi_{i}\otimes\chi_{j}^{c}) \sim_{E;K} (2\pi i)^{\frac{l(l-1)}{2}}\prod\limits_{j=1}^{l}p(\widecheck{\chi_{j}},1)^{2j-l-1}.
\end{equation}

Since $(2s_{j}-n'-2j)+(2j-l-1)=2s_{j}-n'-l-1=2s_{j}-n$, we get finally
 \begin{equation}\nonumber
 p(\Pi^{\#})\sim_{E(\Pi)E;K} \Omega(\Pi^{\#}_{\infty})p(\Pi')\Omega(\Pi'_{\infty})^{-1}(2\pi i)^{rl+\frac{l(l-1)}{2}}\prod\limits_{j=1}^{l}P^{(t_{j})}(\Pi')\prod\limits_{j=1}^{l}p(\widecheck{\chi_{j}},1)^{2s_{j}-n}.\end{equation}

\section{Calculate the arithmetic automorphic periods and conclude, the simplest case}

Since $\Pi$ and $\Pi'$ are  cuspidal, we may apply Corollary \ref{Whittaker period local formula} and get: 
\begin{eqnarray}\label{ppi}
&&p(\Pi)\sim_{E(\Pi);K} Z(\Pi_{\infty})\prod\limits_{i=1}^{n-1}P^{(i)}(\Pi) 
\\\label{ppi0}
&\text{and }&p(\Pi')\sim_{E(\Pi');K} Z(\Pi'_{\infty})\prod\limits_{k=1}^{n'-1}P^{(k)}(\Pi').
\end{eqnarray}

Therefore, the right hand side of equation (\ref{simplest case starting point})
\begin{eqnarray}
&&p(\Pi)p(\Pi^{\#})p(m,\Pi_{\infty},\Pi^{\#}_{\infty})\nonumber \\
&\sim_{E(\Pi)E(\Pi')E;K} &Z(\Pi_{\infty})\Omega(\Pi^{\#}_{\infty})Z(\Pi'_{\infty})\Omega(\Pi'_{\infty})^{-1}p(m,\Pi_{\infty},\Pi^{\#}_{\infty})(2\pi i)^{rl+\frac{l(l-1)}{2}}\times \nonumber \\
&&\prod\limits_{j=1}^{l}p(\widecheck{\chi_{j}},1)^{2s_{j}-n}
 \prod\limits_{i=1}^{n-1}P^{(i)}(\Pi)\prod\limits_{k=1}^{n'-1}P^{(k)}(\Pi')\prod\limits_{j=1}^{l}P^{(t_{j})}(\Pi')\nonumber.
\end{eqnarray}

\paragraph{\textbf{Archimedean factors:}} 
\text{}
Recall that by lemma \ref{threeproduct}, we have
\begin{equation}\nonumber
Z(\Pi_{\infty})\Omega(\Pi^{\#}_{\infty})p(m,\Pi_{\infty},\Pi^{\#}_{\infty})\sim_{E(\Pi)E(\Pi')E;K} (2\pi i)^{n(n-1)(m+\frac{1}{2})-\frac{(n-1)(n-2)}{2}}.\end{equation}

By corollary \ref{divise}, we know 
\begin{equation}\nonumber
Z(\Pi'_{\infty})\Omega(\Pi'_{\infty})^{-1}\sim _{E(\Pi');K} (2\pi i)^{\frac{n'(n'-1)}{2}}.\end{equation}

Therefore $Z(\Pi_{\infty})\Omega(\Pi^{\#}_{\infty})Z(\Pi'_{\infty})\Omega(\Pi'_{\infty})^{-1}p(m,\Pi_{\infty},\Pi^{\#}_{\infty})(2\pi i)^{n'l+\frac{l(l-1)}{2}}$ 
\begin{equation}\nonumber
\sim_{E(\Pi)E(\Pi')E;K} (2\pi i)^{n(n-1)(m+\frac{1}{2})-\frac{(n-1)(n-2)}{2}+\frac{n'(n'-1)}{2}+n'l+\frac{l(l-1)}{2}}.\end{equation}

Note that $n-1=l+n'$ and hence ${n-1 \choose 2}={l\choose 2}+ln'+{n'\choose 2}$, we obtain that 

\begin{equation}\nonumber
Z(\Pi_{\infty})\Omega(\Pi^{\#}_{\infty})Z(\Pi'_{\infty})\Omega(\Pi'_{\infty})^{-1}p(m,\Pi_{\infty},\Pi^{\#}_{\infty})(2\pi i)^{n'l+\frac{l(l-1)}{2}} \end{equation}

\begin{equation}\nonumber
\sim_{E(\Pi)E(\Pi');K} (2\pi i)^{n(n-1)(m+\frac{1}{2})}.\end{equation}

\paragraph{\textbf{Arithmetic automorphic periods:}} 
\text{}

At last, we have to determine the value of $t_{j}=\#\{1\leq i\leq n'\mid b_{i}-k_{j}<0\}$ for $1\leq j\leq l$. 

For fixed $1\leq k\leq n'-1$, from the equation (\ref{mainrelation}), we see that the number of $1\leq j\leq l$ such that $t_{j}=k$ is $w(k)-w(k+1)-1$, the number of $1\leq j\leq l$ such that $t_{j}=r$ is $w(n')-1$, and the number of $1\leq j\leq l$ such that $t_{j}=0$ is $n-1-w(1)$.\\

For example, we see $t_{1}=t_{2}=\cdots=t_{w(n'-1)-1}=n'$, $t_{w(n')}=\cdots=t_{w(n'-1)-2}=n'-1$, $\cdots$, $t_{w(n'+2-i)-i+2}=\cdots =t_{w(n'+1-i)-1}=n'-i+1$, $\cdots$, and $t_{w(1)-n'+1}=\cdots=t_{l}=0$.\\

We then deduce that
\begin{eqnarray}
\prod\limits_{k=1}^{n'-1}P^{(k)}(\Pi')\prod\limits_{j=1}^{l}P^{(t_{j})}(\Pi')&=&\prod\limits_{k=1}^{n'-1}P^{(k)}(\Pi')^{w(k)-w(k+1)}P^{(0)}(\Pi')^{n-1-w(1)}P^{(n')}(\Pi')^{w(n')-1}\nonumber \\
&\sim_{E(\Pi')}&\prod\limits_{k=1}^{n'-1}P^{(k)}(\Pi')^{w(k)-w(k+1)}P^{(0)}(\Pi')^{n-w(1)}P^{(n')}(\Pi')^{w(n')}\nonumber \\
&=&\prod\limits_{k=0}^{n'}P^{(k)}(\Pi')^{sp(k,\Pi';\Pi)}\nonumber 
\end{eqnarray}

by the fact that $P^{0}(\Pi')\times P^{(n')}(\Pi')\sim_{E(\Pi')}1$ and equation (\ref{split number}).\\

Finally, we get $p(\Pi)p(\Pi^{\#})p(m,\Pi_{\infty},\Pi^{\#}_{\infty})$

 \begin{equation}\nonumber
\sim_{E(\Pi)E(\Pi')E;K}(2\pi i)^{n(n-1)(m+\frac{1}{2})}
 \prod\limits_{j=1}^{l}p(\widecheck{\chi_{j}},1)^{2s_{j}-n} \prod\limits_{i=1}^{n-1}P^{(i)}(\Pi)\prod\limits_{k=0}^{n'}P^{(k)}(\Pi')^{sp(k,\Pi';\Pi)}.\end{equation}
 
\bigskip

\paragraph{Final conclusion, simplest case:}
\text{}

When $L(\cfrac{1}{2}+m, \Pi\times \Pi^{\#})\neq 0$, we have that
\begin{eqnarray}\nonumber
L(\cfrac{1}{2}+m, \Pi\times \Pi')(2\pi i)^{(m+\frac{1}{2})nl}
\prod\limits_{i=1}^{n-1}P^{(i)}(\Pi)\prod\limits_{k=1}^{n'}P^{(w(k))}(\Pi)^{-1}\prod\limits_{j=1}^{l}p(\widecheck{\chi_{j}},1)^{2s_{j}-n}\\
\nonumber\sim_{E(\Pi)E(\Pi')E;K}
(2\pi i)^{n(n-1)(m+\frac{1}{2})}\prod\limits_{j=1}^{l}p(\widecheck{\chi_{j}},1)^{2s_{j}-n} \prod\limits_{i=1}^{n-1}P^{(i)}(\Pi)\prod\limits_{k=0}^{n'}P^{(k)}(\Pi')^{sp(k,\Pi';\Pi)}.\end{eqnarray}

We deduce that 
\begin{equation}\nonumber
L(\cfrac{1}{2}+m, \Pi\times \Pi')\sim_{E(\Pi)E(\Pi');K}(2\pi i)^{(m+\frac{1}{2})nn'}\prod\limits_{k=1}^{n'}P^{(w(k))}(\Pi)\prod\limits_{k=0}^{n'}P^{(k)}(\Pi')^{sp(k,\Pi';\Pi)}.\end{equation}

We can read from (\ref{relationn*r}) that for $0\leq i\leq n$, $sp(i,\Pi;\Pi')=0$ unless $i\in\{w(k)\mid 1\leq k\leq n'\}$. Moreover, if $i\in\{w(k)\mid 1\leq k\leq n'\}$ then $sp(i,\Pi;\Pi')=1$. We can then write the above formula in a symmetric way:
\begin{equation}\nonumber
L(\cfrac{1}{2}+m, \Pi\times \Pi')\sim_{E(\Pi)E(\Pi');K}(2\pi i)^{(m+\frac{1}{2})nn'}\prod\limits_{i=0}^{n}P^{(i)}(\Pi)^{sp(i,\Pi;\Pi')}\prod\limits_{k=0}^{n'}P^{(k)}(\Pi')^{sp(k,\Pi';\Pi)}.\end{equation}

\begin{rem}
If $L(\cfrac{1}{2}+m, \Pi\times \Pi')=0$, then the above formula is automatically true. Otherwise the condition $L(\cfrac{1}{2}+m, \Pi\times \Pi^{\#})\neq 0$ is equivalent to that $L(\cfrac{1}{2}+m, \Pi\times \chi_{j})\neq 0$ for all $1\leq j\leq l$. For $m\geq 1$, we can always choose $k_{j}$ and $\chi_{j}$ such that the above is true, see Section $3$ of \cite{harrisunitaryperiod}. For $m=0$, we don't know how to prove it at the moment. We will assume this is true henceforth.
\end{rem}

\section{Settings, the general cases}
Let $n>r$ be arbitrary integers. We still want to apply the previous strategy to get special values of $L$-function for $\Pi\times \Pi'$. But if we take $\Pi^{\#}$ to be Langlands sum of $\Pi'$ and some algebraic Hecke characters, it may be no longer algebraic. For example, if $n-1\nequiv n' (\text{mod }2)$, we know the Langlands parameters of $\Pi'$ are in $\Z+\frac{n'-1}{2}$. But the Langlands parameters of an algebraic representation of $GL_{n-1}$ should be in $\Z+\frac{n-1}{2}=\Z+\frac{n'}{2}$. In order to fix this, we will tensor the character $||\cdot||_{\AK}^{-\frac{1}{2}}\psi$, a Hecke character of infinity type $(\frac{1}{2},-\frac{1}{2})$, when necessary. \\

When $n-1\equiv r (\text{mod }2)$, we write $T_{1}=0$ and we will expand $\Pi'$ to an algebraic representation of $GL_{n-1}$ as previously. When $n-1\nequiv r (\text{mod }2)$, we write $T_{1}=\frac{1}{2}$ and we will expand $\Pi'\times ||\cdot||_{\AK}^{-\frac{1}{2}}\psi  $ to an algebraic representation of $GL_{n-1}$. In both cases, we assume the pair $(\Pi,\Pi')$ is \textbf{in good position}, namely, \begin{eqnarray}\label{n*r good position}
&&\text{ each }b_{i}+T_{1} \text{ are included in one of the intervals }]-a_{j+1},-a_{j}[, 1\leq j\leq n-1\nonumber\\
&&\text{ and each such interval contains at most one }b_{i}.
\end{eqnarray} 

Let $w(1)>w(2)>\cdots>w(n)$ be the integers such that \begin{equation}\label{wi}-a_{n+1-w(i)}>b_{n'+1-i}+T_{1}>-a_{n-w(i)} \end{equation} for all $1\leq i \leq n'$.\\

Let $\chi_{1},\chi_{2},\cdots,\chi_{l}$ be conjugate self-dual algebraic Hecke characters of $\AK$ of infinity type $z^{k_{1}}\overline{z}^{-k_{1}},z^{k_{2}}\overline{z}^{-k_{2}},\cdots z^{k_{l}}\overline{z}^{-k_{l}}$ respectively. These characters will help us expand $\Pi'$ or $\Pi'\otimes ||\cdot||_{\AK}^{-\frac{1}{2}}\psi $ to an algebraic representation of $GL_{n-1}$. Similarly, we will tensor them by $||\cdot||_{\AK}^{-\frac{1}{2}}\psi $ if $n \nequiv 0 (\text{mod } 2)$ to settle the parity issue. We write $T_{2}=\frac{1}{2}$ in this case and $0$ otherwise.\\

 We assume that $k_1+T_{2}>k_2+T_{2}>\cdots>k_{l}+T_{2}$ and lie in different intervals $]-a_{j+1},-a_{j}[$ which doesn't contain any of $b_{i}+T_{1}$.\\ 
 
More precisely, we have

\begin{eqnarray}
k_{1}+T_{2} > k_{2} +T_{2}>\cdots >k_{w(n')-1}+T_{2}>&b_{1}+T_{1}&>\nonumber \\
>k_{w(n')} +T_{2}>k_{w(n')+1}+T_{2}>\cdots >k_{w(n'-1)-2}+T_{2} >&b_{2}+T_{1}&>\nonumber \\
\cdots& &\nonumber \\
k_{w(n'+2-i)-i+2}+T_{2}>k_{w(n'+2-i)-i+3}+T_{2}>\cdots>k_{w(n'+1-i)-i}+T_{2}>&b_{i}+T_{1}&>\nonumber \\
\cdots& &\nonumber\\
k_{w(2)-n'+2}+T_{2}>k_{w(2)-n'+3}+T_{2}>\cdots k_{w(1)-n'}+T_{2}>&b_{n'}+T_{1}&>\nonumber\\
k_{w(1)-n'+1}+T_{2}>k_{w(1)-n'+2}+T_{2}>\cdots>k_{l}+T_{2}& &
\end{eqnarray}
and the above $l+k=n-1$ numbers lie in the gaps between the $n$ numbers $-a_{n}>-a_{n-1}>\cdots>-a_{1}$. Note the above $n-1$ numbers are in $\Z+\cfrac{n}{2}$ when $a_{i}\in \Z+\cfrac{n-1}{2}$ for all $1\leq i\leq n$.

There are four cases:

\begin{enumerate}[label=(\Alph*)]
\item $n$ is even and $n'$ is odd, then $T_{1}=0$ and $T_{2}=0$. We set $\Pi^{\#}=\Pi'\boxplus\chi_{1}\boxplus \chi_{2}\boxplus\cdots\boxplus \chi_{l}$ as in previous sections. 
\item $n$ is even and $n'$ is even, then $T_{1}=\frac{1}{2}$ and $T_{2}=0$. We set $\Pi^{\#}=(\Pi'\otimes ||\cdot||_{\AK}^{-\frac{1}{2}}\psi  )\boxplus\chi_{1}\boxplus \chi_{2}\boxplus\cdots\boxplus \chi_{l}$. 
\item $n$ is odd and $n'$ is even, then $T_{1}=0$ and $T_{2}=\frac{1}{2}$. We set $\Pi^{\#}=\Pi'  \boxplus(\chi_{1}\otimes ||\cdot||_{\AK}^{-\frac{1}{2}}\psi ) \boxplus (\chi_{2}\otimes ||\cdot||_{\AK}^{-\frac{1}{2}}\psi )\boxplus\cdots\boxplus (\chi_{l}\otimes ||\cdot||_{\AK}^{-\frac{1}{2}}\psi )$.
\item $n$ is odd and $n'$ is odd, then  $T_{1}=\frac{1}{2}$ and $T_{2}=\frac{1}{2}$. We set $\Pi^{\#}=(\Pi'\boxplus\chi_{1}\boxplus \chi_{2}\boxplus\cdots\boxplus \chi_{l})\otimes ||\cdot||_{\AK}^{-\frac{1}{2}}\psi $.
\end{enumerate}

In all cases, $\Pi^{\#}$ is a generic cohomological conjugate self-dual automorphic representation of $GL_{n-1}(\AK)$ and Proposition \ref{Whittaker period theorem CM} gives us that if $m+\frac{1}{2}$ is critical for $\Pi\times \Pi^{\#}$, then 
\begin{equation}\label{general case starting point}
L(\frac{1}{2}+m,\Pi\times \Pi^{\#})\sim_{E(\Pi)E(\Pi^{\#});K} p(\Pi)p(\Pi^{\#})p(m,\Pi_{\infty},\Pi^{\#}_{\infty}).\end{equation}

Again, we shall simplify both sides of this equation.

\section{Simplify the left hand side, general cases} 
\text{}

For the left hand side of equation (\ref{general case starting point}), we know by construction that:
\begin{enumerate}[label=(\Alph*)]
\item $L(\frac{1}{2}+m,\Pi\times \Pi^{\#})=L(\cfrac{1}{2}+m, \Pi\times \Pi')\prod\limits_{j=1}^{l}L(\frac{1}{2}+m,\Pi\otimes \chi_{j})$

\item $\begin{array}{lcl}
L(\frac{1}{2}+m,\Pi\times \Pi^{\#})&=&L(\cfrac{1}{2}+m, \Pi\times (\Pi'\otimes ||\cdot||_{\AK}^{-\frac{1}{2}}\psi ))\prod\limits_{j=1}^{l}L(\frac{1}{2}+m,\Pi\otimes \chi_{j})\\
&=& L(m, \Pi\times (\Pi'\otimes \psi ))\prod\limits_{j=1}^{l}L(\frac{1}{2}+m,\Pi\otimes \chi_{j})
\end{array}$

\item $\begin{array}{lcl}
L(\frac{1}{2}+m,\Pi\times \Pi^{\#})&=&L(\cfrac{1}{2}+m, \Pi\times \Pi')\prod\limits_{j=1}^{l}L(\frac{1}{2}+m,\Pi\otimes (\chi_{j}\otimes ||\cdot||_{\AK}^{-\frac{1}{2}}\psi) )\\
&=&L(\cfrac{1}{2}+m, \Pi\times \Pi')\prod\limits_{j=1}^{l}L(m,\Pi\otimes (\chi_{j}\otimes\psi) )

\end{array}$

\item $L(\frac{1}{2}+m,\Pi\times \Pi^{\#})=L(\cfrac{1}{2}, \Pi\times(\Pi'\otimes \psi ))\prod\limits_{j=1}^{l}L(m,\Pi\otimes (\chi_{j}\otimes \psi) )$

\end{enumerate}

\bigskip

We set $s_{j}=\#\{1\leq i\leq n \mid k_{j}+T_{2}<-a_{i}\}=j+\#\{1\leq i\leq n'\mid b_{i}+T_{1}>k_{j}+T_{2}\}$ and $t_{j}=\#\{1\leq i\leq n'\mid (b_{i}+T_{1})-(k_{j}+T_{2})<0\}$ as before. Recall that $s_{j}+t_{j}=n'+j$ for all $1\leq j\leq l$.

If $n$ is even (case (A) and (B)), we have for all $1\leq j\leq l$:

\begin{equation}\nonumber
L(\frac{1}{2}+m,\Pi\otimes \chi_{j})\sim_{E(\Pi)E(\chi_{j});K}(2\pi i)^{(m+\frac{1}{2})n}P^{(s_{j})}(\Pi)
p(\widecheck{\chi_{j}},1)^{2s_{j}-n}.\end{equation}

If $n$ is odd (case (C) and (D)), we have for all $1\leq j\leq l$:

\begin{equation}\nonumber
L(m,\Pi\otimes (\chi_{j})\otimes \psi)\sim_{E(\Pi)E(\chi_{j});K}(2\pi i)^{mn}P^{(s_{j})}(\Pi)
p(\widecheck{\chi_{j}},1)^{2s_{j}-n} p(\widecheck{\psi},1)^{s_{j}}p(\widecheck{\psi},\iota)^{n-s_{j}} .\end{equation}

Therefore for cases (A) and (B), we have
\begin{eqnarray}\nonumber
\prod\limits_{j=1}^{l}L(\frac{1}{2}+m,\Pi\otimes \chi_{j})\sim_{E(\Pi)E;K}(2\pi i)^{(m+\frac{1}{2})nl}\prod\limits_{k=1}^{n-1}P^{(k)}(\Pi)\prod\limits_{k=1}^{n'}P^{(w(k))}(\Pi)^{-1}\prod\limits_{j=1}^{l}p(\widecheck{\chi_{j}},1)^{2s_{j}-n}.
\end{eqnarray}

For cases (C) and (D), we put $s:=\sum\limits_{j=1}^{l}s_{j}$ and then we have:
\begin{eqnarray}\nonumber&
\prod\limits_{j=1}^{l}L(m,\Pi\otimes (\chi_{j}\otimes \psi))\sim_{E(\Pi)EE(\psi);K}&\\\nonumber
&(2\pi i)^{mnl}\times
\prod\limits_{k=1}^{n-1}P^{(k)}(\Pi)\prod\limits_{k=1}^{n'}P^{(w(k))}(\Pi)^{-1}\prod\limits_{j=1}^{l}p(\widecheck{\chi_{j}},1)^{2s_{j}-n} p(\widecheck{\psi},1)^{s}p(\widecheck{\psi},\iota)^{nl-s}&
\end{eqnarray}

\section{Simplify the right hand side, general cases}
\paragraph{\textbf{Calculate } $p(\Pi^{\#})$:} 
\text{}
Apply Corollary \ref{Shahidi Whittaker}, we get
\begin{enumerate}[label=(\Alph*)]
\item $p(\Pi^{\#})\sim_{E(\Pi^{\#});K} \Omega(\Pi^{\#}_{\infty})p(\Pi')\Omega(\Pi'_{\infty})^{-1}\prod\limits_{1\leq j\leq l}L(1,\Pi'\otimes \chi_{j}^{c})\prod\limits_{1\leq i<j\leq l}L(1,\chi_{i}\otimes \chi_{j}^{c})$

\item $p(\Pi^{\#})\sim_{E(\Pi^{\#});K} \Omega(\Pi^{\#}_{\infty})p(\Pi')\Omega(\Pi'_{\infty})^{-1}\prod\limits_{1\leq j\leq l}L(1,(\Pi'\otimes ||\cdot||_{\AK}^{-\frac{1}{2}}\psi)\otimes \chi_{j}^{c})\prod\limits_{1\leq i<j\leq l}L(1,\chi_{i}\otimes \chi_{j}^{c})$

\item $p(\Pi^{\#})\sim_{E(\Pi^{\#});K} \Omega(\Pi^{\#}_{\infty})p(\Pi')\Omega(\Pi'_{\infty})^{-1}\prod\limits_{1\leq j\leq l}L(1,\Pi'\otimes (\chi_{j}\otimes ||\cdot||_{\AK}^{-\frac{1}{2}}\psi)^{c})\prod\limits_{1\leq i<j\leq l}L(1,\chi_{i}\otimes \chi_{j}^{c})$
\item $p(\Pi^{\#})\sim_{E(\Pi^{\#});K} \Omega(\Pi^{\#}_{\infty})p(\Pi')\Omega(\Pi'_{\infty})^{-1}\prod\limits_{1\leq j\leq l}L(1,\Pi'\otimes \chi_{j}^{c})\prod\limits_{1\leq i<j\leq l}L(1,\chi_{i}\otimes \chi_{j}^{c})$
\end{enumerate}

Here we have used that: \begin{lem}
If $\eta$ is a conjugate self-dual Hecke character then:
\begin{equation}\nonumber
\cfrac{p(\Pi'\otimes \eta)}{\Omega((\Pi'\otimes \eta)_{\infty})}\sim_{E(\Pi')E(\eta);K}\cfrac{p(\Pi')}{\Omega(\Pi'_{\infty})}.\end{equation} \end{lem}
\begin{dem}
By Corollary \ref{Whittaker period local formula}, we have:
\begin{equation}
p(\Pi'\otimes \eta)\sim_{E(\Pi')E(\eta);K} Z((\Pi'\otimes \eta)_{\infty})\prod\limits_{1\leq i\leq n'-1} P^{(i)}(\Pi'\otimes \eta).
\end{equation}

By the definition of arithmetic automorphic period (c.f. Definition-Lemma \ref{definable periods}), we know $P^{(i)}(\Pi'\otimes \eta)\sim_{E(\Pi')E(\eta);K} p(\widecheck{\eta},1)^{i}p(\widecheck{\eta},\iota)^{n-i}$. The latter is equivalent to $p(\widecheck{\eta},1)^{2i-n}$ since $\eta$ is conjugate self-dual.

We see that:
\begin{equation}
\prod\limits_{1\leq i\leq n}P^{(i)}(\Pi'\otimes \eta)\sim_{E(\Pi')E(\eta);K} \prod\limits_{1\leq i\leq n}[P^{(i)}(\Pi')p(\widecheck{\eta},1)^{2i-n}]\sim_{E(\Pi')E(\eta);K} \prod\limits_{1\leq i\leq n}P^{(i)}(\Pi').
\end{equation}

By Corollary \ref{Whittaker period local formula}, This will imply that:
\begin{equation}\nonumber
\cfrac{p(\Pi'\otimes \eta)}{Z((\Pi'\otimes \eta)_{\infty})}\sim_{E(\Pi')E(\eta);K}\cfrac{p(\Pi')}{Z(\Pi'_{\infty})}.\end{equation}

But we know by Corollary \ref{divise} that $Z(\Pi'_{\infty})\sim_{E(\Pi'_{\infty};K)} (2\pi i)^{\frac{n'(n'-1)}{2}}\Omega(\Pi'_{\infty})$ and a similar formula for $(\Pi'\otimes\eta)_{\infty}$. The lemma then follows.

\end{dem}

\bigskip

By Theorem \ref{main theorem CM}, for all $1\leq j\leq l$, we have
\begin{equation}\nonumber
L(1,\Pi'\otimes \chi_{j}^{c})\sim_{E(\Pi')E(\chi_{j});K} (2\pi i)^{n'}P^{(t_{j})}(\Pi')p(\widecheck{\chi_{j}},1)^{n'-2t_{j}}.\end{equation}

Similarly, we have
\begin{eqnarray}\nonumber
&&L(1,(\Pi'\otimes ||\cdot||_{\AK}^{-\frac{1}{2}}\psi)\otimes \chi_{j}^{c})=L(\frac{1}{2},\Pi'\otimes (\psi \chi_{j}^{c}))\\ \nonumber
&\sim_{E(\Pi')E(\chi_{j})E(\psi);K} &(2\pi i)^{\frac{n'}{2}}P^{(t_{j})}(\Pi')p(\widecheck{\chi_{j}},1)^{n'-2t_{j}}p(\widecheck{\psi},1)^{t_{j}}p(\widecheck{\psi},\iota)^{n'-t_{j}};\end{eqnarray}

\begin{eqnarray}\nonumber
&\text{ and }& L(1,\Pi'\otimes (\chi_{i}\otimes ||\cdot||_{\AK}^{-\frac{1}{2}}\psi)^{c})=L(\frac{1}{2},\Pi'\otimes (\chi_{i}\otimes\psi)^{c})\\
\nonumber
&\sim_{E(\Pi')E(\chi_{j})E(\psi);K} &(2\pi i)^{\frac{n'}{2}}P^{(t_{j})}(\Pi')p(\widecheck{\chi_{j}},1)^{n'-2t_{j}}p(\widecheck{\psi},1)^{n'-t_{j}}p(\widecheck{\psi},\iota)^{t_{j}}.\end{eqnarray}

Along with equation (\ref{chitimeschi}), we get 

\begin{enumerate}[label=(\Alph*)]
\item 
and (D): $p(\Pi^{\#})\sim_{E(\Pi')EE(\psi);K} \Omega(\Pi^{\#}_{\infty})p(\Pi')\Omega(\Pi'_{\infty})^{-1}(2\pi i)^{n'l+\frac{l(l-1)}{2}}\times$
\begin{equation}\nonumber
\prod\limits_{j=1}^{l}P^{(t_{j})}(\Pi')\prod\limits_{j=1}^{l}p(\widecheck{\chi_{j}},1)^{2s_{j}-n}\end{equation}

\item:
$p(\Pi^{\#})\sim_{E(\Pi')EE(\psi);K} \Omega(\Pi^{\#}_{\infty})p(\Pi')\Omega(\Pi'_{\infty})^{-1}(2\pi i)^{\frac{n'l}{2}+\frac{l(l-1)}{2}}\times$
\begin{equation}\nonumber
\prod\limits_{j=1}^{l}P^{(t_{j})}(\Pi')\prod\limits_{j=1}^{l}p(\widecheck{\chi_{j}},1)^{2s_{j}-n}p(\widecheck{\psi},1)^{t}p(\widecheck{\psi},\iota)^{n'l-t}\end{equation}

\item:
$p(\Pi^{\#})\sim_{E(\Pi')EE(\psi);K} \Omega(\Pi^{\#}_{\infty})p(\Pi')\Omega(\Pi'_{\infty})^{-1}(2\pi i)^{\frac{n'l}{2}+\frac{l(l-1)}{2}}\times$
\begin{equation}\nonumber
\prod\limits_{j=1}^{l}P^{(t_{j})}(\Pi')\prod\limits_{j=1}^{l}p(\widecheck{\chi_{j}},1)^{2s_{j}-n}p(\widecheck{\psi},1)^{n'l-t}p(\widecheck{\psi},\iota)^{t}\end{equation}
\end{enumerate}
where $t=\sum\limits_{j=1}^{l}t_{j}=\sum\limits_{j=1}^{l}(n'+j-s_{j})=n'l+\cfrac{l(l+1)}{2}-s$.

We then apply equations (\ref{ppi}), (\ref{ppi0}) and Lemma \ref{threeproduct}, Corollary \ref{divise} to get:
\begin{enumerate}[label=(\Alph*)]
\item 
$p(\Pi)p(\Pi^{\#})p(m,\Pi_{\infty},\Pi^{\#}_{\infty})
 \sim_{E(\Pi)E(\Pi')E;K}
(2\pi i)^{n(n-1)(m+\frac{1}{2})}
\times$
\begin{equation}\nonumber
\prod\limits_{j=1}^{l}p(\widecheck{\chi_{j}},1)^{2s_{j}-n} \prod\limits_{i=1}^{n-1}P^{(i)}(\Pi)\prod\limits_{k=0}^{n'}P^{(k)}(\Pi')^{sp(k,\Pi';\Pi)}.\end{equation}

\item 
$p(\Pi)p(\Pi^{\#})p(m,\Pi_{\infty},\Pi^{\#}_{\infty})
 \sim_{E(\Pi)E(\Pi')EE(\psi);K}
(2\pi i)^{n(n-1)(m+\frac{1}{2})-\frac{n'l}{2}}
\times$
\begin{equation}\nonumber
\prod\limits_{j=1}^{l}p(\widecheck{\chi_{j}},1)^{2s_{j}-n} p(\widecheck{\psi},1)^{t}p(\widecheck{\psi},\iota)^{n'l-t}\prod\limits_{i=1}^{n-1}P^{(i)}(\Pi)\prod\limits_{k=0}^{n'}P^{(k)}(\Pi')^{sp(k,\Pi'\otimes \psi;\Pi)}.\end{equation}

\item 
$p(\Pi)p(\Pi^{\#})p(m,\Pi_{\infty},\Pi^{\#}_{\infty})
 \sim_{E(\Pi)E(\Pi')EE(\psi);K}
(2\pi i)^{n(n-1)(m+\frac{1}{2})-\frac{n'l}{2}}
\times$
\begin{equation}\nonumber
\prod\limits_{j=1}^{l}p(\widecheck{\chi_{j}},1)^{2s_{j}-n}p(\widecheck{\psi},1)^{n'l-t}p(\widecheck{\psi},\iota)^{t} \prod\limits_{i=1}^{n-1}P^{(i)}(\Pi)\prod\limits_{k=0}^{n'}P^{(k)}(\Pi')^{sp(k,\Pi'\otimes \psi;\Pi)}.\end{equation}

\item 
$p(\Pi)p(\Pi^{\#})p(m,\Pi_{\infty},\Pi^{\#}_{\infty})
 \sim_{E(\Pi)E(\Pi')EE(\psi);K}
(2\pi i)^{n(n-1)(m+\frac{1}{2})}
\times$
\begin{equation}\nonumber
\prod\limits_{j=1}^{l}p(\widecheck{\chi_{j}},1)^{2s_{j}-n} \prod\limits_{i=1}^{n-1}P^{(i)}(\Pi)\prod\limits_{k=0}^{n'}P^{(k)}(\Pi')^{sp(k,\Pi';\Pi)}.\end{equation}
\end{enumerate}

\section{Compare both sides, general cases}
At first, observe that 
\begin{equation}\nonumber
p(\widecheck{\psi},1)p(\widecheck{\psi},\iota)\sim_{E(\psi);K} p(\widecheck{\psi},1)p(\widecheck{\psi^{c}},1)\sim_{E(\psi);K} p(\widecheck{\psi\psi^{c}},1)\sim_{E(\psi);K} p(||\cdot||_{\AK}^{-1},1)\sim_{E(\psi)} 2\pi i.\end{equation}
We can then conclude:
\begin{enumerate}[label=(\Alph*)]
\item 
$L(\cfrac{1}{2}+m, \Pi\times \Pi')\sim_{E(\Pi)E(\Pi');K}$

\begin{equation}\nonumber
(2\pi i)^{(m+\frac{1}{2})nn'}\prod\limits_{k=1}^{n'}P^{(w(k))}(\Pi)\prod\limits_{k=0}^{n'}P^{(k)}(\Pi')^{sp(k,\Pi';\Pi)}\end{equation}

\item Since $n(n-1)(m+\frac{1}{2})-\frac{n'l}{2}-(m+\frac{1}{2})nl=(m+\frac{1}{2})n(n-1-l)-\frac{n'l}{2}=(m+\frac{1}{2})nn'-\frac{n'l}{2}=mnn'+\frac{nn'}{2}-\frac{n'l}{2}$, we have
\begin{eqnarray}
L(m, \Pi\times (\Pi'\otimes\psi))\sim_{E(\Pi)E(\Pi')E(\psi);K}(2\pi i)^{mnn'+\frac{nn'}{2}-\frac{n'l}{2}}\times\\\nonumber
\prod\limits_{k=1}^{n'}P^{(w(k))}(\Pi)\prod\limits_{k=0}^{n'}P^{(k)}(\Pi')^{sp(k,\Pi'\otimes \psi;\Pi)}p(\widecheck{\psi},1)^{t}p(\widecheck{\psi},\iota)^{n'l-t}\end{eqnarray}

Since $(2\pi i)^{\frac{nn'}{2}-\frac{n'l}{2}}\sim_{E(\psi)} p(\widecheck{\psi},1)^{\frac{nn'}{2}-\frac{n'l}{2}}p(\widecheck{\psi},\iota)^{\frac{nn'}{2}-\frac{n'l}{2}}$,  and 

\begin{eqnarray}
\frac{nn'}{2}-\frac{n'l}{2}+t&=&\frac{nn'}{2}-\frac{n'l}{2}+(n'l+\frac{l(l+1)}{2}-s)=\frac{nn'}{2}+\frac{n'l}{2}+\frac{l(l+1)}{2}-s\nonumber \\
&=&\frac{nn'}{2}+\frac{(n'+l+1)l}{2}-s=\frac{nn'}{2}+\frac{nl}{2}-s\nonumber \\
&=&\frac{n(n'+l)}{2}-s=\frac{n(n-1)}{2}-s;\nonumber 
\end{eqnarray}

\begin{eqnarray}
 \frac{nn'}{2}-\frac{n'l}{2}+n'l-t&=&\frac{nn'}{2}-\frac{n'l}{2}+n'l-(n'l+\frac{l(l+1)}{2}-s)=s+\frac{nn'}{2}-\frac{(n'+l+1)l}{2}\nonumber \\
&=&s+nn'-\frac{nn'}{2}-\frac{nl}{2}=s+nn'-\frac{n(n-1)}{2}\nonumber 
\end{eqnarray}

We get $L(m, \Pi\times (\Pi'\otimes\psi))\sim_{E(\Pi)E(\Pi')E(\psi);K}(2\pi i)^{mnn'}\times$

\begin{equation}\nonumber
\prod\limits_{k=1}^{n'}P^{(w(k))}(\Pi)\prod\limits_{k=0}^{n'}P^{(k)}(\Pi')^{sp(k,\Pi';\Pi)}p(\widecheck{\psi},1)^{\frac{n(n-1)}{2}-s}p(\widecheck{\psi},\iota)^{s+nn'-\frac{n(n-1)}{2}}\end{equation}

\item Since $n(n-1)(m+\frac{1}{2})-\frac{n'l}{2}-mnl=n(n-1)(m+\frac{1}{2})-\frac{n'l}{2}-(m+\frac{1}{2})nl+\frac{nl}{2}=(m+\frac{1}{2})nn'+\frac{nl}{2}-\frac{n'l}{2}$, we have
\begin{eqnarray}
L(\cfrac{1}{2}+m, \Pi\times \Pi')\sim_{E(\Pi)E(\Pi')E(\psi);K}(2\pi i)^{(m+\frac{1}{2})nn'+\frac{nl}{2}-\frac{n'l}{2}}\times\\\nonumber
\prod\limits_{k=1}^{n'}P^{(w(k))}(\Pi)\prod\limits_{k=0}^{n'}P^{(k)}(\Pi')^{sp(k,\Pi';\Pi)}p(\widecheck{\psi},1)^{n'l-t-s}p(\widecheck{\psi},\iota)^{t+s-nl}
\end{eqnarray}

Moreover, we know $t+s=n'l+\frac{l(l+1)}{2}$, we have $2(t+s)=2n'l+(l+1)l=n'l+(n'+l+1)l=n'l+nl$. Thus $n'l-t-s=t+s-nl=\frac{n'l}{2}-\frac{nl}{2}$. We then get $p(\widecheck{\psi},1)^{n'l-t-s}p(\widecheck{\psi},\iota)^{t+s-nl}=p(\widecheck{\psi\otimes \psi^{c}},1)^{\frac{n'l}{2}-\frac{nl}{2}}=(2\pi i)^{\frac{n'l}{2}-\frac{nl}{2}}.$ 

Therefore:
\begin{equation}\nonumber
L(\cfrac{1}{2}+m, \Pi\times \Pi')\sim_{E(\Pi)E(\Pi');K}(2\pi i)^{(m+\frac{1}{2})nn'}\prod\limits_{k=1}^{n'}P^{(w(k))}(\Pi)\prod\limits_{k=0}^{n'}P^{(k)}(\Pi')^{sp(k,\Pi';\Pi)}.\end{equation}

\item Similarly, since  $n(n-1)(m+\frac{1}{2})-mnl=n(n-1)m+\frac{n(n-1)}{2}-mnl=mnn'+\frac{n(n-1)}{2}$, we have
\begin{eqnarray}
L(m, \Pi\times (\Pi'\otimes\psi)) \sim_{E(\Pi)E(\Pi')E(\psi);K}(2\pi i)^{mnn'}\times\\\nonumber
\prod\limits_{k=1}^{n'}P^{(w(k))}(\Pi)\prod\limits_{k=0}^{n'}P^{(k)}(\Pi')^{sp(k,\Pi'\otimes \psi;\Pi)}p(\widecheck{\psi},1)^{\frac{n(n-1)}{2}-s}p(\widecheck{\psi},\iota)^{s+nn'-\frac{n(n-1)}{2}}.
\end{eqnarray}

It is easy to verify that $s-nl+\frac{n(n-1)}{2}=s-nl+n(n-1)-\frac{n(n-1)}{2}=s+nn'-\frac{n(n-1)}{2}$.

\end{enumerate}

\section{Final conclusion: general cases}
Before concluding, we notice that in case (B) or (D), 
\begin{equation}\nonumber
s=\sum\limits_{1\leq j\leq n-1}s_{j}=\sum\limits_{j=1}^{n-1}j-\sum\limits_{j=1}^{n'}w(j)=\frac{n(n-1)}{2}-\sum\limits_{j=1}^{n'}w(j).
\end{equation}
 Recall that $w(j)=\sum\limits_{k=j}^{n'}sp(k,\Pi'\otimes \psi;\Pi) \text{ for all  }1\leq k\leq n'$ by (\ref{spandw}). Therefore:
\begin{eqnarray}\label{sumofw}
\frac{n(n-1)}{2}-s=\sum\limits_{j=1}^{n'}w(j)&=&\sum\limits_{j=1}^{n'}\sum\limits_{j\leq k\leq n'}sp(j,\Pi'\otimes\psi;\Pi) =\sum\limits_{k=1}^{n'}k*sp(k,\Pi'\otimes\psi;\Pi)\nonumber \\
&=&\sum\limits_{k=0}^{n'}k*sp(k,\Pi' \otimes\psi;\Pi);
\end{eqnarray}

\begin{eqnarray}
\text{and }s+nn'-\frac{n(n-1)}{2}&=&nn'-\sum\limits_{k=0}^{n'}k*sp(k,\Pi'\otimes\psi;\Pi)\nonumber \\
&=&r\sum\limits_{k=0}^{n'}sp(k,\Pi'\otimes\psi;\Pi)-\sum\limits_{k=0}^{n'}j*sp(k,\Pi'\otimes\psi;\Pi)\nonumber \\
&=&\sum\limits_{k=0}^{n'}(n'-k)sp(k,\Pi'\otimes\psi;\Pi)\nonumber 
\end{eqnarray}
by Lemma \ref{split lemma} which says that $\sum\limits_{k=0}^{n'}sp(k,\Pi'\otimes\psi;\Pi)=n$.

Therefore, we get
\begin{eqnarray}
&&\prod\limits_{k=0}^{n'}P^{(k)}(\Pi')^{sp(k,\Pi'\otimes \psi;\Pi)}p(\widecheck{\psi},1)^{\frac{n(n-1)}{2}-s}p(\widecheck{\psi},\iota)^{s+nn'-\frac{n(n-1)}{2}}\nonumber\\
\nonumber
&\sim_{E(\Pi')E(\psi)}&
\prod\limits_{k=0}^{n'}P^{(k)}(\Pi')^{sp(k,\Pi'\otimes \psi;\Pi)}
p(\widecheck{\psi},1)^{\sum\limits_{k=0}^{n'}k*sp(k,\Pi' \otimes\psi;\Pi)}p(\widecheck{\psi},\iota)^{\sum\limits_{k=0}^{n'}(n'-k)sp(k,\Pi'\otimes\psi;\Pi)}\\
\nonumber &\sim_{E(\Pi')E(\psi)}& \prod\limits_{k=0}^{n'}
\left(P^{(k)}(\Pi')p(\widecheck{\psi},1)^{k}p(\widecheck{\psi},\iota)^{n'-k}\right)^{sp(k,\Pi'\otimes \psi;\Pi)}.
\end{eqnarray}

Recall that $P^{(k)}(\Pi'\otimes \psi):=P^{(k)}(\Pi')p(\widecheck{\psi},1)^{k}p(\widecheck{\psi},\iota)^{n'-k}$ by definition, we obtain that:

\begin{thm}\label{n*rmaintheorem}
Let $n>n'$ be two positive integers. Let $K$ be a quadratic imaginary field.
Let $\Pi$ and $\Pi'$ be cuspidal representations of $GL_{n}$ and $GL_{n'}$ respectively which are very regular, cohomological, conjugate self-dual and supercuspidal at at least two finite split places.
We assume that $(\Pi,\Pi')$ is in good position in the sense of Definition \ref{definition good position}.

\begin{enumerate}[label=(\roman*)]
\item If $n\nequiv n' (\text{mod }2)$, then for any critical value $m+\frac{1}{2}$ for $\Pi\otimes \Pi'$ such that $m\geq 1$, or $m\geq 0$ along with a non-vanishing condition, we have:
\begin{equation}\nonumber
L(\cfrac{1}{2}+m, \Pi\times \Pi')\sim_{E(\Pi)E(\Pi');K}
(2\pi i)^{(m+\frac{1}{2})nn'}\prod\limits_{i=0}^{n}P^{(i)}(\Pi)^{sp(i,\Pi;\Pi')}\prod\limits_{k=0}^{n'}P^{(k)}(\Pi')^{sp(k,\Pi';\Pi)}.\end{equation}

\item If $n\equiv n' (\text{mod }2)$, then for any critical value $m$ for $\Pi\otimes \Pi'$ such that $m\geq 1$, or $m\geq 0$ along with a non-vanishing condition, we have:
\begin{equation}\nonumber
L(m, \Pi\times (\Pi'\otimes\psi)) \sim_{E(\Pi)E(\Pi')E(\psi);K}(2\pi i)^{mnn'}\prod\limits_{i=0}^{n}P^{(i)}(\Pi)^{sp(i,\Pi;\Pi'\otimes \psi)}\prod\limits_{k=0}^{n'}P^{(k)}(\Pi'\otimes \psi)^{sp(k,\Pi' \otimes \psi;\Pi)}.\end{equation}

\end{enumerate}

\end{thm}

\chapter{Special values at $1$ of $L$-functions for automorphic pairs over quadratic imaginary fields}
\section{Settings}
\text{}

Let $r_{1}$ and $r_{2}$ be two positive integers.

Let $\Pi_{1}$ and $\Pi_{2}$ be two cuspidal representations of $GL_{r_{1}}(\AK)$ and $GL_{r_{2}}(\AK)$ respectively which has definable arithmetic automorphic periods. We assume they are also conjugate self-dual.\\

We write the infinity type of $\Pi_{1}$ (resp. $\Pi_{2}$) by $(z^{b_{j}}\overline{z}^{-b_{j}})_{1\leq j\leq r_{1}}$ (resp. $(z^{c_{k}}\overline{z}^{-c_{k}})_{1\leq k\leq r_{2}}$).  We see that $b_{j}\in \Z+\frac{r_{1}-1}{2}$ for all $1\leq j\leq r_{1}$ (resp. $c_{k}\in \Z+\frac{r_{2}-1}{2}$ for all $1\leq k\leq r_{2}$).\\

\begin{enumerate}[label=(\Alph*)]
\item If $r_{1}\equiv r_{2}\equiv 0(\text{ mod }2)$, we write $\Pi^{\#}=\Pi_{1}\boxplus \Pi_{2}^{c}$. We define $T_{3}=T_{4}=0$.
\item If $r_{1}\equiv r_{2}\equiv 1(\text{ mod }2)$, we write $\Pi^{\#}=(\Pi_{1}\otimes ||\cdot||_{\AK}^{-\frac{1}{2}}\psi ) \boxplus (\Pi_{2}^{c}\otimes ||\cdot||_{\AK}^{-\frac{1}{2}}\psi)$. We define $T_{3}=T_{4}=\frac{1}{2}$.
\item If $r_{1}\nequiv r_{2}(\text{mod }2)$, we may assume that $r_{1}$ is even and $r_{2}$ is odd. We write $\Pi^{\#}=(\Pi_{1}\otimes ||\cdot||_{\AK}^{-\frac{1}{2}}\psi ) \boxplus \Pi_{2}^{c}$. We define $T_{3}=\frac{1}{2}$ and $T_{4}=0$.
\end{enumerate}

It is easy to see that $\Pi^{\#}$ is an algebraic generic representation of $GL_{r_{1}+r_{2}}(\AK)$ with infinity type $(z^{b_{j}+T_{3}}\overline{z}^{-b_{j}-T_{3}}, z^{-c_{k}+T_{4}}\overline{z}^{c_{k}-T_{4}})_{1\leq j\leq r_{1},1\leq k\leq r_{2}}$. \\

We assume that $\Pi^{\#}$ is regular, i.e. for any $1\leq j\leq r_{1}$ and any $1\leq k\leq r_{2}$, we have $b_{j}+T_{3}\neq -c_{k}+T_{4}$.\\

Set $n=r_{1}+r_{2}+1$. We see that $\{b_{j}+T_{3}\mid 1\leq j\leq r_{1}\}\cup \{-c_{k}+T_{4}\mid 1\leq k\leq r_{2}\}$ are $n-1$ different numbers in $\Z+\frac{n-2}{2}$. We take $a_{1}>a_{2}>\cdots>a_{n}\in \Z+\frac{n-1}{2}$ such that the $n-1$ numbers above are in different gaps between $\{a_{i}\mid 1\leq i\leq n\}$. Let $\Pi$ be a cuspidal conjugate self-dual representation of $GL_{n}(\AK)$ which has arithmetic automorphic periods and infinity type $(z^{a_{i}}\overline{z}^{-a_{i}})$. \\

Our method also requires $\Pi$ to be $3$-regular. To guarantee this, we assume that \begin{equation}\label{r1*r2 very regular}
|(b_{j}+T_{3})-(-c_{k}+T_{4})|\geq 3  \text{ for all }1\leq j\leq r_{1},1\leq k\leq r_{2}.\end{equation}
In this case, we say the pair $(\Pi_{1},\Pi_{2})$ is \textbf{very regular}. We can then take $a_{i}$ as above such that $1+\frac{1}{2}$ is critical for $\Pi\otimes \Pi^{\#}$. Moreover, results in \cite{harrisunitaryperiod} show the existence of $\Pi$ as above, such that $L(1+\frac{1}{2},\Pi\otimes \Pi^{\#})\neq 0$. \\

We fix such $\Pi$ and $m=1$, then $m+\frac{1}{2}$ is critical for $\Pi\times \Pi^{\#}$ and moreover

\begin{equation}\label{application}
L(\frac{1}{2}+m,\Pi\times \Pi^{\#})\sim_{E(\Pi)E(\Pi^{\#});K} p(\Pi)p(\Pi^{\#})p(m,\Pi_{\infty},\Pi^{\#}_{\infty})
\end{equation}
with both sides non zero.\\

In the end of this section, let us show some simple facts on the split index. We can read from the construction of $a_{i}$ that 

\begin{equation}\nonumber
sp(j,\Pi_{1}\otimes \psi^{2T_{3}};\Pi)=sp(j,\Pi_{1}\otimes \psi^{2T_{3}};\Pi_{2}\otimes (\psi)^{2T_{4}})+1\text{ for all }0\leq j\leq r_{1}\end{equation}

\begin{eqnarray}
\text{ and similarly, }&&sp(j,\Pi_{2}^{c}\otimes \psi^{2T_{4}};\Pi)=sp(j,(\Pi_{2}\otimes (\psi^{c})^{2T_{4}})^{c};(\Pi_{1}\otimes \psi^{2T_{3}})^{c})+1\nonumber\\
&=&sp(r_{2}-j,\Pi_{2}\otimes (\psi^{c})^{2T_{4}});\Pi_{1}\otimes \psi^{2T_{3}})+1\text{ for all }0\leq j\leq r_{2}\nonumber
\end{eqnarray}

Here we have used Lemma \ref{split lemma}. \\

Moreover, for each $1\leq i\leq n-1$, one of $sp(i,\Pi;\Pi_{1}\otimes (\psi^{c})^{2T_{3}})$ and $sp(i,\Pi;\Pi_{2}^{c}\otimes \psi^{2T_{4}})$ is $1$ and another is $0$. We also know that $sp(0,\Pi;\Pi_{1}\otimes \psi^{2T_{3}})=sp(0,\Pi;\Pi_{2}^{c}\otimes \psi^{2T_{4}})=0$ and $sp(n,\Pi;\Pi_{1}\otimes \psi^{2T_{3}})=sp(n,\Pi;\Pi_{2}^{c}\otimes \psi^{2T_{3}})=0$.

\section{Simplify the left hand side}
We are going to simply the left hand side of equation (\ref{application}) now.
\begin{enumerate}[label=(\Alph*)]
\item In this case we have $L(m+\frac{1}{2},\Pi\times \Pi^{\#})=L(m+\frac{1}{2},\Pi\times \Pi_{1})\times L(m+\frac{1}{2},\Pi\times \Pi_{2}^{c})$.

By Theorem \ref{n*rmaintheorem}, we know that\begin{eqnarray} 
&L(\cfrac{1}{2}+m, \Pi\times \Pi_{1})\sim_{E(\Pi)E(\Pi_{1});K}& \\\nonumber
&(2\pi i)^{(m+\frac{1}{2})nr_{1}}\prod\limits_{i=0}^{n}P^{(i)}(\Pi)^{sp(i,\Pi;\Pi_{1})}\prod\limits_{j=0}^{r_{1}}P^{(j)}(\Pi_{1})^{sp(j,\Pi_{1};\Pi)}\\\nonumber
\text{and similarly} &L(\cfrac{1}{2}+m, \Pi\times \Pi_{2}^{c})\sim_{E(\Pi)E(\Pi_{2});K} &\\\nonumber
&(2\pi i)^{(m+\frac{1}{2})nr_{2}}\prod\limits_{i=0}^{n}P^{(i)}(\Pi)^{sp(i,\Pi;\Pi_{2}^{c})}\prod\limits_{k=0}^{r_{2}}P^{(k)}(\Pi_{2}^{c})^{sp(k,\Pi_{2}^{c};\Pi)}&.\end{eqnarray}

Therefore, since $sp(i,\Pi;\Pi_{1})+sp(i,\Pi;\Pi_{2}^{c})=1$ for all $1\leq i\leq n-1$, we obtain that
\begin{eqnarray}\label{lhsapplication1}
&&L(m+\frac{1}{2},\Pi\times \Pi^{\#})\\
&\sim_{E(\Pi)E(\Pi)E(\Pi_{2});K} &(2\pi i)^{(m+\frac{1}{2})n(n-1)} \prod\limits_{i=0}^{n}P^{(i)}(\Pi)^{sp(i,\Pi;\Pi_{1})+sp(i,\Pi;\Pi_{2}^{c})}\nonumber \\
&&
\prod\limits_{j=0}^{r_{1}}P^{(j)}(\Pi_{1})^{sp(j,\Pi_{1};\Pi)}
\prod\limits_{k=0}^{r_{2}}P^{(k)}(\Pi_{2})^{sp(k,\Pi_{2}^{c};\Pi)}
\nonumber\\
&\sim_{E(\Pi)E(\Pi)E(\Pi_{2});K} &(2\pi i)^{(m+\frac{1}{2})n(n-1)}\prod\limits_{i=1}^{n-1}P^{(i)}(\Pi)\nonumber \\\nonumber
&&\prod\limits_{j=0}^{r_{1}}P^{(j)}(\Pi_{1})^{sp(j,\Pi_{1};\Pi)}\prod\limits_{k=0}^{r_{2}}P^{(k)}(\Pi_{2}^{c})^{sp(k,\Pi_{2}^{c};\Pi)}.
\end{eqnarray}

\item In this case, we have $L(m+\frac{1}{2},\Pi\times \Pi^{\#})=L(m,\Pi\times (\Pi_{1}\otimes\psi))\times L(m,\Pi\times (\Pi_{2}^{c}\otimes\psi))$.

Applying the second part of Theorem \ref{n*rmaintheorem}, we have \begin{eqnarray}
&L(m+\frac{1}{2},\Pi\times \Pi^{\#}) \sim_{E(\Pi)E(\Pi)E(\Pi_{2});K}&\\
\nonumber&
(2\pi i)^{mn(n-1)} \prod\limits_{i=1}^{n-1}P^{(i)}(\Pi)
\prod\limits_{j=0}^{r_{1}}P^{(j)}(\Pi_{1})^{sp(j,\Pi_{1}\otimes\psi;\Pi)}
\prod\limits_{k=0}^{r_{2}}P^{(k)}(\Pi_{2}^{c})^{sp(k,\Pi_{2}^{c}\otimes \psi;\Pi)}\times&\\\nonumber&
p(\widecheck{\psi},1)^{\sum\limits_{j=0}^{r_{1}}j*sp(j,\Pi_{1}\otimes \psi;\Pi)+\sum\limits_{k=0}^{r_{2}}k*sp(k,\Pi_{2}^{c}\otimes \psi;\Pi)}
p(\widecheck{\psi},\iota)^{\sum\limits_{j=0}^{n'}(r_{1}-j)*sp(j,\Pi_{1}\otimes \psi;\Pi)+\sum\limits_{k=0}^{r_{2}}(r_{2}-k)*sp(k,\Pi_{2}^{c}\otimes \psi;\Pi)}.&
\end{eqnarray}

\begin{lem} We have:
\begin{eqnarray}
&&\sum\limits_{j=0}^{r_{1}}j*sp(j,\Pi_{1}\otimes \psi;\Pi)+\sum\limits_{k=0}^{r_{2}}k*sp(k,\Pi_{2}^{c}\otimes \psi;\Pi)\nonumber \\
&=&\sum\limits_{j=0}^{n'}(r_{1}-j)*sp(j,\Pi_{1}\otimes \psi;\Pi)+\sum\limits_{k=0}^{r_{2}}(r_{2}-k)*sp(k,\Pi_{2}^{c}\otimes \psi;\Pi)\nonumber\\
&=&\frac{n(n-1)}{2}\nonumber\end{eqnarray}

\end{lem}

\begin{dem}
We set $w(j,\Pi_{1}\otimes \psi;\Pi)$, $1\leq j\leq r_{1}$ (resp. $w(k,\Pi_{2}^{c}\otimes \psi;\Pi)$, $1\leq k\leq r_{2}$) to be the index $w(j)$ for the pair $(\Pi,\Pi_{1}\otimes \psi)$ (resp. $(\Pi,\Pi_{2}^{c}\otimes \psi)$)) as in (\ref{relationn*r}). We see from (\ref{sumofw}) that $\sum\limits_{j=0}^{r_{1}}j*sp(j,\Pi_{1}\otimes \psi;\Pi)=\sum\limits_{j=1}^{r_{1}}w(j,\Pi_{1}\otimes \psi;\Pi)$ and $\sum\limits_{k=0}^{r_{2}}k*sp(k,\Pi'\otimes \psi;\Pi)=\sum\limits_{k=1}^{r_{1}}w(k,\Pi_{2}^{c}\otimes \psi;\Pi)$.

Recall that $w(j,\Pi_{1}\otimes \psi;\Pi)$ (resp. $w(k,\Pi_{2}^{c}\otimes \psi;\Pi)$) is the position of the infinity type of $\Pi_{1}\otimes \psi$ (resp. $\Pi_{2}^{c}\otimes \psi$) in the gaps of the infinity type of $\Pi$. It is easy to see that the $n-1$ numbers $w(j,\Pi_{1}\otimes \psi;\Pi)$, $w(k,\Pi_{2}^{c}\otimes \psi;\Pi)$ for $1\leq j\leq r_{1}$ and $1\leq k\leq r_{2}$ runs over $1,2,\cdots, n-1$. We then deduce the first formula of the lemma. 

The second follows easily from the first one.
\end{dem}

From the lemma we see that 
\begin{eqnarray}
&p(\widecheck{\psi},1)^{\sum\limits_{j=0}^{r_{1}}j*sp(j,\Pi_{1}\otimes \psi;\Pi)+\sum\limits_{k=0}^{r_{2}}k*sp(k,\Pi_{2}^{c}\otimes \psi;\Pi)}p(\widecheck{\psi},\iota)^{\sum\limits_{j=0}^{n'}(r_{1}-j)*sp(j,\Pi_{1}\otimes \psi;\Pi)+\sum\limits_{k=0}^{r_{2}}(r_{2}-k)*sp(k,\Pi_{2}^{c}\otimes \psi;\Pi)}\nonumber&\\&
\sim_{E(\psi);K} (2\pi i)^{\frac{n(n-1)}{2}}.&
\end{eqnarray}

We thus obtain that
\begin{eqnarray}\label{lhsapplication2}
&L(m+\frac{1}{2},\Pi\times \Pi^{\#})\sim_{E(\Pi)E(\Pi)E(\Pi_{2});K} (2\pi i)^{(m+\frac{1}{2})n(n-1)}\prod\limits_{i=1}^{n-1}P^{(i)}(\Pi)\nonumber &\\
&\prod\limits_{j=0}^{r_{1}}P^{(j)}(\Pi_{1})^{sp(j,\Pi_{1}\otimes\psi;\Pi)}\prod\limits_{k=0}^{r_{2}}P^{(k)}(\Pi_{2}^{c})^{sp(k,\Pi_{2}^{c}\otimes\psi;\Pi)}.&
\end{eqnarray}

\item In this case, we have $L(m+\frac{1}{2},\Pi\times \Pi^{\#})=L(m,\Pi\times (\Pi_{1}\otimes \psi))\times L(m+\frac{1}{2},\Pi\times \Pi_{2}^{c})$.

Similarly, we get:
\begin{eqnarray}\label{lhsapplication3}
&L(m+\frac{1}{2},\Pi\times \Pi^{\#}) \sim_{E(\Pi)E(\Pi)E(\Pi_{2});K} (2\pi i)^{(m+\frac{1}{2})n(n-1)-\frac{nr_{1}}{2}} \prod\limits_{i=1}^{n-1}P^{(i)}(\Pi)&\nonumber\\
&\prod\limits_{j=0}^{r_{1}}P^{(j)}(\Pi_{1})^{sp(j,\Pi_{1}\otimes \psi;\Pi)}
\prod\limits_{k=0}^{r_{2}}P^{(k)}(\Pi_{2}^{c})^{sp(k,\Pi_{2}^{c};\Pi)}p(\widecheck{\psi},1)^{\sum\limits_{j=0}^{r_{1}}j*sp(j,\Pi_{1}\otimes \psi;\Pi)}
p(\widecheck{\psi},\iota)^{\sum\limits_{j=0}^{n'}(r_{1}-j)*sp(j,\Pi_{1}\otimes \psi;\Pi)}.&
\end{eqnarray}

\end{enumerate}

\section{Simplify the right hand side}
By Corollary \ref{Shahidi Whittaker} and Corollary \ref{divise}, for cases (A) and (B), we have:
\begin{eqnarray}
p(\Pi^{\#})&\sim_{E(\Pi^{\#});K}&\Omega(\Pi^{\#}_{\infty})p(\Pi_{1})\Omega(\Pi_{1,\infty})^{-1}p(\Pi_{2})\Omega(\Pi_{2,\infty})^{-1}L(1,\Pi_{1}\times\Pi_{2})\nonumber\\
&\sim_{E(\Pi^{\#});K}&\Omega(\Pi^{\#}_{\infty})Z(\Pi_{1,\infty})\Omega(\Pi_{1,\infty})^{-1}Z(\Pi_{2,\infty})\Omega(\Pi_{2,\infty})^{-1}L(1,\Pi_{1}\times\Pi_{2})\times\nonumber\\
&& \prod\limits_{j=1}^{r_{1}-1}P^{(j)}(\Pi_{1})\prod\limits_{k=1}^{r_{2}-1}P^{(k)}(\Pi_{2}^{c})\nonumber\\
&\sim_{E(\Pi^{\#});K}& (2\pi i)^{\frac{(r_{1}-1)r_{1}}{2}+\frac{(r_{2}-1)r_{2}}{2}}\Omega(\Pi^{\#}_{\infty})L(1,\Pi_{1}\times\Pi_{2})\prod\limits_{j=1}^{r_{1}-1}P^{(j)}(\Pi_{1})\prod\limits_{k=1}^{r_{2}-1}P^{(k)}(\Pi_{2}^{c}).\nonumber
\end{eqnarray}

Therefore, for cases (A) and (B), we obtain that:
\begin{eqnarray}\label{rhsapplication}
&&p(\Pi)p(\Pi^{\#})p(m,\Pi_{\infty},\Pi^{\#}_{\infty})\nonumber\\
&\sim_{E(\Pi^{\#});K} &(2\pi i)^{\frac{(r_{1}-1)r_{1}}{2}+\frac{(r_{2}-1)r_{2}}{2}}\Omega(\Pi^{\#}_{\infty})Z(\Pi_{\infty})p(m,\Pi_{\infty},\Pi^{\#}_{\infty})\times\nonumber\\
&&L(1,\Pi_{1}\times\Pi_{2})\prod\limits_{i=1}^{n-1}P^{(i)}(\Pi)\prod\limits_{j=1}^{r_{1}-1}P^{(j)}(\Pi_{1})\prod\limits_{k=1}^{r_{2}-1}P^{(k)}(\Pi_{2}^{c})\nonumber\\
&\sim_{E(\Pi^{\#});K}& (2\pi i)^{n(n-1)(m+\frac{1}{2})-\frac{n(n-1)}{2}+\frac{(r_{1}-1)r_{1}}{2}+\frac{(r_{2}-1)r_{2}}{2}}L(1,\Pi_{1}\times\Pi_{2})\times\nonumber \\
&&\prod\limits_{i=1}^{n-1}P^{(i)}(\Pi)\prod\limits_{j=1}^{r_{1}-1}P^{(j)}(\Pi_{1})\prod\limits_{k=1}^{r_{2}-1}P^{(k)}(\Pi_{2}^{c})\nonumber \\
&\sim_{E(\Pi^{\#});K}& (2\pi i)^{n(n-1)(m+\frac{1}{2})-r_{1}r_{2}}L(1,\Pi_{1}\times\Pi_{2})\prod\limits_{i=1}^{n-1}P^{(i)}(\Pi)\times\nonumber\\
&&\prod\limits_{j=0}^{r_{1}}P^{(j)}(\Pi_{1})\prod\limits_{k=0}^{r_{2}}P^{(k)}(\Pi_{2}^{c})
\end{eqnarray}
We have used Lemma \ref{threeproduct}, the fact that ${n-1\choose 2}={r_{1}+r_{2}\choose 2}={r_{1}\choose 2}+{r_{2} \choose 2}+r_{1}r_{2}$ and also the fact that $P^{(0)}(\Pi_{1})P^{(r_{1})}(\Pi_{1})\sim_{E(\Pi_{1})} 1$, $P^{(0)}(\Pi_{2}^{c})P^{(r_{2})}(\Pi_{2}^{c})\sim_{E(\Pi_{2})} 1$.

For case (C), we only need to change $L(1,\Pi_{1}\times \Pi_{2})$ to $L(\frac{1}{2},(\Pi_{1}\otimes\psi)\times \Pi_{2})$ in the above formula.

\section{Final conclusion}
Comparing (\ref{lhsapplication1}) and (\ref{rhsapplication}), we get for case (A):
\begin{eqnarray}
L(1,\Pi_{1}\times \Pi_{2})&\sim_{E(\Pi_{1})E(\Pi_{2});K}&(2\pi i)^{r_{1}r_{2}}\prod\limits_{j=0}^{r_{1}}P^{(j)}(\Pi_{1})^{sp(j,\Pi_{1};\Pi)-1}
\prod\limits_{k=0}^{r_{2}}P^{(k)}(\Pi_{2}^{c})^{sp(k,\Pi_{2}^{c};\Pi)-1}\nonumber\\
&\sim_{E(\Pi_{1})E(\Pi_{2});K}&(2\pi i)^{r_{1}r_{2}}\prod\limits_{j=0}^{r_{1}}P^{(j)}(\Pi_{1})^{sp(j,\Pi_{1};\Pi_{2})}
\prod\limits_{k=0}^{r_{2}}P^{(k)}(\Pi_{2}^{c})^{sp(k,\Pi_{2}^{c};\Pi_{1}^{c})}\nonumber\\
&\sim_{E(\Pi_{1})E(\Pi_{2});K}&(2\pi i)^{r_{1}r_{2}}\prod\limits_{j=0}^{r_{1}}P^{(j)}(\Pi_{1})^{sp(j,\Pi_{1};\Pi_{2})}
\prod\limits_{k=0}^{r_{2}}P^{(r_{2}-k)}(\Pi_{2})^{sp(r_{2}-k,\Pi_{2};\Pi_{1})}\nonumber\\
&\sim_{E(\Pi_{1})E(\Pi_{2});K}&(2\pi i)^{r_{1}r_{2}}\prod\limits_{j=0}^{r_{1}}P^{(j)}(\Pi_{1})^{sp(j,\Pi_{1};\Pi_{2})}
\prod\limits_{k=0}^{r_{2}}P^{(k)}(\Pi_{2})^{sp(k,\Pi_{2};\Pi_{1})}\nonumber.
\end{eqnarray}

Comparing (\ref{lhsapplication2}) and (\ref{rhsapplication}), we get for case (B):
\begin{eqnarray}
L(1,\Pi_{1}\times \Pi_{2})&\sim_{E(\Pi_{1})E(\Pi_{2});K}&(2\pi i)^{r_{1}r_{2}}\prod\limits_{j=0}^{r_{1}}P^{(j)}(\Pi_{1})^{sp(j,\Pi_{1}\otimes\psi;\Pi_{2}\otimes\psi^{c})}
\prod\limits_{k=0}^{r_{2}}P^{(k)}(\Pi_{2})^{sp(k,\Pi_{2}\otimes\psi^{c};\Pi_{1}\otimes\psi)}\nonumber\\
&\sim_{E(\Pi_{1})E(\Pi_{2});K}&(2\pi i)^{r_{1}r_{2}}\prod\limits_{j=0}^{r_{1}}P^{(j)}(\Pi_{1})^{sp(j,\Pi_{1};\Pi_{2})}
\prod\limits_{k=0}^{r_{2}}P^{(k)}(\Pi_{2})^{sp(k,\Pi_{2};\Pi_{1})}\nonumber.
\end{eqnarray}

Here we have used that $sp(j,\Pi_{1}\otimes\psi;\Pi_{2}\otimes\psi^{c})=sp(j,\Pi_{1}\otimes\psi;\Pi_{2}\otimes\psi^{-1})=sp(j,\Pi_{1};\Pi_{2})$ by Lemma \ref{split lemma}.

Similarly, for case $(C)$, comparing (\ref{lhsapplication3}) and (\ref{rhsapplication}), we obtain that:
\begin{eqnarray}
&&L(\frac{1}{2},(\Pi_{1}\otimes\psi)\times \Pi_{2})\nonumber\\
&\sim_{E(\Pi_{1})E(\Pi_{2})E(\psi);K}&(2\pi i)^{r_{1}r_{2}-\frac{nr_{1}}{2}}\prod\limits_{j=0}^{r_{1}}P^{(j)}(\Pi_{1})^{sp(j,\Pi_{1}\otimes\psi;\Pi_{2})}
\prod\limits_{k=0}^{r_{2}}P^{(k)}(\Pi_{2})^{sp(k,\Pi_{2};\Pi_{1}\otimes\psi)}\times\nonumber\\
&&p(\widecheck{\psi},1)^{\sum\limits_{j=0}^{r_{1}}j*(sp(j,\Pi_{1}\otimes\psi;\Pi_{2})+1)}
p(\widecheck{\psi},\iota)^{\sum\limits_{j=0}^{r_{1}}(r_{1}-j)*(sp(j,\Pi_{1}\otimes\psi;\Pi_{2}+1))}\nonumber\\\nonumber
&\sim_{E(\Pi_{1})E(\Pi_{2})E(\psi);K}&(2\pi i)^{\frac{r_{1}r_{2}}{2}-\frac{r_{1}(r_{1}+1)}{2}}
\prod\limits_{j=0}^{r_{1}}P^{(j)}(\Pi_{1})^{sp(j,\Pi_{1}\otimes\psi;\Pi_{2})}
\prod\limits_{k=0}^{r_{2}}P^{(k)}(\Pi_{2})^{sp(k,\Pi_{2};\Pi_{1}\otimes\psi)}\times\\
&&p(\widecheck{\psi},1)^{\sum\limits_{j=0}^{r_{1}}j*sp(j,\Pi_{1}\otimes\psi;\Pi_{2})+\frac{r_{1}(r_{1}+1)}{2}}
p(\widecheck{\psi},\iota)^{\sum\limits_{j=0}^{r_{1}}(r_{1}-j)*sp(j,\Pi_{1}\otimes\psi;\Pi_{2})+\frac{r_{1}(r_{1}+1)}{2}}\nonumber\\
&\sim_{E(\Pi_{1})E(\Pi_{2})E(\psi);K}&(2\pi i)^{\frac{r_{1}r_{2}}{2}}\prod\limits_{j=0}^{r_{1}}P^{(j)}(\Pi_{1})^{sp(j,\Pi_{1}\otimes\psi;\Pi_{2})}
\prod\limits_{k=0}^{r_{2}}P^{(k)}(\Pi_{2})^{sp(k,\Pi_{2};\Pi_{1}\otimes\psi)}\times\nonumber\\
&&p(\widecheck{\psi},1)^{\sum\limits_{j=0}^{r_{1}}j*sp(j,\Pi_{1}\otimes\psi;\Pi_{2})}
p(\widecheck{\psi},\iota)^{\sum\limits_{j=0}^{r_{1}}(r_{1}-j)*sp(j,\Pi_{1}\otimes\psi;\Pi_{2})}\nonumber\\
&\sim_{E(\Pi_{1})E(\Pi_{2})E(\psi);K}&(2\pi i)^{\frac{r_{1}r_{2}}{2}}\prod\limits_{j=0}^{r_{1}}P^{(j)}(\Pi_{1}\otimes \psi)^{sp(j,\Pi_{1}\otimes\psi;\Pi_{2})}\prod\limits_{k=0}^{r_{2}}P^{(k)}(\Pi_{2})^{sp(k,\Pi_{2};\Pi_{1}\otimes\psi)}\nonumber.
\end{eqnarray}
The last step is deduced by definition of $P^{(*)(\Pi_{1}\otimes \psi)}$ (c.f. Definition-Lemma \ref{definable periods}).

\begin{thm}
Let $r_{1}$ and $r_{2}$ be two positive integers.
Let $\Pi_{1}$ and $\Pi_{2}$ be two cuspidal representations of $GL_{r_{1}}(\AK)$ and $GL_{r_{2}}(\AK)$ respectively which are very regular, cohomological, conjugate self-dual and supercuspidal at at least two finite split places.
Assume that the pair $\Pi_{1},\Pi_{2}$ is very regular in the sense of (\ref{r1*r2 very regular}).
\begin{enumerate}[label=(\roman*)]
\item If $r_{1}\equiv r_{2} (\text{mod }2)$, then $L(1,\Pi_{1}\times \Pi_{2})\sim_{E(\Pi_{1})E(\Pi_{2});K}$
 \begin{equation}
 (2\pi i)^{r_{1}r_{2}}\prod\limits_{j=0}^{r_{1}}P^{(j)}(\Pi_{1})^{sp(j,\Pi_{1};\Pi_{2})}
\prod\limits_{k=0}^{r_{2}}P^{(k)}(\Pi_{2})^{sp(k,\Pi_{2};\Pi_{1})}.\nonumber
\end{equation}
\item If $r_{1}\nequiv r_{2} (\text{mod }2)$, then $L(\frac{1}{2},(\Pi_{1}\otimes\psi)\times \Pi_{2}) \sim_{E(\Pi_{1})E(\Pi_{2})E(\psi);K}$
\begin{eqnarray}
(2\pi i)^{\frac{r_{1}r_{2}}{2}}\prod\limits_{j=0}^{r_{1}}P^{(j)}(\Pi_{1}\otimes \psi)^{sp(j,\Pi_{1}\otimes\psi;\Pi_{2})}\prod\limits_{k=0}^{r_{2}}P^{(k)}(\Pi_{2})^{sp(k,\Pi_{2};\Pi_{1}\otimes\psi)}\nonumber.
\end{eqnarray}

\end{enumerate}

\end{thm}

\bibliography{bibfile}
\bibliographystyle{alpha}

\end{document}